\documentstyle[12pt,epic,eepic]{article}
\title{{\bf On $\alpha$-Induction, Chiral Generators
and Modular Invariants for Subfactors}}

\author{{\sc Jens B\"ockenhauer} and {\sc David E. Evans}\\
School of Mathematics\\
University of Wales, Cardiff\\
PO Box 926, Senghennydd Road\\
Cardiff CF2 4YH, Wales, U.K.\\
e-mail: {\tt BockenhauerJM@cf.ac.uk, EvansDE@cf.ac.uk}  \\
\vphantom{X}\\
{\sc Yasuyuki Kawahigashi}\\
Department of Mathematical Sciences\\
University of Tokyo, Komaba, Tokyo, 153-8914, JAPAN\\
e-mail: {\tt yasuyuki@ms.u-tokyo.ac.jp}}
\date{June 23, 1999}

\begin{document}
\maketitle

\input amssym.def
\newsymbol\varnothing 203F
\newsymbol\rtimes 226F
\def\emptyset{\varnothing}

\def\Ad            {{\rm{Ad}}}
\def\Aut           {{\rm{Aut}}}
\def\bbC           {\Bbb{C}}
\def\bbM           {\Bbb{M}}
\def\bbN           {\Bbb{N}}
\def\bbNo          {\Bbb{N}_0}
\def\bbR           {\Bbb{R}}
\def\bbT           {\Bbb{T}}
\def\bbZ           {\Bbb{Z}}
\def\be            {\begin{equation}}
\def\bearl         {\begin{array}{l}}
\def\bearll        {\begin{array}{ll}}
\def\bearlll       {\begin{array}{lll}}
\def\bearrl        {\begin{array}{rl}}
\def\bea           {\begin{eqnarray}}
\def\beaa          {\begin{eqnarray*}}
\def\bfe           {{\bf1}}
\def\can           {\gamma}
\def\canr          {\theta}
\def\cA            {{\cal{A}}}
\def\cC            {{\cal{C}}}
\def\cCA           {\frak{C}}
\def\cD            {{\cal{D}}}
\def\cE            {{\cal{E}}}
\def\cF            {{\cal{F}}}
\def\cG            {{\cal{G}}}
\def\cH            {{\cal{H}}}
\def\cJ            {{\cal{J}}}
\def\cK            {{\cal{K}}}
\def\cL            {{\cal{L}}}
\def\cM            {{\cal{M}}}
\def\cN            {{\cal{N}}}
\def\cO            {{\cal{O}}}
\def\cP            {{\cal{P}}}
\def\cR            {{\cal{R}}}
\def\cS            {{\cal{S}}}
\def\cT            {{\cal{T}}}
\def\cV            {{\cal{V}}}
\def\cW            {{\cal{W}}}
\def\cX            {{\cal{X}}}
\def\cY            {{\cal{Y}}}
\def\cZ            {{\cal{Z}}}
\newcommand\co[1]  {\bar{{#1}}}
\def\dim           {{\rm{dim}}}
\newcommand\del[2] {\delta_{{#1},{#2}}}
\def\E             {{\rm{e}}}
\def\ee            {\end{equation}}
\def\eear          {\end{array}}
\def\eea           {\end{eqnarray}}
\def\eeaa          {\end{eqnarray*}}
\def\End           {{\rm{End}}}
\newcommand\eps[2] {\varepsilon({#1},{#2})}
\newcommand\epsm[2]{\varepsilon^-({#1},{#2})}
\newcommand\epsp[2]{\varepsilon^+({#1},{#2})}
\newcommand\epspm[2]{\varepsilon^\pm({#1},{#2})}
\newcommand\epsmp[2]{\varepsilon^\mp({#1},{#2})}
\newcommand\Eps[2] {{\cal E}({#1},{#2})}
\newcommand\Epsm[2]{{\cal E}^-({#1},{#2})}
\newcommand\Epsp[2]{{\cal E}^+({#1},{#2})}
\newcommand\Epspm[2]{{\cal E}^\pm({#1},{#2})}
\newcommand\Epsmp[2]{{\cal E}^\mp({#1},{#2})}
\newcommand\erf[1] {Eq.\ (\ref{#1})}
\def\fA            {\frak{A}}
\def\fAg           {\frak{A}}
\def\fB            {\frak{B}}
\def\fF            {\frak{F}}
\def\fFg           {\frak{F}}
\def\fh            {\frak{h}}
\def\fW            {\frak{W}}
\def\GLZ           {\it{GL}\,(2;\bbC)}
\def\Hom           {{\rm{Hom}}}
\def\I             {{\rm{i}}}
\def\id            {{\rm{id}}}
\def\iotab         {{\co\iota}}
\def\Jz            {\cal{J}_z}
\def\lan           {\langle}
\def\lab           {{\co \lambda}}
\newcommand\laend[2]{\lambda_{{#1},{#2}}}
\def\LG            {{\it{LG}}}
\def\LH            {{\it{LH}}}
\def\LIG           {{\it{L}_I\it{G}}}
\def\LIcG          {{\it{L}_{\Ic}\it{G}}}
\def\LIH           {{\it{L}_I\it{H}}}
\def\LIcSUn        {{\it{L}_{I'}\it{SU}(n)}}
\def\LISUn         {{\it{L}_I\it{SU}(n)}}
\def\LISUk         {{\it{L}_I\it{SU}(k)}}
\def\LISUnk        {{\it{L}_I\it{SU}(nk)}}
\newcommand\ls[1]  {[\lambda_{{#1}}]}
\def\LSE           {L^2(S^1)}
\def\LSn           {L^2(S^1;\Bbb{C}^n)}
\def\LSUd          {{\it{LSU}(3)}}
\def\LSUk          {{\it{LSU}(k)}}
\def\LSUn          {{\it{LSU}(n)}}
\def\LSUnk         {{\it{LSU}(nk)}}
\def\LSUz          {{\it{LSU}(2)}}
\def\Mor           {{\rm{Mor}}}
\def\mub           {{\co{\mu}}}
\def\MXN           {{}_M {\cal X}_N}
\def\MXM           {{}_M {\cal X}_M}
\newcommand\N[3]   {N_{{#1},{#2}}^{{#3}}}
\def\Nres          {\tilde{N}}
\def\NXN           {{}_N {\cal X}_N}
\def\NXM           {{}_N {\cal X}_M}
\def\nub           {\overline{\nu}}
\def\oto           {=0,1,2,\ldots}
\def\pio           {\pi_0}
\def\PSLZ          {{\it{PSL}}(2;\bbZ)}
\def\PSU           {{\it{PSU}}(1,1)}
\def\reso          {|_{\cH_0}}
\def\Sect          {{\rm{Sect}}}
\def\SLC           {{\it{SL}}(2;\bbC)}
\def\SLnC          {{\it{SL}}(n;\bbC)}
\def\SLZ           {{\it{SL}}(2;\bbZ)}
\def\Ssys          {{\Sigma(\sys)}}
\def\SOf           {{\it{SO}}(5)}
\def\son           {\frak{so}(N)}
\def\SON           {{\it{SO}}(N)}
\def\SUd           {{\it{SU}}(3)}
\def\SUk           {{\it{SU}}(k)}
\def\SUm           {{\it{SU}}(m)}
\def\SUn           {{\it{SU}}(n)}
\def\SUnk          {{\it{SU}}(nk)}
\def\SUz           {{\it{SU}}(2)}
\def\SUzk          {{\it{SU}}(2k)}
\def\SUf           {{\it{SU}}(4)}
\def\sud           {\frak{su}(3)}
\def\sun           {\frak{su}(n)}
\def\suz           {\frak{su}(2)}
\def\sudh          {\widehat{\frak{su}}(3)}
\def\sunh          {\widehat{\frak{su}}(n)}
\def\suzh          {\widehat{\frak{su}}(2)}
\def\sys           {{\Delta}}
\def\tr            {{\rm{tr}}}
\def\Tr            {{\rm{Tr}}}
\def\Un            {\it{U}(n)}
\newcommand\V[3]   {V_{{#1};{#2}}^{#3}}
\def\Vir           {\frak{Vir}}

%%%%%%%%%%%%%%%%%%%%%%%%%%%%%%%%%%%%%%%%%%%%%%%%%%%%%%%%%%%%%%%%%%%%%%%%%%%
                         %Text Definitions%

\def\qed{{\unskip\nobreak\hfil\penalty50
\hskip2em\hbox{}\nobreak\hfil  $\Box$      %\rm Q.E.D.
\parfillskip=0pt \finalhyphendemerits=0\par}\medskip}
\def\proof{\trivlist \item[\hskip \labelsep{\bf Proof.\ }]}
\def\endproof{\null\hfill\qed\endtrivlist}

\def\equi{\sim}
\def\isom{\cong}
\def\ti{\tilde}
\def\lan{\langle}
\def\ran{\rangle}

\def\a{\alpha}
\def\de{\delta}
\def\e{\varepsilon}
\def\ga{\gamma}
\def\Ga{\Gamma}
\def\la{\lambda}
\def\La{\Lambda}
\def\th{\theta}
\def\om{\omega}
\def\Om{\Omega}
  
\newcommand\dta{\begin{picture}(12,10)
\thicklines
\path(2,4)(6,8)(10,4)(2,4)(6,0)(10,4)
\end{picture}}

\newcommand\dtap{\begin{picture}(12,10)\thicklines
\path(6,8)(10,4)(6,0)(6,8)(2,4)(6,0)\end{picture}}

%%%%%%%%%%%%%%%%%%%%%%%%%%%%%%%%%%%%%%%%%%%%%%%%%%%%%%%%%%%%%%%%%%%

\def\thinlines{\allinethickness{0.3pt}}
\def\thicklines{\allinethickness{1.0pt}}
\def\Thicklines{\allinethickness{2.0pt}}

%%%%%%%%%%%%%%%%%%%%%%%%%%%%%%%%%%%%%%%%%%%%%%%%%%%%%%%%%%%%%%%%%%%

\newtheorem{theorem}{Theorem}[section]
\newtheorem{lemma}[theorem]{Lemma}
\newtheorem{conjecture}[theorem]{Conjecture}
\newtheorem{corollary}[theorem]{Corollary}
\newtheorem{definition}[theorem]{Definition}
\newtheorem{assumption}[theorem]{Assumption}
\newtheorem{proposition}[theorem]{Proposition}
\newtheorem{remark}[theorem]{Remark}
\newtheorem{example}[theorem]{Example}

\begin{abstract}
We consider a type III subfactor
$N\subset M$ of finite index with a finite system of
braided $N$-$N$ morphisms which includes the irreducible
constituents of the dual canonical endomorphism.
We apply $\a$-induction and, developing further some
ideas of Ocneanu, we define chiral generators
for the double triangle algebra. Using a new concept
of intertwining braiding fusion relations, we show that the
chiral generators can be naturally identified with the
$\a$-induced sectors. A matrix $Z$ is defined and shown
to commute with the S- and T-matrices arising from the
braiding. If the braiding is non-degenerate, then $Z$
is a ``modular invariant mass matrix'' in the usual sense
of conformal field theory. We show that in that case the
fusion rule algebra of the dual system of $M$-$M$ morphisms
is generated by the images of both kinds of $\a$-induction,
and that the structural information about its irreducible
representations is encoded in the mass matrix $Z$.
Our analysis sheds further light on the connection between
(the classifications of) modular invariants and subfactors,
and we will construct and analyze modular invariants from
$\SUn_k$ loop group subfactors in a forthcoming
publication, including the treatment of all
$\SUz_k$ modular invariants.
\end{abstract}

\newpage

\tableofcontents

\section{Introduction}
\label{sec-intro}
It is a surprising fact that a series of at first sight
unrelated phenomena in mathematics and physics are
governed by the scheme of A-D-E Dynkin diagrams, such
as simple Lie algebras, finite subgroups of $\SLC$,
simple singularities of complex surfaces, quivers of
finite type, modular invariant partition functions of
$\SUz$ WZW models and subfactors of Jones index less than four.
Though a good understanding of the interrelations has not yet
been achieved, this coincidence indicates that there are deep
connections between these different fields which even seem to
go beyond the A-D-E governed cases, e.g.\ finite subgroups
of $\SLnC$, modular invariants of $\SUn$ WZW models, or (certain)
$\SUn_k$ subfactors of larger index. This paper is addressed to
the relation between the (classifications of) modular
invariants in conformal field theory and subfactors in
operator algebras.

In rational (chiral) conformal field theory
one deals with a chiral algebra which possesses a certain
finite spectrum of representations (or superselection sectors)
$\pi_\la$ acting on a Hilbert space $\cH_\la$. Its characters
$\chi_\la(\tau)=\tr_{\cH_\la}(\E^{2\pi\I\tau(L_0-c/24)})$,
$\rm{Im}(\tau)>0$, $L_0$ being the conformal Hamiltonian
and $c$ the central charge, transform unitarily under
``reparametrization of the torus'',
i.e.\ there are matrices $S$ and $T$ such that
\[ \chi_\la(- 1/\tau) = \sum_\mu S_{\la,\mu} \chi_\mu(\tau)\,,
\qquad \chi_\la(\tau+1) = \sum_\mu T_{\la,\mu} \chi_\mu(\tau)\,, \]
which are the generators of a unitary
representation of the (double cover of the) modular group
$\SLZ$ in which $T$ is diagonal.\footnote{More precisely,
for current algebras the characters depend also on other variables
than $\tau$, corresponding to Cartan subalgebra generators which
are omitted here for simplicity. But these variables are responsible
that one is in general dealing with the whole group $\SLZ$ rather
than $\PSLZ$.}
In order to classify conformal field theories, in particular
extensions in a certain sense of a given theory, one searches
for modular invariant partition functions
$Z(\tau)=Z(-1/\tau)=Z(\tau+1)$ of the form
\[ Z(\tau)= \sum_{\la,\mu} Z_{\la,\mu}
\chi_\la(\tau) \chi_\mu(\tau)^* \,,\]
where
\be
Z_{\la,\mu} = 0,1,2,\ldots \,, \qquad\qquad Z_{0,0}=1 \,.
\label{massmatZ}
\ee
Here the label ``0'' refers to the ``vacuum'' representation,
and the condition $Z_{0,0}=1$ reflects the physical concept of
uniqueness of the vacuum state. The matrix $Z$ arising this way is
called a modular invariant mass matrix. Mathematically speaking,
the problem can be rephrased like this: Find all the matrices $Z$
in the commutant of the unitary representation of $\SLZ$
defined by $S$ and $T$ subject to the conditions in \erf{massmatZ}.
In this paper we study this mathematical problem in
the subfactor context. We start with a von Neumann algebra,
more precisely a factor $N$ endowed with a system of braided
endomorphisms. Such a braiding defines matrices $S$ and $T$ which
provide a unitary representation of $\SLZ$ if it is non-degenerate.
We then study embeddings $N\subset M$ in larger factors
$M$ which are in a certain sense compatible with the braided
system of endomorphisms. We show that such an embedding
$N\subset M$ determines a modular invariant mass matrix
in exactly the sense specified above.

Longo and Rehren have studied nets of subfactors and
defined a useful formula to extend a localized transportable
endomorphism of the smaller to the larger observable algebra,
realizing a suggestion in \cite{R}. 
Xu \cite{X1,X2} has worked on essentially the same
construction applied to subfactors arising from
conformal inclusions with the loop group construction of
A. Wassermann \cite{W2}. Two of us systematically analyzed
the Longo-Rehren extension for nets of subfactors
on $S^1$ \cite{BE1,BE3}. As sectors, a reciprocity
between extension and restriction of localized
transportable endomorphisms was established, analogous
to the induction-restriction machinery of group representations,
and therefore the extension was called $\a$-induction in order
to avoid confusion with the different sector induction. It was also
noticed in \cite{BE1} that the extended endomorphisms leave
local algebras invariant and hence $\a$-induction
can also be considered as a map which takes certain
endomorphisms of a local subfactor to endomorphisms of the
embedding factor. This theory was applied to nets arising from
conformal field theory models in \cite{BE2,BE3}, and it
was shown that for all type I modular invariants of
$\SUz$ respectively $\SUd$ there are associated nets of subfactors
and in turn $\a$-induction gives rise to fusion graphs.
In fact it was shown that that these graphs are the
A-D-E Dynkin diagrams respectively their generalizations
of \cite{DZ1,DZ2}, and this is no accident: The homomorphism
property of $\a$-induction relates the spectrum of the
fusion graphs to the non-zero diagonal entries of the
modular invariant mass matrix.

A few months after the work of Longo-Rehren, Ocneanu presented
his theory of ``quantum symmetries'' of Coxeter graphs and gave
lectures \cite{O7} one year later.
He introduced a notion of a ``double triangle
algebra'' and defined elements $p_j^\pm$ which we refer to as
``chiral generators'' as they were not specifically named there.
Ocneanu's analysis has much in common with work of
Xu \cite{X1} and two of us \cite{BE2,BE3} about subfactors of
type E$_6$, E$_8$ and D$_{{\rm{even}}}$.
The reason for this is that the same structures are studied from
different viewpoints, as we will outline in this paper.

We start with a fairly general setting which admits both
constructions, $\a$-induction as well as Ocneanu's double
triangle algebras and chiral generators. Namely,
we consider a type III subfactor
$N\subset M$ of finite index with a finite system
of $N$-$N$ morphisms which includes the irreducible
constituents of the dual canonical endomorphism.
(A ``system of morphisms'' means essentially
that, as sectors, the morphisms form a closed
algebra under the sector ``fusion'' product,
see Definition \ref{system} below.)
Therefore the subfactor is in particular forced to
have finite depth. The inclusion structure associates
to the $N$-$N$ system automatically $N$-$M$, $M$-$N$
and $M$-$M$ systems. The typical situation is that
the system of $M$-$M$ morphisms is the ``unknown part''
of the theory. As an easy reformulation of Ocneanu's idea from
his work on Goodman-de la Harpe-Jones subfactors
associated with Dynkin diagrams one can define the
double triangle algebra for such a setting, and it
provides a powerful tool to gain information about
the ``unknown part'' from the ``known part'' of the theory.
Namely, the double triangle algebra is a direct sum of
intertwiner spaces equipped with two
different product structures, and its center $\cZ_h$
with respect to the ``horizontal product'' turns out
to be isomorphic to the (in general non-commutative)
fusion rule algebra of the $M$-$M$ system when endowed
with the ``vertical product''. This kind of duality is
the subfactor analogue to the group algebra with its
pointwise and convolution products.

Under the assumption that the $N$-$N$ system is braided
there is automatically the notion of $\a$-induction, which
extends $N$-$N$ to (possibly reducible) $M$-$M$ morphisms.
(This notion does not even depend on the finite depth condition.)
The braiding provides powerful tools to analyze the structure
of the center $\cZ_h$ at the same time, and the analysis is
most conveniently carried out with a graphical intertwiner
calculus which will be explained in detail in this paper.
Besides the standard ``braiding fusion symmetries'' for
wire diagrams representing intertwiners of the braided
$N$-$N$ morphisms, we show that the theory of $\a$-induction
gives rise naturally to an extended symmetry
which we call ``intertwining braiding fusion relations''.
This reduces all graphical manipulations representing
the relations between intertwiners to easily visible purely
topological moves, and it allows us to work
without the ``sliding moves along walls'' involving
``quantum $6j$-symbols for subfactors'' which
are the main technical tool in \cite{O7}.
With a braiding on the $N$-$N$ system we can define
chiral generators $p_\la^\pm$ in the center $\cZ_h$,
and our notion essentially coincides with Ocneanu's
definition of elements $p_j^\pm$ given graphically
in his A-D-E setup. We show that
the decomposition of the $p_\la^\pm$'s into minimal central
projections in $\cZ_h$ corresponds exactly to the sector
decomposition of the $\a$-induced sectors $[\a_\la^\pm]$,
and therefore they can be naturally identified.

As shown by Rehren \cite{R0}, a system of braided endomorphisms
gives rise to S- and T-matrices which provide a unitary
representation of the modular group $\SLZ$ whenever the
braiding is non-degenerate.
(Relations between modular S- and T-matrices and
braiding data are also discussed in \cite{MS3,FG,FRS2}.)
In terms of $\a$-induction we define a matrix $Z$ with entries
$Z_{\la,\mu}=\langle \a_\la^+,\a_\mu^- \rangle$
for $N$-$N$ morphisms $\la,\mu$, where the brackets denote
the dimension of the intertwiner space $\Hom(\a_\la^+,\a_\mu^-)$.
As it corresponds to the ``vacuum'' in physical applications,
we use the label ``0'' for the identity morphism $\id_N$, and
hence our matrix $Z$ satisfies the conditions in \erf{massmatZ},
where now $Z_{0,0}=1$ is just the factor property of $M$.
We show that $Z$ commutes with $S$ and $T$ and therefore
$Z$ is a ``modular invariant mass matrix'' in the sense of
conformal field theory if the braiding is non-degenerate.
In fact, the non-degenerate case is the most interesting one,
as in the $\SUn_k$ examples in conformal field theory.
We apply an argument of Ocneanu to our situation
to show that in that case, due to the identification
with chiral generators, both kinds of
$\a$-induction together generate the whole $M$-$M$ fusion rule
algebra. Moreover, the essential information about its
representation theory (or equivalently, about the decomposition
of the center $\cZ_h$ with the vertical product into simple
matrix algebras) is then encoded in the mass matrix $Z$:
We show that the irreducible representations of the
$M$-$M$ fusion rule algebra are labelled by pairs
$\la,\mu$ with $Z_{\la,\mu}\neq 0$, and that their dimensions
are given exactly by the number $Z_{\la,\mu}$. Consequently,
the $M$-$M$ fusion rules are then commutative if and only if
all $Z_{\la,\mu}\in\{0,1\}$. An analogous result has been
claimed by Ocneanu for his A-D-E setting related to the modular
invariant mass matrices of the $\SUz$ WZW models of \cite{CIZ,Kt}.
He has his own geometric construction of modular invariants sketched
in the lectures but not included in the lecture notes \cite{O7}.
Our construction is different and based on the results of
\cite{BE3}, and it shows that the structural results do not
depend on the very special properties of Dynkin diagrams and
hold in a far more general context.
We also analyze the representation of the $M$-$M$ fusion rule
algebra arising from its left action on $M$-$N$ sectors.
As corollaries of our analysis we find that the number
of $N$-$M$ (or $M$-$N$) morphisms is given by the trace
$\tr(Z)$, whereas the number of $M$-$M$ morphisms is
given by $\tr(Z \,{}^{{\rm t}} \! Z)$.

In a forthcoming publication we will further analyze and apply
our construction to subfactors constructed by means of the
level $k$ positive energy representations of the
$\SUn$ loop group theory. For these examples, the braiding
is always non-degenerate and, moreover, the S- and T-matrices
are the modular matrices performing the character
transformations of the corresponding $\SUn_k$ WZW theory.
Therefore the construction of
braided subfactors\footnote{We remark that our short-hand
notion of a ``braided subfactor'' meaning a subfactor for
which Assumptions \ref{set-fin} and \ref{set-braid} below
hold is different from the notion used in \cite{L2.5}.}
for these models yields non-diagonal modular
invariants $Z$. E.g.\ for $\SUz_k$ one can construct
the subfactors in terms of local loop groups which recover
the A-D-E modular invariants of \cite{CIZ,Kt}.
In our setting also the ``type II'' or
``non-blockdiagonal'' invariants can be treated
by dropping the chiral locality condition.
(The chiral locality condition, expressing local
commutativity of the extended chiral theory in the
formulation of nets of subfactors \cite{LR}, implies
``$\alpha\sigma$-reciprocity'' \cite{BE1} which in turn
forces the modular invariant to be of type I. Detailed
explanation and non-local examples will be provided
in \cite{BEK2}.) Thus this paper extends the known
results on conformal inclusions \cite{X1,X2,BE2,BE3} and
simple current extensions \cite{BE2,BE3} of $\SUn_k$,
and it generally relates (the classification of) modular
invariants to (non-degenerately) braided subfactors.
Furthermore our results prove two conjectures by
two of us \cite[Conj.\ 7.1 \& 7.2]{BE3}.

This paper is organized as follows.
In Sect.\ \ref{sec-morsecbraid} we review some basic facts about
morphisms, intertwiners, sectors and braidings, and we
reformulate Rehren's result about S- and T-matrices arising from
superselection sectors in our context of braided factors.
In Sect.\ \ref{sec-graphcal} we establish the graphical
methods for the intertwiner calculus we use in this paper.
The abstract mathematical structure underlying the
basic graphical calculus (Subsect.\ \ref{sec-basgraphcal})
is ``strict monoidal $C^*$-categories'' \cite{DR}.
Graphical methods for calculations involving fusion and
braiding have been used in various publications, see
e.g.\ \cite{MS1,KR,Wi,FK,FG,Kf,J2}. However,
for our purposes it turns out to be extremely important to
handle normalization factors with special care, and to the
best of our knowledge, a comprehensive exposition which applies
to our framework has not been published somewhere. So we
work out a ``rotation covariant'' intertwiner calculus here,
based on a formulation of Frobenius reciprocity by Izumi \cite{I2}.
We then define $\a$-induction for braided subfactors
and use it to extend our graphical calculus conveniently.
In Sect.\ \ref{sec-dta} we present the double triangle algebra
and analyze its properties.
In Sect.\ \ref{sec-aicpmodinv} we present our version of
Ocneanu's graphical notion of chiral generators,
and we show that it can be naturally identified with the
$\a$-induced sectors. We then define the ``mass matrix'' $Z$
and show that it commutes with the S- and T-matrices of the
$N$-$N$ system. Assuming now that the braiding is non-degenerate,
we show that the $M$-$M$ fusion rule algebra is generated by
the images of the two kinds ($+$ and $-$) of $\a$-induction.
In Sect.\ \ref{sec-repMM} we decompose $\cZ_h$ with the vertical
product into simple matrix algebras which is equivalent to the
determination of all the irreducible representations of the
$M$-$M$ fusion rule algebra, and we show that their dimensions
are given by the entries of the modular invariant mass matrix.
Then we analyze the representation arising from the left
action on $M$-$N$ sectors.
In Sect.\ \ref{sec-concl} we finally conclude this paper with
general remarks and comments and an outlook to the applications
to subfactors arising from conformal field theory which will
be treated in \cite{BEK2}.

\section{Preliminaries}
\label{sec-morsecbraid}

\subsection{Morphisms and sectors}
\label{sec-morsec}

For our purposes it turns out to be convenient to make use of
the formulation of sectors between different factors. We follow
here (up to minor notational changes) Izumi's presentation
\cite{I2,I3} based on Longo's sector theory \cite{L2}.
Let $A$, $B$ be infinite factors. We denote by
$\Mor(A,B)$ the set of unital $\ast$-homomorphisms from
$A$ to $B$. We also  denote $\End(A)=\Mor(A,A)$, the
set of unital $\ast$-endomorphisms.
For $\rho\in\Mor(A,B)$ we define the statistical
dimension $d_\rho=[B:\rho(A)]^{1/2}$, where $[B:\rho(A)]$ is the
minimal index \cite{J,Ko}. A morphism $\rho\in\Mor(A,B)$ is
called irreducible if the subfactor $\rho(A)\subset B$ is
irreducible, i.e.\ if the relative commutant
$\rho(A)'\cap B$ consists only of scalar multiples of
the identity in $B$. Two morphisms $\rho,\rho'\in\Mor(A,B)$
are called equivalent if there exists a unitary $u\in B$
such that $\rho'(a)=u\rho(a)u^*$ for all $a\in A$.
We denote by $\Sect(A,B)$ the quotient of $\Mor(A,B)$ by
unitary equivalence, and we call its elements $B$-$A$ sectors.
Similar to the case $A=B$, $\Sect(A,B)$ has the natural
operations, sums and products:  For $\rho_1,\rho_2\in\Mor(A,B)$
choose generators $t_1,t_2\in B$ of a Cuntz algebra $\cO_2$,
i.e.\ such that $t_i^*t_j=\del ij \bfe$ and
$t_1t_1^*+t_2t_2^*=\bfe$. Define $\rho\in\Mor(A,B)$ by
putting $\rho(a)=t_1\rho_1(a)t_1^*+t_2\rho_2(a)t_2^*$
for all $a\in A$, and then the sum of sectors is defined
as $[\rho_1]\oplus[\rho_2]=[\rho]$. The product of sectors
comes from the composition of endomorphisms,
$[\rho_1][\rho_2]=[\rho_1\circ\rho_2]$.
We often omit the composition symbol ``$\circ$'', so
$[\rho_1][\rho_2]=[\rho_1\rho_2]$.
The statistical dimension is an invariant for sectors
(i.e.\ equivalent morphisms have equal dimension) and
is additive and multiplicative with respect to these
operations. Moreover, for $[\rho]\in\Sect(A,B)$ there
is a unique conjugate sector $\overline{[\rho]}\in\Sect(B,A)$
such that, if $[\rho]$ is irreducible,
$\overline{[\rho]}$ is irreducible as well and
$\overline{[\rho]}\times[\rho]$ contains the
identity sector $[\id_A]$ and $[\rho]\times\overline{[\rho]}$
contains $[\id_B]$ precisely once. We choose a representative
endomorphism of $\overline{[\rho]}$ and denote it
naturally by $\co\rho$, thus $[\co\rho]=\overline{[\rho]}$.
For conjugates we have $d_{\co\rho}=d_\rho$. As for
bimodules one may decorate $B$-$A$ sectors $[\rho]$
with suffixes, ${}_B[\rho]_A$, and then we can
multiply ${}_B[\rho]_A \times {}_A[\sigma]_B$ but
not, for instance, ${}_B[\rho]_A$ with itself.
For $\rho,\tau\in\Mor(A,B)$ we denote
\[ \Hom (\rho,\tau) = \{ t\in B :
t\, \rho(a) = \tau (a) \, t \,,\,\,\, a\in A \} \]
and
\[ \langle \rho, \tau \rangle = \dim\,\Hom (\rho,\tau) \,.\]
If $[\rho]=[\rho_1]\oplus[\rho_2]$ then
\[ \langle \rho, \tau \rangle = \langle \rho_1, \tau \rangle
+ \langle \rho_2, \tau \rangle \,. \]
Note that if $\rho$ is irreducible then for $t,t'\in\Hom (\rho,\tau)$
it follows that $t^*t'$ is a scalar and then putting
\be   t^*t'=\langle t,t' \rangle \bfe_B
\label{inprint}
\ee
defines an inner product on $\Hom (\rho,\tau)$. One often calls
$\Hom (\rho,\tau)$ a ``Hilbert space of isometries'' in this case.

If $\rho\in\Mor(A,B)$ with $d_\rho<\infty$ then
$\co\rho\in\Mor(B,A)$ is a conjugate if there are
isometries $r_\rho\in\Hom(\id_A,\co\rho\rho)$ and
${\co r}_\rho\in\Hom(\id_B,\rho\co\rho)$ such that
\[ \rho(r_\rho)^* {\co r}_\rho = d_\rho^{-1} \bfe_B \qquad
\mbox{and} \qquad \co\rho ({\co r}_\rho)^* r_\rho
= d_\rho^{-1}\bfe_A \,, \]
and in the case that $\rho$ is irreducible such isometries
$r_\rho$ and ${\co r}_\rho$ are unique
up to a common phase. If $C$ is another factor and
$\sigma\in\Mor(C,A)$ and $\tau\in\Mor(C,B)$ are
morphisms with finite statistical
dimensions $d_\sigma,d_\tau<\infty$, and conjugate
morphisms $\co\sigma\in\Mor(A,C)$, $\co\tau\in\Mor(B,C)$,
respectively, then the
``left and right Frobenius reciprocity maps'',
\[ \bearll
\cL_\rho: \Hom(\tau,\rho\sigma) \longrightarrow
\Hom(\sigma,\co\rho\tau)\,,\qquad & t \longmapsto
\sqrt{\displaystyle\frac{d_\rho d_\sigma}{d_\tau}} \,
\co\rho(t)^* r_\rho \,,\\[.8em]
\cR_\rho: \Hom(\co\sigma,\co\tau\rho) \longrightarrow
\Hom(\co\tau,\co\sigma\co\rho)\,,& s \longmapsto
\sqrt{\displaystyle\frac{d_\rho d_\tau}{d_\sigma}} \,
s^*\co\tau({\co r}_\rho)\,,
\eear \]
are anti-linear (vector space) isomorphisms with inverses
\[ \bearll
\cL_\rho^{-1}: \Hom(\sigma,\co\rho\tau) \longrightarrow
\Hom(\tau,\rho\sigma)\,,\qquad & x \longmapsto
\sqrt{\displaystyle\frac{d_\rho d_\tau}{d_\sigma}} \,
\rho(x)^* {\co r}_\rho \,,\\[.8em]
\cR_\rho^{-1}: \Hom(\co\tau,\co\sigma\co\rho) \longrightarrow
\Hom(\co\sigma,\co\tau\rho)\,,& y \longmapsto
\sqrt{\displaystyle\frac{d_\rho d_\sigma}{d_\tau}} \,
y^*\co\sigma(r_\rho)\,,
\eear \]
respectively \cite{I2}.
(See also \cite[Sect.\ 5]{FG} and \cite[App.\ A]{FRS2} 
for such formulae arising from superselection sectors.)
Hence we have in particular
Frobenius reciprocity \cite{I2,L3},
\[ \langle \tau,\rho\sigma \rangle =
\langle \co\rho\tau, \sigma \rangle =
\langle \co\rho, \sigma\co\tau \rangle =
\langle \co\sigma\co\rho, \co\tau \rangle =
\langle \co\sigma, \co\tau\rho \rangle =
\langle \tau\co\sigma,\rho \rangle \,. \]
If $\tau$ and $\sigma$ are irreducible then the
Frobenius reciprocity maps are even (anti-linearly) isometric:
With the inner products as in \erf{inprint}
on the above intertwiner spaces we have
$\langle t,t' \rangle = \langle \cL_\rho(t'), \cL_\rho(t) \rangle$
for $t,t'\in\Hom(\tau,\rho\sigma)$ and similarly
$\langle s,s' \rangle = \langle \cR_\rho(s'), \cR_\rho(s) \rangle$
for $s,s'\in\Hom(\co\sigma,\co\tau\rho)$.

The map $\phi_\rho:B\rightarrow A$ defined by
\[ \phi_\rho(b) = r_\rho^* \, \co\rho (b) \,
r_\rho \,, \qquad b\in B \]
is completely positive, normal, unital $\phi_\rho(\bfe_B)=\bfe_A$
and satisfies
\[ \phi_\rho (\rho(a_1)b\rho(a_2)) = a_1 \phi_\rho(b) a_2 \,,
\qquad a_1,a_2\in A \,, \quad b\in B\,.\]
The map is called the (unique) standard left inverse.
The minimal conditional expectation for the subfactor
$\rho(A)\subset B$ is given by $E_\rho=\rho\circ\phi_\rho$.
Let now $\rho,\sigma,\tau$ as above be irreducible with standard left
inverses $\phi_\rho,\phi_\sigma,\phi_\tau$, respectively, and
let $t\in\Hom(\tau,\rho\sigma)$ be non-zero. Then
$\phi_\rho (tt^*)\in\Hom(\sigma,\sigma)$ is a positive scalar
and $\tilde{E}_\tau:B\rightarrow \tau(C)$ given by
$\rho\circ\phi_\rho (tt^*)\tilde{E}_\tau(b)
=\tau\circ\phi_\sigma\circ\phi_\rho (tbt^*)$
for all $b\in B$ is a conditional expectation
for the subfactor $\tau(C)\subset B$.
Since conditional expectations for irreducible subfactors
are unique we conclude that
\[ \phi_\tau(b) \, E_\rho (tt^*) =
\phi_\sigma\circ\phi_\rho  (tbt^*)\,,\qquad b\in B \]
holds for any $t\in\Hom(\tau,\rho\sigma)$.
Moreover, $t^*t'$ is a scalar for any $t,t'\in\Hom(\tau,\rho\sigma)$,
$t^*t'=\langle t,t' \rangle \bfe_B$,
and so is $\cL_\rho(t)^*\cL_\rho(t')$, in fact
\[ \langle t,t' \rangle \bfe_A
= \langle \cL_\rho(t'), \cL_\rho(t) \rangle \bfe_A \equiv
\cL_\rho(t')^*\cL_\rho(t) = \frac{d_\rho d_\sigma}{d_\tau}
\, r_\rho^* \co\rho(t't^*) r_\rho \]
and this is
\be
\phi_\rho (t't^*) = \frac{d_\tau}{d_\rho d_\sigma} \,
\langle t,t' \rangle \bfe_A \,.
\label{phiddd}
\ee

Now let $N\subset M$ be an infinite subfactor of finite
index. Let $\can\in\End(M)$ be a canonical endomorphism from
$M$ into $N$ and $\canr=\can|_N\in\End(N)$. By $\iota\in\Mor(N,M)$
we denote the injection map, $\iota(n)=n\in M$, $n\in N$.
Then $d_\iota=[M:N]^{1/2}$, and
a conjugate $\iotab\in\Mor(M,N)$ is given by
$\iotab(m)=\can(m)\in N$, $m\in M$. (These formulae
could in fact be used to define the canonical and
dual canonical endomorphism.) Note that
$\can=\iota\iotab$ and $\canr=\iotab\iota$, and
there are isometries $w\equiv r_\iota\in\Hom(\id_N,\canr)$
and $v\equiv {\co r}_\iota\in\Hom(\id_M,\can)$ such that
$w^*v  = \can(v ^*)w = [M:N]^{-1/2}\bfe$. Moreover, we have
the pointwise equality $M=Nv$, and for each $m\in M$ the
decomposition $m=nv$ yields a unique element $n\in N$.
Explicitly, $n=[M:N]^{1/2}w^*\can(m)$.

Now let us consider a single factor $A$ and its sectors.
For a set of irreducible sectors which is closed under
conjugation and irreducible decomposition of products
(a ``sector basis'' in the notation of \cite{BE1,BE2,BE3}
in the case that the set is finite) it is often useful
to choose one representative endomorphism for each sector.

\begin{definition}
\label{system}
{\rm We call a subset $\sys\subset\End(A)$ a
{\sl system of endomorphisms} if it
satisfies the following properties.
\begin{enumerate}
\item Each $\la\in \sys$ is irreducible and has finite statistical
dimension.
\item Different elements in $\sys$ are inequivalent, i.e.\
different as sectors.
\item $\id_A\in \sys$.
\item For any $\la\in\sys$, we have a morphism
$\bar\la\in\sys$ such that $[\bar\la]$ is the conjugate sector
of $[\la]$.
\item $\sys$ is closed under composition and subsequent
irreducible decomposition, i.e.\ for any $\la,\mu\in\sys$
we have non-negative integers $N_{\la,\mu}^\nu$ with
$[\la][\mu]=\sum_{\nu\in\sys} N_{\la,\mu}^\nu [\nu]$ as sectors.
\end{enumerate}
}\end{definition}
Note that we do not assume finiteness of $\sys$ in this definition.
The numbers $N_{\la\mu}^\nu=\langle\la\mu,\nu\rangle$ are called
fusion coefficients. Frobenius reciprocity now reads
$N_{\la,\mu}^\nu=N_{\co\la,\nu}^\mu=N_{\nu,\co\mu}^\la$, and
associativity of the sector product yields
$\sum_{\mu\in\sys} N_{\la,\mu}^\nu N_{\rho,\sigma}^\mu
=\sum_{\tau\in\sys} N_{\la,\rho}^\tau N_{\tau,\sigma}^\nu$.
The additivity and multiplicativity of the statistical
dimension with respect to sector sums and products implies
$\sum_{\nu\in\sys} N_{\la,\mu}^\nu d_\nu
= d_\la d_\mu$, $\la,\mu,\nu\in\sys$.
Defining matrices $N_\mu$ with entries
$(N_\mu)_{\la,\nu}=N_{\la,\mu}^\nu$ gives
$N_{\co\mu}$ as the transpose of $N_\mu$ and defines the
``regular representation'' of the sector products,
$N_\la N_\mu = \sum_{\nu\in\sys} N_{\la,\mu}^\nu N_\nu$,
and the statistical dimension can be regarded as a
one-dimensional representation or as a simultaneous
eigenvector of all matrices $N_\mu$ with eigenvalues
$d_\mu$ ($\la,\mu,\nu\in\sys$).

\subsection{Braided endomorphisms}
\label{braid-prelim}

Let $A$ again be an infinite factor and $\sys$ a system of
endomorphisms of $A$. In general the sector products are not
commutative. If the sectors commute, then a
``systematic choice of unitary intertwiners''
in each space $\Hom(\la\mu,\mu\la)$, $\la,\mu\in\sys$, is
called a braiding (which need not exist in general).
To be more precise, we give the following

\begin{definition}
\label{braid}{\rm
We say that a system $\sys$ of endomorphisms is {\sl braided}
if for any pair $\la,\mu\in\sys$ there is a unitary operator
$\eps\lambda\mu\in\Hom(\la\mu,\mu\la)$ subject to
initial conditions
\be
\eps {\id_A}\mu=\eps\lambda{\id_A}=\bfe \,,
\label{ini}
\ee
and whenever $t\in \Hom (\lambda,\mu\nu)$ we have
the braiding fusion equations (BFE's)
\be
\begin{array}{rl}
\rho(t) \, \eps \lambda \rho 
&=\,\, \eps \mu \rho \, \mu (\eps \nu \rho) \, t \,, \\[.4em]
t \, \eps \rho \lambda
&=\,\, \mu (\eps \rho \nu) \, \eps \rho \mu \, \rho (t) \,, \\[.4em]
\rho(t)^* \, \eps \mu\rho \, \mu( \eps \nu\rho )
&=\,\, \eps \la\rho \, t^* \,, \\[.4em]
t^* \, \mu( \eps \rho\nu ) \, \eps \rho\mu
&= \,\, \eps \rho\la \, \rho(t)^* \,,
\eear
\label{BFE}
\ee
for any $\la,\mu,\nu\in\sys$.
}\end{definition}
The unitaries $\eps\la\mu$ are called {\sl braiding operators}
(or {\sl statistics operators}).
Note that a braiding $\varepsilon\equiv\varepsilon^+$ always
comes along with another ``opposite'' braiding $\varepsilon^-$,
namely operators $\epsm \la\mu = (\epsp \mu\la )^*$,
$\epsp \mu\la \equiv \eps \mu\la$,
satisfy the same relations. The unitaries $\epsp \la\mu$
and $\epsm \la\mu$ are different in general but may coincide
for some $\lambda$, $\mu$. Later we will also use the following
notion of non-degeneracy of a braiding (cf.\ \cite{R0}). 

\begin{definition}
\label{non-deg}{\rm
We say that a braiding $\e$ on a system of endomorphisms
$\sys$ is {\sl non-degenerate}, if the following condition is
satisfied: If some morphism $\la\in\sys$
satisfies $\e^+(\la,\mu)=\e^-(\la,\mu)$
for all morphisms $\mu\in\sys$,
then we have $\la=\id_A$.
}\end{definition}

We may also extend a given braiding from $\sys$ in a well
defined manner to all equivalent and sum endomorphisms as follows.
We denote by $\Sigma(\sys)$ the set of all endomorphisms
$\la,\rho\in\End(A)$ given as
$\la(a)=\sum_{i=1}^n t_i\la_i(a)t_i^*$ and
$\rho(a)=\sum_{j=1}^m s_j\rho_j(a)s_j^*$ for all $a\in A$,
where $t_i\in A$, $i=1,2\dots,n$, and $s_j\in A$, $j=1,2,\dots,m$,
are Cuntz algebra generators, i.e.\ $t_i^*t_k=\del ik\bfe$ and
$\sum_{i=1}^n t_it_i^*=\bfe$, and similarly
$s_j^*s_l=\del jl\bfe$ and $\sum_{j=1}^m s_js_j^*=\bfe$, and 
$\la_i,\rho_j\in\sys$. (Here $n,m\ge 1$.) For $\la,\rho$
as above we put
\be
\eps \la\rho = \sum_{i=1}^n \sum_{j=1}^m s_j \rho_j (t_i) \,
\eps {\la_i}{\rho_j} \, \la_i(s_j^*) t_i^* \,,
\label{esum}
\ee
and one can check that this definition is independent
of the ambiguities in the choice
of isometries $t_i\in\Hom(\la_i,\la)$ and
$s_j\in\Hom(\rho_j,\rho)$.
Note that in the case $n=m=1$ this reads
\be
\eps {\Ad(u)\circ\la}{\Ad(q)\circ\rho} = q \rho (u) \,
\eps \la\rho \, \la(q^*) u^* 
\label{eun}
\ee
with some unitaries $u,q\in A$. Then for any
sum endomorphisms $\la,\mu,\rho\in\Sigma(\sys)$ the
BFE's (\ref{BFE}) hold as well or,
alternatively, we have the naturality equations
\be
\rho(t) \, \eps \la \rho = \eps \mu \rho \, t \,,
\qquad t \, \eps \rho \la = \eps \rho \mu \, \rho (t)
\label{nat}
\ee
whenever $t\in \Hom(\la,\mu)$. Using decompositions of products
$\la\mu$, $\la,\mu\in\Sigma(\sys)$ one can then easily show by
use of the BFE's that
\be
\eps {\lambda\mu}\rho = \eps \lambda\rho \,
\lambda (\eps \mu\rho) \,, \qquad
\eps \lambda{\mu\rho} = \mu (\eps \lambda\rho) \,
\eps \lambda\mu \,.
\label{comp}
\ee
By plugging this in \erf{nat} we find that BFE's hold for
endomorphisms in $\Sigma(\sys)$ as well and \erf{nat} yields
for $\eps \la\mu\in\Hom(\la\mu,\mu\la)$ the braid relation
(or ``Yang-Baxter equation'')
\be
\rho(\eps\lambda\mu) \, \eps\lambda\rho \, \lambda(\eps\mu\rho) =
\eps\mu\rho \, \mu(\eps\lambda\rho) \, \eps\lambda\mu \,.
\label{YBE}
\ee

Now let $\sys$ be a braided system of endomorphisms and let
$\rho,\co\rho\in\sys$ be conjugate morphisms.
Denote by $r\equiv r_\rho\in\Hom(\id_A,\bar\rho\rho)$ and
$\co r\equiv {\co r}_\rho \in\Hom(\id_A,\rho\bar\rho)$
isometries such that
\[ \rho (r)^* \co r = \co\rho (\co r)^* r = d_\rho^{-1} \bfe \,,\]
which are then unique up to a common phase.\footnote{If $\rho$
is not self-conjugate then we may choose $r_{\co\rho}={\co r}_\rho$
and ${\co r}_{\co\rho}=r_\rho$. However, if $\rho$ is
self-conjugate, $\rho=\co\rho$, we do not have
$r_\rho=\co r_{\rho}$ in general. This is only true for
so-called ``real'' sectors, and for ``pseudo-real''
sectors we have $r_\rho=-\co r_{\rho}$.}
Note that
$\eps {\co\rho}\rho ^* \co r\in\Hom(\id_A,\bar\rho\rho)$ is an
isometry and hence $\eps {\co\rho}\rho ^* \co r = \om_\rho r$
for some phase $\om_\rho\in\bbT$ which is called the
{\sl statistics phase} and is obviously independent of the common phase
of $r$ and $\co r$. In fact $\om_\rho$ is even independent of
the choice of $\rho$ and $\co\rho$ within their sectors: If
$\rho'=\Ad\; u \circ \rho$ and $\co\rho'=\Ad\; \co u \circ \co\rho$
for some unitaries $u,\co u\in A$,
then it is easy to see that isometries
$r'=\co u \co\rho(u)r\in\Hom(\id_A,\bar\rho '\rho')$ and
$\co r'=u\rho(\co u)\co r\in\Hom(\id_A,\rho'\bar\rho')$ also fulfill
$\rho (r')^* \co r' = \co\rho (\co r')^* r' = d_\rho^{-1} \bfe$.
Now the braiding operator transforms as
$\eps {\co\rho'}{\rho'} =
u\rho(\co u)\eps{\co\rho}\rho \co\rho(u)^*\co u ^*$ and hence
\[ \eps {\co\rho'}{\rho'} ^* \co r'
= \co u\co\rho(u)\eps{\co\rho}\rho ^*\co r = \om_\rho r' \,.\]
The statistics phase can also be obtained by
\[ \phi_\rho(\eps\rho\rho) = r^* \co\rho (\eps\rho\rho)r
= \om_\rho d_\rho^{-1} \bfe \,.\]
(The number $\om_\rho d_\rho^{-1}$ is usually called the
{\sl statistics parameter}.)
This is obtained from the initial condition and the BFE:
\[ \rho(r) = \rho(r) \eps {\id_A}\rho = \eps {\co\rho}\rho
\co\rho (\eps\rho\rho) r \,, \]
but since $r^*\eps {\co\rho}\rho ^*=\om_\rho\co r ^*$ we obtain
\[ r^* \co\rho (\eps\rho\rho)r = r^*  \eps {\co\rho}\rho ^* \rho(r)
= \om_\rho \co r ^* \rho (r)
= \om_\rho d_\rho^{-1} \bfe \,.\]
Moreover we have $\om_\rho=\om_{\co\rho}$. This can be seen as
follows. We have
\[ r=r\eps \rho{\id_A} = \co\rho(\eps\rho\rho)\eps\rho{\co\rho}
\rho(r) \,,\]
hence $r^* \co\rho(\eps\rho\rho)=\rho(r)^*\eps\rho{\co\rho} ^*$, thus
\[ \om_\rho d_\rho^{-1} \bfe = \rho(r)^*\eps\rho{\co\rho} ^*r
= \om_{\co\rho}\rho(r)^* \co r = \om_{\co\rho} d_\rho^{-1} \,, \]
since $\eps \rho{\co\rho} ^* r = \om_{\co\rho} \co r$ by
definition. Therefore we have
$\om_\rho r^*=\co r^* \eps \rho{\co\rho} ^*$. Another
application of the BFE yields
$\eps\rho\rho \rho(\co r) = \rho(\eps\rho{\co\rho})^* \co r$,
hence we have
\[ \rho(\co r)^* \eps\rho\rho \rho (\co r) =
\rho (\co r)^* \rho (\eps \rho{\co \rho})^* \co r =
\om_\rho \rho (r)^* \co r = \om_\rho d_\rho^{-1} \bfe \,.\]

Now let $\la,\mu,\nu\in\sys$. Let
$r\equiv r_\la\in\Hom(\id_A,\co\la\la)$
and $\co r\equiv {\co r}_\la\in\Hom(\id_A,\la\co\la)$
be isometries such that
$\la(r)^*\co r=\co\la(\co r)^*r = d_\la^{-1}\bfe$. Let
$t,t'\in\Hom(\la,\mu\nu)$. Recall that
$\phi_\mu (t't^*)=d_\la d_\mu^{-1} d_\nu^{-1} t^*t'\in\Hom(\la,\la)$
is a scalar. We can now compute
\[ \bearll
\om_\lambda d_\mu^{-1} d_\nu^{-1}\, t^*t'
&= \om_\lambda d_\lambda^{-1}  \, \phi_\nu\circ\phi_\mu (t't^*)
 = \phi_\nu\circ\phi_\mu (t' \la(\co r)^* \eps \la\la
   \la (\co r) t^*) \\[.4em]
&= \co r ^* \, \phi_\nu\circ\phi_\mu (t' \eps \la\la t^*) \, \co r 
 = \co r ^* \, \phi_\nu\circ\phi_\mu ( \eps\la{\mu\nu} \la(t')t^*)\,
   \co r \\[.4em]
&=  \co r ^* \, \phi_\nu\circ\phi_\mu ( \eps\la{\mu\nu} t^*) \, t' \co r
 = \co r ^* t^* \, \phi_\nu\circ\phi_\mu ( \eps{\mu\nu}{\mu\nu}) \,
   t' \co r \\[.4em]
&= \co r ^* t^* \, \phi_\nu\circ\phi_\mu ( \mu(\eps\mu\nu)\mu^2
   (\eps\nu\nu)    \eps\mu\mu \mu(\eps\nu\mu))\, t' \co r \\[.4em]
&= \om_\mu d_\mu^{-1} \, \co r ^* t^* \, \phi_\nu ( \eps\mu\nu
   \mu (\eps\nu\nu) \eps\nu\mu ) \, t' \co r \\[.4em]
&= \om_\mu d_\mu^{-1} \, \co r ^* t^* \, \phi_\nu ( \nu (\eps\nu\mu
   \eps\nu\nu \nu (\eps\mu\nu )\,  t' \co r \\[.4em]
&= \om_\mu \om_\nu d_\mu^{-1}  d_\nu^{-1}\,  \co r ^* t^* \,
   \eps\nu\mu \eps\mu\nu \, t' \co r \\[.4em]
&= \om_\mu \om_\nu d_\mu^{-1}  d_\nu^{-1}\, t^* 
   \eps\nu\mu \eps\mu\nu \, t' \,,
\eear \]
where we finally could omit the $\co r$'s since
$t^* \eps\nu\mu \eps\mu\nu t'\in\Hom(\la,\la)$ is a scalar.
As $\eps\nu\mu \eps\mu\nu t'\in\Hom(\la,\mu\nu)$ we find
$\om_\la \langle t,t'\rangle = \om_\mu \om_\nu \langle t,
\eps\nu\mu \eps\mu\nu t' \rangle$ for any
$t,t'\in\Hom(\la,\mu\nu)$, and therefore we arrive at the
important relation
\be
\eps\nu\mu \eps\mu\nu \, t
= \frac{\om_\la}{\om_\mu \om_\nu} \, t  \qquad
\mbox{for all} \qquad t\in\Hom(\la,\mu\nu)\,.
\label{mondiag}
\ee

Decomposing $[\mu\nu]$ in all irreducible sectors $[\la]$
and choosing for each
$\lambda\in\sys$ some orthonormal bases of intertwiners
$t_{\la;i}\in\Hom(\la,\mu\nu)$, $i=1,2,\dots,N_{\mu,\nu}^\la$,
where $N_{\mu,\nu}^\la = \langle \la,\mu\nu\rangle$ as usual,
we have
$\sum_{\la\in\sys}\sum_i t_{\la;i}t_{\la;i}^*=\bfe$,
and therefore we find by Eqs.\ (\ref{phiddd}) and (\ref{mondiag}),
\[ \phi_\mu (\eps\nu\mu \eps\mu\nu)^* = \phi_\mu \left(
\eps\nu\mu \eps\mu\nu \sum_{\la\in\sys}
\sum_i t_{\la;i}t_{\la;i}^* \right)^*
= \sum_{\la\in\sys} \frac{\om_\mu \om_\nu}{\om_\la}
N_{\mu,\nu}^\la \frac{d_\la}{d_\mu d_\nu} \bfe \,. \]
One then defines a matrix $Y$ in terms of these numbers
\cite{R0} (see also \cite{FG,FRS2}):
\be
Y_{\mu,\nu} = \sum_{\la\in\sys}
\frac{\om_\mu \om_\nu}{\om_\la}
N_{\mu,\nu}^\la \, d_\la \,, \qquad \mu,\nu\in\sys \,,
\label{Ydef}
\ee
i.e.
$d_\mu d_\nu \phi_\mu (\eps\nu\mu \eps\mu\nu)^* = Y_{\mu,\nu}\bfe$.
Then one has
\[Y_{\la,\mu}=Y_{\mu,\la}=Y_{\co\la,\mu}^*=Y_{\co\la,\co\mu}\,.\]
The first equality is obvious from \erf{Ydef}, so we only need to show 
$Y_{\la,\mu}=(Y_{\co\la,\mu})^*$. In fact, applying the BFE again
yields $\co\la(\eps\la\mu) r_\la=\eps{\co\la}\mu ^* \mu(r_\la)$
and $r_\la^* \co\la (\eps\mu\la) = \mu(r_\la)^* \eps\mu{\co\la} ^*$.
Hence
\[ \bearll
Y_{\la,\mu} \bfe &= \phi_\mu (Y_{\la,\mu}) = d_\la d_\mu
(r_\mu^* \co\mu (r_\la^* \co\la (\eps\mu\la \eps\la\mu)
r_\la)r_\mu)^* \\[.4em]
&= d_\la d_\mu (r_\la^* r_\mu^* \co\mu(\eps \mu{\co\la} ^*
\eps {\co\la}\mu ^*) r_\mu r_\la)^* = (r_\la ^* Y_{\co\la,\mu} r_\la)^*
= (Y_{\co\la,\mu})^* \bfe \,.
\eear \]
Moreover, we have
\[ Y_{\nu,\rho}Y_{\mu,\rho}
=d_\rho \sum_\la N_{\mu,\nu}^\la Y_{\rho,\la}\,, \]
since
\[ \bearll
Y_{\nu,\rho}Y_{\mu,\rho} \bfe &= d_\rho^2 d_\mu d_\nu \, \phi_\nu
(\eps\rho\nu \phi_\mu(\eps\rho\mu\eps\mu\rho)\eps\nu\rho )^* \\[.4em]
&= d_\rho^2 d_\mu d_\nu \, \phi_\nu\circ\phi_\mu (\mu(\eps\rho\nu)
\eps\rho\mu \eps\mu\rho \mu(\eps\nu\rho) )^* \\[.4em]
&= d_\rho^2 d_\mu d_\nu \sum_\la \sum_i \phi_\nu\circ\phi_\mu
(\eps\rho{\mu\nu} \rho(t_{\la;i}t_{\la;i}^*)
\eps{\mu\nu}\rho )^* \\[.4em]
&=  d_\rho^2 d_\mu d_\nu \sum_\la \sum_i \phi_\nu\circ\phi_\mu
(t_{\la;i} \eps\rho\la \eps\la\rho t_{\la;i}^*)^* \\[.4em]
&= d_\rho^2 d_\mu d_\nu \sum_\la \sum_i
\phi_\mu (t_{\la;i}t_{\la;i}^*)^*
\phi_\la (\eps\rho\la \eps\la\rho)^* =
d_\rho \sum_\la N_{\mu,\nu}^\la Y_{\rho,\la} \bfe \,.
\eear \]

{}From now on we assume that the system $\sys$ is finite.
We define the complex number
\[ z_\sys= \sum_{\la\in\sys} d_\la^2 \om_\la \,, \]
and if $z_\sys \neq 0$ we put
$c =  4 \arg (z_\sys)/ \pi$.
Note that the $c$ is here only defined mod $8$
and we may make a choice. 
Let $C$ be the conjugation matrix with entries
$C_{\la,\mu}=\del \la{\co\mu}$. Clearly, $C=C^*=C^{-1}$. 
We then have the following

\begin{proposition}
\label{ST0}
Let $\sys$ be finite system of endomorphisms with
$z_\sys\neq 0$. Then S- and T-matrices defined by
\[ S_{\la,\mu} =  |z_\sys|^{-1} \, Y_{\la,\mu} \,, \qquad
T_{\la,\mu} = \E^{-\pi\I c/12} \, \om_\la \, \del \la\mu \,,
\qquad \la,\mu\in\sys \,, \]
obey the partial Verlinde modular algebra
$TSTST=S$, $CTC=T$, $CSC=S$ and $T^*T=\bfe$.
\end{proposition}

To prove the proposition, we simply compute
\[ \bearll
\sum_\mu \om_\la Y_{\la,\mu} \om_\mu Y_{\mu,\nu} \om_\nu
&= \om_\la \om_\nu \sum_\mu \om_\mu Y_{\la,\co\mu}^*
Y_{\nu,\co\mu}^* 
= \om_\la \om_\nu \sum_{\mu,\sigma} \om_\mu d_\mu
  N_{\la,\nu}^\sigma  Y_{\co\mu,\sigma}^* \\[.4em]
&= \om_\la \om_\nu \sum_{\mu,\rho,\sigma} \om_\mu d_\mu
  N_{\la,\nu}^\sigma  N_{\co\mu,\sigma}^\rho
\frac{\om_\rho}{\om_\mu \om_\sigma} d_\rho
= \om_\la \om_\nu \sum_{\rho,\sigma} d_\rho^2 d_\sigma
  N_{\la,\nu}^\sigma \frac{\om_\rho}{\om_\sigma} \\[.4em]
&= Y_{\la,\nu} \sum_\rho d_\rho^2 \om_\rho
= Y_{\la,\nu} z_\sys \,,
\eear \]
hence
$TSTST = \E^{-\pi\I c/4}|z_\sys|^{-1} S z_\sys=S$.
The remaining relations $CTC=T$, $CSC=S$ and $T^*T=\bfe$
are obvious.

We define {\sl weight vectors} $y^\la$ with components
$y^\la_\mu=Y_{\la,\mu}$ and {\sl statistics characters}
$\chi_\la:\sys\rightarrow\bbC$ with evaluations
$\chi_\la(\mu)=d_\la^{-1}Y_{\la,\mu}$, $\la,\mu\in\sys$.
We have seen that the weight vectors $y^\la$ are simultaneous
eigenvectors of the fusion matrices $N_\mu$ with eigenvalues
$\chi_\la(\mu)$, $N_\mu y^\la = \chi_\la(\mu) y^\la$.
Hence we obtain by computing inner products,
\[ \chi_\mu(\rho) \langle y^\la,y^\mu \rangle =
\langle y^\la , N_\rho y^\mu \rangle =
\langle N_{\co\rho} y^\la,y^\mu \rangle =
\chi_\la(\co\rho)^* \langle y^\la,y^\mu \rangle =
\chi_\la(\rho) \langle y^\la,y^\mu \rangle \,.\]
Therefore the eigenvectors are either orthogonal,
$\langle y^\la,y^\mu \rangle =0$, or parallel,
$d_\mu y^\la = d_\la y^\mu$ since then the characters
are equal, $\chi_\la=\chi_\mu$. It is obvious that if
some $\la\in\sys$ is degenerate, i.e.\ has trivial
monodromy with all other $\mu\in\sys$, then $y^\la$ is parallel
to the vector $y^0$. (Here and later we use the
label ``$0$'' for the identity $\id_A\in\sys$.)
Note that we have $y^0_\mu=d_\mu$, and then
$Y_{\la,\mu}=d_\la d_\mu$. Conversely, if $y^\la$ is
parallel to $y^0$ we have seen that then necessarily
$Y_{\la,\mu}=d_\la d_\mu$, hence
\[ Y_{\la,\mu} = \sum_{\rho\in\sys}
\frac{\om_\la \om_\mu}{\om_\rho}
N_{\la,\mu}^\rho \, d_\rho = d_\la d_\mu =
\sum_{\rho\in\sys} N_{\la,\mu}^\rho \, d_\rho \,,
\qquad \mu\in\sys \,, \]
and this is clearly only possible if all the eigenvalues
$\om_\la \om_\mu \om^{-1}_\rho$ of the
monodromy are trivial, i.e.\ if $\la$ is degenerate.
We conclude that a braiding on $\sys$ is non-degenerate
if and only if $\langle y^\la , y^0 \rangle = \del \la 0 w$,
where $w=\sum_{\la\in\sys} d_\la^2$ is the {\sl global index}.
We now arrive at Rehren's result \cite{R0}.

\begin{theorem}
\label{ST}
The following conditions are equivalent
for a finite braided system of endomorphisms $\sys$:
\begin{enumerate}
\item The braiding on $\sys$ is non-degenerate.
\item We have $w=|z_\sys|^2$ and the matrices $S$ and $T$
obey the full Verlinde modular algebra
\[ S^*S=T^*T=\bfe \,,\quad (ST)^3=S^2 = C \,, \quad CTC=T \,, \]
moreover $S$ diagonalizes the fusion rules (Verlinde formula):
\[ N_{\la,\mu}^\nu = \sum_{\rho\in\sys}
\frac{S_{\la,\rho}S_{\mu,\rho}
S_{\nu,\rho}^*}{S_{0,\rho}} \,. \]
\end{enumerate}
\end{theorem}

Note that the implication ${\it 2.}\Rightarrow{\it 1.}$
is trivial since invertibility of $S$ implies that there
is no vector $y^\la$ parallel $y^0$. So let us assume that
the braiding is non-degenerate:
$\langle y^\la , y^0 \rangle = \del \la 0 w$ for all $\la\in\sys$.
Then we can first check
\[ \bearll
w &= \sum_\mu \langle y^0,y^\mu \rangle
d_\mu \om_\mu^{-1} = \sum_{\mu,\nu} d_\nu
 Y_{\mu,\nu} d_\mu \om_\mu^{-1} = \sum_{\mu,\nu,\la}
d_\nu \frac{\om_\mu\om_\nu}{\om_\la} N_{\mu,\nu}^\la
d_\la d_\mu \om_\mu^{-1} \\[.4em]
&= \sum_{\mu,\nu,\la}
d_\la d_\nu \frac{\om_\nu}{\om_\la} N_{\co\nu,\la}^\mu
d_\mu = \sum_{\la,\nu} d_\la^2 \om_\la^{-1}
d_\nu^2 \om_\nu \,,
\eear \]
thus
$w=\left|\sum_{\la\in\sys} d_\la^2 \om_\la\right|^2\equiv |z_\sys|^2$.
Next we compute
\[ \langle y^\la,y^\mu \rangle = \sum_\rho
Y_{\la,\rho}^* Y_{\mu,\rho} = \sum_{\rho,\nu}
N_{\co\la,\mu}^\nu Y_{\rho,\nu} d_\rho =
\sum_\nu N_{\co\la,\mu}^\nu \langle y^0,y^\nu \rangle
= N_{\co\la,\mu}^0 w = \del \la\mu w \,, \]
hence $S^*S=\bfe$. Similarly we observe that
$\sum_\rho Y_{\la,\rho} Y_{\mu,\rho}=
\sum_\rho Y_{\co\la,\rho}^* Y_{\mu,\rho}=
\del {\co\la}\mu w$, giving $S^2=C$ which obviously commutes
with $T$.
Finally we check
\[ \sum_\rho \frac{S_{\la,\rho}S_{\mu,\rho}
S_{\nu,\rho}^*}{S_{0,\rho}} = w^{-1} \sum_\rho
\frac{Y_{\la,\rho}Y_{\mu,\rho} Y_{\nu,\rho}^*}{d_\rho}
= w^{-1} \sum_{\rho,\sigma} N_{\la,\mu}^\sigma Y_{\rho,\sigma}
Y_{\nu,\rho}^* = \sum_\sigma N_{\la,\mu}^\sigma \del \nu\sigma
= N_{\la,\mu}^\nu \,,\]
proving the Verlinde identity.

\begin{corollary}
\label{sltzaction}
If the braiding on $\sys$ is non-degenerate, then the
matrix $S$ and the diagonal matrix $T$ are the images
$S=U(\cS)$ and $T=U(\cT)$ of canonical generators
\[ \cS = \left(\begin{array}{cc} 0 & -1 \\ 1 & 0
\end{array}\right), \qquad
\cT = \left(\begin{array}{cc} 1 & 1 \\ 0 & 1
\end{array}\right), \]
in a unitary representation $U$ of the modular
group\footnote{In the literature the name ``modular
group'' is often reserved for $\PSLZ=\SLZ/\bbZ_2$
rather than $\SLZ$. Clearly, we obtain a
representation of $\PSLZ$ whenever the charge
conjugation is trivial, $C=\bfe$.}
$\SLZ$ with dimension $|\sys|$,
the cardinality of $\sys$.
\end{corollary}

\section{Graphical Intertwiner Calculus}
\label{sec-graphcal}

\subsection{Basic graphical intertwiner calculus}
\label{sec-basgraphcal}

\begin{figure}[htb]
\begin{center}
\unitlength 0.6mm
\begin{picture}(60,30)
\thinlines
%\path(5,10)(5,20)(55,20)(55,10)(5,10)
\put(5,10){\dashbox{2}(50,10){$x$}}
\put(10,30){\vector(0,-1){10}}
\put(25,30){\vector(0,-1){10}}
\put(50,30){\vector(0,-1){10}}
\put(10,10){\vector(0,-1){10}}
\put(25,10){\vector(0,-1){10}}
\put(50,10){\vector(0,-1){10}}
\put(3,25){\makebox(0,0){$\la_1$}}
\put(18,25){\makebox(0,0){$\la_2$}}
\put(57,25){\makebox(0,0){$\la_n$}}
\put(3,5){\makebox(0,0){$\mu_1$}}
\put(18,5){\makebox(0,0){$\mu_2$}}
\put(57,5){\makebox(0,0){$\mu_m$}}
%\put(30,15){\makebox(0,0){$x$}}
\put(37.5,25){\makebox(0,0){$\cdots$}}
\put(37.5,5){\makebox(0,0){$\cdots$}}
\end{picture}
\end{center}
\caption{An intertwiner $x$}
\label{boxx}
\end{figure}
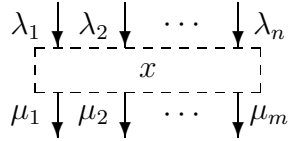

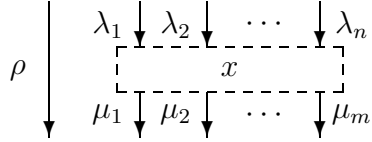
\begin{figure}[htb]
\begin{center}
\unitlength 0.6mm
\begin{picture}(80,30)
\thinlines
%\path(25,10)(25,20)(75,20)(75,10)(25,10)
\put(25,10){\dashbox{2}(50,10){$x$}}
\put(10,30){\vector(0,-1){30}}
\put(30,30){\vector(0,-1){10}}
\put(45,30){\vector(0,-1){10}}
\put(70,30){\vector(0,-1){10}}
\put(30,10){\vector(0,-1){10}}
\put(45,10){\vector(0,-1){10}}
\put(70,10){\vector(0,-1){10}}
\put(3,15){\makebox(0,0){$\rho$}}
\put(23,25){\makebox(0,0){$\la_1$}}
\put(38,25){\makebox(0,0){$\la_2$}}
\put(77,25){\makebox(0,0){$\la_n$}}
\put(23,5){\makebox(0,0){$\mu_1$}}
\put(38,5){\makebox(0,0){$\mu_2$}}
\put(77,5){\makebox(0,0){$\mu_m$}}
%\put(50,15){\makebox(0,0){$x$}}
\put(57.5,25){\makebox(0,0){$\cdots$}}
\put(57.5,5){\makebox(0,0){$\cdots$}}
\end{picture}
\end{center}
\caption{The intertwiner $\rho(x)$}
\label{rhoboxx}
\end{figure}

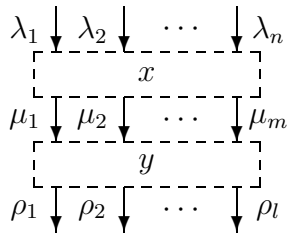
\begin{figure}[htb]
\begin{center}
\unitlength 0.6mm
\begin{picture}(80,50)
\thinlines
\put(5,30){\dashbox{2}(50,10){$x$}}
\put(5,10){\dashbox{2}(50,10){$y$}}
\put(10,50){\vector(0,-1){10}}
\put(25,50){\vector(0,-1){10}}
\put(50,50){\vector(0,-1){10}}
\put(10,30){\vector(0,-1){10}}
\put(25,30){\vector(0,-1){10}}
\put(50,30){\vector(0,-1){10}}
\put(10,10){\vector(0,-1){10}}
\put(25,10){\vector(0,-1){10}}
\put(50,10){\vector(0,-1){10}}
\put(3,45){\makebox(0,0){$\la_1$}}
\put(18,45){\makebox(0,0){$\la_2$}}
\put(57,45){\makebox(0,0){$\la_n$}}
\put(3,25){\makebox(0,0){$\mu_1$}}
\put(18,25){\makebox(0,0){$\mu_2$}}
\put(57,25){\makebox(0,0){$\mu_m$}}
\put(3,5){\makebox(0,0){$\rho_1$}}
\put(18,5){\makebox(0,0){$\rho_2$}}
\put(57,5){\makebox(0,0){$\rho_l$}}
\put(37.5,45){\makebox(0,0){$\cdots$}}
\put(37.5,25){\makebox(0,0){$\cdots$}}
\put(37.5,5){\makebox(0,0){$\cdots$}}
\end{picture}
\end{center}
\caption{Product $yx$ of diagrammatically
         composable intertwiners $x$ and $y$}
\label{x-top-y}
\end{figure}

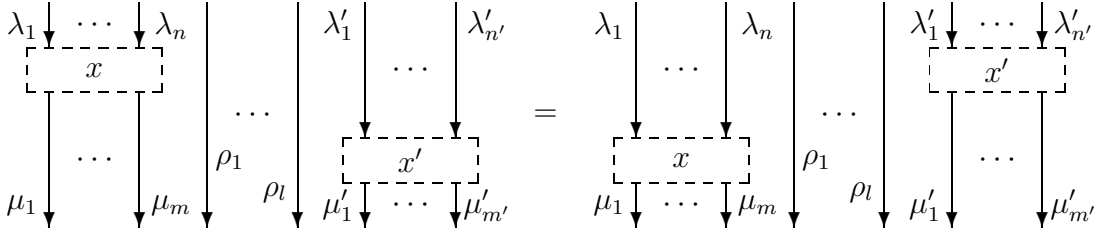
\begin{figure}[htb]
\begin{center}
\unitlength 0.6mm
\begin{picture}(240,50)
\thinlines
\put(5,30){\dashbox{2}(30,10){$x$}}
\put(75,10){\dashbox{2}(30,10){$x'$}}
\put(10,50){\vector(0,-1){10}}
\put(30,50){\vector(0,-1){10}}
\put(10,30){\vector(0,-1){30}}
\put(30,30){\vector(0,-1){30}}
\put(45,50){\vector(0,-1){50}}
\put(65,50){\vector(0,-1){50}}
\put(80,50){\vector(0,-1){30}}
\put(100,50){\vector(0,-1){30}}
\put(80,10){\vector(0,-1){10}}
\put(100,10){\vector(0,-1){10}}
\put(20,45){\makebox(0,0){$\cdots$}}
\put(20,15){\makebox(0,0){$\cdots$}}
\put(55,25){\makebox(0,0){$\cdots$}}
\put(90,35){\makebox(0,0){$\cdots$}}
\put(90,5){\makebox(0,0){$\cdots$}}
\put(4,45){\makebox(0,0){$\la_1$}}
\put(37,45){\makebox(0,0){$\la_n$}}
\put(4,5){\makebox(0,0){$\mu_1$}}
\put(37,5){\makebox(0,0){$\mu_m$}}
\put(50,15){\makebox(0,0){$\rho_1$}}
\put(60,8){\makebox(0,0){$\rho_l$}}
\put(74,45){\makebox(0,0){$\la_1'$}}
\put(107,45){\makebox(0,0){$\la_{n'}'$}}
\put(74,5){\makebox(0,0){$\mu_1'$}}
\put(107,5){\makebox(0,0){$\mu_{m'}'$}}
\put(120,25){\makebox(0,0){$=$}}
\put(135,10){\dashbox{2}(30,10){$x$}}
\put(205,30){\dashbox{2}(30,10){$x'$}}
\put(140,50){\vector(0,-1){30}}
\put(160,50){\vector(0,-1){30}}
\put(140,10){\vector(0,-1){10}}
\put(160,10){\vector(0,-1){10}}
\put(175,50){\vector(0,-1){50}}
\put(195,50){\vector(0,-1){50}}
\put(210,50){\vector(0,-1){10}}
\put(230,50){\vector(0,-1){10}}
\put(210,30){\vector(0,-1){30}}
\put(230,30){\vector(0,-1){30}}
\put(150,35){\makebox(0,0){$\cdots$}}
\put(150,5){\makebox(0,0){$\cdots$}}
\put(185,25){\makebox(0,0){$\cdots$}}
\put(220,45){\makebox(0,0){$\cdots$}}
\put(220,15){\makebox(0,0){$\cdots$}}
\put(134,45){\makebox(0,0){$\la_1$}}
\put(167,45){\makebox(0,0){$\la_n$}}
\put(134,5){\makebox(0,0){$\mu_1$}}
\put(167,5){\makebox(0,0){$\mu_m$}}
\put(180,15){\makebox(0,0){$\rho_1$}}
\put(190,8){\makebox(0,0){$\rho_l$}}
\put(204,45){\makebox(0,0){$\la_1'$}}
\put(237,45){\makebox(0,0){$\la_{n'}'$}}
\put(204,5){\makebox(0,0){$\mu_1'$}}
\put(237,5){\makebox(0,0){$\mu_{m'}'$}}
\end{picture}
\end{center}
\caption{Vertical translation intertwiners $x$ and $x'$}
\label{xx'}
\end{figure}

We now introduce our conventions to represent and
manipulate intertwiners graphically. We consider a braided
system of endomorphisms $\sys\subset\End(A)$ with $A$ a
type III factor. Essentially we represent intertwiners
by ``wire diagrams'' where the (oriented) wires represent
endomorphisms $\la\in\sys$. This works as follows.
For an intertwiner
$x\in\Hom(\la_1\la_2\cdots\la_n,\mu_1\mu_2\cdots\mu_m)$
we draw a (dashed) box with $n$ (downward) incoming wires
labelled by $\la_1,\dots,\la_n$ and $m$ (downward) outgoing
wires $\mu_1,\dots,\mu_m$ as in Fig.\ \ref{boxx},
$\la_i,\mu_j\in\sys$.
%
% intertwiner as box
Therefore the diagrammatic representation of $x$ does not
only specify it as an operator, it even specifies the intertwiner
space it is considered  to belong to. (Note that the same operator
can belong to different intertwiner spaces as e.g.\ the identity
operator belongs to any $\Hom(\la,\la)$ with $\la$ varying.)
If a morphism $\rho\in\sys$ is applied to $x$, then
$\rho(x)\in\Hom(\rho\la_1\la_2\cdots\la_n,\rho\mu_1\mu_2\cdots\mu_m)$
is represented graphically by adding a straight wire on the left
as in Fig.\ \ref{rhoboxx}.
%
% rho applied to an intertwiner as box
Reflecting the fact that $x$ can also
be considered as an intertwiner in
$\Hom(\la_1\la_2\cdots\la_n\rho,\mu_1\mu_2\cdots\mu_m\rho)$
we can always add (or remove) a straight wire on the right
as in Fig.\ \ref{boxxrho}
%
% an intertwiner as box with added right wire rho
\begin{figure}[htb]
\begin{center}
\unitlength 0.6mm
\begin{picture}(80,30)
\thinlines
%\path(5,10)(5,20)(55,20)(55,10)(5,10)
\put(5,10){\dashbox{2}(50,10){$x$}}
\put(10,30){\vector(0,-1){10}}
\put(25,30){\vector(0,-1){10}}
\put(50,30){\vector(0,-1){10}}
\put(10,10){\vector(0,-1){10}}
\put(25,10){\vector(0,-1){10}}
\put(50,10){\vector(0,-1){10}}
\put(70,30){\vector(0,-1){30}}
\put(3,25){\makebox(0,0){$\la_1$}}
\put(18,25){\makebox(0,0){$\la_2$}}
\put(57,25){\makebox(0,0){$\la_n$}}
\put(3,5){\makebox(0,0){$\mu_1$}}
\put(18,5){\makebox(0,0){$\mu_2$}}
\put(57,5){\makebox(0,0){$\mu_m$}}
\put(77,15){\makebox(0,0){$\rho$}}
%\put(30,15){\makebox(0,0){$x$}}
\put(37.5,25){\makebox(0,0){$\cdots$}}
\put(37.5,5){\makebox(0,0){$\cdots$}}
\end{picture}
\end{center}
\caption{The intertwiner $x$}
\label{boxxrho}
\end{figure}
without changing the intertwiner as an operator.
We say that intertwiners
$x\in\Hom(\la_1\la_2\cdots\la_n,\mu_1\mu_2\cdots\mu_m)$
and
$y\in\Hom(\nu_1\nu_2\cdots\nu_k,\rho_1\rho_2\cdots\rho_l)$,
$\rho_j\in\sys$,
are {\sl diagrammatically composable} if $m=k$ and
$\mu_i=\nu_i$ for all $i=1,2,\dots,m$. Then the composed
intertwiner
$yx\in\Hom(\la_1\la_2\cdots\la_n,\rho_1\rho_2\cdots\rho_l)$
is represented graphically by putting the wire diagram
for $x$ on top of that for $y$ as in Fig.\ \ref{x-top-y}.
We also call this graphical procedure composition of
wire diagrams.
%
% intertwiner composition x on top of y
Sometimes diagrammatic composability may be achieved by
adding or removing straight wires on the right.
Now let also
$x'\in\Hom(\la_1'\la_2'\cdots\la_{n'}',\mu_1'\mu_2'\cdots\mu_{m'}')$
with $\la_i',\mu_j'\in\sys$. The intertwining property of $x$
yields the identity
$\mu_1\mu_2\cdots\mu_m \rho_1\rho_2\cdots\rho_l(x')x
= x\la_1\la_2\cdots\la_n \rho_1\rho_2\cdots\rho_l(x')$,
and this is diagrammatically given in Fig.\ \ref{xx'}.
%
% vertical translation of x and x'
Thus we have some freedom in translating intertwiner boxes vertically
without actually changing the represented intertwiner.

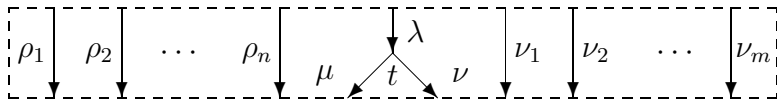
\begin{figure}[htb]
\begin{center}
\unitlength 0.6mm
\begin{picture}(170,20)
\thinlines
\put(0,0){\dashbox{2}(170,20){}}
\put(10,20){\vector(0,-1){20}}
\put(25,20){\vector(0,-1){20}}
\put(60,20){\vector(0,-1){20}}
\put(85,20){\vector(0,-1){10}}
\put(85,10){\vector(-1,-1){10}}
\put(85,10){\vector(1,-1){10}}
\put(110,20){\vector(0,-1){20}}
\put(125,20){\vector(0,-1){20}}
\put(160,20){\vector(0,-1){20}}
\put(5,10){\makebox(0,0){$\rho_1$}}
\put(20,10){\makebox(0,0){$\rho_2$}}
\put(55,10){\makebox(0,0){$\rho_n$}}
\put(90,15){\makebox(0,0){$\la$}}
\put(70,5){\makebox(0,0){$\mu$}}
\put(100,5){\makebox(0,0){$\nu$}}
\put(85,5){\makebox(0,0){$t$}}
\put(115,10){\makebox(0,0){$\nu_1$}}
\put(130,10){\makebox(0,0){$\nu_2$}}
\put(165,10){\makebox(0,0){$\nu_m$}}
\put(37.5,10){\makebox(0,0){$\cdots$}}
\put(147.5,10){\makebox(0,0){$\cdots$}}
\end{picture}
\end{center}
\caption{$\protect\sqrt[4]{\frac{d_\mu d_\nu}{d_\la}}
\rho_1\rho_2\cdots\rho_n(t)$
where $t\in\Hom(\la,\mu\nu)$ is an isometry}
\label{rrrt}
\end{figure}

\begin{figure}[htb]
\begin{center}
\unitlength 0.6mm
\begin{picture}(170,20)
\thinlines
\put(0,0){\dashbox{2}(170,20){}}
\put(10,20){\vector(0,-1){20}}
\put(25,20){\vector(0,-1){20}}
\put(60,20){\vector(0,-1){20}}
\put(75,20){\vector(1,-1){10}}
\put(95,20){\vector(-1,-1){10}}
\put(85,10){\vector(0,-1){10}}
\put(110,20){\vector(0,-1){20}}
\put(125,20){\vector(0,-1){20}}
\put(160,20){\vector(0,-1){20}}
\put(5,10){\makebox(0,0){$\rho_1$}}
\put(20,10){\makebox(0,0){$\rho_2$}}
\put(55,10){\makebox(0,0){$\rho_n$}}
\put(90,5){\makebox(0,0){$\la$}}
\put(70,15){\makebox(0,0){$\mu$}}
\put(100,15){\makebox(0,0){$\nu$}}
\put(87,16){\makebox(0,0){$t^*$}}
\put(115,10){\makebox(0,0){$\nu_1$}}
\put(130,10){\makebox(0,0){$\nu_2$}}
\put(165,10){\makebox(0,0){$\nu_m$}}
\put(37.5,10){\makebox(0,0){$\cdots$}}
\put(147.5,10){\makebox(0,0){$\cdots$}}
\end{picture}
\end{center}
\caption{$\protect\sqrt[4]{\frac{d_\mu d_\nu}{d_\la}}
\rho_1\rho_2\cdots\rho_n(t)^*$
where $t\in\Hom(\la,\mu\nu)$ is an isometry}
\label{rrrt*}
\end{figure}
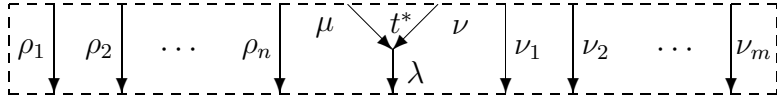

The intertwiners we consider are (sums over)
compositions of {\sl elementary intertwiners}
arising from the unitary braiding operators
$\eps\la\mu\in\Hom(\la\mu,\mu\la)$
and isometries $t\in\Hom(\la,\mu\nu)$. The wire diagrams
and boxes we are dealing with are therefore compositions
of ``elementary boxes'' representing the elementary intertwiners.
We now have to introduce some normalization convention. First,
the identity intertwiner $\bfe\equiv\bfe_A$ is naturally
given by the ``trivial box'' with only straight
wires of arbitrary labels. The next elementary
intertwiner is $\rho_1\rho_2\cdots\rho_n(\eps\la\mu)$
for which we draw a box as in Fig.\ \ref{rrrepslm}
%
% \rho_1\rho_2\cdots\rho_n(\eps\la\mu)
\begin{figure}[htb]
\begin{center}
\unitlength 0.6mm
\begin{picture}(170,20)
\thinlines
\put(0,0){\dashbox{2}(170,20){}}
\put(10,20){\vector(0,-1){20}}
\put(25,20){\vector(0,-1){20}}
\put(60,20){\vector(0,-1){20}}
\put(75,20){\vector(1,-1){20}}
\put(95,20){\line(-1,-1){8}}
\put(83,8){\vector(-1,-1){8}}
\put(110,20){\vector(0,-1){20}}
\put(125,20){\vector(0,-1){20}}
\put(160,20){\vector(0,-1){20}}
\put(5,10){\makebox(0,0){$\rho_1$}}
\put(20,10){\makebox(0,0){$\rho_2$}}
\put(55,10){\makebox(0,0){$\rho_n$}}
\put(100,5){\makebox(0,0){$\la$}}
\put(70,5){\makebox(0,0){$\mu$}}
\put(115,10){\makebox(0,0){$\nu_1$}}
\put(130,10){\makebox(0,0){$\nu_2$}}
\put(165,10){\makebox(0,0){$\nu_m$}}
\put(37.5,10){\makebox(0,0){$\cdots$}}
\put(147.5,10){\makebox(0,0){$\cdots$}}
\end{picture}
\end{center}
\caption{$\rho_1\rho_2\cdots\rho_n(\eps\la\mu)$}
\label{rrrepslm}
\end{figure}
where the arbitrary labels $\nu_1,\dots,\nu_m$ are irrelevant
and may be omitted. Similarly, the box of Fig.\ \ref{rrrt}
%
% \rho_1\rho_2\cdots\rho_n(t)
represents the elementary intertwiner 
$d_\mu^{1/4} d_\nu^{1/4} d_\la^{-1/4}\rho_1\rho_2\cdots\rho_n(t)$,
where $t\in\Hom(\la,\mu\nu)$ is an isometry. We label the
trivalent vertex in the box by $t$ since $\Hom(\la,\mu\nu)$ may
be more than one-dimensional and so we have to specify the
intertwiner.
(Note that there would still be an ambiguity of a phase for
the choice of an isometry even if $\Hom(\la,\mu\nu)$ is only
one-dimensional.)
Finally, the elementary intertwiners
$\eps\la\mu ^*=\epsm\mu\la$ and 
$d_\mu^{1/4} d_\nu^{1/4} d_\la^{-1/4}\rho_1\rho_2\cdots\rho_n(t)^*$
are represented by Figs.\ \ref{rrrepslm*}
%
% \rho_1\rho_2\cdots\rho_n(\eps\la\mu)^*
\begin{figure}[htb]
\begin{center}
\unitlength 0.6mm
\begin{picture}(170,20)
\thinlines
\put(0,0){\dashbox{2}(170,20){}}
\put(10,20){\vector(0,-1){20}}
\put(25,20){\vector(0,-1){20}}
\put(60,20){\vector(0,-1){20}}
\put(75,20){\line(1,-1){8}}
\put(87,8){\vector(1,-1){8}}
\put(95,20){\vector(-1,-1){20}}
\put(110,20){\vector(0,-1){20}}
\put(125,20){\vector(0,-1){20}}
\put(160,20){\vector(0,-1){20}}
\put(5,10){\makebox(0,0){$\rho_1$}}
\put(20,10){\makebox(0,0){$\rho_2$}}
\put(55,10){\makebox(0,0){$\rho_n$}}
\put(100,5){\makebox(0,0){$\mu$}}
\put(70,5){\makebox(0,0){$\la$}}
\put(115,10){\makebox(0,0){$\nu_1$}}
\put(130,10){\makebox(0,0){$\nu_2$}}
\put(165,10){\makebox(0,0){$\nu_m$}}
\put(37.5,10){\makebox(0,0){$\cdots$}}
\put(147.5,10){\makebox(0,0){$\cdots$}}
\end{picture}
\end{center}
\caption{$\rho_1\rho_2\cdots\rho_n(\eps\la\mu)^*
         =\rho_1\rho_2\cdots\rho_n(\epsm\mu\la)$}
\label{rrrepslm*}
\end{figure}
and \ref{rrrt*}, i.e.\ they are obtained from the original
boxes in Figs.\ \ref{rrrepslm} and \ref{rrrt} by vertical
reflection and inversion of all the arrows. 
%
% \rho_1\rho_2\cdots\rho_n(t)^*
Note that $\e\equiv\e^+$ represents overcrossing
and $\e^-$ undercrossing of wires.
We will consider intertwiners which are products
of diagrammatically composable elementary intertwiners.
In terms of wire diagrams we are correspondingly dealing
with compositions of elementary boxes of Figs.\ \ref{rrrepslm},
\ref{rrrt}, \ref{rrrepslm*}, \ref{rrrt*} so that the wires
with the same labels (and orientations) can and will be glued
together in parallel and then we finally forget about the
boundaries of the (dashed) boxes. Therefore, if a wire
diagram represents some intertwiner $x$ then $x^*$ is
represented by the diagram obtained by vertical reflection
and reversing all the arrows. Note that our resulting wire
diagrams are then composed only from straight lines, over- and
undercrossings (in X-shape) and trivalent vertices (in Y-shape
or inverted Y-shape).

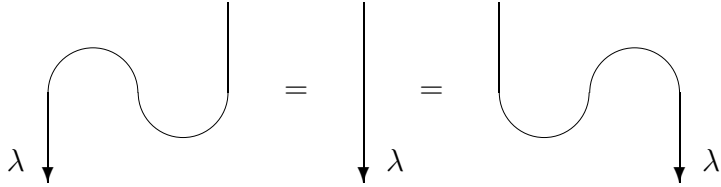
\begin{figure}[htb]
\begin{center}
\unitlength 0.6mm
\begin{picture}(170,40)
\thinlines
\put(8,5){\makebox(0,0){$\la$}}
\put(15,20){\vector(0,-1){20}}
\put(25,20){\arc{20}{3.142}{0}}
\put(45,20){\arc{20}{0}{3.142}}
\put(55,40){\line(0,-1){20}}
\put(70,20){\makebox(0,0){$=$}}
\put(85,40){\vector(0,-1){40}}
\put(92,5){\makebox(0,0){$\la$}}
\put(100,20){\makebox(0,0){$=$}}
\put(115,40){\line(0,-1){20}}
\put(125,20){\arc{20}{0}{3.142}}
\put(145,20){\arc{20}{3.142}{0}}
\put(155,20){\vector(0,-1){20}}
\put(162,5){\makebox(0,0){$\la$}}
\end{picture}
\end{center}
\caption{A topological invariance for intertwiners represented
        by wire diagrams}
\label{isoinv1}
\end{figure}

\begin{figure}[htb]
\begin{center}
\unitlength 0.6mm
\begin{picture}(80,20)
\thinlines
\put(22,19.84){\vector(1,0){0}}
\put(20,10){\circle{20}}
\put(5,10){\makebox(0,0){$\la$}}
\put(40,10){\makebox(0,0){$=$}}
\put(58,19.84){\vector(-1,0){0}}
\put(60,10){\circle{20}}
\put(75,10){\makebox(0,0){$\la$}}
\end{picture}
\end{center}
\caption{Wire diagrams for the statistical dimension $d_\la$}
\label{dla}
\end{figure}
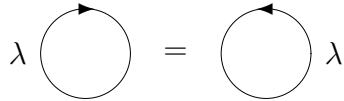

\begin{figure}[htb]
\begin{center}
\unitlength 0.6mm
\begin{picture}(130,40)
\thinlines
\put(0,20){\arc{50}{5.356}{0.927}}
\put(40,20){\arc{50}{2.214}{2.418}}
\put(40,20){\arc{50}{2.578}{3.705}}
\put(40,20){\arc{50}{3.865}{4.069}}
\put(8,5){\makebox(0,0){$\la$}}
\put(15,0){\vector(-4,-3){0}}
\put(25,0){\vector(4,-3){0}}
\put(32,5){\makebox(0,0){$\mu$}}
\put(45,20){\makebox(0,0){$=$}}
\put(60,40){\vector(0,-1){40}}
\put(70,40){\vector(0,-1){40}}
\put(53,5){\makebox(0,0){$\la$}}
\put(77,5){\makebox(0,0){$\mu$}}
\put(85,20){\makebox(0,0){$=$}}
\put(90,20){\arc{50}{5.356}{5.560}}
\put(90,20){\arc{50}{5.720}{0.564}}
\put(90,20){\arc{50}{0.724}{0.927}}
\put(130,20){\arc{50}{2.214}{4.069}}
\put(98,5){\makebox(0,0){$\la$}}
\put(105,0){\vector(-4,-3){0}}
\put(115,0){\vector(4,-3){0}}
\put(122,5){\makebox(0,0){$\mu$}}
\end{picture}
\end{center}
\caption{Unitarity of braiding operators as a vertical
Reidemeister move of type II}
\label{Reid2}
\end{figure}
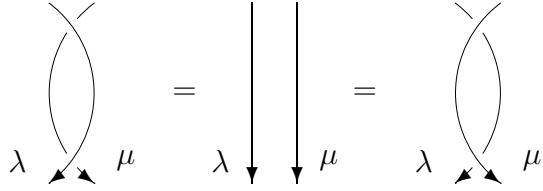

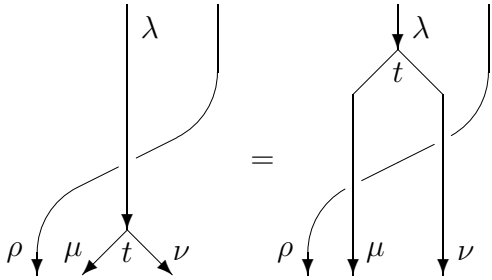
\begin{figure}[htb]
\begin{center}
\unitlength 0.6mm
\begin{picture}(120,60)
\thinlines
\put(26.180,5){\arc{32.361}{3.142}{4.249}}
\put(10,5){\vector(0,-1){5}}
\put(28,24){\line(-2,-1){9.1}}   %11.1 for 11.056
\put(33.820,45){\arc{32.361}{0}{1.107}}
\put(50,60){\line(0,-1){15}}
\put(32,26){\line(2,1){9.1}}
\put(30,60){\line(0,-1){35}}
\put(30,25){\vector(0,-1){15}}
\put(30,10){\vector(1,-1){10}}
\put(30,10){\vector(-1,-1){10}}
\put(30,5){\makebox(0,0){$t$}}
\put(35,55){\makebox(0,0){$\la$}}
\put(18,5){\makebox(0,0){$\mu$}}
\put(42,5){\makebox(0,0){$\nu$}}
\put(5,5){\makebox(0,0){$\rho$}}
\put(60,25){\makebox(0,0){$=$}}
\put(86.180,5){\arc{32.361}{3.142}{4.199}}    %4.199 for 4.249
\put(70,5){\vector(0,-1){5}}
\put(90,25){\line(-2,-1){8.1}}   %11.1 for 11.056
\put(93.820,45){\arc{32.361}{0}{1.057}} %1.057 for 1.107
\put(110,60){\line(0,-1){15}}
\put(90,25){\line(2,1){8.1}}
\put(90,60){\vector(0,-1){10}}
\put(90,50){\line(1,-1){10}}
\put(90,50){\line(-1,-1){10}}
\put(90,45){\makebox(0,0){$t$}}
\put(80,40){\line(0,-1){20}}
\put(100,40){\line(0,-1){10}}
\put(80,10){\line(0,1){10}}
\put(100,10){\line(0,1){20}}
\put(80,10){\vector(0,-1){10}}
\put(100,10){\vector(0,-1){10}}
\put(95,55){\makebox(0,0){$\la$}}
\put(85,5){\makebox(0,0){$\mu$}}
\put(105,5){\makebox(0,0){$\nu$}}
\put(65,5){\makebox(0,0){$\rho$}}
\end{picture}
\end{center}
\caption{The first braiding fusion equation}
\label{wireBFE1}
\end{figure}

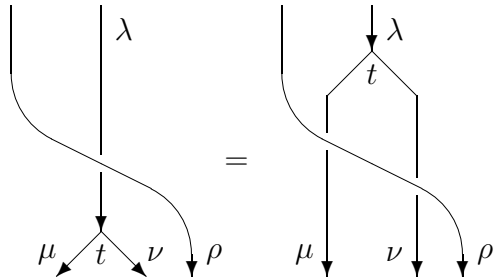
\begin{figure}[htb]
\begin{center}
\unitlength 0.6mm
\begin{picture}(120,60)
\thinlines
\put(26.180,45){\arc{32.361}{2.034}{3.142}}
\put(10,60){\line(0,-1){15}}
\put(30,25){\line(-2,1){11.1}}   %11.1 for 11.056
\put(33.820,5){\arc{32.361}{5.176}{6.283}}
\put(50,5){\vector(0,-1){5}}
\put(30,25){\line(2,-1){11.1}}
\put(30,60){\line(0,-1){33}}
\put(30,23){\vector(0,-1){13}}
\put(30,10){\vector(1,-1){10}}
\put(30,10){\vector(-1,-1){10}}
\put(30,5){\makebox(0,0){$t$}}
\put(35,55){\makebox(0,0){$\la$}}
\put(18,5){\makebox(0,0){$\mu$}}
\put(42,5){\makebox(0,0){$\nu$}}
\put(55,5){\makebox(0,0){$\rho$}}
\put(60,25){\makebox(0,0){$=$}}
\put(86.180,45){\arc{32.361}{2.034}{3.142}}
\put(70,60){\line(0,-1){15}}
\put(90,25){\line(-2,1){11.1}}   %11.1 for 11.056
\put(93.820,5){\arc{32.361}{5.176}{6.283}}
\put(110,5){\vector(0,-1){5}}
\put(90,25){\line(2,-1){11.1}}
\put(90,60){\vector(0,-1){10}}
\put(90,50){\line(1,-1){10}}
\put(90,50){\line(-1,-1){10}}
\put(90,45){\makebox(0,0){$t$}}
\put(80,40){\line(0,-1){8}}
\put(100,40){\line(0,-1){18}}
\put(80,10){\line(0,1){18}}
\put(100,10){\line(0,1){8}}
\put(80,10){\vector(0,-1){10}}
\put(100,10){\vector(0,-1){10}}
\put(95,55){\makebox(0,0){$\la$}}
\put(75,5){\makebox(0,0){$\mu$}}
\put(95,5){\makebox(0,0){$\nu$}}
\put(115,5){\makebox(0,0){$\rho$}}
\end{picture}
\end{center}
\caption{The second braiding fusion equation}
\label{wireBFE2}
\end{figure}

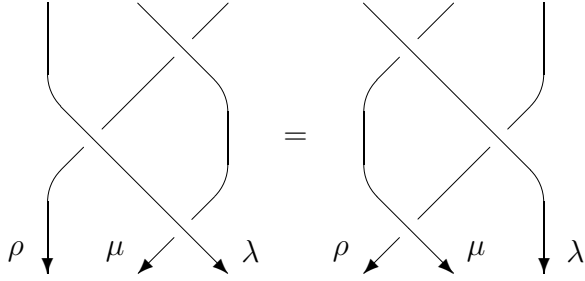
\begin{figure}[htb]
\begin{center}
\unitlength 0.6mm
\begin{picture}(140,60)
\thinlines
\put(15,60){\line(0,-1){15.858}}
\put(25,30){\vector(1,-1){30}}
\put(25,30){\line(-1,1){7.071}}
\put(35,60){\line(1,-1){17.071}}
\put(55,30){\line(0,1){5.858}}
\put(55,30){\line(0,-1){5.858}}
\put(47,12){\line(1,1){5.071}}
\put(43,8){\vector(-1,-1){8}}
\put(55,60){\line(-1,-1){8}}
\put(43,48){\line(-1,-1){16}}
\put(23,28){\line(-1,-1){5.071}}
\put(15,15.858){\vector(0,-1){15.858}}
\put(25,44.142){\arc{20}{2.356}{3.142}}
\put(45,35.858){\arc{20}{5.498}{0}}
\put(45,24.142){\arc{20}{0}{0.785}}
\put(25,15.858){\arc{20}{3.142}{3.927}}
\put(8,5){\makebox(0,0){$\rho$}}
\put(30,5){\makebox(0,0){$\mu$}}
\put(60,5){\makebox(0,0){$\la$}}
\put(70,30){\makebox(0,0){$=$}}
\put(85,60){\line(1,-1){37.071}}
\put(125,15.858){\vector(0,-1){15.858}}
\put(105,60){\line(-1,-1){8}}
\put(93,48){\line(-1,-1){5.071}}
\put(85,30){\line(0,1){5.858}}
\put(85,30){\line(0,-1){5.858}}
\put(95,10){\vector(1,-1){10}}
\put(95,10){\line(-1,1){7.071}}
\put(125,60){\line(0,-1){15.858}}
\put(117,32){\line(1,1){5.071}}
\put(113,28){\line(-1,-1){16}}
\put(93,8){\vector(-1,-1){8}}
\put(95,24.142){\arc{20}{2.356}{3.142}}
\put(115,15.858){\arc{20}{5.498}{0}}
\put(115,44.142){\arc{20}{0}{0.785}}
\put(95,35.858){\arc{20}{3.142}{3.927}}
\put(80,5){\makebox(0,0){$\rho$}}
\put(110,5){\makebox(0,0){$\mu$}}
\put(132,5){\makebox(0,0){$\la$}}
\end{picture}
\end{center}
\caption{The braid relation as a vertical Reidemeister move of type III}
\label{Reid3}
\end{figure}

So far, we have considered only wires with downward orientation.
We now introduce also the reversed orientation in terms of conjugation
as follows: Reversing the orientation of an arrow on a wire
changes its label $\la$ to $\co\la$. Also we will usually omit
drawing a wire labelled by $\id\equiv\id_A$. For each $\lambda\in\sys$
we fix (the common phase of) isometries
$r_\la\in\Hom(\id,\co\la\la)$ and
${\co r}_\la\in\Hom(\id,\la\co\la)$ such that
$\la(r_\la)^*{\co r}_\la=\co\la({\co r}_\la)^*r_\la=d_\la^{-1}\bfe$
and in turn for $\sqrt{d_\la} r_\la$ we draw one of the
equivalent diagrams in Fig.\ \ref{risom}.
%
% diagram for r
\begin{figure}[htb]
\begin{center}
\unitlength 0.6mm
\begin{picture}(145,20)
\thinlines
\put(25,20){\vector(0,-1){10}}
\put(25,10){\vector(-1,-1){10}}
\put(25,10){\vector(1,-1){10}}
\put(30,15){\makebox(0,0){$\id$}}
\put(10,5){\makebox(0,0){$\co\la$}}
\put(40,5){\makebox(0,0){$\la$}}
\put(50,10){\makebox(0,0){$=$}}
\put(75,10){\vector(-1,-1){10}}
\put(75,10){\vector(1,-1){10}}
\put(60,5){\makebox(0,0){$\co\la$}}
\put(90,5){\makebox(0,0){$\la$}}
\put(100,10){\makebox(0,0){$=$}}
\put(120,2){\arc{20}{3.142}{0}}
\put(110,0){\line(0,1){2}}
\put(130,0){\line(0,1){2}}
\put(135,5){\makebox(0,0){$\la$}}
\put(130,0){\vector(0,-1){0}}
\end{picture}
\end{center}
\caption{Wire diagrams for $\protect\sqrt{d_\la} r_\la$}
\label{risom}
\end{figure}
So the normalized isometries and their adjoints appear in wire
diagrams as ``caps'' and ``cups'', respectively.
The point is that with our normalization convention, the
relation $\la(r_\la)^*{\co r}_\la=d_\la^{-1}\bfe$ (and its
adjoint) gives a {\sl topological invariance}
for intertwiners represented by wire diagrams,
displayed in Fig.\ \ref{isoinv1}.
%
% first isotopic invariance
Note that then the wire diagrams in Fig.\ \ref{dla}
%
% statistical dimension
represent the scalar $d_\la$ (i.e.\ the intertwiner
$d_\la \bfe\in\Hom(\id,\id)$). Also note the
``vertical Reidemeister move of type II'' in Fig.\ \ref{Reid2}
%
% Reidemeister move II
is just the unitarity condition
$\eps\la\mu ^* \eps\la\mu=\bfe=\eps\mu\la\eps\mu\la ^*$.
The BFE's yield another topological invariance,
see Fig.\ \ref{wireBFE1} for the first equation
%
% BFE (1)
and Fig.\ \ref{wireBFE2} for the second equation.
%
% BFE (2)
The third and fourth equations are obtained similarly by use
of the co-isometry $t^*$; we leave it as an exercise to the reader
to draw the corresponding wire diagrams. Up to conjugation they
can also be obtained by changing over- to undercrossings in
Figs.\ \ref{wireBFE1} and \ref{wireBFE2}. Finally, the braid
relation, \erf{YBE}, represents graphically a
vertical Reidemeister move of type III, presented in Fig.\ \ref{Reid3}.
%
% Reidemeister III
The topological invariance gives us the freedom to write
down the intertwiner algebraically from a given wire diagram:
We can deform the wire diagram by finite sequences of
the above moves and then split it in
elementary wire diagrams --- in whatever way we
decompose the wire diagrams into horizontal slices
of elementary intertwiners,
we always obtain the same intertwiner due to our topological
invariance identities.

\begin{figure}[htb]
\begin{center}
\unitlength 0.6mm
\begin{picture}(120,40)
\thinlines
\put(10,20){\arc{20}{0.785}{5.498}}
\put(24.142,20){\line(-1,1){7.071}}
\put(22.142,18){\line(-1,-1){5.071}}
\put(24.142,20){\line(1,-1){7.071}}
\put(26.142,22){\line(1,1){5.071}}
\put(24.142,34.142){\arc{20}{0}{0.785}}
\put(24.142,5.858){\arc{20}{5.498}{6.283}}
\put(34.142,40){\line(0,-1){5.858}}
\put(34.142,5.858){\vector(0,-1){5.858}}
\put(41.142,5){\makebox(0,0){$\la$}}
\put(60,20){\makebox(0,0){$=$}}
\put(110,20){\arc{20}{3.927}{2.357}}
\put(95.858,20){\line(1,-1){7.071}}
\put(97.858,22){\line(1,1){5.071}}
\put(95.858,20){\line(-1,1){7.071}}
\put(93.858,18){\line(-1,-1){5.071}}
\put(95.858,34.142){\arc{20}{2.357}{3.142}}
\put(95.858,5.858){\arc{20}{3.142}{3.927}}
\put(85.858,34.142){\vector(0,1){5.858}}
\put(85.858,0){\line(0,1){5.858}}
\put(78.858,35){\makebox(0,0){$\la$}}
\end{picture}
\end{center}
\caption{Statistics phase $\om_\la$ as a ``twist''}
\label{statph}
\end{figure}
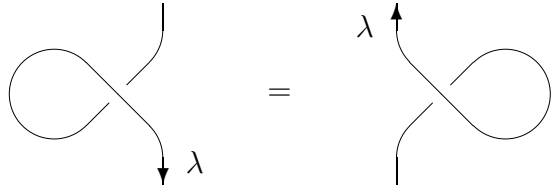

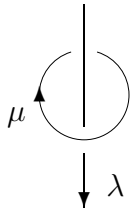
\begin{figure}[htb]
\begin{center}
\unitlength 0.6mm
\begin{picture}(35,45)
\thinlines
\put(20,25){\arc{20}{5.012}{4.412}}
\put(20,45){\line(0,-1){27}}
\put(20,12){\vector(0,-1){12}}
\put(10.12,27){\vector(0,1){0}}
\put(5,20){\makebox(0,0){$\mu$}}
\put(27,5){\makebox(0,0){$\la$}}
\end{picture}
\end{center}
\caption{Rehren's statistics character $\chi_\la$
         evaluated on $\mu$: $\chi_\la(\mu)$}
\label{statch}
\end{figure}

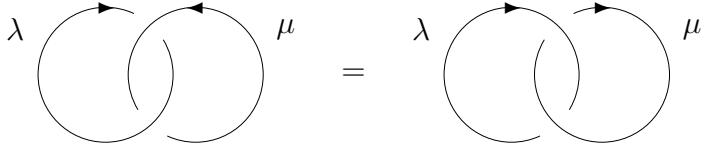
\begin{figure}[htb]
\begin{center}
\unitlength 0.6mm
\begin{picture}(160,32)
\thinlines
\put(25,15){\arc{30}{5.742}{5.142}}
\put(45,15){\arc{30}{2.601}{2.001}}
\put(27,29.8){\vector(1,0){0}}
\put(43,29.8){\vector(-1,0){0}}
\put(5,25){\makebox(0,0){$\la$}}
\put(65,25){\makebox(0,0){$\mu$}}
\put(80,15){\makebox(0,0){$=$}}
\put(115,15){\arc{30}{1.141}{0.541}}
\put(135,15){\arc{30}{4.283}{3.683}}
\put(117,29.8){\vector(1,0){0}}
\put(137,29.8){\vector(1,0){0}}
\put(95,25){\makebox(0,0){$\la$}}
\put(155,25){\makebox(0,0){$\mu$}}
\end{picture}
\end{center}
\caption{Matrix element $Y_{\la,\mu}$ of Rehren's Y-matrix
         as a ``Hopf link''}
\label{Ymatrix}
\end{figure}

Next we recall that we can write the statistics phase $\om_\la$
as the intertwiner $d_\la r_\la^* \co\la(\eps\la\la)r_\la$.
Therefore we obtain for $\om_\la$ the wire diagram on the
left-hand side of Fig.\ \ref{statph}.
%
% statistics phase
The diagram on the right-hand side expresses that $\om_\la$
can also be obtained as
$d_\la \co\la(r_\la)^*\eps{\co\la}{\co\la} \co\la (r_\la)$.
Note that we obtain the complex conjugate $\om_\la^*$
by exchanging over- and undercrossings.
Similarly, we recall that we can write a matrix element
$Y_{\la,\mu}=Y_{\mu,\la}$ of Rehren's Y-matrix
as $d_\la d_\mu \phi_\mu(\eps \la\mu \eps \mu\la)^*
= d_\la d_\mu r_\mu^* \co\mu (\epsm \la\mu \epsm \mu\la) r_\mu$.
Dividing by $d_\la$ we obtain $\chi_\la(\mu)$, the
statistics character $\chi_\la$ evaluated on $\mu$,
represented graphically by the wire diagram
in Fig.\ \ref{statch}.
%
% statistics character
We have drawn the circle $\mu$ symmetrically relative
to the straight wire $\la$ because it does not make a difference
whether we put the ``caps'' and ``cups'' for the isometry
$r_\mu$ and its conjugate $r_\mu^*$ on the left or on the
right due to the braiding fusion relations.
As it is a scalar, we can write
$Y_{\la,\mu}={\co r}_\mu^*Y_{\la,\mu}{\co r}_\mu$
and therefore its expression
$d_\la d_\mu {\co r}_\mu^*r_\la^* \co\la
(\epsm \mu\la \epsm \la\mu) r_\la {\co r}_\mu$ yields
exactly the ``Hopf link'' as the wire diagram
for the matrix element $Y_{\la,\mu}$, given by the
left-hand side of Fig.\ \ref{Ymatrix}.
%
% Y-matrix
The equality to the right-hand side is just the
relation $Y_{\la,\mu}=Y_{\la,\co\mu}^*$ together
with the prescription of representing conjugates.
Recall that if $\sys$ is finite then the Y-matrix
differs from the S-matrix
just by an overall normalization factor $\sqrt w$,
where $w$ is the global index.

Often we consider intertwiners which are sums over
intertwiners represented by the same wire diagram but
the sum runs over one or more of the labels. Then
we simply write the sum symbol in front of the diagram,
we may similarly insert scalar factors.
Now recall that for finite $\sys$ the non-degeneracy of
the braiding is encoded in the orthogonality relation
$\langle y^0,y^\la \rangle = \del \la 0 w$.
In terms of the statistics characters this reads
$\sum_\mu d_\mu \chi_\la(\mu) =d_\la^{-1}\del \la 0 w=\del \la 0 w$.
Graphically this can be represented as in Fig. \ref{ort0}.
%
% non-degeneracy and the "killing ring"
\begin{figure}[htb]
\begin{center}
\unitlength 0.6mm
\begin{picture}(195,45)
\thinlines
\put(10,22){\makebox(0,0){$\displaystyle\sum_{\mu\in\sys}\,\,d_\mu$}}
\put(45,25){\arc{20}{5.012}{4.412}}
\put(45,45){\line(0,-1){27}}
\put(45,12){\vector(0,-1){12}}
\put(35.12,27){\vector(0,1){0}}
\put(30,20){\makebox(0,0){$\mu$}}
\put(52,5){\makebox(0,0){$\la$}}
\put(70,24){\makebox(0,0){$=$}}
\put(90,22){\makebox(0,0){$\displaystyle\sum_{\mu\in\sys}$}}
\put(120,25){\arc{20}{0}{6.283}}
\put(110.12,27){\vector(0,1){0}}
\put(105,20){\makebox(0,0){$\mu$}}
\put(155,25){\arc{20}{5.012}{4.412}}
\put(155,45){\line(0,-1){27}}
\put(155,12){\vector(0,-1){12}}
\put(145.12,27){\vector(0,1){0}}
\put(140,20){\makebox(0,0){$\mu$}}
\put(162,5){\makebox(0,0){$\la$}}
\put(185,24){\makebox(0,0){$=\,\,\del \la 0 \, w$}}
\end{picture}
\end{center}
\caption{Orthogonality relation for a non-degenerate braiding
         (``killing ring'')}
\label{ort0}
\end{figure}
This kind of (graphical) relation has also been used more recently
in \cite{TW,O6,KL} and was called a ``killing ring'' in \cite{O6}.

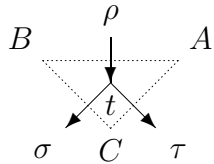
\begin{figure}[htb]
\begin{center}
\unitlength 0.6mm
\begin{picture}(60,40)
\thinlines
\dottedline{1}(15,25)(30,10)(45,25)(15,25)
\put(30,30){\vector(0,-1){10}}
\put(30,20){\vector(-1,-1){10}}
\put(30,20){\vector(1,-1){10}}
\put(10,30){\makebox(0,0){$B$}}
\put(30,5){\makebox(0,0){$C$}}
\put(50,30){\makebox(0,0){$A$}}
\put(15,5){\makebox(0,0){$\sigma$}}
\put(30,35){\makebox(0,0){$\rho$}}
\put(45,5){\makebox(0,0){$\tau$}}
\put(30,15){\makebox(0,0){$t$}}
\end{picture}
\end{center}
\caption{The intertwiner
         $\protect\sqrt[4]{\frac{d_\sigma d_\tau}{d_\rho}}t$
         as a triangle}
\label{tABC}
\end{figure}

Wire diagrams can also be used for intertwiners of morphisms
between different factors. Let $A,B,C$ infinite factors,
$\rho\in\Mor(A,B)$, $\sigma\in\Mor(C,B)$,
$\tau\in\Mor(A,C)$ irreducible morphisms and
$t\in\Hom(\rho,\sigma\tau)$ an isometry. Then Fig.\ \ref{tABC}
%
% triangle ABC
represents the intertwiner
$d_\sigma^{1/4} d_\tau^{1/4} d_\rho^{-1/4} t$.
Similarly we can draw a picture using a co-isometry.
Along the lines of the previous paragraphs, we can
similarly build up larger wire diagrams out of trivalent
vertices involving different factors. We do
not need the triangles with corners labelled by factors as
we can also label the regions between the wires. 
So far we do not have a meaningful way to cross wires with
differently labelled regions left and right, but all the
arguments listed above which do not involve braidings
can be used for intertwiners of morphisms between different
factors exactly as proceeded above. Moreover, the diagrams
may also involve wires where left and right regions are
labelled by the same factor, i.e.\ these wires correspond to
{\em endo}morphisms of some factor which may well
form a braided system, and then one may have crossings for
those wires.

\subsection{Frobenius reciprocity and rotations}
\label{graph2-prelim}

\begin{figure}[htb]
\begin{center}
\unitlength 0.6mm
\begin{picture}(150,50)
\thinlines
\put(10,20){\dashbox{2}(40,20){$t$}}
\put(30,50){\vector(0,-1){10}}
\put(20,20){\vector(0,-1){20}}
\put(40,20){\line(0,-1){5}}
\put(60,15){\vector(0,1){35}}
\put(50,15){\arc{20}{0}{3.142}}
\put(15,5){\makebox(0,0){$\rho$}}
\put(25,45){\makebox(0,0){$\tau$}}
\put(65,45){\makebox(0,0){$\sigma$}}
\put(80,30){\makebox(0,0){$=$}}
\put(100,20){\dashbox{2}(40,20)
{$\sqrt{\frac{d_\tau}{d_\rho}} \cR_\sigma(t)^*$}}
\put(110,50){\vector(0,-1){10}}
\put(130,40){\vector(0,1){10}}
\put(120,20){\vector(0,-1){20}}
\put(125,5){\makebox(0,0){$\rho$}}
\put(105,45){\makebox(0,0){$\tau$}}
\put(135,45){\makebox(0,0){$\sigma$}}
\end{picture}
\end{center}
\caption{Right Frobenius reciprocity for an intertwiner
$t\in\Hom(\tau,\rho\sigma)$}
\label{RFRt}
\end{figure}
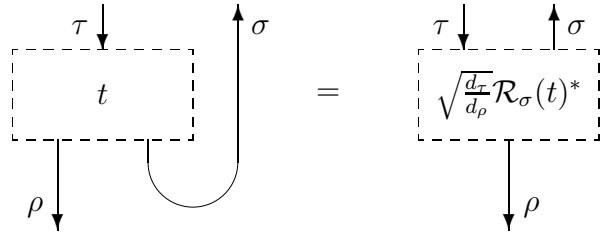

Let $A,B,C$ be infinite factors, $\rho\in\Mor(A,B)$,
$\tau\in\Mor(C,B)$, $\sigma\in\Mor(C,A)$
morphisms with finite statistical dimensions
$d_\rho,d_\tau,d_\sigma<\infty$, respectively, and let
$t\in\Hom(\tau,\rho\sigma)$. Then
\[ \cL_\rho(t) = \sqrt \frac{d_\rho d_\sigma}{d_\tau}
\bar\rho(t)^* r_\rho \in \Hom(\sigma,\co\rho \tau) \]
and
\[ \cR_\sigma(t) = \sqrt \frac{d_\rho d_\sigma}{d_\tau}
t^* \rho({\co r}_\sigma) \in \Hom(\rho,\tau \co\sigma) \]
are the images under left and right Frobenius maps.
Displaying the intertwiners $d_\rho^{1/2} r_\rho^* \bar\rho(t)$
and $d_\sigma^{1/2} \rho({\co r}_\sigma)^*t$ graphically yields
the identities in Figs.\ \ref{LFRt}
%
% Left Frobenius
%
\begin{figure}[htb]
\begin{center}
\unitlength 0.6mm
\begin{picture}(150,50)
\thinlines
\put(20,20){\dashbox{2}(40,20){$t$}}
\put(40,50){\vector(0,-1){10}}
\put(50,20){\vector(0,-1){20}}
\put(30,20){\line(0,-1){5}}
\put(10,15){\vector(0,1){35}}
\put(20,15){\arc{20}{0}{3.142}}
\put(5,45){\makebox(0,0){$\rho$}}
\put(45,45){\makebox(0,0){$\tau$}}
\put(55,5){\makebox(0,0){$\sigma$}}
\put(80,30){\makebox(0,0){$=$}}
\put(100,20){\dashbox{2}(40,20)
{$\sqrt{\frac{d_\tau}{d_\sigma}} \cL_\rho(t)^*$}}
\put(110,40){\vector(0,1){10}}
\put(130,50){\vector(0,-1){10}}
\put(120,20){\vector(0,-1){20}}
\put(105,45){\makebox(0,0){$\rho$}}
\put(135,45){\makebox(0,0){$\tau$}}
\put(125,5){\makebox(0,0){$\sigma$}}
\end{picture}
\end{center}
\caption{Left Frobenius reciprocity for an intertwiner
$t\in\Hom(\tau,\rho\sigma)$}
\label{LFRt}
\end{figure}
and \ref{RFRt}, respectively.
%
% Right Frobenius
%
These morphisms need not be irreducible. Taking them as products,
we may replace any of them by bundles of wires. We call the
linear isomorphisms $t\mapsto d_\rho^{1/2} r_\rho^* \bar\rho(t)$
and $t\mapsto d_\sigma^{1/2} \rho({\co r}_\sigma)^*t$ the
left and right Frobenius rotations.

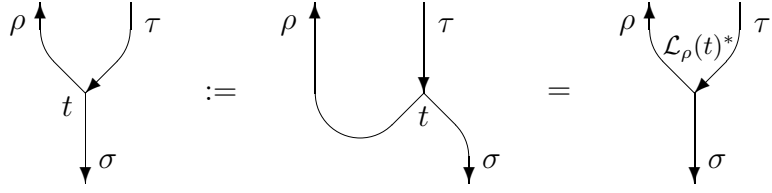
\begin{figure}[htb]
\begin{center}
\unitlength 0.6mm
\begin{picture}(170,40)
\thinlines
\put(20,20){\vector(0,-1){20}}
\put(12.929,27.071){\line(1,-1){7.071}}
\put(27.071,27.071){\vector(-1,-1){7.071}}
\put(20,34.142){\arc{20}{2.356}{3.142}}
\put(20,34.142){\arc{20}{0}{0.785}}
\put(10,34.142){\vector(0,1){5.858}}
\put(30,40){\line(0,-1){5.858}}
\put(5,35){\makebox(0,0){$\rho$}}
\put(35,35){\makebox(0,0){$\tau$}}
\put(25,5){\makebox(0,0){$\sigma$}}
\put(16,17){\makebox(0,0){$t$}}
\put(50,20){\makebox(0,0){$:=$}}
\put(70.858,20){\vector(0,1){20}}
\put(80.858,20){\arc{20}{0.785}{3.142}}
\put(95,40){\vector(0,-1){20}}
\put(95,20){\line(-1,-1){7.071}}
\put(95,20){\line(1,-1){7.071}}
\put(95,5.858){\arc{20}{5.498}{0}}
\put(105,5.858){\vector(0,-1){5.858}}
\put(65,35){\makebox(0,0){$\rho$}}
\put(100,35){\makebox(0,0){$\tau$}}
\put(110,5){\makebox(0,0){$\sigma$}}
\put(95,15){\makebox(0,0){$t$}}
\put(125,20){\makebox(0,0){$=$}}
\put(155,20){\vector(0,-1){20}}
\put(147.929,27.071){\line(1,-1){7.071}}
\put(162.071,27.071){\vector(-1,-1){7.071}}
\put(155,34.142){\arc{20}{2.356}{3.142}}
\put(155,34.142){\arc{20}{0}{0.785}}
\put(145,34.142){\vector(0,1){5.858}}
\put(165,40){\line(0,-1){5.858}}
\put(140,35){\makebox(0,0){$\rho$}}
\put(170,35){\makebox(0,0){$\tau$}}
\put(160,5){\makebox(0,0){$\sigma$}}
\put(156,30){\makebox(0,0){{\footnotesize $\cL_\rho(t)^*$}}}
\end{picture}
\end{center}
\caption{Left Frobenius reciprocity for a trivalent vertex
labelled by an isometry}
\label{Ltight}
\end{figure}

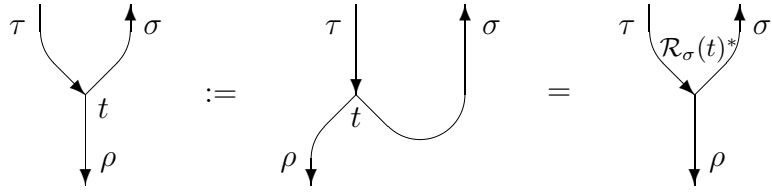
\begin{figure}[htb]
\begin{center}
\unitlength 0.6mm
\begin{picture}(170,40)
\thinlines
\put(20,20){\vector(0,-1){20}}
\put(12.929,27.071){\vector(1,-1){7.071}}
\put(27.071,27.071){\line(-1,-1){7.071}}
\put(20,34.142){\arc{20}{2.356}{3.142}}
\put(20,34.142){\arc{20}{0}{0.785}}
\put(10,40){\line(0,-1){5.858}}
\put(30,34.142){\vector(0,1){5.858}}
\put(5,35){\makebox(0,0){$\tau$}}
\put(35,35){\makebox(0,0){$\sigma$}}
\put(25,5){\makebox(0,0){$\rho$}}
\put(24,17){\makebox(0,0){$t$}}
\put(50,20){\makebox(0,0){$:=$}}
\put(104.142,20){\vector(0,1){20}}
\put(94.142,20){\arc{20}{0}{2.356}}
\put(80,40){\vector(0,-1){20}}
\put(80,20){\line(-1,-1){7.071}}
\put(80,20){\line(1,-1){7.071}}
\put(80,5.858){\arc{20}{3.142}{3.927}}
\put(70,5.858){\vector(0,-1){5.858}}
\put(110,35){\makebox(0,0){$\sigma$}}
\put(75,35){\makebox(0,0){$\tau$}}
\put(65,5){\makebox(0,0){$\rho$}}
\put(80,15){\makebox(0,0){$t$}}
\put(125,20){\makebox(0,0){$=$}}
\put(155,20){\vector(0,-1){20}}
\put(147.929,27.071){\vector(1,-1){7.071}}
\put(162.071,27.071){\line(-1,-1){7.071}}
\put(155,34.142){\arc{20}{2.356}{3.142}}
\put(155,34.142){\arc{20}{0}{0.785}}
\put(145,40){\line(0,-1){5.858}}
\put(165,34.142){\vector(0,1){5.858}}
\put(140,35){\makebox(0,0){$\tau$}}
\put(170,35){\makebox(0,0){$\sigma$}}
\put(160,5){\makebox(0,0){$\rho$}}
\put(156,30){\makebox(0,0){{\footnotesize $\cR_\sigma(t)^*$}}}
\end{picture}
\end{center}
\caption{Right Frobenius reciprocity for a trivalent vertex
labelled by an isometry}
\label{Rtight}
\end{figure}

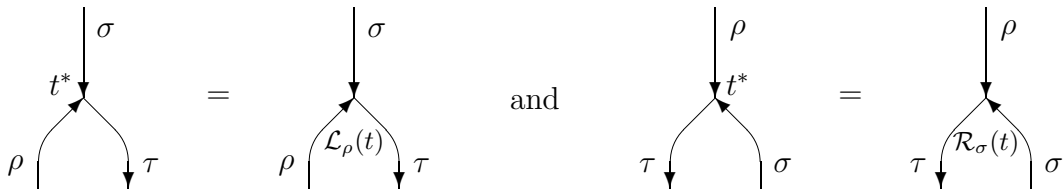
\begin{figure}[htb]
\begin{center}
\unitlength 0.6mm
\begin{picture}(240,40)
\thinlines
\put(20,40){\vector(0,-1){20}}
\put(12.929,12.929){\vector(1,1){7.071}}
\put(27.071,12.929){\line(-1,1){7.071}}
\put(20,5.858){\arc{20}{3.142}{3.927}}
\put(20,5.858){\arc{20}{5.498}{0}}
\put(10,5.858){\line(0,-1){5.858}}
\put(30,5.858){\vector(0,-1){5.858}}
\put(5,5){\makebox(0,0){$\rho$}}
\put(35,5){\makebox(0,0){$\tau$}}
\put(25,35){\makebox(0,0){$\sigma$}}
\put(15,23){\makebox(0,0){$t^*$}}
\put(50,20){\makebox(0,0){$=$}}
\put(80,40){\vector(0,-1){20}}
\put(72.929,12.929){\vector(1,1){7.071}}
\put(87.071,12.929){\line(-1,1){7.071}}
\put(80,5.858){\arc{20}{3.142}{3.927}}
\put(80,5.858){\arc{20}{5.498}{0}}
\put(70,5.858){\line(0,-1){5.858}}
\put(90,5.858){\vector(0,-1){5.858}}
\put(65,5){\makebox(0,0){$\rho$}}
\put(95,5){\makebox(0,0){$\tau$}}
\put(85,35){\makebox(0,0){$\sigma$}}
\put(80,10){\makebox(0,0){{\footnotesize $\cL_\rho(t)$}}}
\put(120,20){\makebox(0,0){and}}
\put(160,40){\vector(0,-1){20}}
\put(152.929,12.929){\line(1,1){7.071}}
\put(167.071,12.929){\vector(-1,1){7.071}}
\put(160,5.858){\arc{20}{3.142}{3.927}}
\put(160,5.858){\arc{20}{5.498}{0}}
\put(150,5.858){\vector(0,-1){5.858}}
\put(170,5.858){\line(0,-1){5.858}}
\put(145,5){\makebox(0,0){$\tau$}}
\put(175,5){\makebox(0,0){$\sigma$}}
\put(165,35){\makebox(0,0){$\rho$}}
\put(165,23){\makebox(0,0){$t^*$}}
\put(190,20){\makebox(0,0){$=$}}
\put(220,40){\vector(0,-1){20}}
\put(212.929,12.929){\line(1,1){7.071}}
\put(227.071,12.929){\vector(-1,1){7.071}}
\put(220,5.858){\arc{20}{3.142}{3.927}}
\put(220,5.858){\arc{20}{5.498}{0}}
\put(210,5.858){\vector(0,-1){5.858}}
\put(230,5.858){\line(0,-1){5.858}}
\put(205,5){\makebox(0,0){$\tau$}}
\put(235,5){\makebox(0,0){$\sigma$}}
\put(225,35){\makebox(0,0){$\rho$}}
\put(220,10){\makebox(0,0){{\footnotesize $\cR_\sigma(t)$}}}
\end{picture}
\end{center}
\caption{Frobenius reciprocity for a trivalent vertex labelled
by a co-isometry}
\label{RLtight*}
\end{figure}

Now let us assume that $t$ is isometric and labels a
trivalent vertex of wires corresponding to irreducible
morphisms $\rho,\tau,\sigma$.
With the above ``transformation law'' we then have
the identity of Fig.\ \ref{Ltight}, where the first equality is just
a definition which gives us some prescription of ``tightening''
wires at trivalent vertices.
%
% Left Frobenius for trivalent vertex with isometry
%
In fact, the label $\cL_\rho(t)^*$ of the trivalent vertex
makes sense since it is a co-isometry: Due to irreducibility
of $\tau$ and $\sigma$, the map $t\mapsto\cL_\rho(t)^*$ is
isometric. Similarly, we get Fig.\ \ref{Rtight} (using
irreducibility of $\tau$ and $\rho$).
%
% Right Frobenius for trivalent vertex with isometry
%
Hence the prefactor in Figs.\ \ref{LFRt} and \ref{RFRt}
is just such that it transforms isometries with natural
normalization prefactors into co-isometries with natural
normalization prefactors and, by taking adjoints,
the other way round which gives the graphical identities given in
Fig.\ \ref{RLtight*}.
%
% Frobenius for trivalent vertex with co-isometry
%
We may now use the replacement prescription three times,
beginning with a trivalent vertex labelled by an isometry
$t\in\Hom(\tau,\rho\sigma)$ and proceeding in a clockwise
direction. Then we end up with a co-isometry
$\Theta(t)^*\in\Hom(\co\sigma\co\rho,\co\tau)$ in the
corner where we originally had the label $t$. In fact,
\[ \Theta(t) = \cR_\rho(\cL_\tau(\cR_\sigma(t))) =
\sqrt{d_\rho d_\sigma d_\tau}\, r_\tau^* \co\tau
(t^*\rho({\co r}_\sigma) {\co r}_\rho) \,. \]
Similarly we can go in the counter-clockwise direction and
then we obtain $\tilde\Theta(t)^*\in\Hom(\co\sigma\co\rho,\co\tau)$,
where
\[ \tilde\Theta(t) = \cL_\sigma(\cR_\tau(\cL_\rho(t))) =
\sqrt{d_\rho d_\sigma d_\tau}\, \co\sigma\co\rho({\co r}_\tau^* t^*)
\co\sigma(r_\rho) r_\sigma \,, \]
and in order to establish a well-defined rotation procedure
we have to show that $\Theta(t)=\tilde\Theta(t)$. Now
\[ \bearll
\Theta(t)^* \tilde\Theta(t) &= \sqrt{d_\rho d_\sigma d_\tau} \,
\co\tau({\co r}_\tau^* t^*) \Theta(t)^*
\co\sigma(r_\rho) r_\sigma \\[.4em]
&= d_\rho d_\sigma d_\tau  \, \co\tau({\co r}_\tau^* t^*) 
\co\tau ({\co r}_\rho^* \rho({\co r}_\sigma^*)t)
\co\tau\tau(\co\sigma(r_\rho) r_\sigma)r_\tau \\[.4em]
&= d_\rho d_\sigma d_\tau  \, \co\tau({\co r}_\tau^* t^* 
{\co r}_\rho^* \rho({\co r}_\sigma^*)
\rho\sigma(\co\sigma(r_\rho) r_\sigma) t) r_\tau \\[.4em]
&= d_\rho d_\tau  \, \co\tau({\co r}_\tau^* t^* 
{\co r}_\rho^* \rho(r_\rho) t) r_\tau
= d_\tau  \, \co\tau({\co r}_\tau^*) r_\tau = \bfe \,,
\eear \]
hence
$(\Theta(t)-\tilde\Theta(t))^*(\Theta(t)-\tilde\Theta(t))=0$,
i.e.\ $\Theta(t)=\tilde\Theta(t)$. Thus a
trivalent vertex labelled with an isometry
$t\in\Hom(\tau,\rho\sigma)$
can equivalently be labelled with a co-isometry
$\Theta(t)^*\in\Hom(\co\sigma\co\rho,\co\tau)$.
So here we have established some ``rotation invariance''
of trivalent vertices (in standard inverted Y-shape or Y-shape)
with a replacement prescription for the rotated labelling
(co-) isometries.

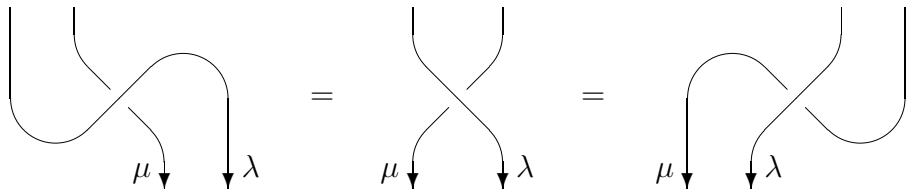
\begin{figure}[htb]
\begin{center}
\unitlength 0.6mm
\begin{picture}(220,40)
\thinlines
\put(10.858,40){\line(0,-1){20}}
\put(20.858,20){\arc{20}{0.785}{3.142}}
\put(35,20){\line(-1,-1){7.071}}
\put(35,20){\line(1,1){7.071}}
\put(49.142,20){\arc{20}{3.926}{0}}
\put(59.142,20){\vector(0,-1){20}}
\put(25,40){\line(0,-1){5.858}}
\put(33,22){\line(-1,1){5.071}}
\put(37,18){\line(1,-1){5.071}}
\put(35,34.142){\arc{20}{2.356}{3.142}}
\put(45,5.858){\vector(0,-1){5.858}}
\put(35,5.858){\arc{20}{5.498}{0}}
\put(40,3.5){\makebox(0,0){$\mu$}}
\put(64.142,5){\makebox(0,0){$\la$}}
\put(80,20){\makebox(0,0){$=$}}
\put(100,40){\line(0,-1){5.858}}
\put(110,20){\line(-1,1){7.071}}
\put(110,20){\line(1,-1){7.071}}
\put(110,34.142){\arc{20}{2.356}{3.142}}
\put(120,5.858){\vector(0,-1){5.858}}
\put(110,5.858){\arc{20}{5.498}{0}}
\put(120,40){\line(0,-1){5.858}}
\put(112,22){\line(1,1){5.071}}
\put(108,18){\line(-1,-1){5.071}}
\put(110,34.142){\arc{20}{0}{0.785}}
\put(100,5.858){\vector(0,-1){5.858}}
\put(110,5.858){\arc{20}{3.142}{3.927}}
\put(95,3.5){\makebox(0,0){$\mu$}}
\put(125,5){\makebox(0,0){$\la$}}
\put(140,20){\makebox(0,0){$=$}}
\put(160.858,20){\vector(0,-1){20}}
\put(170.858,20){\arc{20}{3.142}{5.498}}
\put(183,22){\line(-1,1){5.071}}
\put(187,18){\line(1,-1){5.071}}
\put(199.142,20){\arc{20}{0}{2.356}}
\put(209.142,40){\line(0,-1){20}}
\put(185,20){\line(-1,-1){7.071}}
\put(185,20){\line(1,1){7.071}}
\put(195,40){\line(0,-1){5.858}}
\put(185,34.142){\arc{20}{0}{0.785}}
\put(175,5.858){\vector(0,-1){5.858}}
\put(185,5.858){\arc{20}{3.142}{3.927}}
\put(155.858,3.5){\makebox(0,0){$\mu$}}
\put(180,5){\makebox(0,0){$\la$}}
\end{picture}
\end{center}
\caption{Rotation of crossings}
\label{rotcross}
\end{figure}

Next we turn to the rotation of crossings when we have a braiding.
Assume we have a braided system of endomorphisms
$\Delta\ni\la,\mu,\nu$ of some factor $A$. From the BFE we obtain
$r_\la = \bar\la (\epsmp \mu\la) \epsmp \mu{\bar\la} \mu(r_\la)$.
Applying $\la$ and multiplying by
$d_\la \epspm\la\mu {\co r}_\la^*$ from the left yields
\be
\epspm \la\mu = d_\la {\co r}_\la^* \la
(\epsmp \mu{\bar\la}) \la\mu(r_\la)\,.
\label{rotcreq}
\ee
The BFE yields similarly
$\la({\co r}_\mu)= \epspm \mu\la \mu(\epsmp{\co\mu}\la) {\co r}_\mu$,
and by multiplying with $d_\mu \mu\la(r_\mu^*) \epspm\la\mu$ from
the left we obtain
\[ \epspm \la\mu = d_\mu \mu\la(r_\mu^*)\mu
(\epsm{\co\mu}\la){\co r}_\mu \,,\]
and therefore we have the graphical identity
given in Fig.\ \ref{rotcross}, here displayed only
for overcrossings. 
%
% Rotation of crossings
%
Then this procedure can even be iterated so that we obtain
arbitrarily twisted crossings. Note that for the rotation
of crossings we do not need any relabelling prescription
as this is encoded in the BFE's.

We now turn to the discussion of ``abstract pictures''
which admit different intertwiner interpretations according to
Frobenius rotations. Let $A_1,A_2,...,A_\ell$ be factors equipped
with sets $\sys_{i,j}\subset\Mor(A_i,A_j)$, $i,j=1,2,...,\ell$,
of irreducible, pairwise inequivalent morphisms with finite index
such that $\bigsqcup_{i,j} \sys_{i,j}$ is
closed under conjugation and irreducible decomposition
of products (whenever composable) as sectors, and in particular
each $\sys_{i,i}$ is a system of endomorphisms.
Some of the systems $\sys_{i,i}$ may be braided.

We now consider ``labelled knotted graphs'' of the following
form. On a finite connected and simply connected
region in the plane we have a finite
number of wires (i.e.\ images of piecewise $C^\infty$ maps from
the unit interval into the region). Within the region
there is a finite number of trivalent vertices (i.e.\ common
endpoints of three wires) and crossings of two wires, and
for the latter there
is a notion of over- and undercrossing (i.e.\ for each
crossing there is one wire ``on top of the other'').
If wires are not closed
(i.e.\ if their two endpoints do not coincide)
then they are only allowed to have trivalent vertices or distinguished
points on the boundary of the region as their endpoints.
The wires meet each other only at the trivalent vertices and
crossings, and they are directed and labelled by the
morphisms in $\bigsqcup_{i,j} \sys_{i,j}$ subject to the
following rules. Crossings are only possible
for wires with labelling morphisms in some $\sys_{i,i}$
with braiding. Furthermore it must be possible to associate
the factors $A_i$ to the free regions between
the wires such that any wire labelled by some
$\rho\in\sys_{i,j}$ has the ``source'' factor $A_i$ on its
left and the ``range'' factor $A_j$ on its right
relative to the  orientation
(composition compatibility). We identify graphs which are
transformed into each other by inversion of the orientation
of a wire and simultaneous replacement of its label, say
$\rho\in\sys_{i,j}$, by the representative conjugate
morphism $\bar\rho\in\sys_{j,i}$.
Finally, the trivalent vertices are labelled either by isometric or
co-isometric intertwiners which are associated locally
to one corner region of the trivalent vertex as follows.
If $\tau\in\sys_{i,j}$, $\rho\in\sys_{k,j}$, $\sigma\in\sys_{i,k}$
label the three wires of a trivalent vertex, $\tau$ is entering
and, following counter-clockwise, $\rho$ and $\sigma$ are
outgoing (as e.g.\ the trivalent vertex in Figs.\ \ref{Ltight}
and \ref{Rtight}, possibly up to isotopy and rotation), then in
the local corner region opposite to $\tau$ the label must either
be an isometry $t\in\Hom(\tau,\rho\sigma)$ or a co-isometry
$s^*\in\Hom(\co\sigma\co\rho,\co\tau)$. If the wires at a
trivalent vertex have orientation different from this, the rule
can be derived from the previous case by reversing orientations
and simultaneous relabelling by conjugate morphisms.

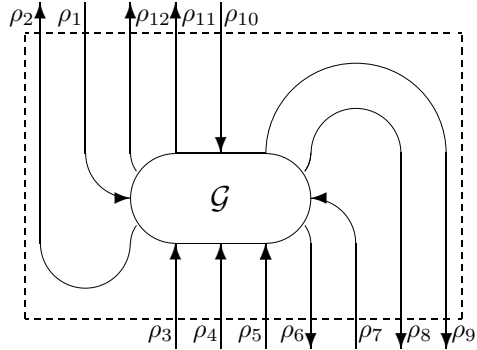
\begin{figure}[htb]
\begin{center}
\unitlength 0.4mm
\begin{picture}(150,115)
\thinlines
\put(5,10){\dashbox{2}(145,95){}}
\put(55,50){\arc{30}{1.571}{4.712}}
\put(85,50){\arc{30}{4.712}{1.571}}
\put(55,35){\line(1,0){30}}
\put(55,65){\line(1,0){30}}
\put(115,65){\arc{60}{3.142}{0}}
\put(40,65){\arc{30}{1.571}{3.142}}
\put(40,65){\vector(0,1){50}}
\put(25,65){\line(0,1){50}}
\put(100,35){\arc{30}{4.712}{0}}
\put(55,65){\vector(0,1){50}}
\put(70,115){\vector(0,-1){50}}
\put(55,0){\vector(0,1){35}}
\put(70,0){\vector(0,1){35}}
\put(85,0){\vector(0,1){35}}
\put(50,65){\arc{20}{2.51}{3.142}}
\put(90,65){\arc{20}{0}{0.632}}
\put(90,35){\arc{20}{5.652}{0}}
\put(50,35){\arc{20}{3.142}{3.773}}
\put(115,65){\arc{30}{3.142}{0}}
\put(90,65){\arc{20}{0}{0.632}}
\put(25,35){\arc{30}{0}{3.142}}
\put(85,0){\vector(0,1){35}}
\put(10,35){\vector(0,1){80}}
\put(100,35){\vector(0,-1){35}}
\put(115,0){\line(0,1){35}}
\put(130,65){\vector(0,-1){65}}
\put(145,65){\vector(0,-1){65}}
\put(40,50){\vector(1,0){0}}
\put(100,50){\vector(-1,0){0}}
\put(70,50){\makebox(0,0){$\cG$}}
\put(4,110){\makebox(0,0){{\footnotesize $\rho_2$}}}
\put(20,110){\makebox(0,0){{\footnotesize $\rho_1$}}}
\put(47.5,110){\makebox(0,0){{\footnotesize $\rho_{12}$}}}
\put(63,110){\makebox(0,0){{\footnotesize $\rho_{11}$}}}
\put(77,110){\makebox(0,0){{\footnotesize $\rho_{10}$}}}
\put(151,5){\makebox(0,0){{\footnotesize $\rho_9$}}}
\put(136,5){\makebox(0,0){{\footnotesize $\rho_8$}}}
\put(120,5){\makebox(0,0){{\footnotesize $\rho_7$}}}
\put(50,5){\makebox(0,0){{\footnotesize $\rho_3$}}}
\put(65,5){\makebox(0,0){{\footnotesize $\rho_4$}}}
\put(80,5){\makebox(0,0){{\footnotesize $\rho_5$}}}
\put(94,5){\makebox(0,0){{\footnotesize $\rho_6$}}}
\end{picture}
\end{center}
\caption{A Frobenius annulus surrounding $\cG$}
\label{annul}
\end{figure}

Now let $\cG$ be such a labelled knotted graph as above.
To interpret $\cG$ as an intertwiner, we may put it
in some ``Frobenius annulus'' as shown in Fig.\ \ref{annul}
for an example.\footnote{Our notion of a Frobenius annulus
is inspired by the annular invariance used in Jones'
definition of a ``general planar algebra'' \cite{J2}.}
%
% Frobenius annulus
%
A Frobenius annulus has labelled wires inside such that each
of them meets an open end of a wire of $\cG$ at one endpoint
(labelled by $\rho_1$,...,$\rho_{12}$ in our example),
matching the label and orientation of this wire,
and this way all the open ends of the wires of $\cG$ are either
connected to the top or bottom of the outside square boundary
of the annulus. No crossings or trivalent vertices are allowed in
the annulus, but it may contain cups or caps. Gluing the wires
together and forgetting about the boundary of $\cG$
and the annulus, we will read the result as a wire diagram
and therefore the annulus corresponds to a ``Frobenius choice'',
deciding whether we will get a certain intertwiner or its image
by  certain Frobenius rotations,
cf.\ Figs.\ \ref{LFRt} and \ref{RFRt} (and their adjoints).

Reading vertically downwards, we may now have the problem that
on a finite number of horizontal levels a finite
number of singular points of crossings, trivalent vertices,
cups and caps are exactly on the same level (or ``height'')
so that we cannot time slice the diagram into stripes containing
only one elementary intertwiner. Also some wires
may have pieces going exactly horizontally.
We now allow to make small vertical
translations such that these crossings and trivalent vertices
are put on slightly different levels and all wires obtain
piecewise slopes, without letting wires
touch or producing new crossings, but we may possibly produce
some new cups or caps. In the latter case we can always arrange
it so that even each new cup or cap appears on a distinct level.
The trivalent vertices and crossings may not be in ``standard form'',
i.e.\ in Y- or inverted Y shape respectively X-shape.
In an ``$\epsilon$-neighborhood'' of a trivalent vertex,
we now bend the wires so that the angles are arranged
in standard form. Similarly we modify the crossings
to bend them into an X-shape.
Using for labels at trivalent vertices our replacement prescription
by Frobenius reciprocity, we can obtain isometries as labels
for trivalent vertices in inverted Y-shape, located on the bottom
corner region, and co-isometries as labels for trivalent vertices in
Y-shape, located on the top corner region.

Again, these topological moves are allowed to produce at most
new cups or caps, all on different levels so that the resulting
diagram can be time sliced into stripes of elementary diagrams.
Clearly, this procedure of deforming a labelled knotted graph
in a Frobenius annulus into a regular wire diagram is highly
ambiguous. However, the ambiguities in the above procedures
are irrelevant: The ambiguities arising from the production
of slopes of wires and different levels
of certain elementary intertwiners are irrelevant due to the
topological invariance of Fig.\ \ref{isoinv1} and the
freedom of translating intertwiners vertically
as shown in Fig.\ \ref{xx'},
and the ambiguities arising from rotations of the elementary
intertwiners are irrelevant due to the rotation invariance
of trivalent vertices and crossings, as we have established in
Figs.\ \ref{Ltight}, \ref{Rtight}, \ref{RLtight*} and
\ref{rotcross}.

Now let $\cG_1$ and $\cG_2$ be two labelled knotted graphs
as above which are defined on the same (connected, simply
connected) region in the plane
and have the same entering and outgoing wires at the same
points with the same orientation, i.e.\ they have coinciding
open ends so that they fit in the same Frobenius annuli.
When embedded in some Frobenius annulus it may now happen
that the corresponding intertwiners are the same, even if
$\cG_1$ and $\cG_2$ are different. Because of the isomorphism
property of Frobenius rotations it is clear that then
$\cG_1$ and $\cG_2$ yield the same intertwiner through
embedding in any Frobenius annulus. We can write down
sufficient conditions for such equality in terms of some
``regular isotopy'': For given $\cG_1$ and $\cG_2$ as above
choose a Frobenius annulus and regularize the pictures
into two wire diagrams $\cW_1$ and $\cW_2$, respectively.
We call $\cG_1$ and $\cG_2$ regularly isotopic if $\cW_1$
can be transformed into $\cW_2$ by the following list of moves:
\begin{enumerate}
\item Reversing orientation of some wires with simultaneous relabelling
by conjugate morphisms,
\item any horizontal translations of elementary intertwiners
which may change slopes of wires but which do not let the wires
meet or involve cups or caps,
\item vertical translations of elementary intertwiners as in Fig.\ \ref{xx'},
\item topological moves as in Fig.\ \ref{isoinv1},
\item rotations of trivalent vertices and their labels
as in Figs.\ \ref{Ltight}, \ref{Rtight}, \ref{RLtight*},
\item and for wires corresponding to a braided system
$\sys_{i,i}$ we additionally admit
\begin{enumerate}
\item vertical Reidemeister moves of type II as in Fig.\ \ref{Reid2},
\item moving crossings over and under trivalent vertices,
cups and caps according to the BFE's
(cf.\ Figs.\ \ref{wireBFE1} and \ref{wireBFE2} for the
first two relations),
\item vertical Reidemeister moves of type III for crossings
(cf.\ Fig.\ \ref{Reid3} for overcrossings),
\item rotations of crossings
(cf.\ Fig. \ref{rotcross} for overcrossings).
\end{enumerate}
\end{enumerate}
Thus the ambiguity in the regularization procedure means in
particular that from one graph we can only obtain wire diagrams
that can be transformed into each other by these moves.
It is easy to see that regular isotopy is an equivalence
relation for knotted labelled graphs. Moreover, for
closed labelled knotted graphs (i.e.\ without open ends)
which are then embedded in a trivial annulus, the local
rotation invariance of the elementary intertwiners ends
up in a total rotation invariance: We can rotate the
picture freely, the rotated graph is always regularly isotopic
to the original one and we will always end up with the same
scalar (times $\bfe_{A_i}$, where $A_i$ is the factor
associated to the outside region).\footnote{For a
single kind of wire corresponding to a braided system,
this invariance is similar to the complex number-valued
regular isotopy invariant of knotted graphs obtained
in \cite{MO}.}

\begin{figure}[htb]
\begin{center}
\unitlength 0.6mm
\begin{picture}(160,30)
\thinlines
\put(5,15){\makebox(0,0){$d_\rho$}}
\put(15,10){\dashbox{2}(10,10){$x$}}
\put(20,30){\vector(0,-1){10}}
\put(20,10){\vector(0,-1){10}}
\put(25,26){\makebox(0,0){$\rho$}}
\put(25,4){\makebox(0,0){$\rho$}}
\put(40,15){\makebox(0,0){$=$}}
\put(50,10){\dashbox{2}(10,10){$x$}}
\put(65,20){\arc{20}{3.142}{0}}
\put(65,10){\arc{20}{0}{3.142}}
\put(75,10){\line(0,1){10}}
\put(75,17){\vector(0,1){0}}
\put(79,5){\makebox(0,0){$\rho$}}
\put(105,15){\makebox(0,0){$\longleftrightarrow$}}
\put(135,10){\line(0,1){10}}
\put(135,17){\vector(0,1){0}}
\put(145,20){\arc{20}{3.142}{0}}
\put(145,10){\arc{20}{0}{3.142}}
\put(150,10){\dashbox{2}(10,10){$x$}}
\put(131,5){\makebox(0,0){$\rho$}}
\end{picture}
\end{center}
\caption{Two intertwiners of the same scalar value}
\label{cut}
\end{figure}
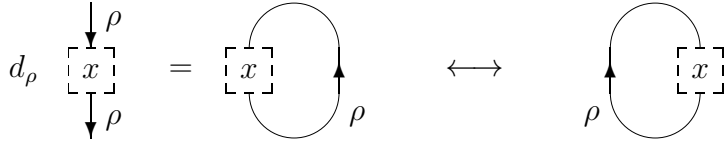

Let us finally consider an intertwiner $x\in\Hom(\rho,\rho)$
with $\rho\in\Mor(A,B)$ irreducible. Then clearly $x$ is a
scalar: $x=\xi\bfe_B$, $\xi\in\bbC$. Hence we have the identity
$d_\rho \xi\bfe_B \equiv d_\rho x
= d_\rho {\co r}_\rho^* x {\co r}_\rho$,
and this is graphically the left-hand side in
Fig.\ \ref{cut}.
%
% Cutting a wire
On the other hand, application of the left inverse yields
$d_\rho \phi_\rho(x)=d_\rho r_\rho^*\co\rho(x)r_\rho=d_\rho\xi\bfe_A$,
which is a different intertwiner of the same scalar value,
and it is represented graphically by the right-hand side in
Fig.\ \ref{cut}. Thus the left- and right-hand side in
Fig.\ \ref{cut} represent the same scalar. If we consider
closed wire diagrams and are only interested in the scalars
they represent, then we therefore have a
``regular isotopy on the 2-sphere''.

\subsection{$\a$-Induction for braided subfactors}
\label{sec-aindbsf}

We now consider $\a$-induction of \cite{BE1,BE2,BE3} in the
setting of braided subfactors. Here we work with a 
type III subfactor $N\subset M$,
equipped with a braided system $\sys\subset\End(N)$ in the
sense of Definition \ref{system} such that for the injection
map $\iota:N\to M$, the sector $[\co\iota\iota]$ decomposes
into a finite sum of sectors of morphisms in $\sys$.
(Here $\co\iota$ denotes any choice of a representative morphism
for the conjugate sector of $[\iota]$.) Note that since elements
in $\sys$ have by definition finite statistical dimension,
it follows that the injection map has finite statistical
dimension and thus the subfactor $N\subset M$ has finite index.
But also note that we did neither assume the finite depth
condition on $N\subset M$ (we did not assume finiteness of $\sys$)
nor non-degeneracy of the braiding at this point.
As usual, we denote the canonical endomorphism
$\iota\co\iota\in\End(M)$ by $\can=\iota\co\iota$,
the dual canonical endomorphism $\co\iota\iota\in\End(N)$
by $\canr=\co\iota\iota$ and ``canonical'' isometries
by $v\in M$ and $w\in N$, more precisely, we have
$v\in\Hom(\id_M,\can)$ and $w\in\Hom(\id_N,\canr)$
such that $w^*v=\can(v^*)w=[M:N]^{-1/2}\bfe$. Recall
that we have pointwise equality $M=Nv$.

With a braiding $\e$ on $\sys$ and its
extension to $\Sigma(\sys)$ as in Subsection \ref{braid-prelim}
we can define the $\a$-induced $\a^\pm_\la$ for $\la\in\Ssys$
exactly as in \cite{LR,BE1}, namely we define
\[ \a_\la^\pm = \co\iota^{\,-1} \circ \Ad (\epspm \lambda\canr)
\circ \lambda \circ \co\iota \,. \]
Then $\a^+_\la$ and $\a^-_\la$ are morphisms in $\Mor(M,M)$
with the properties $\a^\pm_\la\circ\iota=\iota\circ\la$,
$\a^\pm_\la(v)=\epspm \la\canr ^* v$,
$\a_{\la\mu}^\pm=\a_\la^\pm \a_\mu^\pm$ if also
$\mu\in\Ssys$, and clearly $\a_{{\rm{id}}_N}^\pm={{\rm{id}}}_M$.
Note that the first property yields immediately
$d_{\a_\la^\pm}=d_\la$ by the multiplicativity of
the minimal index \cite{L2.5}. We also obtain
easily that $\overline{\a_\la^\pm}=\a_{\co\la}^\pm$, since
we obtain $r_\la=\epspm \canr{\co\la\la}\canr(r_\la)$,
and similarly
${\co r}_\la=\epspm \canr {\la\co\la} \canr({\co r}_\la)$
easily from \erf{nat}. Multiplying both relations by $v$ from
the right yields $r_\la v= \a_{\co\la}^\pm \a_\la^\pm(v) r_\la$
and ${\co r}_\la v= \a_\la^\pm \a_{\co\la}^\pm(v) {\co r}_\la$,
hence $r_\la\in\Hom(\id_M,\a_{\co\la}^\pm\a_\la^\pm)$,
${\co r}_\la\in\Hom(\id_M,\a_\la^\pm\a_{\co\la}^\pm)$
as $M=Nv$, thus we can put $R_{\a_\la^\pm}=\iota(r_\la)$,
${\co R}_{\a_\la^\pm}=\iota({\co r}_\la)$ as R-isometries
for the $\a$-induced morphisms, i.e.\
$\overline{\a_\la^\pm}=\a_{\co\la}^\pm$.
Note also that the definition of $\a^\pm_\la$ does not
depend on the choice of the representative morphism
$\co\iota$ for the conjugate sector of $[\iota]$ due to the
transformation properties of the braiding operators, \erf{eun}.

Though the local net structure for $N(I)\subset M(I)$
is assumed in \cite{LR,BE1}, we need only an
assumption of a braiding for the definition of
$\a^\pm_\la$. We, however, have to be careful,
because we do not assume the chiral locality condition
$\eps \canr\canr \can(v) = \can(v)$ in this paper.
(The name ``chiral locality'' is motivated from the
treatment of extensions of chiral observables in conformal
field theory in the setting of nets of subfactors
\cite{LR}, where the extended net is shown to satisfy
local commutativity if and only if the condition
$\eps \canr\canr \can(v) = \can(v)$ is met
\cite[Thm.\ 4.9]{LR}.)
Some theorems in \cite{BE1,BE2,BE3} do depend on the
chiral locality condition and are {\sl not} true in this more
general setting of $\a$-induction. Namely, with
$\eps \canr\canr \can(v) = \can(v)$ it was easily derived
\cite[Lemma 3.5]{BE1} by using the BFE that then
$\Hom(\a_\la^\pm,\a_\mu^\pm)=\Hom(\iota\la,\iota\mu)$
for $\la,\mu\in\Ssys$. As a surprising corollary
(cf.\ \cite[Cor.\ 3.6]{BE1}) one found
by putting $\la=\mu=\id_N$ that $\iota$, thus the subfactor
$N\subset M$, was irreducible which had not been assumed.
Another corollary was then the ``main formula''
\cite[Thm.\ 3.9]{BE1}, giving
$\langle\a_\la^\pm,\a_\mu^\pm\rangle=\langle\iota\la,\iota\mu\rangle
=\langle\canr\la,\mu\rangle$ by Frobenius reciprocity.
(Moreover, in the framework of {\it nets of} subfactors
$\cN\subset\cM$,
where the braidings arise from the transportability of
localized endomorphisms, a certain reciprocity formula
$\langle \a_\la^\pm,\beta \rangle=\langle \la,\sigma_\beta\rangle$,
called ``$\alpha\sigma$-reciprocity'', between localized
transportable endomorphisms $\la$ and $\beta$ of the smaller
respectively the larger net was established; here
$\sigma$-restriction is essentially
$\sigma_\beta=\co\iota\beta\iota$.)
Without chiral locality, these results are in general not
true: The subfactor $N\subset M$ is neither forced to be
irreducible, nor does the main formula hold, however, we
always have the inequality
$\langle\a_\la^\pm,\a_\mu^\pm\rangle\le\langle\canr\la,\mu\rangle$,
since only the ``$\ge$'' part of the proof of
\cite[Thm.\ 3.9]{BE1} uses chiral locality.

It is a simple application of the braiding fusion equation
and does not involve chiral locality that for
$\la,\mu,\nu\in\Ssys$ we have the 
(equivalent) relations \cite[Lemma 3.25]{BE1}
\be
\a^\mp_\rho(Q) \epspm \la\rho = \epspm \mu\rho Q \,,\qquad
Q \epspm \rho\la = \epspm \rho\mu \a^\pm_\rho (Q)
\label{anat1}
\ee
whenever $Q\in\Hom(\iota\la,\iota\mu)$.

Let $a\in\Mor(M,N)$ be such that $[a]$ is a subsector
of $[\mu\bar\iota]$ for some $\mu\in\Ssys$. Hence $a\iota\in\Ssys$.
Similarly, let $\co b \in\Mor(N,M)$ be such that
$[\bar b]$ is a subsector of $[\iota{\co\nu}]$
for some ${\co\nu}\in\Ssys$.
If $T\in\Hom(\bar b,\iota{\co\nu})$ is an isometry we put
\[ \Epspm \la {\bar b} = T^* \epspm \la{\co\nu}
\a^\pm_\la (T) \,, \qquad
\Epspm {\bar b}\la = (\Epsmp \la {\bar b})^* \,. \]
Note that the definition is independent of the choice of $T$
and ${\co\nu}$ in the following sense: If also
$S\in\Hom(\bar b,\iota{\co\tau})$
is an isometry for some ${\co\tau}\in\Ssys$ then
$ST^*\in\Hom(\iota{\co\nu},\iota{\co\tau})$ and therefore
\[ \Epspm \la {\bar b} = S^*ST^*
\epspm \la{\co\nu} \a^\pm_\la (T)
= S^* \epspm \la{{\co\tau}} \a^\pm_\la (ST^*T)
= S^* \epspm \la{{\co\tau}} \a^\pm_\la (S) \,. \]
Similarly one easily checks that $\Epspm \la {\bar b}$ is unitary.

\begin{proposition}
\label{stat}
Let $\la\in\Ssys$, let $a\in\Mor(M,N)$ be such that $[a]$ is
a subsector of $[\mu\bar\iota]$ for some $\mu\in\Ssys$ and
let $\co b \in\Mor(N,M)$ be such that $[\bar b]$ is a
subsector of $[\iota{\co\nu}]$ for some ${\co\nu}\in\Ssys$.
Then we have
\be 
\epspm \la {a\iota} \in \Hom(\la a, a \a_\la^\pm) \,,\qquad
\Epspm \la{\bar b} \in \Hom(\a^\pm_\la \bar b , \bar b\lambda) \,.
\label{genstat}
\ee
\end{proposition}

\begin{proof}
The first relation in \erf{genstat} is trivial on $N$, so we
only need to show it for $v$ since $M=Nv$. Note that
$a(v)\in\Hom(a\iota,a\iota\canr)$, therefore \erf{BFE}
yields
\[ a(v) \epspm \la {a\iota} = a\iota(\epspm \la\canr)
\epspm \la {a\iota} \la(a(v)) \,, \]
hence
\[ a \circ\a^\pm_\la (v) = a\iota(\epspm \la\canr ^*) a(v) =
\epspm \la {a\iota} \la (a(v))  \epspm \la {a\iota} ^* =
 \Ad\;\epspm \la{a\iota}\circ\la\circ a (v) \,. \]
For the second relation we use the fact that
$TT^*\in\Hom(\iota{\co\nu},\iota{\co\nu})$ for
$T\in\Hom(\bar b,\iota{\co\nu})$:
\[ \bearll
\Epspm \la {\bar b} \a_\la^\pm \bar b(n)
&= T^* \epspm \la{\co\nu} \a_\la^\pm(TT^*{\co\nu}(n)T)
= T^* \epspm \la{\co\nu} \la{\co\nu}(n) \a^\pm_\la (T) \\[.4em]
&= T^* {\co\nu}\la(n) \epspm \la{\co\nu} \a^\pm_\la (T)
= \bar b \lambda (n) \Epspm \la {\bar b}
\eear \]
for all $n\in N$.
\end{proof}

Due to Prop.\ \ref{stat} we can now draw the pictures in
Fig.\ \ref{NMbraid} for the operators $\epspm \la {a\iota}$
and $\Epspm \la {\bar b}$.
%
% N-M braidings
%
\begin{figure}[htb]
\begin{center}
\unitlength 0.6mm
\begin{picture}(220,40)
\thinlines
\put(10,40){\line(0,-1){5.858}}
\put(16.5,23.5){\vector(1,-1){0}}
\put(20,20){\line(-1,1){7.071}}
\put(20,34.142){\arc{20}{2.356}{3.142}}
\Thicklines
\put(20,20){\line(1,-1){7.071}}
\put(30,0){\line(0,1){5.858}}
\put(30,0){\vector(0,-1){0}}
\put(20,5.858){\arc{20}{5.498}{0}}
\thicklines
\put(30,40){\line(0,-1){5.858}}
\put(22,22){\line(1,1){5.071}}
\put(18,18){\line(-1,-1){5.071}}
\put(20,34.142){\arc{20}{0}{0.785}}
\put(10,0){\line(0,1){5.858}}
\put(10,0){\vector(0,-1){0}}
\put(20,5.858){\arc{20}{3.142}{3.927}}
\put(5,3.5){\makebox(0,0){$a$}}
\put(36,4){\makebox(0,0){$\a_\la^+$}}
\put(5,36){\makebox(0,0){$\la$}}
\put(35,35){\makebox(0,0){$a$}}
\put(50,20){\makebox(0,0){$;$}}
\thinlines
\put(70,40){\line(0,-1){5.858}}
\put(76.5,23.5){\vector(1,-1){0}}
\put(78,22){\line(-1,1){5.071}}
\put(80,34.142){\arc{20}{2.356}{3.142}}
\Thicklines
\put(82,18){\line(1,-1){5.071}}
\put(90,0){\line(0,1){5.858}}
\put(90,0){\vector(0,-1){0}}
\put(80,5.858){\arc{20}{5.498}{0}}
\thicklines
\put(90,40){\line(0,-1){5.858}}
\put(80,20){\line(1,1){7.071}}
\put(80,20){\line(-1,-1){7.071}}
\put(80,34.142){\arc{20}{0}{0.785}}
\put(70,0){\line(0,1){5.858}}
\put(70,0){\vector(0,-1){0}}
\put(80,5.858){\arc{20}{3.142}{3.927}}
\put(65,3.5){\makebox(0,0){$a$}}
\put(96,4){\makebox(0,0){$\a_\la^-$}}
\put(65,36){\makebox(0,0){$\la$}}
\put(95,35){\makebox(0,0){$a$}}
\put(110,20){\makebox(0,0){$;$}}
\Thicklines
\put(130,40){\line(0,-1){5.858}}
\put(136.5,23.5){\vector(1,-1){0}}
\put(140,20){\line(-1,1){7.071}}
\put(140,34.142){\arc{20}{2.356}{3.142}}
\thinlines
\put(140,20){\line(1,-1){7.071}}
\put(150,0){\line(0,1){5.858}}
\put(150,0){\vector(0,-1){0}}
\put(140,5.858){\arc{20}{5.498}{0}}
\thicklines
\put(150,40){\line(0,-1){5.858}}
\put(150,40){\vector(0,1){0}}
\put(142,22){\line(1,1){5.071}}
\put(138,18){\line(-1,-1){5.071}}
\put(140,34.142){\arc{20}{0}{0.785}}
\put(130,0){\line(0,1){5.858}}
\put(140,5.858){\arc{20}{3.142}{3.927}}
\put(125,4){\makebox(0,0){$b$}}
\put(156,4){\makebox(0,0){$\la$}}
\put(125,35){\makebox(0,0){$\a_\la^+$}}
\put(155,35.5){\makebox(0,0){$b$}}
\put(170,20){\makebox(0,0){$;$}}
\Thicklines
\put(190,40){\line(0,-1){5.858}}
\put(196.5,23.5){\vector(1,-1){0}}
\put(198,22){\line(-1,1){5.071}}
\put(200,34.142){\arc{20}{2.356}{3.142}}
\thinlines
\put(202,18){\line(1,-1){5.071}}
\put(210,0){\line(0,1){5.858}}
\put(210,0){\vector(0,-1){0}}
\put(200,5.858){\arc{20}{5.498}{0}}
\thicklines
\put(210,40){\line(0,-1){5.858}}
\put(210,40){\vector(0,1){0}}
\put(200,20){\line(1,1){7.071}}
\put(200,20){\line(-1,-1){7.071}}
\put(200,34.142){\arc{20}{0}{0.785}}
\put(190,0){\line(0,1){5.858}}
\put(200,5.858){\arc{20}{3.142}{3.927}}
\put(185,4){\makebox(0,0){$b$}}
\put(216,4){\makebox(0,0){$\la$}}
\put(185,35){\makebox(0,0){$\a_\la^-$}}
\put(215,35.5){\makebox(0,0){$b$}}
\end{picture}
\end{center}
\caption{Wire diagrams for $\epsp \la {a\iota}$, $\epsm \la {a\iota}$,
$\Epsp \la {\bar b}$, $\Epsm \la {\bar b}$, respectively}
\label{NMbraid}
\end{figure}
The pictures for their conjugates $\epsmp {a\iota}\la$ and
$\Epsmp {\bar b}\la$ are as usual obtained by horizontal
reflection and inversion of arrows of the pictures in
Fig.\ \ref{NMbraid}.

\begin{lemma}
\label{Enat}
Let $\co a,\co b\in\Mor(M,N)$ be such
that $[\co a]$ and $[\co b]$ are subsectors
of $[\iota\co\mu]$ and $[\iota{\co\nu}]$ for some
$\co\mu,{\co\nu}\in\Ssys$, respectively.
Whenever $Y\in\Hom(\co a,\co b)$ we have
\[ \a_\rho^\mp(Y) \, \Epspm {\co a}\rho =
\Epspm {\co b}\rho \,  Y \,,\qquad
Y \, \Epspm \rho{\co a} =
\Epspm \rho {\co b} \, \a_\rho^\pm (Y) \,. \]
\end{lemma}

\begin{proof}
Let $S\in\Hom(\co a,\iota\co\mu)$ and
$T\in\Hom(\co b,\iota\co\nu)$
be isometries. Then
$\Epspm{\co a}\rho=\a_\rho^\mp(S)^* \epspm {\co\mu}\rho S$
and $\Epspm \rho{\co b}= T^* \epspm \rho{\co\nu} \a_\rho^\pm(T)$.
Now $TYS^*\in\Hom(\iota\co\mu,\iota\co\nu)$. Inserting this
in \erf{anat1} yields the statement.
\end{proof}

In order to establish a symmetry for
``moving crossings over trivalent vertices''
we can now state the following

\begin{proposition}
\label{intbfe}
Let $\la,\rho\in\Ssys$, let $a,b\in\Mor(M,N)$ be such
that $[a]$ and $[b]$ are subsectors
of $[\mu\bar\iota]$ and $[\nu\bar\iota]$ for some $\mu,\nu\in\Ssys$
and let $\co a,\co b \in\Mor(N,M)$ be conjugates, respectively.
Whenever $t\in\Hom(\la,a\bar b)$, $x\in\Hom(a,\la b)$
and $Y\in\Hom(\co a,\co b \la)$,
we have the intertwining braiding fusion equations
(IBFE's):
\begin{eqnarray}
\rho(t) \, \epspm \lambda\rho 
&=& \epspm {a\iota}\rho \, a(\Epspm {\bar b}\rho)
\, t \,, \label{ibfe1} \\
t \, \epspm \rho\lambda
&=& a (\Epspm \rho{\bar b}) \, \epspm \rho{a\iota} \,
\rho (t) \,, \label{ibfe2} \\
\rho(x) \, \epspm {a\iota} \rho 
&=& \epspm \la\rho \, \la(\epspm {b\iota}\rho)
\, x \,, \label{ibfe3} \\
x \, \epspm \rho {a\iota}
&=& \la (\epspm \rho{b\iota}) \, \epspm \rho\la \,
\rho (x) \,,\label{ibfe4} \\
\a_\rho^\mp(Y) \, \Epspm {\co a}\rho 
&=& \Epspm {\co b}\rho \, \co b(\epspm \la\rho)
\, Y \,, \label{ibfe5} \\
Y \, \Epspm \rho{\co a}
&=& \bar b (\epspm \rho\la) \, \Epspm \rho{\bar b} \,
\a_\rho^\pm (Y) \label{ibfe6} \,.
\end{eqnarray}
\end{proposition}

\begin{proof}
Since $[\co b]$ must be a subsector of $[\iota\co\nu]$
for $\co\nu\in\Ssys$ a conjugate of $\nu$, there is
an isometry $T\in\Hom(\bar b,\iota{\co\nu})$. Note that then
$a(T)\in\Hom(a\bar b,a\iota{\co\nu})$. Hence by naturality and
Proposition \ref{stat} we compute
\[ \bearll \epspm \rho{a\bar b} 
&= a(T^*) \epspm \rho{a\iota{\co\nu}} \rho a (T)
 = a(T^*) a (\epspm \rho{\co\nu}) \epspm \rho{a\iota}
   \rho a(T)  \\[.4em]
&= a(T^*) a (\epspm \rho {\co\nu}) a\a_\rho^\pm(T)
   \epspm \rho{a\iota}
 = a( \Epspm \rho{\bar b}) \epspm \rho{a\iota} \,,
\eear \]
and hence also $\epspm {a\bar b}\rho = \epspm {a\iota}\rho
a(\Epspm {\bar b}\rho$. We also obtain
$\epspm {\la b\iota}\rho=\epspm\la\rho \la(\epspm{b\iota}\rho)$
and $\epspm \rho{\la b\iota}=\la(\epspm \rho{b\iota}) \epspm \rho\la$
by \erf{comp}. Note that $x\in\Hom(a\iota,\la b\iota)$ by
restriction. Eqs.\ (\ref{ibfe1})-(\ref{ibfe4}) follow now
by naturality, \erf{nat}. Next, we note that
$T\in\Hom(\bar b \la,\iota{\co\nu}\la)$, and hence
$\Epspm \rho{\co b \la}=T^* \epspm \rho{\co\nu\la} \a_\rho^\pm(T)$.
Therefore
\[ \bearll
 \Epspm \rho{\co b \la} &= T^* \co\nu(\epspm \rho\la)
 \epspm\rho{\co\nu}
\a_\rho^\pm(T) = \co b(\epspm \rho\la) T^* \epspm \rho{\co\nu}
\a_\rho^\pm(T) \\[.4em]
& =  \co b(\epspm \rho\la) \Epspm \rho{\co b} \,,
\eear \]
and hence also $\Epspm {\co b \la}\rho = \Epspm {\co b}\rho
\co b (\epspm \la\rho)$. Now Eqs.\ (\ref{ibfe5}) and (\ref{ibfe6})
follow from Lemma \ref{Enat}.
\end{proof}

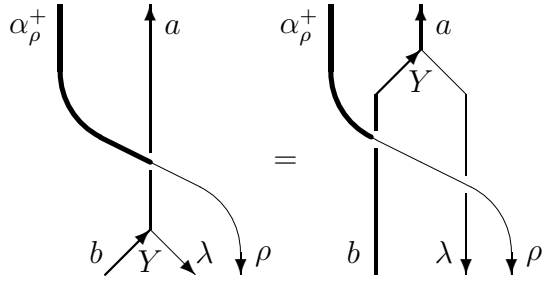
\begin{figure}[htb]
\begin{center}
\unitlength 0.6mm
\begin{picture}(120,60)
\Thicklines
\put(26.180,45){\arc{32.361}{2.034}{3.142}}
\put(10,60){\line(0,-1){15}}
\put(30,25){\line(-2,1){11.1}}   %11.1 for 11.056
\thinlines
\put(33.820,5){\arc{32.361}{5.176}{6.283}}
\put(50,5){\vector(0,-1){5}}
\put(30,25){\line(2,-1){11.1}}
\thicklines
\put(30,27){\vector(0,1){33}}
\put(30,23){\line(0,-1){13}}
\put(20,0){\vector(1,1){10}}
\thinlines
\put(30,10){\vector(1,-1){10}}
\put(3,55){\makebox(0,0){$\a_\rho^+$}}
\put(30,3){\makebox(0,0){$Y$}}
\put(35,55){\makebox(0,0){$a$}}
\put(18,5){\makebox(0,0){$b$}}
\put(42,5){\makebox(0,0){$\la$}}
\put(55,4){\makebox(0,0){$\rho$}}
\put(60,25){\makebox(0,0){$=$}}
\Thicklines
\put(86.180,45){\arc{32.361}{2.034}{3.142}}
\put(70,60){\line(0,-1){15}}
\thinlines
\put(93.820,5){\arc{32.361}{5.176}{6.283}}
\put(110,5){\vector(0,-1){5}}
\put(80,30){\line(-2,1){1.1}}
\put(100,20){\line(2,-1){1.1}}
\thicklines
\put(90,50){\vector(0,1){10}}
\put(80,40){\vector(1,1){10}}
\put(80,40){\line(0,-1){8}}
\put(80,10){\line(0,1){18}}
\put(80,10){\line(0,-1){10}}
\thinlines
\put(100,40){\line(-1,1){10}}
\put(100,40){\line(0,-1){18}}
\put(100,10){\line(0,1){8}}
\put(100,10){\vector(0,-1){10}}
\put(90,25){\line(-2,1){10}}
\put(90,25){\line(2,-1){10}}
\put(90,43){\makebox(0,0){$Y$}}
\put(95,55){\makebox(0,0){$a$}}
\put(75,5){\makebox(0,0){$b$}}
\put(95,5){\makebox(0,0){$\la$}}
\put(115,4){\makebox(0,0){$\rho$}}
\put(63,55){\makebox(0,0){$\a_\rho^+$}}
\end{picture}
\end{center}
\caption{The sixth intertwining braiding fusion equation (overcrossings)}
\label{IBFE2}
\end{figure}

These IBFE's can be nicely visualized in diagrams. We display
\erf{ibfe1} in Fig.\ \ref{IBFE1}
%
% IBFE (1)
\begin{figure}[htb]
\begin{center}
\unitlength 0.6mm
\begin{picture}(120,60)
\thinlines
\put(26.180,5){\arc{32.361}{3.142}{4.249}}
\put(10,5){\vector(0,-1){5}}
\put(28,24){\line(-2,-1){9.1}}
\put(33.820,45){\arc{32.361}{0}{1.107}}
\put(50,60){\line(0,-1){15}}
\put(32,26){\line(2,1){9.1}}
\put(30,60){\line(0,-1){35}}
\put(30,25){\vector(0,-1){15}}
\thicklines
\put(40,0){\vector(-1,1){10}}
\put(30,10){\vector(-1,-1){10}}
\put(30,5){\makebox(0,0){$t$}}
\put(35,55){\makebox(0,0){$\la$}}
\put(18,4.5){\makebox(0,0){$a$}}
\put(42,5){\makebox(0,0){$b$}}
\put(5,3.5){\makebox(0,0){$\rho$}}
\put(60,25){\makebox(0,0){$=$}}
\thinlines
\put(86.180,5){\arc{32.361}{3.142}{4.199}}    %4.199 for 4.249
\put(70,5){\vector(0,-1){5}}
\put(93.820,45){\arc{32.361}{0}{1.057}} %1.057 for 1.107
\put(110,60){\line(0,-1){15}}
\put(90,60){\vector(0,-1){10}}
\thicklines
\put(100,40){\vector(-1,1){10}}
\put(90,50){\line(-1,-1){10}}
\put(90,45){\makebox(0,0){$t$}}
\put(80,40){\line(0,-1){20}}
\put(100,40){\line(0,-1){10}}
\put(80,10){\line(0,1){10}}
\put(100,10){\line(0,1){20}}
\put(80,10){\vector(0,-1){10}}
\put(100,10){\line(0,-1){10}}
\Thicklines
\put(90,25){\line(2,1){8}}
\put(90,25){\line(-2,-1){8}}
\put(95,55){\makebox(0,0){$\la$}}
\put(85,4.5){\makebox(0,0){$a$}}
\put(105,5){\makebox(0,0){$b$}}
\put(65,3.5){\makebox(0,0){$\rho$}}
\put(90,19){\makebox(0,0){$\a_\rho^-$}}
\end{picture}
\end{center}
\caption{The first intertwining braiding fusion equation (overcrossings)}
\label{IBFE1}
\end{figure}
and \erf{ibfe6} in Fig.\ \ref{IBFE2}, both for overcrossings.
%
% IBFE (2)
We leave the remaining diagrams as a straightforward exercise
to the reader. Note that the IBFE's give us the freedom to
move wires with label $\rho$ and $\a_\rho^\pm$ freely over
trivalent vertices which involve one $N$-$N$ wire and two
$N$-$M$ wires. Unitarity of operators $\Epspm \la{\co b}$
yields a ``vertical Reidemeister move of type II'' similar
to Fig.\ \ref{Reid2}. We can now also easily elaborate the
rotation behavior of mixed crossings displayed in
Fig.\ \ref{NMbraid} (and consequently their conjugates).
Crucial for this is the fact that
$R_{\a_\la^\pm}=\iota(r_\la)\equiv r_\la$ and
${\co R}_{\a_\la^\pm}=\iota({\co r}_\la)\equiv {\co r}_\la$
can be used as R-isometries for the $\a$-induced morphisms as
$R_{\a_\la^\pm}\in\Hom(\id_M,
\overline{\a^\pm_\la}\a^\pm_\la)$
and
${\co R}_{\a_\la^\pm}\in\Hom
(\id_M,\a^\pm_\la\overline{\a^\pm_\la})$
satisfy
$\a^\pm_\la(R_{\a_\la^\pm})^*{\co R}_{\a_\la^\pm}
=d_\la^{-1}\bfe_M$ and
$\overline{\a^\pm_\la}({\co R}_{\a_\la^\pm})^*
\co R_{\a_\la^\pm} =d_\la^{-1}\bfe_M$ and
$d_{\a^\pm_\la}=d_\la$. First we notice that
we have
\[ \epspm \la{a\iota} = d_\la \, {\co r}_\la \,
\la(\epsmp {a\iota}{\co\la}) \, \la a(r_\la) \]
by \erf{rotcreq}. Now let $R_a\in\Hom(\id_M,\co a a)$
and ${\co r}_a\in\Hom(\id_N, a\co a)$ be isometries
such that $a(R_a)^*{\co r}_a=d_a^{-1} \bfe_N$ and
$\co a({\co r}_a)^*R_a=d_a^{-1}$, and otherwise we keep
the notations as in Prop.\ \ref{intbfe}. From
\erf{ibfe2} we obtain
$a(\Epsmp {\co a}\la){\co r}_a = \epspm \la{a\iota} \la({\co r}_a)$.
Hence we have
\[ \bearll
\epspm \la{a\iota} &= d_a \, \epspm \la{a\iota}\,
\la a(R_a)^* \la({\co r}_a) = d_a \, a\a_\la^\pm (R_a)^* \,
\epspm \la{a\iota} \, \la({\co r}_a) \\[.4em]
&=  d_a \, a\a_\la^\pm (R_a)^* \,
a(\Epsmp {\co a}\la) \, {\co r}_a \,.
\eear \]
Next we compute, using again \erf{rotcreq},
\[ \bearll
\Epspm \la{\co b} &= T^*\,\epspm \la{\co\nu} \, \a_\la^\pm(T)
= d_\la \, T^*\,{\co r}_\la \la(\epsmp {\co\nu}{\co\la})
\la\co\nu(r_\la) \, \a_\la^\pm (T) \\[.4em]
&= d_\la \, {\co r}_\la^* \a_\la^\pm( \a_{\co\la}^\pm(T)^*
\epsmp {\co\nu}{\co\la}T) \a_\la^\pm \co b(r_\la)
=  d_\la \, {\co r}_\la^* \a_\la^\pm( \Epsmp {\co b}{\co\la})
\a_\la^\pm \co b(r_\la) \,.
\eear \]
Finally, as \erf{ibfe2} yields
${\co r}_a^* a(\Epspm \la{\co a}=\la({\co r}_a)^*
\epsmp {a\iota}\la$, we obtain
\[ \Epspm \la{\co a}= d_a \, \co a({\co r}_a)^* \, \co a a
(\Epspm \la{\co a}) R_a = d_a \, \co a\la({\co r}_a)^* \,
\co a(\epsmp {a\iota}\la)  R_a \,. \]
Drawing for $R_{\a_\la^\pm}=\iota(r_\la)$ and
${\co R}_{\a_\la^\pm}=\iota({\co r}_\la)$ caps
of the wires $\a_\la^\pm$,
these relations yield graphically the analogues of
Fig.\ \ref{rotcross}. We conclude that we can
include the crossings of Fig.\ \ref{NMbraid} consistently
in our ``rotation covariant'' graphical framework.

\section{Double Triangle Algebras for Subfactors}
\label{sec-dta}

We now formulate Ocneanu's construction \cite{O7}
for a subfactor with finite index and finite depth rather
than for bi-unitary connections and bimodules arising from
Goodman-de la Harpe-Jones subfactors associated to
A-D-E Dynkin diagrams in order to apply it in a more
general context. From now on we work with $N\subset M$
satisfying the following

\begin{assumption}
\label{set-fin}{\rm
Let $N\subset M$ be a type III subfactor with finite index.
We assume that we have a system of endomorphisms
$\NXN\subset\Mor(N,N)\equiv\End(N)$ in the sense of
Definition \ref{system} such that for the injection map
$\iota:N\to M$, the sector $[\canr]=[\bar\iota \iota]$
decomposes into a sum of sectors of morphisms in $\NXN$.
We choose sets of morphisms $\NXM\subset\Mor(M,N)$,
$\MXN\subset\Mor(N,M)$ and $\MXM\subset\Mor(M,M)\equiv\End(M)$
consisting of representative endomorphisms of irreducible
subsectors of sectors of the form $[\la\co\iota]$,
$[\iota\la]$ and $[\iota\la\co\iota]$, $\la\in\NXN$,
respectively. (We may and do choose $\id_M$ in $\MXM$
as the endomorphism representing the trivial sector.)
We also assume that $\NXN$ is finite. Consequently, the set
$\cX= \NXN \sqcup \NXM \sqcup \MXN \sqcup \MXM$ is finite.
}\end{assumption}
Note that Assumption \ref{set-fin} implies that representative
morphisms for all irreducible sectors appearing in 
decompositions of powers $[\ga^k]$ ($[\canr^k]$)
of Longo's (dual) canonical endomorphism are contained in
$\MXM$ ($\NXN$). In other words, the set $\cX$ contains
at least the morphisms corresponding to the (equivalence classes
of) bimodules arising from this subfactor through
the Jones tower, and therefore
we may call an $\cX$ which does not contain any other morphisms
a {\sl minimal choice}. We conclude that
finiteness of $\NXN$ in Assumption in \ref{set-fin} automatically
implies that the subfactor $N\subset M$ has finite depth.
We used sectors instead of bimodules in view of our
``identification'' of chiral generators with $\a$-induced
sectors below. Therefore we need a sector approach in order to
define $\a$-induction since its definition involves
$\co\iota^{\,-1}$, and hence we work with factors of type III.
(We do not need hyperfiniteness of $M$ for our purposes.)

We now use the graphical calculus presented in Section \ref{sec-graphcal}.
In the graphical method of \cite{O3} (and \cite[Chapter 12]{EK3}),
factors, bimodules (morphisms), and intertwiners are represented with
trivalent vertices, edges, and triangles, respectively, and this is
where the name ``double triangle algebra'' comes from. However,
here (as in \cite{O6,O7}) these three kinds of objects are represented
by regions, wires, and trivalent vertices, respectively, though
the labels for regions are omitted for notational simplicity.

For $\cal X$ in Assumption \ref{set-fin},
we define the {\sl double triangle algebra}
$\dta$ with two multiplications $*_h$ and $*_v$ as follows.
As a linear space, we set
\[ \dta=\bigoplus_{a,b,c,d\in\NXM}
\Hom(a \bar b, c \bar d) \,. \]
This is a finite dimensional complex linear space.
An element in $\dta$ is presented graphically as in
Fig.\  \ref{dta1} under the interpretation in Section
\ref{sec-graphcal} with the convention of reading the diagram from
the top to the bottom.
%
%double triangle algebra (1)
\thinlines
\begin{figure}[htb]
\begin{center}
\unitlength 0.6mm
\begin{picture}(40,40)
\thicklines
\put(10,10){\line(1,0){20}}
\put(10,30){\line(1,0){20}}
\thinlines
\put(20,30){\vector(0,-1){20}}
\put(7,30){\makebox(0,0){$a$}}
\put(33,30){\makebox(0,0){$b$}}
\put(7,10){\makebox(0,0){$c$}}
\put(33,10){\makebox(0,0){$d$}}
\put(17,20){\makebox(0,0){$\la$}}
\put(20,34){\makebox(0,0){$s^*$}}
\put(20,6){\makebox(0,0){$t$}}
\end{picture}
\end{center}
\caption{An element in $\protect\dta$}
\label{dta1}
\end{figure}
(A general element in $\dta$ is a linear combination of this type
of element.) We can interpret the same diagram with the
convention of reading the diagram from the
left to the right or, equivalently, keeping the top-to-bottom
convention but putting the diagram in a suitable Frobenius
annulus. Then the resulting intertwiner is in 
\[ \dtap=\bigoplus_{a,b,c,d\in\NXM}
\Hom(\bar c a,\bar d b) \,. \]
The isomorphism of these two spaces is given by
application of two Frobenius rotations,
and we can use this isomorphism
to identify $\dta$ and $\dtap$.
By our convention of the normalization in Section \ref{sec-graphcal},
the diagram of Fig.\  \ref{dta1} represents an element
$d_a^{1/4} d_b^{1/4} d_c^{1/4} d_d^{1/4} d_\la^{-1/2}ts^*$
in the block
$\Hom(a\bar b, c\bar d)$, where $s\in\Hom (\la,a\bar b)$
and $t\in\Hom(\la, c\bar d)$ are isometries and $\la\in\NXN$.
Similarly we may use elements in $\dta$ which are graphically
represented as in Fig.\  \ref{dta2} with
isometries $S\in\Hom(\beta,\bar c a)$, $T\in\Hom(\beta,\bar d b)$
and $\beta\in\MXM$.
%
%double triangle algebra (2)
\thinlines
\begin{figure}[htb]
\begin{center}
\unitlength 0.6mm
\begin{picture}(40,40)
\thicklines
\put(10,10){\line(0,1){20}}
\put(30,10){\line(0,1){20}}
\Thicklines
\put(10,20){\vector(1,0){20}}
\put(10,33){\makebox(0,0){$a$}}
\put(30,33){\makebox(0,0){$b$}}
\put(10,7){\makebox(0,0){$c$}}
\put(30,7){\makebox(0,0){$d$}}
\put(20,24){\makebox(0,0){$\beta$}}
\put(5,20){\makebox(0,0){$S^*$}}
\put(35,20){\makebox(0,0){$T$}}
\end{picture}
\end{center}
\caption{An element in $\protect\dta$}
\label{dta2}
\end{figure}
Note that elements of the form in Fig.\ \ref{dta1},
or equivalently of the form in Fig.\ \ref{dta2},
span $\dta$ linearly.

Our graphical convention is as follows. We use thin, thick, and
very thick wires for $N$-$N$ morphisms,
$N$-$M$ morphisms, and $M$-$M$ morphisms, respectively,
analogous to the convention \cite{O7}.  We call
them $N$-$N$ wires, and so on.  We label
$N$-$N$ morphisms with Greek letters $\la,\mu,\nu,\dots$,
$N$-$M$ morphisms with Roman letters $a,b,c,d,\dots$, and
$M$-$M$ morphisms with Greek letters $\beta,\beta', \beta'',\dots$.
We orient $N$-$N$ or $M$-$M$ wires but we put no
orientations on $N$-$M$ wires since it is clear from
the context whether we mean an $N$-$M$ morphism
$a$ or an $M$-$N$ morphism $\bar a$.  We simply put
a label $a$ for an unoriented thick wire for both.
Note that, whatever we consider, $\dta$ or $\dtap$,
the same intertwiner (as an operator) may appear in
different blocks of the double triangle algebra, e.g.\
the identity $\id_N$ is an element in any
$\Hom(a \bar b,a \bar b)$, $a,b\in\NXM$. The graphical
notation is particularly useful in order to avoid this kind
of confusion because diagrams as in Figs.\ \ref{dta1}
and \ref{dta2} always specify also the associated
block in addition to the intertwiner as an operator.

The {\sl horizontal product} $*_h$ on $\dta$
is defined as in Fig.\  \ref{prod}.
%
%horizontal product
\thinlines
\begin{figure}[htb]
\begin{center}
\unitlength 0.6mm
\begin{picture}(150,40)
\thicklines
\put(10,10){\line(1,0){20}}
\put(10,30){\line(1,0){20}}
\put(50,10){\line(1,0){20}}
\put(50,30){\line(1,0){20}}
\put(110,10){\line(1,0){30}}
\put(110,30){\line(1,0){30}}
\thinlines
\put(20,30){\vector(0,-1){20}}
\put(60,30){\vector(0,-1){20}}
\put(120,30){\vector(0,-1){20}}
\put(130,30){\vector(0,-1){20}}
\put(7,30){\makebox(0,0){$a$}}
\put(33,30){\makebox(0,0){$b$}}
\put(7,10){\makebox(0,0){$c$}}
\put(33,10){\makebox(0,0){$d$}}
\put(47,30){\makebox(0,0){$a'$}}
\put(73,30){\makebox(0,0){$b'$}}
\put(47,10){\makebox(0,0){$c'$}}
\put(73,10){\makebox(0,0){$d'$}}
\put(40,20){\makebox(0,0){$*_h$}}
\put(93,20){\makebox(0,0){$=\del b{a'} \del d{c'}$}}
\put(107,30){\makebox(0,0){$a$}}
\put(143,30){\makebox(0,0){$b'$}}
\put(107,10){\makebox(0,0){$c$}}
\put(143,10){\makebox(0,0){$d'$}}
\put(125,14){\makebox(0,0){$d$}}
\put(125,26){\makebox(0,0){$b$}}
\put(17,20){\makebox(0,0){$\la$}}
\put(57,20){\makebox(0,0){$\mu$}}
\put(117,20){\makebox(0,0){$\la$}}
\put(133,20){\makebox(0,0){$\mu$}}
\put(20,34){\makebox(0,0){$s^*$}}
\put(20,6){\makebox(0,0){$t$}}
\put(61,34){\makebox(0,0){${s'}^*$}}
\put(60,6){\makebox(0,0){$t'$}}
\put(120,34){\makebox(0,0){$s^*$}}
\put(120,6){\makebox(0,0){$t$}}
\put(132,34){\makebox(0,0){${s'}^*$}}
\put(130,6){\makebox(0,0){$t'$}}
\end{picture}
\end{center}
\caption{The horizontal product $*_h$ on $\protect\dta$}
\label{prod}
\end{figure}
The meaning of the right-hand side is as follows. The
product is by definition zero if the labels of the
open ends of the wires facing each other do not match.
If they match, we glue the wires of the two diagrams
together as in Fig.\  \ref{prod} and interpret it
as an intertwiner. It belongs to the block
of the double triangle algebra which is specified by
the four remaining open ends of the new diagram.
This is a horizontal version of the composition
of intertwiners described in Section \ref{sec-graphcal}.

We also can represent this horizontal product in terms
of elements in Fig.\  \ref{dta2}.  This is described in
Fig.\  \ref{prod2}, because the
convention of Section \ref{sec-graphcal}
means that this product is just the  composition of the
intertwiners in $\dtap$,
and this composition is realized by taking
the inner product of the two intertwiners
in the right-hand side in Fig.\  \ref{prod2}.
%
%Horizontal product (2)
\begin{figure}[htb]
\begin{center}
\unitlength 0.7mm
\begin{picture}(200,40)
\thicklines
\put(10,10){\line(0,1){20}}
\put(30,10){\line(0,1){20}}
\put(60,10){\line(0,1){20}}
\put(80,10){\line(0,1){20}}
\put(170,10){\line(0,1){20}}
\put(190,10){\line(0,1){20}}
\Thicklines
\put(10,20){\line(1,0){20}}
\put(20,20){\vector(1,0){0}}
\put(60,20){\line(1,0){20}}
\put(70,20){\vector(1,0){0}}
\put(170,20){\line(1,0){20}}
\put(180,20){\vector(1,0){0}}
\put(10,7){\makebox(0,0){$c$}}
\put(30,7){\makebox(0,0){$d$}}
\put(60,7){\makebox(0,0){$c'$}}
\put(80,7){\makebox(0,0){$d'$}}
\put(170,7){\makebox(0,0){$c$}}
\put(190,7){\makebox(0,0){$d'$}}
\put(10,33){\makebox(0,0){$a$}}
\put(30,33){\makebox(0,0){$b$}}
\put(60,33){\makebox(0,0){$a'$}}
\put(80,33){\makebox(0,0){$b'$}}
\put(170,33){\makebox(0,0){$a$}}
\put(190,33){\makebox(0,0){$b'$}}
\put(20,24){\makebox(0,0){$\beta$}}
\put(70,24){\makebox(0,0){$\beta'$}}
\put(180,24){\makebox(0,0){$\beta$}}
\put(6,20){\makebox(0,0){$S^*$}}
\put(34,20){\makebox(0,0){$T$}}
\put(56,20){\makebox(0,0){${S'}^*$}}
\put(84,20){\makebox(0,0){$T'$}}
\put(166,20){\makebox(0,0){$S^*$}}
\put(194,20){\makebox(0,0){$T'$}}
\put(45,20){\makebox(0,0){$*_h$}}
\put(125,20){\makebox(0,0){$=\displaystyle\del b{a'} \del d{c'}
\del \beta{\beta'}\sqrt{\frac{d_b d_d}{d_\beta}} \lan S',T \ran$}}
\end{picture}
\end{center}
\caption{The horizontal product presented in another way}
\label{prod2}
\end{figure}

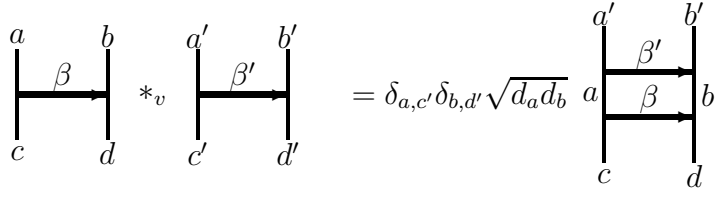
\begin{figure}[htb]
\begin{center}
\unitlength 0.6mm
\begin{picture}(170,50)
\thicklines
\put(10,15){\line(0,1){20}}
\put(30,15){\line(0,1){20}}
\put(50,15){\line(0,1){20}}
\put(70,15){\line(0,1){20}}
\put(140,10){\line(0,1){30}}
\put(160,10){\line(0,1){30}}
\Thicklines
\put(10,25){\vector(1,0){20}}
\put(50,25){\vector(1,0){20}}
\put(140,20){\vector(1,0){20}}
\put(140,30){\vector(1,0){20}}
\put(10,38){\makebox(0,0){$a$}}
\put(30,38){\makebox(0,0){$b$}}
\put(10,12){\makebox(0,0){$c$}}
\put(30,12){\makebox(0,0){$d$}}
\put(50,38){\makebox(0,0){$a'$}}
\put(70,38){\makebox(0,0){$b'$}}
\put(50,12){\makebox(0,0){$c'$}}
\put(70,12){\makebox(0,0){$d'$}}
\put(140,43){\makebox(0,0){$a'$}}
\put(160,43){\makebox(0,0){$b'$}}
\put(140,7){\makebox(0,0){$c$}}
\put(160,7){\makebox(0,0){$d$}}
\put(137,25){\makebox(0,0){$a$}}
\put(163,25){\makebox(0,0){$b$}}
\put(20,29){\makebox(0,0){$\beta$}}
\put(60,29){\makebox(0,0){$\beta'$}}
\put(150,34){\makebox(0,0){$\beta'$}}
\put(150,24){\makebox(0,0){$\beta$}}
\put(40,25){\makebox(0,0){$*_v$}}
\put(108,25){\makebox(0,0){$=\hphantom{}\del a{c'}\del b{d'}
\sqrt{d_a d_b}$}}
\end{picture}
\end{center}
\caption{The vertical product $*_v$ in $\protect\dta$}
\label{prod3}
\end{figure}

We similarly define the {\sl vertical product} $*_v$ on $\dta$
by composing two diagrams vertically, but with extra coefficients
as in Fig.\  \ref{prod3}.
%
%Vertical product
\thinlines
The meaning of the right-hand side is as before. 
Note that the definitions of horizontal and vertical products are not
completely symmetric due to the extra coefficients we chose.
This choice is somewhat arbitrary but it just
turns out to be useful for our purposes. Namely, with this
definition of the products, the minimal central projections of
$(\dta,*_h)$ have simple and useful composition rules with respect
to the vertical product $*_v$, see Theorem \ref{double2} below.
We clearly also have a $*$-structure for the horizontal product
obtained by vertical reflection of the diagram, adjoining labels
for trivalent vertices and reversing orientations of wires.
Analogously, a $*$-structure for the vertical product comes from
horizontal reflection. The basic idea is that the 90-degree
rotation is something like a ``Fourier transform'' which transforms
the two products into each other, similar to the situation of the
group algebra of a finite or compact group. 

For each $\beta,\la,a,b$ we choose orthonormal bases of isometries
$T_{\bar b,a}^{\beta;i} \in\Hom(\beta, \bar b a)$,
$i=1,2,...,N_{\bar b,a}^\beta$, and
$t_{a,\bar b}^{\la;j} \in\Hom(\la, a \bar b)$
$j=1,2,...,N_{a\bar b}^\la$, so that
\be
\sum_{\beta\in\MXM} \sum_{i=1}^{N_{\bar b,a}^\beta}
T_{\bar b,a}^{\beta;i}(T_{\bar b,a}^{\beta;i})^*=\bfe_M
\qquad\mbox{and}\qquad
\sum_{\la\in\NXN} \sum_{j=1}^{N_{a,\bar b}^\la}
t_{a,\bar b}^{\la;j}(t_{a,\bar b}^{\la;j})^*=\bfe_N
\label{totsum}
\ee
for all $a,b\in\NXM$. Then it is easy to see that the elements in
Fig.\ \ref{munit}
%
% matrix unit in dta
\thinlines
\begin{figure}[htb]
\begin{center}
\unitlength 0.6mm
\begin{picture}(232,40)
\thicklines
\put(85,10){\line(0,1){20}}
\put(105,10){\line(0,1){20}}
\Thicklines
\put(85,20){\vector(1,0){20}}
\put(85,33){\makebox(0,0){$a$}}
\put(105,34){\makebox(0,0){$b$}}
\put(85,6){\makebox(0,0){$c$}}
\put(105,7){\makebox(0,0){$d$}}
\put(73,20){\makebox(0,0){$(T_{\bar c,a}^{\beta;i})^*$}}
\put(115,19){\makebox(0,0){$T_{\bar d,b}^{\beta;j}$}}
\put(95,24){\makebox(0,0){$\beta$}}
\put(30,20){\makebox(0,0){$e_{\beta;c,a,i}^{d,b,j}=
\displaystyle\frac{\sqrt{d_\beta}}{\sqrt[4]{d_a d_b d_c d_d}}$}}
\put(125,16){\makebox(0,0){,}}
\thicklines
\put(210,10){\line(1,0){20}}
\put(210,30){\line(1,0){20}}
\thinlines
\put(220,30){\vector(0,-1){20}}
\put(207,33){\makebox(0,0){$a$}}
\put(233,34){\makebox(0,0){$b$}}
\put(207,6){\makebox(0,0){$c$}}
\put(233,7){\makebox(0,0){$d$}}
\put(220,36){\makebox(0,0){$(t_{a, \bar b}^{\la;i})^*$}}
\put(220,4){\makebox(0,0){$t_{c, \bar d}^{\la;j}$}}
\put(224,20){\makebox(0,0){$\la$}}
\put(170,20){\makebox(0,0){$f_{\la;c,d,j}^{a,b,i}=
{\displaystyle\sqrt{\frac{d_\la}{d_a d_b d_c d_d}}}$}}
\end{picture}
\end{center}
\caption{Matrix units $e_{\beta;c,a,i}^{d,b,j}$ for
$(\protect\dta,*_h)$ and $f_{\la;c,d,j}^{a,b,i}$ for
$(\protect\dta,*_v)$}
\label{munit}
\end{figure}
form bases of $\dta$ which constitute complete systems of
matrix units $(\dta, *_h)$ respectively $(\dta,*_v)$.
Thus for each of the two multiplications the double triangle
algebra is a direct sum of full matrix algebras.
The two different bases are transformed into each other by
a unitary transformation with coefficients given by
the $6j$-symbols for subfactors of \cite{O3} (see
\cite[Chapter 12]{EK3} for the basic properties of
``quantum $6j$-symbols''), but this will not be
exploited here.

\begin{definition}\label{defe-beta}
{\rm For each $\beta\in\MXM$ we define an element
$e_\beta=\sum_{a,b,i} e_{\beta;b,a,i}^{b,a,i}\in\dta$.
Graphically, this element is given by the left-hand
side in Fig.\  \ref{e-beta}. We use the convention
shown on the right-hand side in Fig.\  \ref{e-beta}
to represent this element.
}\end{definition}
%
%M-M sector = central proj
\thinlines
\begin{figure}[htb]
\begin{center}
\unitlength 0.6mm
\begin{picture}(173,40)
\thicklines
\put(60,10){\line(0,1){20}}
\put(80,10){\line(0,1){20}}
\Thicklines
\put(60,20){\vector(1,0){20}}
\put(60,33){\makebox(0,0){$a$}}
\put(80,33){\makebox(0,0){$a$}}
\put(60,7){\makebox(0,0){$b$}}
\put(80,7){\makebox(0,0){$b$}}
\put(48,20){\makebox(0,0){$(T_{\bar b,a}^{\beta;i})^*$}}
\put(90,19){\makebox(0,0){$T_{\bar b,a}^{\beta;i}$}}
\put(70,24){\makebox(0,0){$\beta$}}
\put(20,18){\makebox(0,0){$\displaystyle
\sum_{a,b,i} \sqrt{\frac{d_\beta}{d_a d_b}}$}}
\put(120,19){\makebox(0,0){$=:$}}
\thicklines
\put(150,10){\line(0,1){20}}
\put(170,10){\line(0,1){20}}
\Thicklines
\put(150,20){\vector(1,0){20}}
\put(150,33){\makebox(0,0){$a$}}
\put(170,33){\makebox(0,0){$a$}}
\put(150,7){\makebox(0,0){$b$}}
\put(170,7){\makebox(0,0){$b$}}
\thinlines
\put(170,20){\arc{5}{1.57}{4.71}}
\put(150,20){\arc{5}{4.71}{7.85}}
\put(160,24){\makebox(0,0){$\beta$}}
\put(137,16){\makebox(0,0){$\displaystyle \sum_{a,b}$}}
\end{picture}
\end{center}
\caption{The minimal central projection $e_\beta$}
\label{e-beta}
\end{figure}
Due to the summation over $i=1,2,...,N_{\bar b,a}^\beta$,
the definition is independent of the choice
of the intertwiner bases as different orthonormal bases
are related by a unitary matrix.
We will use such a graphical convention whenever we have
a sum over internal ``fusion channels'' of two corresponding
trivalent vertices together with prefactors
which renormalize the trivalent
vertices to isometries. Note that we obtain a prefactor, as
displayed in Fig.\ \ref{trick} for an example, when we turn
around the small arcs at trivalent vertices.
%
%turning around small arcs
\thinlines
\begin{figure}[htb]
\begin{center}
\unitlength 0.6mm
\begin{picture}(180,20)
\thicklines
\put(0,10){\line(1,0){15}}
\put(15,10){\line(1,-1){10}}
\dottedline{2}(25,0)(45,0)
\put(55,10){\line(-1,-1){10}}
\put(55,10){\line(1,0){15}}
\thinlines
\put(15,10){\line(1,1){10}}
\dottedline{2}(25,20)(45,20)
\put(55,10){\line(-1,1){10}}
\put(37,20){\vector(1,0){0}}
\put(5,15){\makebox(0,0){$a$}}
\put(65,15){\makebox(0,0){$a$}}
\put(35,5){\makebox(0,0){$b$}}
\put(35,15){\makebox(0,0){$\la$}}
\put(15,10){\arc{5}{3.142}{0.785}}
\put(55,10){\arc{5}{2.356}{0}}
\put(93,10){\makebox(0,0){$=\;\displaystyle\frac{d_\la}{d_b}$}}
\thicklines
\put(110,10){\line(1,0){15}}
\put(125,10){\line(1,-1){10}}
\dottedline{2}(135,0)(155,0)
\put(165,10){\line(-1,-1){10}}
\put(165,10){\line(1,0){15}}
\thinlines
\put(125,10){\line(1,1){10}}
\dottedline{2}(135,20)(155,20)
\put(165,10){\line(-1,1){10}}
\put(147,20){\vector(1,0){0}}
\put(115,15){\makebox(0,0){$a$}}
\put(175,15){\makebox(0,0){$a$}}
\put(145,5){\makebox(0,0){$b$}}
\put(145,15){\makebox(0,0){$\la$}}
\put(125,10){\arc{5}{5.498}{3.142}}
\put(165,10){\arc{5}{0}{3.927}}
\end{picture}
\end{center}
\caption{Turning around small arcs yields a prefactor}
\label{trick}
\end{figure}
Here the dotted parts mean that there might be expansions as
given in the following lemma or later even be braiding operators
in between; it is just important that the small arcs at
corresponding trivalent vertices denote the same summation over
internal fusion channels.

\begin{lemma}
\label{idexp}
The identity of Fig.\ \ref{exp-id} holds. Analogous identities
hold if $a,b,\beta$ are replaced by wires of other type
(in a compatible way).
\end{lemma}
%
%expansion with sigma
\thinlines
\begin{figure}[htb]
\begin{center}
\unitlength 0.6mm
\begin{picture}(90,40)
\thicklines
\put(10,10){\line(1,0){20}}
\put(10,30){\line(1,0){20}}
\put(60,10){\line(0,1){20}}
\put(80,10){\line(0,1){20}}
\Thicklines
\put(60,20){\vector(1,0){20}}
\put(20,34){\makebox(0,0){$a$}}
\put(20,6){\makebox(0,0){$b$}}
\put(60,33){\makebox(0,0){$a$}}
\put(80,33){\makebox(0,0){$a$}}
\put(60,7){\makebox(0,0){$b$}}
\put(80,7){\makebox(0,0){$b$}}
\put(70,24){\makebox(0,0){$\beta$}}
\thinlines
\put(80,20){\arc{5}{1.57}{4.71}}
\put(60,20){\arc{5}{4.71}{7.85}}
\put(44,18){\makebox(0,0){$\displaystyle =\sum_\beta$}}
\end{picture}
\end{center}
\caption{The identity with expansion using $\beta$}
\label{exp-id}
\end{figure}

\begin{proof}
With the normalization convention as in Fig.\ \ref{e-beta},
this is just the expansion of the identity in \erf{totsum},
and this certainly holds as well using similar expansions with
other intertwiner bases.
\end{proof}

Note that the identity in Fig.\ \ref{exp-id} may, for example, also
appear rotated by 90 degrees as we can put the left- and right-hand sides
in some Frobenius annulus as described in Subsect.\ \ref{graph2-prelim}.

As we have already indicated,
the horizontal product is essentially the composition of
intertwiners in $\dtap$. The main point of the double
triangle algebra is the following.
Suppose we have complete information on the fusion
rules of $N$-$N$, $N$-$M$,
$M$-$N$ morphisms in $\cX$ and their $6j$-symbols.
We can define the algebra $\dta$ in terms of matrix elements
$f_{c,d,j}^{\la;a,b,i}$ and determine their composition with
respect to the horizontal product without any information
of the $M$-$M$ morphisms.  Then we can {\sl find}
$M$-$M$ sectors and determine their fusion rules
by the following theorem which generalizes a result 
for Goodman-de la Harpe-Jones subfactors in \cite{O7}
in a straightforward manner.

\begin{theorem}
\label{double2}
For any $\beta\in\MXM$ the element
$e_\beta\in\dta$ of Definition \ref{defe-beta} is a
minimal central projection with respect to the
horizontal product, and all minimal central projections
arise in this way in a bijective correspondence.
Furthermore, we have\footnote{Note that the fusion coefficients
with dimension prefactors as in \erf{dddN} coincide with the
structure constants used for $C$-algebras \cite{BI}.}
\be
e_\beta *_v e_{\beta'} =
\sum_{{\beta}''\in\MXM} \frac{d_\beta d_{\beta'}}{d_{\beta''}}
N_{\beta,\beta'}^{\beta''} \, e_{\beta''}
\label{dddN}
\ee
for all $\beta,\beta'\in\MXM$.  In particular,
the center $\cZ_h$ of $\dta$ with respect to the
horizontal product is closed under the vertical product.
\end{theorem}

\begin{proof}
That each $e_\beta$ is a minimal central projection and
that all minimal central projections arise in this way
is obvious from the description of the matrix units.
The vertical product $e_\beta *_v e_{\beta'}$ is given
graphically by the left-hand side of Fig\ \ref{v-prod1}.
%
%e_\beta *_v e_\beta' (1)
\begin{figure}[htb]
\begin{center}
\unitlength 0.6mm
\begin{picture}(211,50)
\thinlines
\put(10,23){\makebox(0,0){$\displaystyle \sum_{a,b,c}\; d_b$}}
\put(50,15){\arc{5}{1.571}{4.712}}
\put(30,15){\arc{5}{4.712}{1.571}}
\put(50,35){\arc{5}{1.571}{4.712}}
\put(30,35){\arc{5}{4.712}{1.571}}
\thicklines
\put(30,5){\line(0,1){40}}
\put(50,5){\line(0,1){40}}
\Thicklines
\put(30,15){\line(1,0){20}}
\put(30,35){\line(1,0){20}}
\put(42,15){\vector(1,0){0}}
\put(42,35){\vector(1,0){0}}
\put(26,3){\makebox(0,0){$c$}}
\put(27,25){\makebox(0,0){$b$}}
\put(26,47){\makebox(0,0){$a$}}
\put(54,3){\makebox(0,0){$c$}}
\put(53,25){\makebox(0,0){$b$}}
\put(54,47){\makebox(0,0){$a$}}
\put(40,10){\makebox(0,0){$\beta$}}
\put(40,40){\makebox(0,0){$\beta'$}}
\put(84,21){\makebox(0,0){$=\displaystyle
\sum_{a,b,c,\atop\beta'',\beta''',\beta''''} d_b$}}
\thinlines
\put(125,25){\arc{5}{1.571}{4.712}}
\put(110,25){\arc{5}{4.712}{1.571}}
\put(165,25){\arc{5}{1.571}{4.712}}
\put(150,25){\arc{5}{4.712}{1.571}}
\put(205,25){\arc{5}{1.571}{4.712}}
\put(190,25){\arc{5}{4.712}{1.571}}
\put(180,35){\arc{5}{1.571}{4.712}}
\put(135,35){\arc{5}{4.712}{1.571}}
\put(180,15){\arc{5}{1.571}{4.712}}
\put(135,15){\arc{5}{4.712}{1.571}}
\thicklines
\put(110,5){\line(0,1){40}}
\put(125,10){\line(0,1){30}}
\put(135,10){\line(0,1){30}}
\put(130,40){\arc{10}{3.142}{0}}
\put(130,10){\arc{10}{0}{3.142}}
\put(205,5){\line(0,1){40}}
\put(180,10){\line(0,1){30}}
\put(190,10){\line(0,1){30}}
\put(185,40){\arc{10}{3.142}{0}}
\put(185,10){\arc{10}{0}{3.142}}
\Thicklines
\put(110,25){\line(1,0){15}}
\put(119.5,25){\vector(1,0){0}}
\put(135,15){\line(1,0){10}}
\put(135,35){\line(1,0){10}}
\put(142,15){\vector(1,0){0}}
\put(142,35){\vector(1,0){0}}
\put(145,30){\arc{10}{4.712}{0}}
\put(145,20){\arc{10}{0}{1.571}}
\put(150,20){\line(0,1){10}}
\put(150,25){\line(1,0){15}}
\put(159.5,25){\vector(1,0){0}}
\put(170,30){\arc{10}{3.142}{4.712}}
\put(170,20){\arc{10}{1.571}{3.142}}
\put(165,20){\line(0,1){10}}
\put(190,25){\line(1,0){15}}
\put(199.5,25){\vector(1,0){0}}
\put(170,15){\line(1,0){10}}
\put(170,35){\line(1,0){10}}
\put(177,15){\vector(1,0){0}}
\put(177,35){\vector(1,0){0}}
\put(106,3){\makebox(0,0){$c$}}
\put(106,47){\makebox(0,0){$a$}}
\put(209,3){\makebox(0,0){$c$}}
\put(209,47){\makebox(0,0){$a$}}
\put(123,3){\makebox(0,0){$c$}}
\put(123,47){\makebox(0,0){$a$}}
\put(192,3){\makebox(0,0){$c$}}
\put(192,47){\makebox(0,0){$a$}}
\put(145,10){\makebox(0,0){$\beta$}}
\put(145,40){\makebox(0,0){$\beta'$}}
\put(170,10){\makebox(0,0){$\beta$}}
\put(170,40){\makebox(0,0){$\beta'$}}
\put(117.5,30){\makebox(0,0){$\beta'''$}}
\put(157.5,30){\makebox(0,0){$\beta''$}}
\put(197.5,30){\makebox(0,0){$\beta''''$}}
\put(139,25){\makebox(0,0){$b$}}
\put(176,25){\makebox(0,0){$b$}}
\end{picture}
\end{center}
\caption{The vertical product $e_\beta *_v e_{\beta'}$}
\label{v-prod1}
\end{figure}
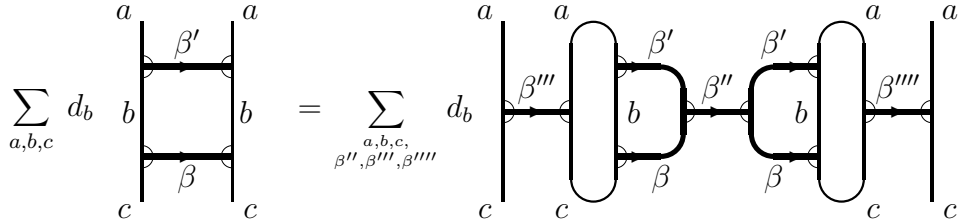
We can use the expansion of Lemma \ref{idexp} for the two parallel
wires $\beta$ and $\beta'$ in the middle. Now note that the horizontal
unit is given by $\bfe_h=\sum_\beta e_\beta$. Therefore, by
multiplying $\bfe_h$ from the left and from the right, we
obtain the diagram on the right-hand side of Fig.\ \ref{v-prod1}.
Reading the diagram from left to right, we observe that
intertwiners in $\Hom(\beta''',\beta'')$ and $\Hom(\beta'',\beta'''')$
are involved here. Hence we first obtain a factor
$\del {\beta'''}{\beta''} \del {\beta''}{\beta''''}$. Next, we can
use the trick of Fig.\ \ref{trick} to turn around the small arcs
at the trivalent vertices involving $a,b,\beta'$. This yields
a factor $d_\beta'/d_b$. This way we see that the diagram on the
right-hand side of Fig.\ \ref{v-prod1} represents the same element
of the $\dta$ as the diagram in. Fig.\ \ref{v-prod2}.
%
%e_\beta *_v e_\beta' (2)
\begin{figure}[htb]
\begin{center}
\unitlength 0.6mm
\begin{picture}(206,50)
\thinlines
\put(14,23){\makebox(0,0){$\displaystyle \sum_{a,b,c,\beta''} d_{\beta'}$}}
\put(40,12.5){\arc{5}{4.712}{1.571}}
\put(60,12.5){\arc{5}{1.571}{4.712}}
\put(75,27.5){\arc{5}{4.712}{1.571}}
\put(165,27.5){\arc{5}{1.571}{4.712}}
\put(90,12.5){\arc{5}{4.712}{1.571}}
\put(150,12.5){\arc{5}{1.571}{4.712}}
\put(110,27.5){\arc{5}{4.712}{1.571}}
\put(130,27.5){\arc{5}{1.571}{4.712}}
\put(180,12.5){\arc{5}{4.712}{1.571}}
\put(200,12.5){\arc{5}{1.571}{4.712}}
\dottedline{4}(65,0)(65,50)(175,50)(175,0)(65,0)
\thicklines
\put(40,5){\line(0,1){40}}
\put(60,10){\line(0,1){5}}
\put(65,10){\arc{10}{1.571}{3.142}}
\put(65,15){\arc{10}{3.142}{4.712}}
\put(65,5){\line(1,0){20}}
\put(65,20){\line(1,0){5}}
\put(75,25){\line(0,1){2.5}}
\put(70,25){\arc{10}{0}{1.571}}
\put(90,10){\line(0,1){12.5}}
\put(85,10){\arc{10}{0}{1.571}}
\put(85,22.5){\arc{10}{4.712}{0}}
\put(75,27.5){\line(1,0){10}}
\put(165,25){\line(0,1){2.5}}
\put(170,25){\arc{10}{1.571}{3.142}}
\put(170,20){\line(1,0){5}}
\put(180,10){\line(0,1){5}}
\put(175,10){\arc{10}{0}{1.571}}
\put(175,15){\arc{10}{4.712}{0}}
\put(155,5){\line(1,0){20}}
\put(155,27.5){\line(1,0){10}}
\put(150,10){\line(0,1){12.5}}
\put(155,10){\arc{10}{1.571}{3.142}}
\put(155,22.5){\arc{10}{3.142}{4.712}}
\put(200,5){\line(0,1){40}}
\Thicklines
\put(40,12.5){\line(1,0){20}}
\put(52,12.5){\vector(1,0){0}}
\put(99.5,12.5){\vector(1,0){0}}
\put(144.5,12.5){\vector(1,0){0}}
\put(122,27.5){\vector(1,0){0}}
\put(87,45){\vector(1,0){0}}
\put(157,45){\vector(1,0){0}}
\put(65,35){\line(1,0){5}}
\put(65,45){\line(1,0){40}}
\put(65,40){\arc{10}{1.571}{4.712}}
\put(75,27.5){\line(0,1){2.5}}
\put(70,30){\arc{10}{4.712}{0}}
\put(105,17.5){\arc{10}{0}{1.571}}
\put(105,40){\arc{10}{4.712}{0}}
\put(110,17.5){\line(0,1){22.5}}
\put(90,12.5){\line(1,0){15}}
\put(110,27.5){\line(1,0){20}}
\put(135,17.5){\arc{10}{1.571}{3.142}}
\put(135,40){\arc{10}{3.142}{4.712}}
\put(130,17.5){\line(0,1){22.5}}
\put(135,12.5){\line(1,0){15}}
\put(170,35){\line(1,0){5}}
\put(135,45){\line(1,0){40}}
\put(175,40){\arc{10}{4.712}{1.571}}
\put(165,27.5){\line(0,1){2.5}}
\put(170,30){\arc{10}{3.142}{4.712}}
\put(180,12.5){\line(1,0){20}}
\put(192,12.5){\vector(1,0){0}}
\put(36,3){\makebox(0,0){$c$}}
\put(36,47){\makebox(0,0){$a$}}
\put(204,3){\makebox(0,0){$c$}}
\put(204,47){\makebox(0,0){$a$}}
\put(75,8){\makebox(0,0){$c$}}
\put(70,16){\makebox(0,0){$a$}}
\put(165,8){\makebox(0,0){$c$}}
\put(170,16){\makebox(0,0){$a$}}
\put(85,22.5){\makebox(0,0){$b$}}
\put(155,22.5){\makebox(0,0){$b$}}
\put(97.5,17.5){\makebox(0,0){$\beta$}}
\put(95,40){\makebox(0,0){$\beta'$}}
\put(142.5,17.5){\makebox(0,0){$\beta$}}
\put(145,40){\makebox(0,0){$\beta'$}}
\put(50,17.5){\makebox(0,0){$\beta''$}}
\put(120,32.5){\makebox(0,0){$\beta''$}}
\put(190,17.5){\makebox(0,0){$\beta''$}}
\end{picture}
\end{center}
\caption{The vertical product $e_\beta *_v e_{\beta'}$}
\label{v-prod2}
\end{figure}
Now let us look at the part of this picture inside
the dotted box. Reading it from the left, this
part can be read for fixed $a$ and $c$ as
$\sum_{i,k} T_i T_{\beta,\beta'}^{\beta'';k}
(T_{\beta,\beta'}^{\beta'';k})^* T_i^*$,
and the sum over $i$ runs over a full orthonormal basis of
isometries $T_i$ in the Hilbert space $\Hom(\beta,\co c a\co{\beta'})$
since we have the summation over $b$. Next we look
at the part inside the dotted box of the diagram in
Fig.\ \ref{v-prod3}.
%
%e_\beta *_v e_\beta' (3)
\begin{figure}[htb]
\begin{center}
\unitlength 0.6mm
\begin{picture}(206,50)
\thinlines
\put(16.5,23){\makebox(0,0){$\displaystyle
\sum_{a,c,\beta'',\beta'''} d_{\beta'}$}}
\put(40,12.5){\arc{5}{4.712}{1.571}}
\put(60,12.5){\arc{5}{1.571}{4.712}}
\put(90,22.5){\arc{5}{4.712}{1.571}}
\put(150,22.5){\arc{5}{1.571}{4.712}}
\put(75,12.5){\arc{5}{4.712}{1.571}}
\put(165,12.5){\arc{5}{1.571}{4.712}}
\put(110,32.5){\arc{5}{4.712}{1.571}}
\put(130,32.5){\arc{5}{1.571}{4.712}}
\put(180,12.5){\arc{5}{4.712}{1.571}}
\put(200,12.5){\arc{5}{1.571}{4.712}}
\dottedline{4}(65,0)(65,50)(175,50)(175,0)(65,0)
\thicklines
\put(40,5){\line(0,1){40}}
\put(60,10){\line(0,1){5}}
\put(65,10){\arc{10}{1.571}{3.142}}
\put(65,15){\arc{10}{3.142}{4.712}}
\put(65,5){\line(1,0){5}}
\put(65,20){\line(1,0){5}}
\put(70,15){\arc{10}{4.712}{0}}
\put(75,10){\line(0,1){5}}
\put(70,10){\arc{10}{0}{1.571}}
\put(170,15){\arc{10}{3.142}{4.712}}
\put(170,20){\line(1,0){5}}
\put(180,10){\line(0,1){5}}
\put(175,10){\arc{10}{0}{1.571}}
\put(175,15){\arc{10}{4.712}{0}}
\put(170,5){\line(1,0){5}}
\put(165,10){\line(0,1){5}}
\put(170,10){\arc{10}{1.571}{3.142}}
\put(200,5){\line(0,1){40}}
\Thicklines
\put(40,12.5){\line(1,0){20}}
\put(52,12.5){\vector(1,0){0}}
\put(99.5,22.5){\vector(1,0){0}}
\put(144.5,22.5){\vector(1,0){0}}
\put(82,12.5){\vector(1,0){0}}
\put(162,12.5){\vector(1,0){0}}
\put(122,32.5){\vector(1,0){0}}
\put(87,45){\vector(1,0){0}}
\put(157,45){\vector(1,0){0}}
\put(65,35){\line(1,0){20}}
\put(65,45){\line(1,0){40}}
\put(65,40){\arc{10}{1.571}{4.712}}
\put(90,17.5){\line(0,1){12.5}}
\put(85,17.5){\arc{10}{0}{1.571}}
\put(85,30){\arc{10}{4.712}{0}}
\put(105,27.5){\arc{10}{0}{1.571}}
\put(105,40){\arc{10}{4.712}{0}}
\put(110,27.5){\line(0,1){12.5}}
\put(75,12.5){\line(1,0){10}}
\put(90,22.5){\line(1,0){15}}
\put(135,22.5){\line(1,0){15}}
\put(110,32.5){\line(1,0){20}}
\put(135,27.5){\arc{10}{1.571}{3.142}}
\put(135,40){\arc{10}{3.142}{4.712}}
\put(130,27.5){\line(0,1){12.5}}
\put(155,12.5){\line(1,0){10}}
\put(155,35){\line(1,0){20}}
\put(135,45){\line(1,0){40}}
\put(175,40){\arc{10}{4.712}{1.571}}
\put(150,17.5){\line(0,1){12.5}}
\put(155,17.5){\arc{10}{1.571}{3.142}}
\put(155,30){\arc{10}{3.142}{4.712}}
\put(180,12.5){\line(1,0){20}}
\put(192,12.5){\vector(1,0){0}}
\put(36,3){\makebox(0,0){$c$}}
\put(36,47){\makebox(0,0){$a$}}
\put(204,3){\makebox(0,0){$c$}}
\put(204,47){\makebox(0,0){$a$}}
\put(69,8){\makebox(0,0){$c$}}
\put(69,17){\makebox(0,0){$a$}}
\put(171,8){\makebox(0,0){$c$}}
\put(171,17){\makebox(0,0){$a$}}
\put(87,6){\makebox(0,0){$\beta'''$}}
\put(153,6){\makebox(0,0){$\beta'''$}}
\put(101,17){\makebox(0,0){$\beta$}}
\put(101,40){\makebox(0,0){$\beta'$}}
\put(139,17){\makebox(0,0){$\beta$}}
\put(139,40){\makebox(0,0){$\beta'$}}
\put(50,17.5){\makebox(0,0){$\beta''$}}
\put(120,37.5){\makebox(0,0){$\beta''$}}
\put(190,17.5){\makebox(0,0){$\beta''$}}
\end{picture}
\end{center}
\caption{The vertical product $e_\beta *_v e_{\beta'}$}
\label{v-prod3}
\end{figure}
Here, since we introduced the sum over $\beta'''$,
the part can be similarly read for fixed $a$ and $c$ as
$\sum_{j,k} S_j T_{\beta,\beta'}^{\beta'';k}
(T_{\beta,\beta'}^{\beta'';k})^* S_j^*$,
where the sum over $j$ runs over another orthonormal basis of
isometries $S_i$ in the Hilbert space
$\Hom(\beta,\co c a\co{\beta'})$.
Since such bases $\{T_i\}$ and $\{S_j\}$ are related by a
unitary matrix transformation
(this is essentially ``unitarity of $6j$-symbols''),
we conclude that the diagrams in Figs.\ \ref{v-prod2} and
\ref{v-prod3} represent the same element in $\dta$.
We now see that we first obtain a factor $\del {\beta''}{\beta'''}$.
Next we can turn around the small arcs at the outer two trivalent
vertices involving $\beta,\beta'$ and $\beta'''=\beta''$ so that
we obtain a factor $d_\beta/d_{\beta''}$. Then, by ``stretching''
the diagram a bit, we can read the diagram for fixed $a,c,\beta''$ as
\[ \bearl
\displaystyle\sum_{i,j,m=1}^{N_{\co c,a}^{\beta''}}
\sum_{k,l=1}^{N_{\beta,\beta'}^{\beta''}}
\frac{d_\beta d_{\beta'}}{d_{\beta''}}
T_{\co c,a}^{\beta'';i} (T_{\co c,a}^{\beta'';i})^*
T_{\co c,a}^{\beta'';j} (T_{\beta,\beta'}^{\beta'';l})^*
T_{\beta,\beta'}^{\beta'';k} (T_{\beta,\beta'}^{\beta'';k})^*
T_{\beta,\beta'}^{\beta'';l} (T_{\co c,a}^{\beta'';j})^*
T_{\co c,a}^{\beta'';m} (T_{\co c,a}^{\beta'';m})^* \\[.4em]
\qquad\qquad\qquad = \displaystyle\sum_{i=1}^{N_{\co c,a}^{\beta''}}
\frac{d_\beta d_{\beta'}}{d_{\beta''}} N_{\beta,\beta'}^{\beta''}\,
T_{\co c,a}^{\beta'';i} (T_{\co c,a}^{\beta'';i})^* \,.
\eear \]
Now proceeding with the summations over $a,c,\beta''$
yields the statement.
\end{proof}

Now consider the vector space with basis elements
$[\beta]$, $\beta\in\MXM$ which we can endow with a product
through
$[\beta][\beta']=\sum_{\beta''}N_{\beta,\beta'}^{\beta''}[\beta'']$.
We call the algebra defined this way the $M$-$M$ fusion rule algebra.
Similarly we define the $N$-$N$ fusion rule algebra using
morphisms in $\NXN$.

\begin{definition}
\label{Phi}{\rm
We define a linear map $\Phi$ from the $M$-$M$ fusion rule algebra
to ${\cal Z}_h$ by linear extension of $\Phi([\beta])=e_\beta/d_\beta$.
}\end{definition}
Theorem \ref{double2} now says that this map $\Phi$
is an isomorphism from the
$M$-$M$ fusion rule algebra onto $(\cZ_h, *_v)$.
Note that $(\cZ_h, *_v)$ is a non-unital subalgebra
of $(\dta, *_v)$. The unit $\bfe_v$ of $(\dta, *_v)$
is given by $\bfe_v=\sum_\la f_\la$, where
$f_\la=\sum_{a,b,j} f_{\la;a,b,j}^{a,b,j}$ whereas
the unit of  $(\cZ_h, *_v)$ is given by $e_0$.

\begin{definition}
\label{traces}{\rm
We define two linear functionals
$\varphi_h$ and $\tau_v$ on $\dta$
corresponding to the two product structures $*_h$ and
$*_v$ by linear extension of
\be
\begin{array}{rl}
\varphi_h (e_{\beta;c,a,i}^{d,b,j}) &=
\del ab \del cd \del ij  \, d_a d_c d_\beta / w^2 \,,\\[.4em]
\tau_v (f_{\la;c,d,j}^{a,b,i}) &=
\del ac \del cd \del ij \,  d_\la \,.
\eear
\label{eqstates}
\ee
}\end{definition}
Applied to an element in Fig.\ \ref{dta1} (Fig.\ \ref{dta2})
the functional $\varphi_h$ ($\tau_v$) can be characterized
graphically as in Fig.\ \ref{hori-st}
%
%The horizontal state
\thinlines
\begin{figure}[htb]
\begin{center}
\unitlength 0.6mm
\begin{picture}(240,40)
\thicklines
\put(30,10){\line(0,1){20}}
\put(50,10){\line(0,1){20}}
\put(135,20){\circle{20}}
\Thicklines
\put(30,20){\vector(1,0){20}}
\put(125,20){\vector(1,0){20}}
\put(30,33){\makebox(0,0){$a$}}
\put(50,33){\makebox(0,0){$b$}}
\put(30,7){\makebox(0,0){$c$}}
\put(50,7){\makebox(0,0){$d$}}
\put(40,24){\makebox(0,0){$\beta$}}
\put(135,33){\makebox(0,0){$a$}}
\put(135,7){\makebox(0,0){$c$}}
\put(135,24){\makebox(0,0){$\beta$}}
\put(10,20){\makebox(0,0){$\varphi_h:$}}
\put(25,20){\makebox(0,0){$S^*$}}
\put(55,20){\makebox(0,0){$T$}}
\put(120,20){\makebox(0,0){$S^*$}}
\put(150,20){\makebox(0,0){$T$}}
\put(87,20){\makebox(0,0){$\longmapsto \; \del ab \del cd \,
\displaystyle\frac{d_a d_c}{w^2}$}}
\put(200,20){\makebox(0,0){$= \del ab \del cd \,
\displaystyle\frac{(d_a d_c)^{3/2}d_\beta^{1/2}}{w^2} \,
\langle S,T \rangle$}}
\end{picture}
\end{center}
\caption{The horizontal functional $\varphi_h$}
\label{hori-st}
\end{figure}
(Fig.\ \ref{vert-st}).
%
%The vertical state
\thinlines
\begin{figure}[htb]
\begin{center}
\unitlength 0.6mm
\begin{picture}(240,40)
\thicklines
\put(30,10){\line(1,0){20}}
\put(30,30){\line(1,0){20}}
\put(145,20){\circle{20}}
\thinlines
\put(40,30){\vector(0,-1){20}}
\put(145,30){\vector(0,-1){20}}
\put(27,30){\makebox(0,0){$a$}}
\put(53,30){\makebox(0,0){$b$}}
\put(27,10){\makebox(0,0){$c$}}
\put(53,10){\makebox(0,0){$d$}}
\put(37,20){\makebox(0,0){$\la$}}
\put(142,20){\makebox(0,0){$\la$}}
\put(40,35){\makebox(0,0){$s^*$}}
\put(40,5){\makebox(0,0){$t$}}
\put(145,35){\makebox(0,0){$s^*$}}
\put(145,5){\makebox(0,0){$t$}}
\put(10,20){\makebox(0,0){$\tau_v:$}}
\put(130,20){\makebox(0,0){$a$}}
\put(160,20){\makebox(0,0){$b$}}
\put(90,20){\makebox(0,0){$\longmapsto \;\del ac \del bd\,
\displaystyle \sqrt{d_a d_b}$}}
\put(205,20){\makebox(0,0){$= \del ab \del cd \,
d_a d_b d_\la^{1/2} \, \langle s,t \rangle$}}
\end{picture}
\end{center}
\caption{The vertical functional $\tau_v$}
\label{vert-st}
\end{figure}
Therefore these functionals correspond to closing the
open ends of a diagram with prefactors as in the middle
part of Figs.\ \ref{hori-st} and \ref{vert-st}.

Recall that the global index of $\NXN$ is given by
$w=\sum_{\la\in\NXN} d_\la^2$. Note that we have sector
decompositions
$[a\iota]=\sum_\la \langle \la,a\iota \rangle [\la]$
and hence
$d_a d_\iota=\sum_\la \langle \la,a\iota \rangle d_\la$
for any $a\in\NXM$.
Using Frobenius reciprocity
$\langle \la,a\iota \rangle = \langle \la\co\iota,a \rangle$
we obtain similarly
$d_\la d_\iota=\sum_a \langle \la,a\iota \rangle d_a$.
Hence
$w=\sum_\la d_\la^2 = \sum_{\la,a} \langle \la,a\iota \rangle
d_\la d_a/d_\iota = \sum_a d_a^2$.
Similarly we obtain
$w=\sum_\beta d_\beta^2$ (cf.\ \cite{O3}).

\begin{lemma}
\label{tr-val}
We have $\varphi_h(e_\beta)=d_\beta^2/w$. In particular, the
functional $\varphi_h$ is a faithful state on $(\dta,*_h)$.
The functional $\tau_v$ is a (un-normalized) faithful
trace on $(\dta,*_v)$.
\end{lemma}

\begin{proof}
By Definition \ref{traces} and Fig.\  \ref{e-beta},
we compute
\[ \varphi_h(e_\beta) = \sum_{a,b\in\NXM} N_{\bar a, b}^\beta
d_a d_b d_\beta w^{-2} = \sum_{a\in\NXM}
\left(\sum_{b\in\NXM} N_{a,\beta}^b d_b\right)
d_\beta d_a w^{-2} = d_\beta^2 w^{-1}\,.\]
Since the horizontal unit $\bfe_h$ is given by
$\bfe_h=\sum_\beta e_\beta$ we find that
$\varphi(\bfe_h)=1$. As $\varphi_h$ sends
off-diagonal matrix units to zero and the diagonal
ones to strictly positive numbers, this proves that
$\varphi_h$ is a faithful state.
Obviously also $\tau_v$ sends off-diagonal matrix units
(with respect to $*_v$) to zero and the diagonal
ones to strictly positive numbers, and hence
it is a strictly positive functional but it is
not normalized. The trace property $\tau_v(xy)=\tau_v(yx)$
is clear from the definition of $\tau_v$ using
matrix units for $x$ and $y$.
\end{proof}

For $\tau_v$ we could have gained analogous properties
as for $\varphi_h$ by replacing the scalar
$d_\la$ in \erf{eqstates} by $d_a d_b d_\la/w^2$
(and by multiplying the scalars in Fig.\ \ref{vert-st}
also by $d_a d_b/w^2$). However, we chose a different
normalization on each matrix unit in  order to turn
$\tau_v$ into a trace on $(\dta,*_v)$. Later we want to study
the center $(\cZ_h,*_v)$ which is, as we have seen, a subalgebra
of $(\dta,*_v)$. Therefore $\tau_v$ provides  a faithful trace on
$(\cZ_h,*_v)$ but it has in general different weightings on its simple
summands. To construct from $\tau_v$ a trace which sends one-dimensional
projections to one will in particular be possible in the case that
$\NXN$ is non-degenerately braided, see Subsect.\ \ref{sec-dfa} below.

This is also the case in the following most basic example of
the double triangle algebra.
Let $N$ be a type III factor and $G$ a finite group acting freely
on $N$.  Consider the subfactor $N\subset N\rtimes G=M$.
Then (with the minimal choice for $\cX$)
the double triangle algebra $\dta$ for this subfactor
is just the group algebra of $G$. That is, the double
triangle algebra is spanned by the group elements linearly.
The horizontal product is given by the group
multiplication. By  Proposition \ref{double2} we conclude
that the minimal central projections in $\dta$
and thus irreducible $M$-$M$ sectors are labelled by
the irreducible representations of $G$.
(Of course, this identification of the $M$-$M$ sectors
is well-known for that example.) 
The functional $\tau_v$ gives the standard trace on the group
algebra, and the vertical product corresponds to the
ordinary tensor product of group representations.

\section{$\a$-Induction, Chiral Generators and
Modular Invariants}
\label{sec-aicpmodinv}

\subsection{Relating $\a$-induction to chiral generators}
\label{sec-iden}

We will now define chiral generators for braided
subfactors and prove that the concepts of $\a$-induction
and chiral generators are essentially the same.
For the rest of this paper deal with the following

\begin{assumption}
\label{set-braid}
{\rm In addition to Assumption \ref{set-fin} we now assume that
the system $\NXN$ is braided.
}\end{assumption}
With the braiding we have now the notion of $\a$-induction
in the sense of Subsect.\ \ref{sec-aindbsf}. From now on
we are also dealing with crossings of $N$-$N$ wires and
mixed crossings introduced in  Subsect.\ \ref{sec-aindbsf}.
We now present chiral generators as our version of a
definition Ocneanu originally introduced for systems of
bimodules arising from A-D-E Dynkin diagrams in \cite{O7}.
The construction of the chiral generator is
similar to the ``Ocneanu projection'' in the
tube algebra \cite{O6} (see also \cite{EK4})
and also related to Izumi's analysis \cite{I3} of the
tube algebra in terms of sectors for the
Longo-Rehren inclusion \cite{LR}.

\begin{definition}
\label{chiral1}{\rm
For any $\la\in\NXN$, we define an element $p^+_\la\in\dta$
by the diagram on the left-hand side of Fig.\  \ref{chiral-proj}
and call it a {\sl chiral generator}.
Similarly, we also define $p^-_\la$ by exchanging over- and
undercrossings.
}\end{definition}
%
%chiral generator
\thinlines
\begin{figure}[htb]
\begin{center}
\unitlength 0.6mm
\begin{picture}(215,40)
\thicklines
\put(10,17){\makebox(0,0){$\displaystyle\sum_{a,b}$}}
\put(30,5){\line(0,1){8}}
\put(30,17){\line(0,1){18}}
\put(50,5){\line(0,1){8}}
\put(50,17){\line(0,1){18}}
\Thicklines
\put(30,15){\line(1,0){20}}
\put(42,15){\vector(1,0){0}}
\thinlines
\put(30,20){\arc{10}{1.571}{4.712}}
\put(50,20){\arc{10}{4.712}{1.571}}
\put(30,25){\arc{5}{1.571}{4.712}}
\put(50,25){\arc{5}{4.712}{1.571}}
\put(25,35){\makebox(0,0){$a$}}
\put(55,35){\makebox(0,0){$a$}}
\put(25,5){\makebox(0,0){$b$}}
\put(55,5){\makebox(0,0){$b$}}
\put(40,21){\makebox(0,0){$\a_\la^+$}}
\put(70,20){\makebox(0,0){$=$}}
\thicklines
\put(80,17){\makebox(0,0){$\displaystyle\sum_{a,b}$}}
\put(100,5){\line(0,1){18}}
\put(100,27){\line(0,1){8}}
\put(120,5){\line(0,1){18}}
\put(120,27){\line(0,1){8}}
\Thicklines
\put(100,25){\line(1,0){20}}
\put(112,25){\vector(1,0){0}}
\thinlines
\put(100,20){\arc{10}{1.571}{4.712}}
\put(120,20){\arc{10}{4.712}{1.571}}
\put(100,15){\arc{5}{1.571}{4.712}}
\put(120,15){\arc{5}{4.712}{1.571}}
\put(95,35){\makebox(0,0){$a$}}
\put(125,35){\makebox(0,0){$a$}}
\put(95,5){\makebox(0,0){$b$}}
\put(125,5){\makebox(0,0){$b$}}
\put(110,19){\makebox(0,0){$\a_\la^+$}}
\put(140,20){\makebox(0,0){$=$}}
\thinlines
\put(150,17){\makebox(0,0){$\displaystyle\sum_{a,b,\nu}$}}
\put(185,10){\line(0,1){7}}
\put(185,30){\line(0,-1){7}}
\put(186,20.5){\vector(1,0){0}}
\put(185,10){\vector(0,-1){0}}
\put(185,30){\arc{20}{0}{3.14}}
\put(185,30){\arc{5}{0}{3.14}}
\put(175,30){\arc{5}{0}{3.14}}
\put(195,30){\arc{5}{0}{3.14}}
\put(185,10){\arc{5}{3.14}{6.28}}
\thicklines
\put(170,10){\line(1,0){30}}
\put(170,30){\line(1,0){30}}
\put(170,35){\arc{10}{1.571}{3.142}}
\put(200,35){\arc{10}{0}{1.571}}
\put(170,5){\arc{10}{3.142}{4.712}}
\put(200,5){\arc{10}{4.712}{0}}
\put(195,22){\makebox(0,0){$\la$}}
\put(181,14){\makebox(0,0){$\nu$}}
\put(160,5){\makebox(0,0){$b$}}
\put(210,5){\makebox(0,0){$b$}}
\put(160,35){\makebox(0,0){$a$}}
\put(210,35){\makebox(0,0){$a$}}
\put(180,34){\makebox(0,0){$b$}}
\put(190,34){\makebox(0,0){$b$}}
\end{picture}
\end{center}
\caption{A chiral generator $p^+_\la$}
\label{chiral-proj}
\end{figure}
Note that we do {\sl not} assume the non-degeneracy
of the braiding for the definition $p^+_\la$.

We obtain the diagram in the middle from the one on the
left-hand side in Fig.\ \ref{chiral-proj} by applying
two IBFE's. This way we obtain two twists in the semi-circular
thin wires which correspond to the label $\la$ but they give
complex conjugate phases so that their effects cancel out.
The diagram on the right-hand side is obtained by
Lemma \ref{idexp} and application of the IBFE,
and this shows that our definition
coincides with Ocneanu's notion given in his setting.

Since $\a_\la^\pm \iota=\iota\la$ we find that each irreducible
subsector $[\beta]$ of $[\a_\la^\pm]$ is the equivalence class
of some $\beta\in\MXM$ if $\la\in\NXN$. Therefore we have the sector
decomposition
$[\a_\la^\pm]=\sum_{\beta\in\MXM}\langle\beta,\a_\la^\pm\rangle [\beta]$,
and we can consider $[\a_\la^\pm]$ as an element of the $M$-$M$
fusion algebra. The relation between the sector decomposition of
$[\a_\la^\pm]$ and the chiral generator is clarified
by the following result.

\begin{theorem}
\label{identify}
For any $\la\in\NXN$, we have
$d_\la^{-1} p^\pm_\la = \sum_{\beta\in\MXM} d_\beta^{-1}
\langle\beta,\a_\la^+\rangle e_\beta$, and consequently
$p^\pm_\la= d_\la \Phi([\a^\pm_\la])$. In particular,
$p^\pm_\la$ is in the center $\cZ_h$.
\end{theorem}

\begin{proof}
We only show the statement for the $+$-sign; the other case is
analogous. First we fix $a,b\in\NXM$ and $\la\in\NXN$.
For each $\beta\in\MXM$ we choose orthonormal bases
of isometries $T_{\co b a}^{\beta;i}\in\Hom(\beta,\co b a)$,
$i=1,2,...,N_{\co b,a}^\beta$, so that
$\sum_{\beta,i} T_{\co b a}^{\beta;i}
(T_{\co b a}^{\beta;i})^* = \bfe_M$.
Using Frobenius reciprocity, we obtain an orthonormal basis of
isometries
$\cL^{-1}_b(T_{\co b a}^{\beta;i})=d_a^{1/2} d_b^{1/2} d_\beta^{-1/2}
b(T_{\co b a}^{\beta;i})^* {\co r}_b\in\Hom(a,b\beta)$.
Here we chose an isometry ${\co r}_b\in\Hom(\id_N,b\co b)$
such that there is an isometry $R_b\in\Hom(\id_M,\co b b)$
subject to relations $b(R_b)^*{\co r}_b=d_b^{-1}\bfe_N$
and $\co b({\co r}_b)^*R_b=d_b^{-1}\bfe_M$, as usual.
Choosing also orthonormal bases of isometries
$V_{\beta;\ell}\in\Hom(\beta,\a_\la^+)$,
$\ell=1,2,...,\langle \beta, \a_\la^+ \rangle$,
for each $\beta\in\MXM$ (so that
$\sum_{\beta,\ell} V_{\beta;\ell} V_{\beta;\ell}^*=\bfe_M$)
we find that
$\{ b(V_{\beta;\ell}) \cL^{-1}_b(T_{\co b a}^{\beta;i}) \}_{\beta,i,\ell}$
gives an orthonormal basis of isometries of
$\Hom(a,b\a_\la^+)$. Finally, using Proposition \ref{stat},
we find that putting
\[ s_{\beta;\ell,i}=\epsp \la{b\iota} ^*
 b(V_{\beta;\ell})\cL^{-1}_b(T_{\co b a}^{\beta;i}) =
\sqrt\frac{d_a d_b}{d_\beta} \epsp \la{b\iota} ^*
b(V_{\beta;\ell} (T_{\co b a}^{\beta;i})^*){\co r}_b \]
defines an orthonormal basis of isometries
$\{s_{\beta;\ell,i}\}_{\beta,i,\ell}$ of $\Hom(a,\la b)$. Then
we have for any $\ell=1,2,...,\langle \beta,\a_\la^+ \rangle$
by the elementary relations for the intertwiners $R_b,{\co r}_b$
the following identity:
\[ \bearll
T_{\co b a}^{\beta;i}(T_{\co b a}^{\beta;i})^*
&= d_b^2 \, \co b ({\co r}_b)^* \,
\co b b(T_{\co b a}^{\beta;i} V_{\beta;\ell}^*) \, R_b R_b^* \,
\co b b(V_{\beta;\ell} (T_{\co b a}^{\beta;i})^*) \,
\co b({\co r}_b) \\[.5em]
&= \displaystyle\frac{d_\beta d_b}{d_a} \, \co b (s_{\beta;\ell,i}
\epsp \la{b\iota} ^*) \,
R_bR_b^* \, \co b( \epsp \la{b\iota} s_{\beta;\ell,i}^* )\,.
\eear \]
The second line yields graphically exactly the diagram in
Fig.\ \ref{inter1} where we read the diagram from the left to
the right in order to interpret it as an intertwiner in $\dtap$.
%
%diagram for TT^*
\thinlines
\begin{figure}[htb]
\begin{center}
\unitlength 0.6mm
\begin{picture}(120,50)
\put(55,25){\line(1,0){10}}
\put(45,35){\vector(1,-1){10}}
\put(85,25){\vector(1,0){10}}
\put(95,25){\line(1,1){10}}
\thicklines
\put(35,5){\line(1,0){20}}
\put(35,35){\line(1,0){10}}
\put(45,35){\line(1,1){10}}
\put(65,35){\line(0,-1){7}}
\put(65,15){\line(0,1){7}}
\put(55,15){\arc{20}{0}{1.57}}
\put(55,35){\arc{20}{4.71}{6.28}}
\put(115,5){\line(-1,0){20}}
\put(115,35){\line(-1,0){10}}
\put(105,35){\line(-1,1){10}}
\put(85,35){\line(0,-1){7}}
\put(85,15){\line(0,1){7}}
\put(95,15){\arc{20}{1.57}{3.14}}
\put(95,35){\arc{20}{3.14}{4.71}}
\Thicklines
\put(65,25){\vector(1,0){10}}
\put(75,25){\line(1,0){10}}
\put(32,5){\makebox(0,0){$b$}}
\put(32,35){\makebox(0,0){$a$}}
\put(65,45){\makebox(0,0){$b$}}
\put(59,21){\makebox(0,0){$\la$}}
\put(55,34){\makebox(0,0){$s_{\beta;\ell,i}$}}
\put(118,5){\makebox(0,0){$b$}}
\put(118,35){\makebox(0,0){$a$}}
\put(85,45){\makebox(0,0){$b$}}
\put(91,21){\makebox(0,0){$\la$}}
\put(75,20){\makebox(0,0){$\a_\la^+$}}
\put(95,34.5){\makebox(0,0){$s_{\beta;\ell,i}^*$}}
\put(13,22){\makebox(0,0){$\displaystyle\frac{d_\beta}
{\sqrt{d_\la d_a d_b}}$}}
\end{picture}
\end{center}
\caption{Diagram for $T_{\co b a}^{\beta;i}(T_{\co b a}^{\beta;i})^*$}
\label{inter1}
\end{figure}
Now let us take on both sides first the summation over
$i=1,2,...,N_{\co b,a}^\beta$. Then the left-hand side gives
exactly the $\Hom(\co b a,\co b a)$ part of $e_\beta$
(in $\dtap$) as defined in Definition \ref{defe-beta}.
Next we divide by $d_\beta$ and we proceed with the
summation over $\ell=1,2,...,\langle\beta,\a_\la^+\rangle$
and $\beta\in\MXM$.
On the left-hand side we obtain the $\Hom(\co b a,\co b a)$ part
of $\sum_\beta d_\beta^{-1} \langle\beta,\a_\la^+\rangle e_\beta$
this way, and this is exactly the $\Hom(\co b a,\co b a)$
part of $\Phi([\a_\la^+])$.
On the right-hand side we now have a summation over the
full basis $\{s_{\beta;\ell,i}\}_{\beta,i,\ell}$ of $\Hom(a,\la b)$.
Therefore we can use the graphical convention of
Fig.\ \ref{e-beta} to put a small semi-circle around
the wire labelled by $\la$ at the two trivalent vertices.
This gives us a factor $\sqrt{d_a d_b/d_\la}$ so that only
a factor $d_\la^{-1}$ remains from the original prefactor
in Fig.\ \ref{inter1}. Thus, by repeating the above procedure
for all $a,b\in\NXM$ and making finally the summation
over $a,b\in\NXM$, we obtain on the left the full
$\Phi([\a_\la^+])$ whereas the right-hand side gives
graphically the diagram in Fig.\ \ref{ident}.
%
%chiral proj = alpha induction
\thinlines
\begin{figure}[htb]
\begin{center}
\unitlength 0.6mm
\begin{picture}(115,50)
\put(50,25){\line(1,0){10}}
\put(40,35){\vector(1,-1){10}}
\put(80,25){\vector(1,0){10}}
\put(90,25){\line(1,1){10}}
\put(40,35){\arc{5}{5.49}{9.42}}
\put(100,35){\arc{5}{0}{3.93}}
\thicklines
\put(30,5){\line(1,0){20}}
\put(30,35){\line(1,0){10}}
\put(40,35){\line(1,1){10}}
\put(60,35){\line(0,-1){7}}
\put(60,15){\line(0,1){7}}
\put(50,15){\arc{20}{0}{1.57}}
\put(50,35){\arc{20}{4.71}{6.28}}
\put(110,5){\line(-1,0){20}}
\put(110,35){\line(-1,0){10}}
\put(100,35){\line(-1,1){10}}
\put(80,35){\line(0,-1){7}}
\put(80,15){\line(0,1){7}}
\put(90,15){\arc{20}{1.57}{3.14}}
\put(90,35){\arc{20}{3.14}{4.71}}
\Thicklines
\put(60,25){\vector(1,0){10}}
\put(70,25){\line(1,0){10}}
\put(27,5){\makebox(0,0){$b$}}
\put(27,35){\makebox(0,0){$a$}}
\put(60,45){\makebox(0,0){$b$}}
\put(54,21){\makebox(0,0){$\la$}}
\put(113,5){\makebox(0,0){$b$}}
\put(113,35){\makebox(0,0){$a$}}
\put(80,45){\makebox(0,0){$b$}}
\put(86,21){\makebox(0,0){$\la$}}
\put(70,20){\makebox(0,0){$\a_\la^+$}}
\put(11,22){\makebox(0,0){$\displaystyle
\sum_{a,b}\; \frac{1}{d_\la}$}}
\end{picture}
\end{center}
\caption{The image $\Phi([\a^+_\la])=\sum_\beta d_\beta^{-1}
\langle\beta,\a_\la^+\rangle e_\beta$}
\label{ident}
\end{figure}
The diagram on the left-hand side in Fig.\ \ref{chiral-proj}
is obtained from Fig.\ \ref{ident}, up to the factor
$d_\la$, by a topological move.
\end{proof}

Note that it was not clear from the definition that
the chiral generators are in the center $\cZ_h$,
but Theorem \ref{identify} proves this centrality as it states
that $p_\la^\pm$ is a linear combination of $e_\beta$'s.
Also note that if $\a^\pm_\la$ is irreducible then $p^\pm_\la$
is a (horizontal) projection, however, if $\a^\pm_\la$
is not irreducible, then $p_\la^\pm$ is a sum over projections
with weight coefficients arising from the nature of the
isomorphism $\Phi$ in Definition \ref{Phi}.

Two of us \cite[Subsection 3.3]{BE3} established a relative
braiding between the two kinds of $\a$-induction,
which holds in a fairly general context. (It does neither
depend on chiral locality nor even on finite depth.)
Theorem \ref{identify} now shows that Ocneanu's relative
braiding \cite{O7} is a special case of the analysis
in \cite[Subsection 3.3]{BE3}.

{}From Theorem \ref{identify} and the homomorphism
property of $\a$-induction \cite[Lemma 3.10]{BE1}, we obtain
immediately the following

\begin{corollary}
\label{homom}
The chiral generators
$p^\pm_\la$ are in $\cZ_h$.  For $\la,\mu\in\NXN$, we have
\[ p^\pm_\la *_v p^\pm_\mu = \sum_{\nu\in\NXN}
\frac{d_\la d_\mu}{d_\nu} N_{\la,\mu}^\nu \, p^\pm_\nu \,. \]
\end{corollary}
Note that this corollary shows that the $M$-$M$ fusion rule
algebra contains two representations of the $N$-$N$ fusion
rule algebra.

\subsection{Modular invariants for braided subfactors}
\label{sec-modular}

We will now show that a notion of ``modular invariant'' arises
naturally for a braided subfactor.
We first note that under Assumption \ref{set-braid}, we have
matrices $Y=(Y_{\la,\mu})$ and $T=(T_{\la,\mu})$ for the system
$\sys=\NXN$ as in Subsection \ref{braid-prelim}. 
We recall that in the case that the braiding is
non-degenerate, the matrix $S=w^{-1/2}Y$ is unitary and the
matrices $S$ and (the diagonal) $T$ obey the Verlinde modular
algebra by Theorem \ref{ST}.
Motivated by the results of \cite{BE3} we now
construct a certain matrix $Z$ commuting with $Y$ and $T$
such that it is a ``modular invariant mass matrix'' in the
usual sense of conformal field theory whenever the braiding
is non-degenerate.

\begin{definition}
\label{modular1}
{\rm For a system $\cal X$ satisfying Assumption \ref{set-braid},
we define a matrix $Z$ with entries
$Z_{\la,\mu}=\langle \a_\la^+,\a_\mu^- \rangle$, $\la,\mu\in\NXN$.
}\end{definition}
As $Z_{\la,\mu}$ is by definition a dimension and since
$\a_{\id_N}^\pm=\id_M$ is irreducible by virtue
of the factor property of $M$, the matrix elements obviously
satisfy the conditions in \erf{massmatZ} for $\la,\mu\in\NXN$,
where the label ``0'' refers as usual to the identity morphism
$\id_N\in\NXN$. We relate the definition of $Z$ to the
chiral generators by the following

\begin{theorem}
\label{modular1+}
We have the identity
\be
Z_{\la,\mu} = \frac w{d_\la d_\mu} \, \varphi_h(p_\la^+ *_h  p_\mu^-) \,,
\qquad \la,\mu\in\NXN \,.
\label{Zvarphipp}
\ee
Therefore the number $Z_{\la,\mu}$ is graphically
represented as in Fig.\ \ref{Zgraph}.
\end{theorem}
%
%Modular invariant matrix
\begin{figure}[htb]
\begin{center}
\unitlength 0.7mm
\begin{picture}(115,50)
\thinlines
\put(34,23){\makebox(0,0){$Z_{\la,\mu} \;= \;\displaystyle
\sum_{b,c} \;\;\frac{d_b d_c}{w d_\la d_\mu}$}}
\put(80,5){\line(0,1){10}}
\put(80,45){\line(0,-1){10}}
\put(100,5){\line(0,1){8}}
\put(100,45){\line(0,-1){8}}
\put(80,45){\arc{5}{0}{3.142}}
\put(100,45){\arc{5}{0}{3.142}}
\put(80,5){\arc{5}{3.142}{0}}
\put(100,5){\arc{5}{3.142}{0}}
\Thicklines
\put(80,15){\line(0,1){20}}
\put(100,17){\line(0,1){16}}
\put(80,23){\vector(0,-1){0}}
\put(100,27){\vector(0,1){0}}
\thicklines
\put(77,5){\line(1,0){26}}
\put(77,45){\line(1,0){26}}
\put(77,15){\line(1,0){1}}
\put(77,35){\line(1,0){1}}
\put(82,15){\line(1,0){21}}
\put(82,35){\line(1,0){21}}
\put(77,40){\arc{10}{1.571}{4.712}}
\put(103,40){\arc{10}{4.712}{1.571}}
\put(77,10){\arc{10}{1.571}{4.712}}
\put(103,10){\arc{10}{4.712}{1.571}}
\put(90,2){\makebox(0,0){$c$}}
\put(90,48){\makebox(0,0){$c$}}
\put(90,12){\makebox(0,0){$b$}}
\put(90,38){\makebox(0,0){$b$}}
\put(74,25.4){\makebox(0,0){$\a_\la^+$}}
\put(107,24.6){\makebox(0,0){$\a_\mu^-$}}
\end{picture}
\end{center}
\caption{Graphical representation of $Z_{\la,\mu}$}
\label{Zgraph}
\end{figure}

\begin{proof}
{}From Theorem \ref{identify} we obtain
\[  \sum_{\beta\in\MXM}\frac 1{d_\beta}
\lan \a^+_\la,\beta \ran e_\beta = \frac 1{d_\la} p_\la^+ \,. \]
Hence
\[  \sum_{\beta\in\MXM}
\frac 1{d_\beta^2} \lan \a^+_\la,\beta \ran
\lan \a^-_\mu,\beta\ran e_\beta
= \frac 1{d_\la d_\mu} \, p_\la^+ *_h p_\mu^-  \,. \]
Application of the horizontal state $\varphi_h$ of Definition
\ref{traces} and multiplication by $w$ yields \erf{Zvarphipp}
since $[\a_\la^+]$ and $[\a_\mu^-]$ decompose into sectors
$[\beta]$ with $\beta\in\MXM$, and by Lemma \ref{tr-val}.
Now the right-hand side of\ \erf{Zvarphipp} is given graphically
by the diagram on the left in Fig.\ \ref{Z1},
%
%proof of Z=silly picture
\begin{figure}[htb]
\begin{center}
\unitlength 0.6mm
\begin{picture}(215,60)
\thicklines
\put(15,30){\makebox(0,0){$\displaystyle\sum_{b,c}\;
\frac{d_b d_c}{w d_\la d_\mu}$}}
\put(40,10){\line(0,1){23}}
\put(40,37){\line(0,1){13}}
\put(60,15){\line(0,1){18}}
\put(60,37){\line(0,1){8}}
\put(45,10){\arc{10}{1.571}{3.142}}
\put(45,50){\arc{10}{3.142}{4.712}}
\put(67.5,15){\arc{15}{0}{3.142}}
\put(67.5,45){\arc{15}{3.142}{0}}
\put(90,10){\arc{10}{0}{1.571}}
\put(90,50){\arc{10}{4.712}{0}}
\put(75,15){\line(0,1){30}}
\put(95,10){\line(0,1){40}}
\put(45,5){\line(1,0){45}}
\put(45,55){\line(1,0){45}}
\Thicklines
\put(40,35){\line(1,0){20}}
\put(52,35){\vector(1,0){0}}
\put(77,35){\line(1,0){16}}
\put(87,35){\vector(1,0){0}}
\thinlines
\put(40,30){\arc{10}{1.571}{4.712}}
\put(60,30){\arc{10}{4.712}{1.571}}
\put(40,25){\arc{5}{1.571}{4.712}}
\put(60,25){\arc{5}{4.712}{1.571}}
\put(75,30){\arc{10}{1.571}{4.212}}
\put(95,30){\arc{10}{5.212}{1.571}}
\put(75,25){\arc{5}{1.571}{4.712}}
\put(95,25){\arc{5}{4.712}{1.571}}
\put(36,55){\makebox(0,0){$c$}}
\put(99,5){\makebox(0,0){$b$}}
\put(67.5,48.5){\makebox(0,0){$c$}}
\put(67.5,12.5){\makebox(0,0){$b$}}
\put(50,30){\makebox(0,0){$\a_\la^+$}}
\put(85,30){\makebox(0,0){$\a_\mu^-$}}
\put(110,30){\makebox(0,0){$=$}}
\thicklines
\put(135,30){\makebox(0,0){$\displaystyle\sum_{b,c}\;
\frac{d_a d_b}{w d_\la d_\mu}$}}
\put(160,10){\line(0,1){15.5}}
\put(160,29.5){\line(0,1){20.5}}
\put(180,15){\line(0,1){10.5}}
\put(180,29.5){\line(0,1){15.5}}
\put(165,10){\arc{10}{1.571}{3.142}}
\put(165,50){\arc{10}{3.142}{4.712}}
\put(187.5,15){\arc{15}{0}{3.142}}
\put(187.5,45){\arc{15}{3.142}{0}}
\put(210,10){\arc{10}{0}{1.571}}
\put(210,50){\arc{10}{4.712}{0}}
\put(195,15){\line(0,1){30}}
\put(215,10){\line(0,1){40}}
\put(165,5){\line(1,0){45}}
\put(165,55){\line(1,0){45}}
\Thicklines
\put(160,27.5){\line(1,0){20}}
\put(172,27.5){\vector(1,0){0}}
\put(162,32.5){\line(1,0){16}}
\put(168,32.5){\vector(-1,0){0}}
\thinlines
\put(160,22.5){\arc{10}{1.571}{4.712}}
\put(180,22.5){\arc{10}{4.712}{1.571}}
\put(160,17.5){\arc{5}{1.571}{4.712}}
\put(180,17.5){\arc{5}{4.712}{1.571}}
\put(160,37.5){\arc{10}{2.071}{4.712}}
\put(180,37.5){\arc{10}{4.712}{1.071}}
\put(160,42.5){\arc{5}{1.571}{4.712}}
\put(180,42.5){\arc{5}{4.712}{1.571}}
\put(156,30){\makebox(0,0){$c$}}
\put(184,30){\makebox(0,0){$c$}}
\put(199,30){\makebox(0,0){$b$}}
\put(211,30){\makebox(0,0){$b$}}
\put(170,22.5){\makebox(0,0){$\a_\la^+$}}
\put(170,38.5){\makebox(0,0){$\a_\mu^-$}}
\end{picture}
\end{center}
\caption{The scalar $w d_\la^{-1} d_\mu^{-1}
\varphi_h(p_\la^+ *_h  p_\mu^-)$}
\label{Z1}
\end{figure}
and we can slide around the trivalent vertices
to obtain the diagram on the right-hand side.
Without changing the scalar value we can now
open the outer wire labelled by $b$ and close it
on the other side, as in Fig.\ \ref{cut}.
This way we obtain the picture in Fig.\ \ref{Zgraph}
up to a 90 degree rotation, but a rotation is irrelevant
for the scalar values.
\end{proof}

We remark that we can apply Lemma \ref{idexp} to
replace the two horizontal wires labelled by $b$
by a summation over a thin wire $\nu$, and this way
we obtain an equivalent diagram from Fig.\ \ref{Zgraph}
for the matrix elements $Z_{\la,\mu}$, which only consists
of thin ($N$-$N$) wires $\la,\mu,\nu$ and thick ($N$-$M$)
wires $b,c$ but which does not involve very thick ($M$-$M$)
wires labelled by $\a$-induced morphisms $\a_\la^+,\a_\mu^-$.

\begin{theorem}
\label{modular2}
The matrix $Z$ of Definition \ref{modular1}
commutes with the matrices $Y$ and $T$ of
the system $\NXN$.
\end{theorem}

\begin{proof}
Using the diagram for the matrix elements $Y_{\nu,\la}$ in
Fig.\ \ref{Ymatrix}, the sum $\sum_\la Y_{\nu,\la}Z_{\la,\mu}$
can be represented by the diagram on the left-hand side of
Fig.\ \ref{[Y,Z]1}.
%
% Z commutes with Y (1)
\begin{figure}[htb]
\begin{center}
\unitlength 0.7mm
\begin{picture}(174,50)
\thinlines
\put(12,23){\makebox(0,0){$\displaystyle
\sum_{b,c,\la} \;\frac{d_b d_c}{w d_\mu}$}}
\put(40,5){\line(0,1){10}}
\put(40,45){\line(0,-1){5.5}}
\put(40,35.5){\line(0,-1){15.5}}
\put(60,5){\line(0,1){8}}
\put(60,45){\line(0,-1){23}}
\put(40,45){\arc{5}{0}{3.142}}
\put(60,45){\arc{5}{0}{3.142}}
\put(40,5){\arc{5}{3.142}{0}}
\put(60,5){\arc{5}{3.142}{0}}
\put(40,32.5){\arc{10}{1.971}{1.171}}
\put(40,22){\vector(0,-1){0}}
\put(60,34){\vector(0,1){0}}
\put(44.8,34){\vector(0,1){0}}
\Thicklines
\put(40,15){\line(0,1){5}}
\put(60,17){\line(0,1){1}}
\thicklines
\put(37,5){\line(1,0){26}}
\put(37,45){\line(1,0){26}}
\put(37,15){\line(1,0){1}}
\put(37,20){\line(1,0){1}}
\put(42,15){\line(1,0){21}}
\put(42,20){\line(1,0){21}}
\put(37,25){\arc{10}{1.571}{3.142}}
\put(63,25){\arc{10}{0}{1.571}}
\put(37,40){\arc{10}{3.142}{4.712}}
\put(63,40){\arc{10}{4.712}{0}}
\put(32,25){\line(0,1){15}}
\put(68,25){\line(0,1){15}}
\put(37,10){\arc{10}{1.571}{4.712}}
\put(63,10){\arc{10}{4.712}{1.571}}
\put(50,2){\makebox(0,0){$c$}}
\put(50,48){\makebox(0,0){$c$}}
\put(50,12){\makebox(0,0){$b$}}
\put(50,24){\makebox(0,0){$b$}}
\put(36,25){\makebox(0,0){$\la$}}
\put(64,24){\makebox(0,0){$\mu$}}
\put(50,32.5){\makebox(0,0){$\nu$}}
\thinlines
\put(105,23){\makebox(0,0){$= \;\displaystyle
\sum_{a,b,c,\la} \;\frac{d_a d_b d_c}{w d_\mu d_\nu}$}}
\put(145,5){\line(0,1){10}}
\put(145,45){\line(0,-1){5.5}}
\put(145,35.5){\line(0,-1){15.5}}
\put(165,5){\line(0,1){8}}
\put(165,45){\line(0,-1){23}}
\put(145,45){\arc{5}{0}{3.142}}
\put(165,45){\arc{5}{0}{3.142}}
\put(145,5){\arc{5}{3.142}{0}}
\put(165,5){\arc{5}{3.142}{0}}
\put(145,32.5){\arc{10}{4.712}{1.171}}
\put(137,37.5){\line(1,0){8}}
\put(137,27.5){\line(1,0){6}}
\put(137,37.5){\arc{5}{4.712}{1.571}}
\put(137,27.5){\arc{5}{4.712}{1.571}}
\put(145,22){\vector(0,-1){0}}
\put(165,34){\vector(0,1){0}}
\put(149.8,34){\vector(0,1){0}}
\Thicklines
\put(145,15){\line(0,1){5}}
\put(165,17){\line(0,1){1}}
\thicklines
\put(142,5){\line(1,0){26}}
\put(142,45){\line(1,0){26}}
\put(142,15){\line(1,0){1}}
\put(142,20){\line(1,0){1}}
\put(147,15){\line(1,0){21}}
\put(147,20){\line(1,0){21}}
\put(142,25){\arc{10}{1.571}{3.142}}
\put(168,25){\arc{10}{0}{1.571}}
\put(142,40){\arc{10}{3.142}{4.712}}
\put(168,40){\arc{10}{4.712}{0}}
\put(137,25){\line(0,1){15}}
\put(173,25){\line(0,1){15}}
\put(142,10){\arc{10}{1.571}{4.712}}
\put(168,10){\arc{10}{4.712}{1.571}}
\put(133,32.5){\makebox(0,0){$a$}}
\put(155,2){\makebox(0,0){$c$}}
\put(155,48){\makebox(0,0){$c$}}
\put(135,46){\makebox(0,0){$b$}}
\put(155,12){\makebox(0,0){$b$}}
\put(155,24){\makebox(0,0){$b$}}
\put(142,24){\makebox(0,0){$\la$}}
\put(169,24){\makebox(0,0){$\mu$}}
\put(155,32.5){\makebox(0,0){$\nu$}}
\end{picture}
\end{center}
\caption{Commutation of $Y$ and $Z$}
\label{[Y,Z]1}
\end{figure}
Using Lemma \ref{idexp} and also the trick to turn around
the small arcs given in Fig.\ \ref{trick}, we
obtain the right-hand side of Fig.\ \ref{[Y,Z]1}. We can
now slide around the lower trivalent vertex of the wire
$\nu$ to obtain the left-hand side of Fig.\ \ref{[Y,Z]2}.
%
% Z commutes with Y (2)
\begin{figure}[htb]
\begin{center}
\unitlength 0.7mm
\begin{picture}(185,50)
\thinlines
\put(16,23){\makebox(0,0){$\displaystyle
\sum_{a,b,c,\la} \;\frac{d_a d_b d_c}{w d_\mu d_\nu}$}}
\put(50,5){\line(0,1){10}}
\put(50,45){\line(0,-1){6}}
\put(50,35){\line(0,-1){15}}
\put(70,5){\line(0,1){8}}
\put(70,45){\line(0,-1){6}}
\put(70,22){\line(0,1){13}}
\put(50,45){\arc{5}{0}{3.142}}
\put(70,45){\arc{5}{0}{3.142}}
\put(50,5){\arc{5}{3.142}{0}}
\put(70,5){\arc{5}{3.142}{0}}
\put(42,37){\line(1,0){36}}
\put(42,37){\arc{5}{4.712}{1.571}}
\put(78,37){\arc{5}{1.571}{4.712}}
\put(50,27){\vector(0,-1){0}}
\put(70,31){\vector(0,1){0}}
\put(58,37){\vector(-1,0){0}}
\Thicklines
\put(50,15){\line(0,1){5}}
\put(70,17){\line(0,1){1}}
\thicklines
\put(47,5){\line(1,0){26}}
\put(47,45){\line(1,0){26}}
\put(47,15){\line(1,0){1}}
\put(47,20){\line(1,0){1}}
\put(52,15){\line(1,0){21}}
\put(52,20){\line(1,0){21}}
\put(47,25){\arc{10}{1.571}{3.142}}
\put(73,25){\arc{10}{0}{1.571}}
\put(47,40){\arc{10}{3.142}{4.712}}
\put(73,40){\arc{10}{4.712}{0}}
\put(42,25){\line(0,1){15}}
\put(78,25){\line(0,1){15}}
\put(47,10){\arc{10}{1.571}{4.712}}
\put(73,10){\arc{10}{4.712}{1.571}}
\put(60,23){\makebox(0,0){$a$}}
\put(60,2){\makebox(0,0){$c$}}
\put(60,48){\makebox(0,0){$c$}}
\put(40,46){\makebox(0,0){$b$}}
\put(80,46){\makebox(0,0){$b$}}
\put(60,12){\makebox(0,0){$b$}}
\put(46,25){\makebox(0,0){$\la$}}
\put(74,24){\makebox(0,0){$\mu$}}
\put(60,32.5){\makebox(0,0){$\nu$}}
\thinlines
\put(110,23){\makebox(0,0){$= \;\displaystyle
\sum_{a,b,c,\rho} \;\frac{d_a d_b d_c}{w d_\mu d_\nu}$}}
\put(172,5){\line(0,1){6}}
\put(172,15){\line(0,1){20}}
\put(172,45){\line(0,-1){6}}
\put(172,45){\arc{5}{0}{3.142}}
\put(172,5){\arc{5}{3.142}{0}}
\put(140,13){\line(1,0){8}}
\put(180,13){\line(-1,0){23}}
\put(140,13){\arc{5}{4.712}{1.571}}
\put(180,13){\arc{5}{1.571}{4.712}}
\put(140,37){\line(1,0){10}}
\put(155,37){\line(1,0){25}}
\put(140,37){\arc{5}{4.712}{1.571}}
\put(180,37){\arc{5}{1.571}{4.712}}
\put(172,27){\vector(0,1){0}}
\put(166,13){\vector(1,0){0}}
\put(162,37){\vector(-1,0){0}}
\Thicklines
\put(150,37){\line(1,0){5}}
\put(152,13){\line(1,0){1}}
\thicklines
\put(160,5){\line(1,0){15}}
\put(160,45){\line(1,0){15}}
\put(145,40){\arc{10}{3.142}{0}}
\put(145,10){\arc{10}{0}{3.142}}
\put(175,40){\arc{10}{4.712}{0}}
\put(175,10){\arc{10}{0}{1.571}}
\put(160,40){\arc{10}{3.142}{4.712}}
\put(160,10){\arc{10}{1.571}{3.142}}
\put(140,10){\line(0,1){30}}
\put(180,10){\line(0,1){30}}
\put(150,10){\line(0,1){25}}
\put(150,39){\line(0,1){1}}
\put(155,10){\line(0,1){25}}
\put(155,39){\line(0,1){1}}
\put(158,24){\makebox(0,0){$c$}}
\put(137,25){\makebox(0,0){$a$}}
\put(184,25){\makebox(0,0){$a$}}
\put(182,46){\makebox(0,0){$b$}}
\put(182,4){\makebox(0,0){$b$}}
\put(164,17){\makebox(0,0){$\rho$}}
\put(147,25){\makebox(0,0){$b$}}
\put(176,24){\makebox(0,0){$\mu$}}
\put(164,33){\makebox(0,0){$\nu$}}
\end{picture}
\end{center}
\caption{Commutation of $Y$ and $Z$}
\label{[Y,Z]2}
\end{figure}
Next, we can use Lemma \ref{idexp} to replace the two parallel
horizontal wires with labels $a$ and $b$ by a summation over
a thin wire $\rho$. Similarly, but the other way round, we can
then use Lemma \ref{idexp} to replace the summation over the
wire with label $\la$ by two straight horizontal wires with
labels $b$ and $c$. This way we obtain the right-hand side of
Fig.\ \ref{[Y,Z]2}. Now it should be clear how to proceed:
We slide around the upper trivalent vertex of the wire $\mu$
counter-clockwise. Then we see that the result gives us
the diagram for $\sum_\rho Z_{\nu,\rho} Y_{\rho,\mu}$,
rotated by 90 degrees. This proves $YZ=ZY$. Next we show
commutativity of $Z$ with $T$. We have to show
$\om_\la Z_{\la,\mu}=Z_{\la,\mu}\om_\mu$. Using the
graphical expression for the statistics phase $\om_\la$
on the left-hand side of Fig.\ \ref{statph}, we can represent
$\om_\la Z_{\la,\mu}$ by the left-hand side of Fig.\ \ref{[T,Z]}.
%
% Z commutes with T
\begin{figure}[htb]
\begin{center}
\unitlength 0.7mm
\begin{picture}(169,50)
\thinlines
\put(13.5,23){\makebox(0,0){$\displaystyle
\sum_{b,c} \;\frac{d_b d_c}{w d_\la d_\mu}$}}
\put(45,5){\line(0,1){10}}
\put(45,45){\line(0,-1){4.5}}
\put(45,37.5){\line(0,-1){17.5}}
\put(55,5){\line(0,1){8}}
\put(55,45){\line(0,-1){23}}
\put(45,45){\arc{5}{0}{3.142}}
\put(55,45){\arc{5}{0}{3.142}}
\put(45,5){\arc{5}{3.142}{0}}
\put(55,5){\arc{5}{3.142}{0}}
\put(41,37.5){\arc{8}{0.9}{6.283}}
\put(45,25){\vector(0,-1){0}}
\put(55,29){\vector(0,1){0}}
\Thicklines
\put(45,15){\line(0,1){5}}
\put(55,17){\line(0,1){1}}
\thicklines
\put(37,5){\line(1,0){26}}
\put(37,45){\line(1,0){26}}
\put(37,15){\line(1,0){6}}
\put(37,20){\line(1,0){6}}
\put(47,15){\line(1,0){16}}
\put(47,20){\line(1,0){16}}
\put(37,25){\arc{10}{1.571}{3.142}}
\put(63,25){\arc{10}{0}{1.571}}
\put(37,40){\arc{10}{3.142}{4.712}}
\put(63,40){\arc{10}{4.712}{0}}
\put(32,25){\line(0,1){15}}
\put(68,25){\line(0,1){15}}
\put(37,10){\arc{10}{1.571}{4.712}}
\put(63,10){\arc{10}{4.712}{1.571}}
\put(50,2){\makebox(0,0){$c$}}
\put(50,48){\makebox(0,0){$c$}}
\put(50,12){\makebox(0,0){$b$}}
\put(50,24){\makebox(0,0){$b$}}
\put(40,28){\makebox(0,0){$\la$}}
\put(60,27){\makebox(0,0){$\mu$}}
\thinlines
\put(105,23){\makebox(0,0){$= \;\displaystyle
\sum_{b,c} \;\frac{d_b d_c}{w d_\la d_\mu}$}}
\put(140,5){\line(0,1){10}}
\put(136,34.5){\line(1,0){2}}
\put(142,34.5){\line(1,0){21}}
\put(140,20){\line(0,1){18.5}}
\put(147,5){\line(0,1){8}}
\put(147,22){\line(0,1){4.5}}
\put(155,20){\line(0,1){6.5}}
\put(163,34.5){\arc{5}{1.571}{4.712}}
\put(155,20){\arc{5}{3.142}{0}}
\put(140,5){\arc{5}{3.142}{0}}
\put(147,5){\arc{5}{3.142}{0}}
\put(136,38.5){\arc{8}{1.571}{6.283}}
\put(151,26.5){\arc{8}{3.142}{0}}
\put(140,25){\vector(0,-1){0}}
\put(147,26){\vector(0,1){0}}
\Thicklines
\put(140,15){\line(0,1){5}}
\put(147,17){\line(0,1){1}}
\thicklines
\put(132,5){\line(1,0){26}}
\put(132,45){\line(1,0){26}}
\put(132,15){\line(1,0){6}}
\put(132,20){\line(1,0){6}}
\put(142,15){\line(1,0){16}}
\put(142,20){\line(1,0){16}}
\put(132,25){\arc{10}{1.571}{3.142}}
\put(158,25){\arc{10}{0}{1.571}}
\put(132,40){\arc{10}{3.142}{4.712}}
\put(158,40){\arc{10}{4.712}{0}}
\put(127,25){\line(0,1){15}}
\put(163,25){\line(0,1){15}}
\put(132,10){\arc{10}{1.571}{4.712}}
\put(158,10){\arc{10}{4.712}{1.571}}
\put(143.5,2){\makebox(0,0){$c$}}
\put(145,48.5){\makebox(0,0){$b$}}
\put(152,12){\makebox(0,0){$b$}}
\put(167,27){\makebox(0,0){$c$}}
\put(135,28){\makebox(0,0){$\la$}}
\put(151,24){\makebox(0,0){$\mu$}}
\end{picture}
\end{center}
\caption{Commutation of $T$ and $Z$}
\label{[T,Z]}
\end{figure}
We now start to rotate the upper oval consisting of the
thick wires $b$ and $c$ in a clockwise direction. This way
we obtain the right-hand side of Fig.\ \ref{[T,Z]}.
It should now be clear that, if we continue rotating
to a full rotation by 360 degrees, then we remove the twist from
the wire $\la$ whereas we obtain a twist in the wire
$\mu$ which is of the type displayed on the right-hand side
of Fig.\ \ref{statph}, thus representing $\om_\mu$.
Hence $TZ=ZT$.
\end{proof}

The following is now immediate by Thm.\ \ref{ST}, which
states that in the non-degenerate case matrices
$S=w^{-1/2}Y$ and $T$ provide a unitary representation
of the modular group $\SLZ$.

\begin{corollary}
\label{modular6}
If the braiding on $\NXN$ is non-degenerate, then the matrix $Z$
defined in Definition \ref{modular1} is a modular invariant mass matrix.
\end{corollary}
In conformal field theory the $\SLZ$ action arises from a
``reparametrization of the torus'', and in the parameter space
$S$ corresponds to a 90 degree rotation and $T$ to twisting
the torus. Note that this action is nicely reflected in
the proof of Thm.\ \ref{modular2}.

\subsection{Generating property of $\a$-induction}
\label{sec-genaind}

We now show that both kinds of $\a$-induction generate the
whole $M$-$M$ fusion rule algebra (or the sector algebra
in our terminology of \cite{BE1,BE2,BE3}) in the case that
the $N$-$N$ system is non-degenerately braided.
That is, from now on we work with the following

\begin{assumption}
\label{set-nondeg}{\rm
In addition to Assumption \ref{set-braid}, we now assume 
that the braiding on $\NXN$ is non-degenerate in the sense of
Definition \ref{non-deg}.
}\end{assumption}
With Assumption \ref{set-nondeg} we can now use the
``killing ring'', the orthogonality relation of
Fig.\ \ref{ort0}, and this turns out to be a powerful
tool in the graphical framework.

The following theorem states in particular that any
minimal central projection $e_\beta$ of $(\dta,*_h)$
appears in the linear decomposition of some
$p_\la^+ *_v p_\mu^-$. Such a generating property of
$p_j^\pm$'s has also been noticed by Ocneanu in the
setting of the lectures \cite{O7}. We can apply his idea
of the proof (which is not included in the notes \cite{O7})
to our situation without essential change.

\begin{theorem}
\label{generating}
Under Assumption \ref{set-nondeg}, we have
$\sum_{\la,\mu\in\NXN} p_\la^+ *_v p_\mu^- = w\bfe_h$
in $\dta$, and consequently
\be
\sum_{\la,\mu\in\NXN} d_\la d_\mu
[\a^+_\la] [\a^-_\mu] = w \sum_{\beta\in\MXM} d_{\beta} [\beta]
\ee
in the $M$-$M$ fusion rule algebra.
In particular, for any $\beta\in\MXM$ the sector
$[\beta]$ is a subsector of $[\a^+_\la][\a^-_\mu]$
for some $\la,\mu\in\NXN$.
\end{theorem}

\begin{proof}
The sum  $\sum_{\la,\mu} p^+_\la *_v p^-_\mu$ is given graphically
by the left-hand side of Fig.\ \ref{gen1}.
%
%proof of generating property (1)
\thinlines
\begin{figure}[htb]
\begin{center}
\unitlength 0.6mm
\begin{picture}(191,50)
\thicklines
\put(13,22){\makebox(0,0){$\displaystyle\sum_{a,b,c,\la,\mu}d_b$}}
\put(40,0){\line(0,1){8}}
\put(40,12){\line(0,1){38}}
\put(60,0){\line(0,1){8}}
\put(60,12){\line(0,1){38}}
\Thicklines
\put(40,10){\line(1,0){20}}
\put(52,10){\vector(1,0){0}}
\put(42,30){\line(1,0){16}}
\put(52,30){\vector(1,0){0}}
\thinlines
\put(40,15){\arc{10}{1.571}{4.712}}
\put(60,15){\arc{10}{4.712}{1.571}}
\put(40,20){\arc{5}{1.571}{4.712}}
\put(60,20){\arc{5}{4.712}{1.571}}
\put(40,35){\arc{10}{2.071}{4.712}}
\put(60,35){\arc{10}{4.712}{1.071}}
\put(40,40){\arc{5}{1.571}{4.712}}
\put(60,40){\arc{5}{4.712}{1.571}}
\put(35,47){\makebox(0,0){$a$}}
\put(65,47){\makebox(0,0){$a$}}
\put(35,25){\makebox(0,0){$b$}}
\put(65,25){\makebox(0,0){$b$}}
\put(35,3){\makebox(0,0){$c$}}
\put(65,3){\makebox(0,0){$c$}}
\put(50,16){\makebox(0,0){$\a_\la^+$}}
\put(50,36){\makebox(0,0){$\a_\mu^-$}}
\put(95,22){\makebox(0,0){$=\;\displaystyle\sum_{a,b,c,\la,\mu,\nu}d_b$}}
\thicklines
\put(130,0){\arc{10}{3.142}{4.712}}
\put(130,50){\arc{10}{1.571}{3.142}}
\put(180,0){\arc{10}{4.712}{0}}
\put(180,50){\arc{10}{0}{1.571}}
\put(130,5){\line(1,0){50}}
\put(130,45){\line(1,0){50}}
\thinlines
\put(155,13){\vector(0,-1){8}}
\put(155,45){\line(0,-1){28}}
\put(155,5){\arc{5}{3.142}{0}}
\put(155,45){\arc{5}{0}{3.142}}
\put(135,45){\line(0,-1){25}}
\put(175,45){\line(0,-1){25}}
\put(140,15){\line(1,0){30}}
\put(140,20){\arc{10}{1.571}{3.142}}
\put(170,20){\arc{10}{0}{1.571}}
\put(135,45){\arc{5}{0}{3.142}}
\put(175,45){\arc{5}{0}{3.142}}
\put(150,15){\vector(1,0){0}}
\put(145,45){\line(0,-1){10}}
\put(165,35){\vector(0,1){10}}
\put(150,30){\line(1,0){3}}
\put(157,30){\line(1,0){3}}
\put(150,35){\arc{10}{1.571}{3.142}}
\put(160,35){\arc{10}{0}{1.571}}
\put(145,45){\arc{5}{0}{3.142}}
\put(165,45){\arc{5}{0}{3.142}}
\put(121,3){\makebox(0,0){$c$}}
\put(189,3){\makebox(0,0){$c$}}
\put(121,47){\makebox(0,0){$a$}}
\put(189,47){\makebox(0,0){$a$}}
\put(150,48){\makebox(0,0){$c$}}
\put(160,48){\makebox(0,0){$c$}}
\put(140,41){\makebox(0,0){$b$}}
\put(170,41){\makebox(0,0){$b$}}
\put(148,20){\makebox(0,0){$\la$}}
\put(163,26){\makebox(0,0){$\mu$}}
\put(159,10){\makebox(0,0){$\nu$}}
\end{picture}
\end{center}
\caption{The sum  $\sum_{\la,\mu} p^+_\la *_v p^-_\mu$}
\label{gen1}
\end{figure}
By using Lemma \ref{idexp} for the two parallel vertical
wires $c$ on the bottom and the IBFE moves we obtain the
right-hand side of Fig.\ \ref{gen1}. For the summation over
the thin wire $\la$ we can use Lemma \ref{idexp} again
to obtain the left-hand side of Fig.\ \ref{gen2}.
%
%proof of generating property (2)
\thinlines
\begin{figure}[htb]
\begin{center}
\unitlength 0.6mm
\begin{picture}(216,50)
\thicklines
\put(13,22){\makebox(0,0){$\displaystyle\sum_{a,b,c,\mu,\nu}d_b$}}
\put(35,50){\line(0,-1){30}}
\put(95,50){\line(0,-1){30}}
\put(40,40){\line(0,-1){15}}
\put(90,40){\line(0,-1){15}}
\put(45,20){\line(1,0){40}}
\put(45,25){\arc{10}{1.571}{3.142}}
\put(85,25){\arc{10}{0}{1.571}}
\put(40,15){\line(1,0){50}}
\put(45,40){\arc{10}{3.142}{4.712}}
\put(85,40){\arc{10}{4.712}{0}}
\put(40,0){\arc{10}{3.142}{4.712}}
\put(40,20){\arc{10}{1.571}{3.142}}
\put(90,0){\arc{10}{4.712}{0}}
\put(90,20){\arc{10}{0}{1.571}}
\put(40,5){\line(1,0){50}}
\put(45,45){\line(1,0){40}}
\Thicklines
\put(65,17){\line(0,1){1}}
\thinlines
\put(65,13){\vector(0,-1){8}}
\put(65,45){\line(0,-1){23}}
\put(65,5){\arc{5}{3.142}{0}}
\put(65,45){\arc{5}{0}{3.142}}
\put(55,45){\line(0,-1){10}}
\put(75,35){\vector(0,1){10}}
\put(60,30){\line(1,0){3}}
\put(67,30){\line(1,0){3}}
\put(60,35){\arc{10}{1.571}{3.142}}
\put(70,35){\arc{10}{0}{1.571}}
\put(55,45){\arc{5}{0}{3.142}}
\put(75,45){\arc{5}{0}{3.142}}
\put(31,3){\makebox(0,0){$c$}}
\put(99,3){\makebox(0,0){$c$}}
\put(31,47){\makebox(0,0){$a$}}
\put(99,47){\makebox(0,0){$a$}}
\put(60,48){\makebox(0,0){$c$}}
\put(70,48){\makebox(0,0){$c$}}
\put(50,25){\makebox(0,0){$b$}}
\put(73,26){\makebox(0,0){$\mu$}}
\put(69,10){\makebox(0,0){$\nu$}}
\thicklines
\put(123,22){\makebox(0,0){$=\displaystyle\sum_{a,b,c,\mu,\nu}d_b$}}
\put(150,50){\line(0,-1){30}}
\put(210,50){\line(0,-1){30}}
\put(155,40){\line(0,-1){15}}
\put(205,40){\line(0,-1){15}}
\put(160,20){\line(1,0){40}}
\put(160,25){\arc{10}{1.571}{3.142}}
\put(200,25){\arc{10}{0}{1.571}}
\put(155,15){\line(1,0){50}}
\put(160,40){\arc{10}{3.142}{4.712}}
\put(200,40){\arc{10}{4.712}{0}}
\put(155,0){\arc{10}{3.142}{4.712}}
\put(155,20){\arc{10}{1.571}{3.142}}
\put(205,0){\arc{10}{4.712}{0}}
\put(205,20){\arc{10}{0}{1.571}}
\put(155,5){\line(1,0){50}}
\put(160,45){\line(1,0){40}}
\Thicklines
\put(180,17){\line(0,1){1}}
\thinlines
\put(180,13){\vector(0,-1){8}}
\put(180,45){\line(0,-1){15}}
\put(180,22){\line(0,1){4}}
\put(180,5){\arc{5}{3.142}{0}}
\put(180,45){\arc{5}{0}{3.142}}
\put(170,45){\line(0,-1){4}}
\put(163,33){\vector(0,1){12}}
\put(175,36){\line(1,0){3}}
\put(182,36){\line(1,0){3}}
\put(175,41){\arc{10}{1.571}{3.142}}
\put(185,32){\arc{8}{4.712}{1.571}}
\put(168,28){\line(1,0){17}}
\put(168,33){\arc{10}{1.571}{3.142}}
\put(170,45){\arc{5}{0}{3.142}}
\put(163,45){\arc{5}{0}{3.142}}
\put(146,3){\makebox(0,0){$c$}}
\put(214,3){\makebox(0,0){$c$}}
\put(146,47){\makebox(0,0){$a$}}
\put(214,47){\makebox(0,0){$a$}}
\put(175,48){\makebox(0,0){$c$}}
\put(185,48){\makebox(0,0){$c$}}
\put(166.5,40){\makebox(0,0){$b$}}
\put(194,31){\makebox(0,0){$\mu$}}
\put(184,10){\makebox(0,0){$\nu$}}
\end{picture}
\end{center}
\caption{The sum  $\sum_{\la,\mu} p^+_\la *_v p^-_\mu$}
\label{gen2}
\end{figure}
Now we can slide around the right trivalent vertex of
the wire $\mu$, and this yields the right-hand side
of Fig.\ \ref{gen2}. Next we can use the trick of
Fig.\ \ref{trick} to turn around the small
arcs from the wire $\mu$ to the wire $b$. This yields
a factor $d_\mu/d_b$. Then we can proceed with the
summation over $b$, using Lemma \ref{idexp} once more,
and this gives us the left-hand side
of Fig.\ \ref{gen3}.
%
%proof of generating property (3)
\thinlines
\begin{figure}[htb]
\begin{center}
\unitlength 0.6mm
\begin{picture}(207,50)
\thicklines
\put(13,22){\makebox(0,0){$\displaystyle\sum_{a,c,\mu,\nu}d_\mu$}}
\put(35,50){\line(0,-1){30}}
\put(95,50){\line(0,-1){30}}
\put(40,40){\line(0,-1){15}}
\put(90,40){\line(0,-1){15}}
\put(45,20){\line(1,0){40}}
\put(45,25){\arc{10}{1.571}{3.142}}
\put(85,25){\arc{10}{0}{1.571}}
\put(40,15){\line(1,0){50}}
\put(45,40){\arc{10}{3.142}{4.712}}
\put(85,40){\arc{10}{4.712}{0}}
\put(40,0){\arc{10}{3.142}{4.712}}
\put(40,20){\arc{10}{1.571}{3.142}}
\put(90,0){\arc{10}{4.712}{0}}
\put(90,20){\arc{10}{0}{1.571}}
\put(40,5){\line(1,0){50}}
\put(45,45){\line(1,0){40}}
\Thicklines
\put(65,17){\line(0,1){1}}
\thinlines
\put(65,13){\vector(0,-1){8}}
\put(65,45){\line(0,-1){15.5}}
\put(65,22){\line(0,1){3.5}}
\put(65,5){\arc{5}{3.142}{0}}
\put(65,45){\arc{5}{0}{3.142}}
\put(60,27.5){\line(1,0){10}}
\put(66,27.5){\vector(-1,0){0}}
\put(60,37.5){\line(1,0){3}}
\put(67,37.5){\line(1,0){3}}
\put(60,32.5){\arc{10}{1.571}{4.7122}}
\put(70,32.5){\arc{10}{4.712}{1.571}}
\put(31,3){\makebox(0,0){$c$}}
\put(99,3){\makebox(0,0){$c$}}
\put(31,47){\makebox(0,0){$a$}}
\put(99,47){\makebox(0,0){$a$}}
\put(50,24){\makebox(0,0){$c$}}
\put(80,31.5){\makebox(0,0){$\mu$}}
\put(69,10){\makebox(0,0){$\nu$}}
\thicklines
\put(125,22){\makebox(0,0){$=\;\displaystyle\sum_{a,c}\;\frac w{d_c}$}}
\put(150,50){\line(0,-1){30}}
\put(200,50){\line(0,-1){30}}
\put(155,40){\line(0,-1){15}}
\put(195,40){\line(0,-1){15}}
\put(160,20){\line(1,0){30}}
\put(160,25){\arc{10}{1.571}{3.142}}
\put(190,25){\arc{10}{0}{1.571}}
\put(155,15){\line(1,0){40}}
\put(160,40){\arc{10}{3.142}{4.712}}
\put(190,40){\arc{10}{4.712}{0}}
\put(155,0){\arc{10}{3.142}{4.712}}
\put(155,20){\arc{10}{1.571}{3.142}}
\put(195,0){\arc{10}{4.712}{0}}
\put(195,20){\arc{10}{0}{1.571}}
\put(155,5){\line(1,0){40}}
\put(160,45){\line(1,0){30}}
\put(146,3){\makebox(0,0){$c$}}
\put(204,3){\makebox(0,0){$c$}}
\put(146,47){\makebox(0,0){$a$}}
\put(204,47){\makebox(0,0){$a$}}
\put(175,25){\makebox(0,0){$c$}}
\end{picture}
\end{center}
\caption{The sum  $\sum_{\la,\mu} p^+_\la *_v p^-_\mu$}
\label{gen3}
\end{figure}
Now we observe that the summation over $\mu$ provides a
killing ring, and hence we obtain a factor $w\del \nu 0$.
The normalization convention for the small
arcs yields another factor $1/d_c$, and hence we get
exactly the right-hand side of Fig.\ \ref{gen3}. The
circular wire $c$ cancels the factor $1/d_c$, and thus
we are left exactly with the global index $w$ times
a summation over two straight horizontal wires, and the
latter is exactly the horizontal unit
$\bfe_h=\sum_\beta e_\beta$. The rest is application of
the isomorphism $\Phi$.
\end{proof}

We remark that the non-degeneracy of the braiding played
an essential role in the proof. In fact there are
counter-examples showing that the generating property
does not hold in general if the braiding is degenerate
(e.g.\ the finite group case discussed in Section 4.2
of \cite{BE1} serves as such an example).

\section{Representations of the $M$-$M$ Fusion Rule Algebra}
\label{sec-repMM}

\subsection{Irreducible representations of the $M$-$M$ fusion rules}
\label{sec-dfa}

We next study in detail the algebra $(\cZ_h,*_v)$ or, equivalently,
the $M$-$M$ fusion rule algebra in the case that
the $N$-$N$ system is non-degenerately braided.
Note that the Assumption \ref{set-braid} implies in particular
that the $N$-$N$ fusion rules algebra is Abelian.
However, the $M$-$M$ fusion rules are in general non-commutative,
and therefore so is the center $(\cZ_h,*_v)$.
We are now going to decompose $(\cZ_h,*_v)$ in simple matrix
algebras. Note that such a decomposition of $(\cZ_h,*_v)$
is equivalent to the determination of the irreducible
representations of the $M$-$M$ fusion rule algebra.

We need some preparation. As in the graphical setting for the double
triangle algebra, we can consider the diagram in Fig.\ \ref{Ombcts}
%
%The vector \Om bcts\la\mu
\begin{figure}[htb]
\begin{center}
\unitlength 0.6mm
\begin{picture}(70,40)
\thicklines
\put(17,0){\line(0,1){15}}
\put(63,0){\line(0,1){15}}
\put(27,15){\arc{20}{3.142}{4.712}}
\put(53,15){\arc{20}{4.712}{0}}
\put(27,15){\arc{10}{1.571}{4.712}}
\put(53,15){\arc{10}{4.712}{1.571}}
\put(27,10){\line(1,0){26}}
\put(27,20){\line(1,0){1}}
\put(32,20){\line(1,0){21}}
\put(27,25){\line(1,0){1}}
\put(32,25){\line(1,0){21}}
\thinlines
\put(30,10){\line(0,1){10}}
\put(30,25){\line(0,1){15}}
\put(50,10){\line(0,1){8}}
\put(50,27){\line(0,1){13}}
\put(30,30.5){\vector(0,-1){0}}
\put(50,34.5){\vector(0,1){0}}
\Thicklines
\put(50,22){\line(0,1){1}}
\put(30,20){\line(0,1){5}}
\put(67,4){\makebox(0,0){$a$}}
\put(40,7){\makebox(0,0){$c$}}
\put(40,16){\makebox(0,0){$b$}}
\put(30,5.8){\makebox(0,0){$t$}}
\put(50,5){\makebox(0,0){$s$}}
\put(25,33.5){\makebox(0,0){$\la$}}
\put(55,32.5){\makebox(0,0){$\mu$}}
\put(3,20){\makebox(0,0){$\displaystyle\sum_a$}}
\end{picture}
\end{center}
\caption{The vector $\Om_{b,c,t,s}^{\la,\mu}\in\cH_{\la,\mu}$}
\label{Ombcts}
\end{figure}
as a vector $\Om_{b,c,t,s}^{\la,\mu}\in\cH_{\la,\mu}$, where
$\cH_{\la,\mu}$ is the vector space
$\cH_{\la,\mu} = \bigoplus_{a\in\NXM} \Hom(\la\co\mu,a\co a)$,
$\la,\mu\in\NXN$.
Here $b,c\in\NXM$, and $t\in\Hom(\la,b\co c)$ and
$s\in\Hom(\co\mu,c\co b)$ are isometries labelling the
two trivalent vertices in Fig.\ \ref{Ombcts}. It is important
to notice that we do not allow coefficients depending on $a$: The
same isometries $t,s$ are used in each block $\Hom(\la\co\mu,a\co a)$
of $\cH_{\la,\mu}$. We next define the subspace
$H_{\la,\mu}\subset\cH_{\la,\mu}$ spanned by such vectors:
\[ H_{\la,\mu} = {\rm span} \{
\Om_{b,c,t,s}^{\la,\mu} \,\, | \,\, b,c\in\NXM\,,\,\,
t\in\Hom(\la,b\co c)\,,\,\,s\in\Hom(\co\mu,c\co b) \} \,. \]
Take two such vectors $\Om_{b,c,t,s}^{\la,\mu}$ and
$\Om_{b',c',t',s'}^{\la,\mu}$. We define an element
$|\Om_{b',c',t',s'}^{\la,\mu}\rangle
\langle\Om_{b,c,t,s}^{\la,\mu}|\in\dta$ by
the diagram in Fig.\ \ref{OmOmdta}.
%
%The element |\Om><\Om| in DTA
\begin{figure}[htb]
\begin{center}
\unitlength 0.7mm
\begin{picture}(70,50)
\thicklines
\put(17,0){\line(0,1){10}}
\put(63,0){\line(0,1){10}}
\put(27,10){\arc{20}{3.142}{4.712}}
\put(53,10){\arc{20}{4.712}{0}}
\put(27,10){\arc{10}{1.571}{4.712}}
\put(53,10){\arc{10}{4.712}{1.571}}
\put(27,40){\arc{20}{1.571}{3.142}}
\put(53,40){\arc{20}{0}{1.571}}
\put(27,40){\arc{10}{1.571}{4.712}}
\put(53,40){\arc{10}{4.712}{1.571}}
\put(27,5){\line(1,0){26}}
\put(27,15){\line(1,0){1}}
\put(32,15){\line(1,0){21}}
\put(27,20){\line(1,0){1}}
\put(32,20){\line(1,0){21}}
\put(17,50){\line(0,-1){10}}
\put(63,50){\line(0,-1){10}}
\put(27,30){\line(1,0){1}}
\put(32,30){\line(1,0){21}}
\put(27,35){\line(1,0){1}}
\put(32,35){\line(1,0){21}}
\put(27,45){\line(1,0){26}}
\thinlines
\put(30,5){\line(0,1){10}}
\put(30,20){\line(0,1){10}}
\put(50,5){\line(0,1){8}}
\put(50,22){\line(0,1){6}}
\put(30,23){\vector(0,-1){0}}
\put(50,27){\vector(0,1){0}}
\put(30,45){\line(0,-1){10}}
\put(50,45){\line(0,-1){8}}
\Thicklines
\put(50,17){\line(0,1){1}}
\put(30,15){\line(0,1){5}}
\put(50,33){\line(0,-1){1}}
\put(30,35){\line(0,-1){5}}
\put(67,4){\makebox(0,0){$a'$}}
\put(40,2){\makebox(0,0){$c'$}}
\put(40,11){\makebox(0,0){$b'$}}
\put(30,1.8){\makebox(0,0){$t'$}}
\put(50,1.8){\makebox(0,0){$s'$}}
\put(67,46){\makebox(0,0){$a$}}
\put(40,48){\makebox(0,0){$c$}}
\put(40,39){\makebox(0,0){$b$}}
\put(30,48.2){\makebox(0,0){$t^*$}}
\put(50,48.2){\makebox(0,0){$s^*$}}
\put(25,25){\makebox(0,0){$\la$}}
\put(55,24){\makebox(0,0){$\mu$}}
\put(3,25){\makebox(0,0){$\displaystyle\sum_{a,a'}$}}
\end{picture}
\end{center}
\caption{The element $|\Om_{b',c',t',s'}^{\la,\mu}\rangle
\langle\Om_{b,c,t,s}^{\la,\mu}|\in\protect\dta$}
\label{OmOmdta}
\end{figure}
(This notation will be justified by Lemma \ref{innerprod}
below.) We now choose orthonormal bases of isometries
$t_{b,\co c}^{\la;i}\in\Hom(\la,b\co c)$, $i=1,2,...,N_{b,\co c}^\la$,
for each $\la,b,c$ and put
$\Om_\xi^{\la,\mu}
=\Om_{b,c,{t_{b,\co c}^{\la;i}},{t_{c,\co b}^{\co\mu;j}}}^{\la,\mu}$
with some multi-index $\xi=(b,c,i,j)$. Varying $\xi$, we obtain a
generating set of $H_{\la,\mu}$ which will, however, in general not be
a basis as the vectors $\Om_\xi^{\la,\mu}$ may be linearly
dependent in $H_{\la,\mu}$. Let $\Phi_j^{\la,\mu}\in H_{\la,\mu}$,
$j=1,2$, any two vectors. We can expand them as
$\Phi_j^{\la,\mu}=\sum_\xi c_j^\xi \Om^{\la,\mu}_\xi$ with
$c_j^\xi\in\bbC$, but note that this expansion is not unique.
We now define an element
$|\Phi_1^{\la,\mu}\rangle\langle\Phi_2^{\la,\mu}|\in\dta$
by
\be |\Phi_1^{\la,\mu}\rangle\langle\Phi_2^{\la,\mu}| = 
\sum_{\xi,\xi'} c_1^\xi (c_2^{\xi'})^*
|\Om_\xi^{\la,\mu}\rangle\langle\Om_{\xi'}^{\la,\mu}| \,,
\label{|P><P|}
\ee
and a scalar
$\langle\Phi_2^{\la,\mu},\Phi_1^{\la,\mu}\rangle\in\bbC$,
\be
 \langle\Phi_2^{\la,\mu},\Phi_1^{\la,\mu}\rangle =
\frac 1{d_\la d_\mu} \, \tau_v
(|\Phi_1^{\la,\mu}\rangle\langle\Phi_2^{\la,\mu}|) \,.
\label{<P,P>}
\ee

\begin{lemma}
\label{innerprod}
\erf{|P><P|} extends to a sesqui-linear map
$H_{\la,\mu}\times H_{\la,\mu}\rightarrow\cZ_h$
which is positive definite: If
$|\Phi^{\la,\mu}\rangle\langle\Phi^{\la,\mu}|=0$
for some $\Phi^{\la,\mu}\in H_{\la,\mu}$ then $\Phi^{\la,\mu}=0$.
Consequently, \erf{<P,P>} defines a scalar product
turning $H_{\la,\mu}$ into a Hilbert space.
\end{lemma}

\begin{proof}
As in particular $\Phi_j\in\cH_{\la,\mu}$, we can write
$\Phi_j=\bigoplus_a (\Phi_j)_a$ with
$(\Phi_j)_a\in\Hom(\la\co\mu,a\co a)$ according to the direct sum
structure of $\cH_{\la,\mu}$, $j=1,2$.
Assume $\Phi_1=0$. Then clearly $(\Phi_1)_a=0$ for all $a$.
Now the $\Hom(a\co a,a' \co{a'})$ part of
$|\Phi_1^{\la,\mu}\rangle\langle\Phi_2^{\la,\mu}|\in\dta$
is given by $(\Phi_1)_{a'} (\Phi_2)_a^*$, hence
$|\Phi_1^{\la,\mu}\rangle\langle\Phi_2^{\la,\mu}|=0$.
A similar argument applies to $\Phi_2$, and hence
the element
$|\Phi_1^{\la,\mu}\rangle\langle\Phi_2^{\la,\mu}|\in\dta$
is independent of the linear expansions of the $\Phi_j$'s.
Therefore \erf{|P><P|} defines a sesqui-linear map
$H_{\la,\mu}\times H_{\la,\mu}\rightarrow\dta$.
Now assume $|\Phi_1^{\la,\mu}\rangle\langle\Phi_1^{\la,\mu}|=0$.
Then in particular $(\Phi_1)_a (\Phi_1)_a^*=0$ for all $a\in\NXM$,
and hence $\Phi_1=0$, proving strict positivity. That the
sesqui-linear form $\langle\cdot,\cdot\rangle$ on $H_{\la,\mu}$
is non-degenerate follows now from positive definiteness of $\tau_v$.
It remains to show that
$|\Phi_1^{\la,\mu}\rangle\langle\Phi_2^{\la,\mu}|\in\cZ_h$.
But this is clear since any element of the form in Fig.\ \ref{dta1}
can be ``pulled through'' the diagram in Fig.\ \ref{OmOmdta}
by using the IBFE's.
\end{proof}

\begin{lemma}
\label{keyrel}
We have the identity in Fig.\ \ref{keyrelpict} for
intertwiners in $\Hom(\la'\co{\mu'},\la\co\mu)$,
$\la,\mu,\la',\mu'\in\NXN$.
\end{lemma}
%
% <\Om,\Om> identity
\begin{figure}[htb]
\begin{center}
\unitlength 0.6mm
\begin{picture}(195,65)
\thicklines
\put(10,30){\makebox(0,0){$\displaystyle\sum_a\; d_a$}}
\put(27,20){\line(0,1){25}}
\put(73,20){\line(0,1){25}}
\put(37,45){\arc{20}{3.142}{4.712}}
\put(63,45){\arc{20}{4.712}{0}}
\put(37,45){\arc{10}{1.571}{4.712}}
\put(63,45){\arc{10}{4.712}{1.571}}
\put(37,20){\arc{20}{1.571}{3.142}}
\put(63,20){\arc{20}{0}{1.571}}
\put(37,20){\arc{10}{1.571}{4.712}}
\put(63,20){\arc{10}{4.712}{1.571}}
\put(37,40){\line(1,0){26}}
\put(37,50){\line(1,0){1}}
\put(42,50){\line(1,0){21}}
\put(37,55){\line(1,0){1}}
\put(42,55){\line(1,0){21}}
\put(37,10){\line(1,0){1}}
\put(42,10){\line(1,0){21}}
\put(37,15){\line(1,0){1}}
\put(42,15){\line(1,0){21}}
\put(37,25){\line(1,0){26}}
\thinlines
\put(40,40){\line(0,1){10}}
\put(40,55){\line(0,1){10}}
\put(60,40){\line(0,1){8}}
\put(60,57){\line(0,1){8}}
\put(40,0){\line(0,1){10}}
\put(60,0){\line(0,1){8}}
\put(40,58){\vector(0,-1){0}}
\put(60,62){\vector(0,1){0}}
\put(40,3){\vector(0,-1){0}}
\put(60,7){\vector(0,1){0}}
\put(40,25){\line(0,-1){10}}
\put(60,25){\line(0,-1){8}}
\Thicklines
\put(60,52){\line(0,1){1}}
\put(40,50){\line(0,1){5}}
\put(60,13){\line(0,-1){1}}
\put(40,15){\line(0,-1){5}}
\put(50,37){\makebox(0,0){$c'$}}
\put(50,46){\makebox(0,0){$b'$}}
\put(40,36.8){\makebox(0,0){$t'$}}
\put(60,36.8){\makebox(0,0){$s'$}}
\put(23,17){\makebox(0,0){$a$}}
\put(50,29){\makebox(0,0){$c$}}
\put(50,19){\makebox(0,0){$b$}}
\put(40,29){\makebox(0,0){$t^*$}}
\put(60,29){\makebox(0,0){$s^*$}}
\put(35,5){\makebox(0,0){$\la$}}
\put(65,4){\makebox(0,0){$\mu$}}
\put(35,60){\makebox(0,0){$\la'$}}
\put(65,60){\makebox(0,0){$\mu'$}}
\put(125,32.5){\makebox(0,0){$=\;\del \la{\la'} \, \del \mu{\mu'}
\; \langle\Om_{b,c,t,s}^{\la,\mu},
\Om_{b',c',t',s'}^{\la,\mu}\rangle$}}
\thinlines
\put(175,0){\line(0,1){65}}
\put(185,0){\line(0,1){65}}
\put(175,30.5){\vector(0,-1){0}}
\put(185,34.5){\vector(0,1){0}}
\put(170,10){\makebox(0,0){$\la$}}
\put(190,10){\makebox(0,0){$\mu$}}
\end{picture}
\end{center}
\caption{An identity in $\Hom(\la'\co{\mu'},\la\co\mu)$}
\label{keyrelpict}
\end{figure}

\begin{proof}
Using Lemma \ref{idexp} we can replace the left-hand side
of Fig.\ \ref{keyrelpict} by the left-hand side of
Fig.\ \ref{proofkey1}.
%
% proof of <\Om,\Om> identity (1)
\begin{figure}[htb]
\begin{center}
\unitlength 0.6mm
\begin{picture}(205,65)
\thicklines
\put(10,30){\makebox(0,0){$\displaystyle\sum_{\nu,a}\; d_a$}}
\put(27,15){\line(0,1){20}}
\put(93,15){\line(0,1){20}}
\put(32,15){\arc{10}{1.571}{3.142}}
\put(88,15){\arc{10}{0}{1.571}}
\put(32,35){\arc{10}{3.142}{4.712}}
\put(88,35){\arc{10}{4.712}{0}}
\put(47,20){\arc{10}{1.571}{4.712}}
\put(73,20){\arc{10}{4.712}{1.571}}
\put(32,40){\line(1,0){56}}
\put(32,10){\line(1,0){16}}
\put(52,10){\line(1,0){36}}
\put(47,15){\line(1,0){1}}
\put(52,15){\line(1,0){21}}
\put(47,25){\line(1,0){26}}
\thinlines
\put(37,40){\line(0,1){5}}
\put(83,40){\line(0,1){5}}
\put(37,40){\arc{5}{3.142}{0}}
\put(83,40){\arc{5}{3.142}{0}}
\put(62,55){\vector(1,0){0}}
\put(47,45){\arc{20}{3.142}{4.712}}
\put(73,45){\arc{20}{4.712}{0}}
\put(47,55){\line(1,0){1}}
\put(52,55){\line(1,0){21}}
\put(50,40){\line(0,1){25}}
\put(70,40){\line(0,1){13}}
\put(70,57){\line(0,1){8}}
\put(50,0){\line(0,1){10}}
\put(70,0){\line(0,1){8}}
\put(50,58){\vector(0,-1){0}}
\put(70,62){\vector(0,1){0}}
\put(50,3){\vector(0,-1){0}}
\put(70,7){\vector(0,1){0}}
\put(50,25){\line(0,-1){10}}
\put(70,25){\line(0,-1){8}}
\Thicklines
\put(70,13){\line(0,-1){1}}
\put(50,15){\line(0,-1){5}}
\put(60,37){\makebox(0,0){$c'$}}
\put(45,44){\makebox(0,0){$b'$}}
\put(75,44){\makebox(0,0){$b'$}}
\put(50,36.8){\makebox(0,0){$t'$}}
\put(70,36.8){\makebox(0,0){$s'$}}
\put(23,17){\makebox(0,0){$a$}}
\put(60,29){\makebox(0,0){$c$}}
\put(60,19){\makebox(0,0){$b$}}
\put(50,29){\makebox(0,0){$t^*$}}
\put(70,29){\makebox(0,0){$s^*$}}
\put(45,5){\makebox(0,0){$\la$}}
\put(75,4){\makebox(0,0){$\mu$}}
\put(45,60){\makebox(0,0){$\la'$}}
\put(75,60){\makebox(0,0){$\mu'$}}
\put(60,60){\makebox(0,0){$\nu$}}
\put(118,30){\makebox(0,0){$=\;\displaystyle\sum_{\nu,\rho,\tau}\; d_\nu$}}
\thicklines
\put(142,20){\line(0,1){25}}
\put(198,20){\line(0,1){25}}
\put(147,20){\arc{10}{1.571}{3.142}}
\put(193,20){\arc{10}{0}{1.571}}
\put(147,45){\arc{10}{3.142}{4.712}}
\put(193,45){\arc{10}{4.712}{0}}
\put(152,32.5){\arc{15}{1.571}{4.712}}
\put(188,32.5){\arc{15}{4.712}{1.571}}
\put(147,50){\line(1,0){13}}
\put(180,50){\line(1,0){13}}
\put(147,15){\line(1,0){6}}
\put(157,15){\line(1,0){3}}
\put(180,15){\line(1,0){13}}
\put(152,25){\line(1,0){1}}
\put(157,25){\line(1,0){3}}
\put(188,25){\line(-1,0){8}}
\put(152,40){\line(1,0){8}}
\put(188,40){\line(-1,0){8}}
\put(160,47.5){\arc{5}{4.712}{0}}
\put(160,22.5){\arc{5}{4.712}{0}}
\put(160,42.5){\arc{5}{0}{1.571}}
\put(160,17.5){\arc{5}{0}{1.571}}
\put(180,47.5){\arc{5}{3.142}{4.712}}
\put(180,22.5){\arc{5}{3.142}{4.712}}
\put(180,42.5){\arc{5}{1.571}{3.142}}
\put(180,17.5){\arc{5}{1.571}{3.142}}
\put(162.5,17.5){\line(0,1){5}}
\put(162.5,42.5){\line(0,1){5}}
\put(177.5,17.5){\line(0,1){5}}
\put(177.5,42.5){\line(0,1){5}}
\thinlines
\put(147,20){\arc{20}{1.571}{3.142}}
\put(193,20){\arc{20}{0}{1.571}}
\put(147,45){\arc{20}{3.142}{4.712}}
\put(193,45){\arc{20}{4.712}{0}}
\put(162.5,45){\line(1,0){15}}
\put(162.5,20){\line(1,0){15}}
\put(172,55){\vector(1,0){0}}
\put(162.5,45){\arc{5}{4.712}{1.571}}
\put(162.5,20){\arc{5}{4.712}{1.571}}
\put(177.5,45){\arc{5}{1.571}{4.712}}
\put(177.5,20){\arc{5}{1.571}{4.712}}
\put(172,45){\vector(1,0){0}}
\put(172,20){\vector(1,0){0}}
\put(147,55){\line(1,0){6}}
\put(157,55){\line(1,0){36}}
\put(147,10){\line(1,0){6}}
\put(157,10){\line(1,0){36}}
\put(137,20){\line(0,1){25}}
\put(203,20){\line(0,1){25}}
\put(155,50){\line(0,1){15}}
\put(185,50){\line(0,1){3}}
\put(185,57){\line(0,1){8}}
\put(155,0){\line(0,1){15}}
\put(185,0){\line(0,1){8}}
\put(185,12){\line(0,1){1}}
\put(155,58){\vector(0,-1){0}}
\put(185,62){\vector(0,1){0}}
\put(155,3){\vector(0,-1){0}}
\put(185,7){\vector(0,1){0}}
\put(155,40){\line(0,-1){15}}
\put(185,40){\line(0,-1){13}}
\Thicklines
\put(185,23){\line(0,-1){6}}
\put(155,25){\line(0,-1){10}}
\put(164,52){\makebox(0,0){$c'$}}
\put(176,52){\makebox(0,0){$c'$}}
\put(159,47){\makebox(0,0){$t'$}}
\put(181,47){\makebox(0,0){$s'$}}
\put(160,37){\makebox(0,0){$c$}}
\put(180,37){\makebox(0,0){$c$}}
\put(160,29){\makebox(0,0){$b$}}
\put(180,29){\makebox(0,0){$b$}}
\put(164,13.5){\makebox(0,0){$b'$}}
\put(176,13.5){\makebox(0,0){$b'$}}
\put(152,44){\makebox(0,0){$t^*$}}
\put(188,44){\makebox(0,0){$s^*$}}
\put(150,5){\makebox(0,0){$\la$}}
\put(190,4){\makebox(0,0){$\mu$}}
\put(150,60){\makebox(0,0){$\la'$}}
\put(190,60){\makebox(0,0){$\mu'$}}
\put(170,60){\makebox(0,0){$\nu$}}
\put(170,41){\makebox(0,0){$\rho$}}
\put(170,24){\makebox(0,0){$\tau$}}
\end{picture}
\end{center}
\caption{The identity in $\Hom(\la'\co{\mu'},\la\co\mu)$}
\label{proofkey1}
\end{figure}
Next we can slide one of the trivalent vertices of
the wire $\nu$ around the wire $a$. Using the
identity of Fig.\ \ref{trick}, we obtain a factor
$d_\nu/d_a$, and we can now proceed with the summation
over $a$, again using Lemma \ref{idexp}. Using also
Lemma \ref{idexp} for the parallel wires $c$, $c'$
as well as $b$ and $b'$, we obtain the right-hand side of
Fig.\ \ref{proofkey1}. Using now Lemma \ref{idexp}
once again for the wires $\rho$, $\tau$, we can pull
the wire $\nu$ over the middle expansion. The summation
over $\nu$ yields a killing ring which disconnects the
picture into two halves, one is an intertwiner in
$\Hom(\la',\la)$ and the other in $\Hom(\co{\mu'},\co\mu)$.
Hence we obtain a factor $\del \la{\la'} \del \mu{\mu'}$,
and we conclude that the left-hand side in Fig.\ \ref{keyrelpict}
represents a scalar intertwiner
$\del \la{\la'} \del \mu{\mu'} \zeta\bfe_N
\in\Hom(\la\co\mu,\la\co\mu)$,
$\zeta\in\bbC$. To compute that scalar, we can start again
on the left-hand side of Fig.\ \ref{keyrelpict}, now putting
$\la'=\la$ and $\mu'=\mu$. The diagram on the left-hand side
of Fig.\ \ref{proofkey2} clearly represents an
intertwiner of the same scalar value $\zeta$.
%
% proof of <\Om,\Om> identity (2)
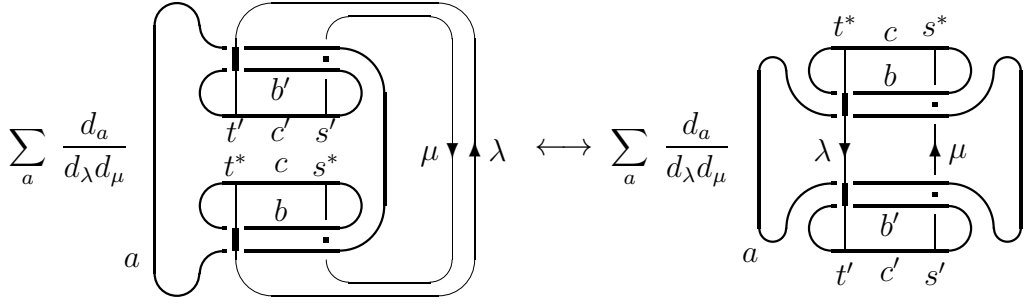
\begin{figure}[htb]
\begin{center}
\unitlength 0.6mm
\begin{picture}(225,65)
\thicklines
\put(13,32.5){\makebox(0,0){$\displaystyle\sum_a\;
\frac{d_a}{d_\la d_\mu}$}}
\put(32,5){\line(0,1){55}}
\put(83,20){\line(0,1){25}}
\put(37,5){\arc{10}{0}{3.142}}
\put(37,60){\arc{10}{3.142}{0}}
\put(47,5){\arc{10}{3.142}{4.712}}
\put(47,60){\arc{10}{1.571}{3.142}}
\put(73,45){\arc{20}{4.712}{0}}
\put(47,45){\arc{10}{1.571}{4.712}}
\put(73,45){\arc{10}{4.712}{1.571}}
\put(73,20){\arc{20}{0}{1.571}}
\put(47,20){\arc{10}{1.571}{4.712}}
\put(73,20){\arc{10}{4.712}{1.571}}
\put(47,40){\line(1,0){26}}
\put(47,50){\line(1,0){1}}
\put(52,50){\line(1,0){21}}
\put(47,55){\line(1,0){1}}
\put(52,55){\line(1,0){21}}
\put(47,10){\line(1,0){1}}
\put(52,10){\line(1,0){21}}
\put(47,15){\line(1,0){1}}
\put(52,15){\line(1,0){21}}
\put(47,25){\line(1,0){26}}
\thinlines
\put(75,57){\arc{10}{3.142}{4.712}}
\put(58,57){\arc{16}{3.142}{4.712}}
\put(75,8){\arc{10}{1.571}{3.142}}
\put(58,8){\arc{16}{1.571}{3.142}}
\put(50,8){\line(0,1){2}}
\put(50,55){\line(0,1){2}}
\put(98,8){\line(0,1){49}}
\put(103,8){\line(0,1){49}}
\put(93,57){\arc{10}{4.712}{0}}
\put(95,57){\arc{16}{4.712}{0}}
\put(93,8){\arc{10}{0}{1.571}}
\put(95,8){\arc{16}{0}{1.571}}
\put(58,0){\line(1,0){37}}
\put(58,65){\line(1,0){37}}
\put(75,3){\line(1,0){18}}
\put(75,62){\line(1,0){18}}
\put(50,40){\line(0,1){10}}
\put(70,40){\line(0,1){8}}
\put(98,30.5){\vector(0,-1){0}}
\put(103,34.5){\vector(0,1){0}}
\put(50,25){\line(0,-1){10}}
\put(70,25){\line(0,-1){8}}
\Thicklines
\put(70,52){\line(0,1){1}}
\put(50,50){\line(0,1){5}}
\put(70,13){\line(0,-1){1}}
\put(50,15){\line(0,-1){5}}
\put(60,37){\makebox(0,0){$c'$}}
\put(60,46){\makebox(0,0){$b'$}}
\put(50,36.8){\makebox(0,0){$t'$}}
\put(70,36.8){\makebox(0,0){$s'$}}
\put(27,8){\makebox(0,0){$a$}}
\put(60,29){\makebox(0,0){$c$}}
\put(60,19){\makebox(0,0){$b$}}
\put(50,29){\makebox(0,0){$t^*$}}
\put(70,29){\makebox(0,0){$s^*$}}
\put(93,31){\makebox(0,0){$\mu$}}
\put(108,32.5){\makebox(0,0){$\la$}}
\put(138,32.5){\makebox(0,0){$\longleftrightarrow\;\displaystyle\sum_a\;
\frac{d_a}{d_\la d_\mu}$}}
\thicklines
\put(166,15){\line(0,1){35}}
\put(224,15){\line(0,1){35}}
\put(169,15){\arc{6}{0}{3.142}}
\put(169,50){\arc{6}{3.142}{0}}
\put(221,15){\arc{6}{0}{3.142}}
\put(221,50){\arc{6}{3.142}{0}}
\put(182,15){\arc{20}{3.142}{4.712}}
\put(208,15){\arc{20}{4.712}{0}}
\put(182,15){\arc{10}{1.571}{4.712}}
\put(208,15){\arc{10}{4.712}{1.571}}
\put(182,50){\arc{20}{1.571}{3.142}}
\put(208,50){\arc{20}{0}{1.571}}
\put(182,50){\arc{10}{1.571}{4.712}}
\put(208,50){\arc{10}{4.712}{1.571}}
\put(182,10){\line(1,0){26}}
\put(182,20){\line(1,0){1}}
\put(187,20){\line(1,0){21}}
\put(182,25){\line(1,0){1}}
\put(187,25){\line(1,0){21}}
\put(182,40){\line(1,0){1}}
\put(187,40){\line(1,0){21}}
\put(182,45){\line(1,0){1}}
\put(187,45){\line(1,0){21}}
\put(182,55){\line(1,0){26}}
\thinlines
\put(185,10){\line(0,1){10}}
\put(185,25){\line(0,1){15}}
\put(205,10){\line(0,1){8}}
\put(205,27){\line(0,1){11}}
\put(185,30.5){\vector(0,-1){0}}
\put(205,34.5){\vector(0,1){0}}
\put(185,55){\line(0,-1){10}}
\put(205,55){\line(0,-1){8}}
\Thicklines
\put(205,22){\line(0,1){1}}
\put(185,20){\line(0,1){5}}
\put(205,43){\line(0,-1){1}}
\put(185,45){\line(0,-1){5}}
\put(164,9){\makebox(0,0){$a$}}
\put(195,6){\makebox(0,0){$c'$}}
\put(195,16){\makebox(0,0){$b'$}}
\put(185,5){\makebox(0,0){$t'$}}
\put(205,5){\makebox(0,0){$s'$}}
\put(195,58){\makebox(0,0){$c$}}
\put(195,49){\makebox(0,0){$b$}}
\put(185,60){\makebox(0,0){$t^*$}}
\put(205,60){\makebox(0,0){$s^*$}}
\put(180,32.5){\makebox(0,0){$\la$}}
\put(210,31.5){\makebox(0,0){$\mu$}}
\end{picture}
\end{center}
\caption{Computation of the scalar $\zeta$}
\label{proofkey2}
\end{figure}
We can now use the move of Fig.\ \ref{cut} which
does not change the scalar value: We open the wire $a$
on the left and close it on the right. The resulting
diagram is regularly isotopic to the diagram on the
right-hand side of Fig.\ \ref{proofkey2}.
Thus we are left with exactly the diagram for
$d_\la^{-1}d_\mu^{-1}\tau_v
(|\Om_{b',c',t',s'}^{\la,\mu}\rangle
\langle\Om_{b,c,t,s}^{\la,\mu}|)$.
This proves the lemma.
\end{proof}

The following is now immediate by the definition of the vertical product.

\begin{corollary}
\label{matunits}
Let $\Phi_j^{\la,\mu}\in H_{\la,\mu}$ and
$\Psi_j^{\la',\mu'}\in H_{\la',\mu'}$, $j=1,2$.
Then we have
\be
|\Phi_1^{\la,\mu}\rangle\langle\Phi_2^{\la,\mu}| *_v
|\Psi_1^{\la',\mu'}\rangle\langle\Psi_2^{\la',\mu'}| \;=\;
\del \la{\la'} \, \del \mu{\mu'} \,
\langle\Phi_2^{\la,\mu},\Psi_1^{\la,\mu}\rangle
\;|\Phi_1^{\la,\mu}\rangle\langle\Psi_2^{\la,\mu}|
\ee
in the double triangle algebra.
\end{corollary}

Whenever $H_{\la,\mu}\neq\{0\}$ we can choose an
orthonormal basis
$\{ E_i^{\la,\mu} \}_{i=1}^{{\rm{dim}} H_{\la,\mu}}$.
Then Lemma \ref{innerprod} and Corollary \ref{matunits} tell us that
$\{\,|E_i^{\la,\mu}\rangle\langle E_j^{\la,\mu}|\,\}_{\la,\mu,i,j}$
forms a set of non-zero matrix units in $(\cZ_h,*_v)$. However,
we do not know yet whether this is a complete set.

\begin{lemma}
\label{miniunit}
Let $\pi_{\la,\mu}(e_\beta)\Om_{b,c,t,s}^{\la,\mu}\in\cH_{\la,\mu}$
denote the vector which is given graphically by the diagram
in Fig.\ \ref{piebOm}, where $\la,\mu\in\NXN$, $b,c\in\NXM$,
and $t\in\Hom(\la,b\co c)$,
$s\in\Hom(\co\mu,c\co b)$ are isometries. Then in fact
$\pi_{\la,\mu}(e_\beta)\Om_{b,c,t,s}^{\la,\mu}\in H_{\la,\mu}$.
\end{lemma}
%
%The vector \pi(e_\beta)\Om bcts\la\mu
\begin{figure}[htb]
\begin{center}
\unitlength 0.6mm
\begin{picture}(86,50)
\thicklines
\put(32,0){\line(0,1){30}}
\put(78,0){\line(0,1){30}}
\put(42,30){\arc{20}{3.142}{4.712}}
\put(68,30){\arc{20}{4.712}{0}}
\put(42,30){\arc{10}{1.571}{4.712}}
\put(68,30){\arc{10}{4.712}{1.571}}
\put(42,25){\line(1,0){26}}
\put(42,35){\line(1,0){1}}
\put(47,35){\line(1,0){21}}
\put(42,40){\line(1,0){1}}
\put(47,40){\line(1,0){21}}
\thinlines
\put(45,25){\line(0,1){10}}
\put(45,40){\line(0,1){10}}
\put(65,25){\line(0,1){8}}
\put(65,42){\line(0,1){8}}
\put(45,43){\vector(0,-1){0}}
\put(65,47){\vector(0,1){0}}
\put(32,10){\arc{5}{4.712}{1.571}}
\put(78,10){\arc{5}{1.571}{4.712}}
\Thicklines
\put(32,10){\line(1,0){46}}
\put(65,37){\line(0,1){1}}
\put(45,35){\line(0,1){5}}
\put(57,10){\vector(1,0){0}}
\put(27,3){\makebox(0,0){$a$}}
\put(82,3){\makebox(0,0){$a$}}
\put(82,22){\makebox(0,0){$a'$}}
\put(55,22){\makebox(0,0){$c$}}
\put(55,31){\makebox(0,0){$b$}}
\put(55,4){\makebox(0,0){$\beta$}}
\put(45,20.8){\makebox(0,0){$t$}}
\put(65,20){\makebox(0,0){$s$}}
\put(40,46){\makebox(0,0){$\la$}}
\put(70,45){\makebox(0,0){$\mu$}}
\put(10,25){\makebox(0,0){$\displaystyle\sum_{a,a'}\,d_{a'}$}}
\end{picture}
\end{center}
\caption{The vector $\pi_{\la,\mu}(e_\beta)
\Om_{b,c,t,s}^{\la,\mu}\in\cH_{\la,\mu}$}
\label{piebOm}
\end{figure}

\begin{proof}
Using Lemma \ref{idexp} and also the trick of
Fig.\ \ref{trick}, we can draw the diagram
on the left-hand side in Fig.\ \ref{proofpieb1} for
$\pi_{\la,\mu}(e_\beta)\Om_{b,c,t,s}^{\la,\mu}$.
%
%proof: the vector \pi(e_\beta)\Om bcts\la\mu (1)
\begin{figure}[htb]
\begin{center}
\unitlength 0.6mm
\begin{picture}(237,55)
\thicklines
\put(11,25){\makebox(0,0){$\displaystyle\sum_{a,a',\nu}\; d_\beta$}}
\put(25,0){\line(0,1){20}}
\put(30,20){\arc{10}{3.142}{4.712}}
\put(30,25){\line(1,0){2.5}}
\put(105,0){\line(0,1){20}}
\put(100,20){\arc{10}{4.712}{0}}
\put(97.5,25){\line(1,0){2.5}}
\put(52.5,15){\line(1,0){25}}
\put(52.5,20){\arc{10}{1.571}{3.142}}
\put(77.5,20){\arc{10}{0}{1.571}}
\put(47.5,20){\line(0,1){10}}
\put(82.5,20){\line(0,1){10}}
\put(42.5,30){\arc{10}{4.712}{0}}
\put(87.5,30){\arc{10}{3.142}{4.712}}
\put(37.5,30){\arc{10}{3.142}{4.712}}
\put(92.5,30){\arc{10}{4.712}{0}}
\put(37.5,35){\line(1,0){5}}
\put(92.5,35){\line(-1,0){5}}
\put(32.5,25){\line(0,1){5}}
\put(97.5,25){\line(0,1){5}}
\Thicklines
\put(32.5,25){\line(1,0){2.5}}
\put(35,20){\arc{10}{4.712}{0}}
\put(40,5){\line(0,1){15}}
\put(45,5){\arc{10}{1.571}{3.142}}
\put(45,0){\line(1,0){40}}
\put(85,5){\arc{10}{0}{1.571}}
\put(90,5){\line(0,1){15}}
\put(95,20){\arc{10}{3.142}{4.712}}
\put(95,25){\line(1,0){2.5}}
\put(67,0){\vector(1,0){0}}
\thinlines
\dottedline{4}(20,20)(110,20)
\put(55,15){\line(0,1){40}}
\put(75,15){\line(0,1){28}}
\put(75,47){\line(0,1){8}}
\put(40,35){\line(0,1){5}}
\put(90,35){\line(0,1){5}}
\put(45,40){\arc{10}{3.142}{4.712}}
\put(85,40){\arc{10}{4.712}{0}}
\put(45,45){\line(1,0){8}}
\put(57,45){\line(1,0){28}}
\put(40,35){\arc{5}{3.142}{0}}
\put(90,35){\arc{5}{3.142}{0}}
\put(32.5,25){\arc{5}{3.142}{0}}
\put(97.5,25){\arc{5}{3.142}{0}}
\put(67,45){\vector(1,0){0}}
\put(55,33){\vector(0,-1){0}}
\put(75,37){\vector(0,1){0}}
\put(20,3){\makebox(0,0){$a$}}
\put(110,3){\makebox(0,0){$a$}}
\put(65,5){\makebox(0,0){$\beta$}}
\put(65,12){\makebox(0,0){$c$}}
\put(55,10){\makebox(0,0){$t$}}
\put(75,10){\makebox(0,0){$s$}}
\put(59,28){\makebox(0,0){$\la$}}
\put(71,26.5){\makebox(0,0){$\mu$}}
\put(65,50){\makebox(0,0){$\nu$}}
\put(34,38){\makebox(0,0){$a'$}}
\put(46,38){\makebox(0,0){$b$}}
\put(84,38){\makebox(0,0){$b$}}
\put(97,38){\makebox(0,0){$a'$}}
\put(130,25){\makebox(0,0){$=\;\displaystyle\sum_{a,a',\nu}\; d_\beta$}}
\thicklines
\put(150,0){\line(0,1){30}}
\put(155,30){\arc{10}{3.142}{4.712}}
\put(230,0){\line(0,1){30}}
\put(225,30){\arc{10}{4.712}{0}}
\put(177.5,15){\line(1,0){25}}
\put(177.5,20){\arc{10}{1.571}{3.142}}
\put(202.5,20){\arc{10}{0}{1.571}}
\put(167.5,20){\arc{10}{4.712}{0}}
\put(212.5,20){\arc{10}{3.142}{4.712}}
\put(220,30){\arc{10}{3.142}{4.712}}
\put(160,30){\arc{10}{4.712}{0}}
\put(155,35){\line(1,0){5}}
\put(225,35){\line(-1,0){5}}
\put(165,25){\line(0,1){5}}
\put(215,25){\line(0,1){5}}
\put(165,25){\line(1,0){2.5}}
\put(212.5,25){\line(1,0){2.5}}
\Thicklines
\put(162.5,25){\line(1,0){2.5}}
\put(162.5,20){\arc{10}{3.142}{4.712}}
\put(157.5,5){\line(0,1){15}}
\put(162.5,5){\arc{10}{1.571}{3.142}}
\put(162.5,0){\line(1,0){55}}
\put(217.5,5){\arc{10}{0}{1.571}}
\put(222.5,5){\line(0,1){15}}
\put(217.5,20){\arc{10}{4.712}{0}}
\put(215,25){\line(1,0){2.5}}
\put(192,0){\vector(1,0){0}}
\thinlines
\dottedline{4}(145,20)(235,20)
\put(180,15){\line(0,1){40}}
\put(200,15){\line(0,1){28}}
\put(200,47){\line(0,1){8}}
\put(157.5,35){\line(0,1){5}}
\put(222.5,35){\line(0,1){5}}
\put(162.5,40){\arc{10}{3.142}{4.712}}
\put(217.5,40){\arc{10}{4.712}{0}}
\put(162.5,45){\line(1,0){15.5}}
\put(182,45){\line(1,0){35.5}}
\put(157.5,35){\arc{5}{3.142}{0}}
\put(222.5,35){\arc{5}{3.142}{0}}
\put(165,25){\arc{5}{3.142}{0}}
\put(215,25){\arc{5}{3.142}{0}}
\put(192,45){\vector(1,0){0}}
\put(180,33){\vector(0,-1){0}}
\put(200,37){\vector(0,1){0}}
\put(145,3){\makebox(0,0){$a$}}
\put(235,3){\makebox(0,0){$a$}}
\put(190,5){\makebox(0,0){$\beta$}}
\put(190,12){\makebox(0,0){$c$}}
\put(180,10){\makebox(0,0){$t$}}
\put(200,10){\makebox(0,0){$s$}}
\put(184,28){\makebox(0,0){$\la$}}
\put(196,26.5){\makebox(0,0){$\mu$}}
\put(190,50){\makebox(0,0){$\nu$}}
\put(165,38){\makebox(0,0){$a'$}}
\put(171,15){\makebox(0,0){$b$}}
\put(209,15){\makebox(0,0){$b$}}
\put(217,38){\makebox(0,0){$a'$}}
\end{picture}
\end{center}
\caption{The vector $\pi_{\la,\mu}(e_\beta)
\Om_{b,c,t,s}^{\la,\mu}\in\cH_{\la,\mu}$}
\label{proofpieb1}
\end{figure}
Now let us look at the part of this picture above
the dotted line. In a suitable Frobenius annulus, this
part can be read for fixed $\nu$ and $a$ as
$\sum_i \la\co\mu(t_i)\epsm\nu{\la\co\mu}t_i^*$,
and the sum runs over a full orthonormal basis of
isometries $t_i$ in the Hilbert space $\Hom(\nu,b\co\beta\co a)$
since we have the summation over $a'$. Next we look
at the part above the dotted line on the right-hand side
of Fig.\ \ref{proofpieb1}. This can be similarly
read for fixed $\nu$ and $a$ as
$\sum_j \la\co\mu(s_j)\epsm\nu{\la\co\mu}s_j^*$,
where the sum runs over another full orthonormal basis of
isometries $s_j\in\Hom(\nu,b\co\beta\co a)$. Since such
bases $\{t_i\}$ and $\{s_j\}$ are related by a
unitary matrix transformation
(this is again just ``unitarity of $6j$-symbols''),
the left- and right-hand side represent
the same vector in $\cH_{\la,\mu}$.
Then, using again Lemma \ref{idexp} and also the trick of
Fig.\ \ref{trick}, we conclude that the vector
$\pi_{\la,\mu}(e_\beta)\Om_{b,c,t,s}^{\la,\mu}$
can be represented by the diagram on the left-hand side of
Fig.\ \ref{proofpieb2}.
%
%proof: the vector \pi(e_\beta)\Om bcts\la\mu (2)
\begin{figure}[htb]
\begin{center}
\unitlength 0.6mm
\begin{picture}(217,50)
\thicklines
\put(32,0){\line(0,1){30}}
\put(98,0){\line(0,1){30}}
\put(42,30){\arc{20}{3.142}{4.712}}
\put(88,30){\arc{20}{4.712}{0}}
\put(42,30){\arc{10}{3.142}{4.712}}
\put(88,30){\arc{10}{4.712}{0}}
\put(37,25){\line(0,1){5}}
\put(93,25){\line(0,1){5}}
\put(42,25){\arc{10}{1.571}{3.142}}
\put(88,25){\arc{10}{0}{1.571}}
\put(42,20){\line(1,0){46}}
\put(42,35){\line(1,0){11}}
\put(57,35){\line(1,0){31}}
\put(42,40){\line(1,0){11}}
\put(57,40){\line(1,0){31}}
\thinlines
\put(55,20){\line(0,1){15}}
\put(55,40){\line(0,1){10}}
\put(75,20){\line(0,1){13}}
\put(75,42){\line(0,1){8}}
\put(55,43){\vector(0,-1){0}}
\put(75,47){\vector(0,1){0}}
\put(48,20){\arc{5}{0}{3.142}}
\put(82,20){\arc{5}{0}{3.142}}
\dottedline{4}(42,0)(88,0)(88,27.5)(42,27.5)(42,0)
\Thicklines
\put(48,20){\line(0,-1){5}}
\put(82,20){\line(0,-1){5}}
\put(58,15){\arc{20}{1.571}{3.142}}
\put(72,15){\arc{20}{0}{1.571}}
\put(58,5){\line(1,0){14}}
\put(75,37){\line(0,1){1}}
\put(55,35){\line(0,1){5}}
\put(67,5){\vector(1,0){0}}
\put(27,3){\makebox(0,0){$a$}}
\put(65,31){\makebox(0,0){$a'$}}
\put(65,17){\makebox(0,0){$c$}}
\put(79,24){\makebox(0,0){$b$}}
\put(51,24){\makebox(0,0){$b$}}
\put(65,9){\makebox(0,0){$\beta$}}
\put(55,15.8){\makebox(0,0){$t$}}
\put(75,15){\makebox(0,0){$s$}}
\put(50,46){\makebox(0,0){$\la$}}
\put(80,45){\makebox(0,0){$\mu$}}
\put(10,25){\makebox(0,0){$\displaystyle\sum_{a,a'}\,d_{a'}$}}
\put(140,12){\makebox(0,0){$\longleftrightarrow
\displaystyle\sum_{c',i,j}\,{\rm{coeff}}_{(c',i,j)}$}}
\thinlines
\dottedline{4}(170,0)(216,0)(216,27.5)(170,27.5)(170,0)
\thicklines
\put(170,20){\line(1,0){46}}
\thinlines
\put(183,20){\line(0,1){7.5}}
\put(203,20){\line(0,1){7.5}}
\put(183,22){\vector(0,-1){0}}
\put(203,25.5){\vector(0,1){0}}
\put(182,11){\makebox(0,0){$t_{a',\co{c'}}^{\la;i}$}}
\put(204,11){\makebox(0,0){$t_{c',\co{a'}}^{\co\mu;j}$}}
\put(175,24){\makebox(0,0){$a'$}}
\put(211,24){\makebox(0,0){$a'$}}
\put(187,24){\makebox(0,0){$\la$}}
\put(199,24){\makebox(0,0){$\mu$}}
\put(193,16.5){\makebox(0,0){$c'$}}
\end{picture}
\end{center}
\caption{The vector $\pi_{\la,\mu}(e_\beta)
\Om_{b,c,t,s}^{\la,\mu}\in\cH_{\la,\mu}$}
\label{proofpieb2}
\end{figure}
Now let us look at the part of the diagram inside the dotted box.
In a suitable Frobenius annulus, this can be interpreted as an
intertwiner in $\Hom(\la\co\mu,a'\co{a'})$. But any element in this
space can be written as a linear combination of elements
constructed from basis isometries $t_{a',\co{c'}}^{\la;i}$,
$t_{c',\co{a'}}^{\co\mu;j}$, as indicated in the dotted box on the
right-hand side of Fig.\ \ref{proofpieb2}. The coefficients in
its linear expansion depend only on $c',i,j$ for fixed $a',\beta,b,c,t,s$,
but certainly not on $a$. This shows that
$\pi_{\la,\mu}(e_\beta)\Om_{b,c,t,s}^{\la,\mu}$
is a linear combination of $\Om_\xi^{\la,\mu}$'s, thus
$\pi_{\la,\mu}(e_\beta)\Om_{b,c,t,s}^{\la,\mu}\in H_{\la,\mu}$.
\end{proof}

The map
$\Om_{b,c,t,s}^{\la,\mu}\mapsto\pi_{\la,\mu}(e_\beta)
\Om_{b,c,t,s}^{\la,\mu}$ defines clearly a linear map
$\pi_{\la,\mu}(e_\beta):H_{\la,\mu}\rightarrow\cH_{\la,\mu}$
since it is just a linear intertwiner multiplication on each
$\Hom(\la\co\mu,a\co a)$ block. From Lemma \ref{miniunit}
we now learn that $\pi_{\la,\mu}(e_\beta)$ is in fact a
linear operator on $H_{\la,\mu}$. With the definition of the
vertical product we now immediately obtain the following

\begin{corollary}
\label{decompieb}
With orthonormal bases
$\{ E_i^{\la,\mu} \}_{i=1}^{{\rm{dim}} H_{\la,\mu}}$
of each $H_{\la,\mu}$ we have
\be
| E_i^{\la,\mu}\rangle\langle E_j^{\la,\mu}| *_v e_\beta *_v
| E_k^{\la',\mu'}\rangle\langle E_l^{\la',\mu'}| \;=\;
\del \la{\la'} \, \del \mu{\mu'} \,
\langle E_j^{\la,\mu}, \pi_{\la,\mu} (e_\beta) E_k^{\la,\mu}\rangle
\; | E_i^{\la,\mu}\rangle\langle\ E_l^{\la,\mu}| \,.
\ee
\end{corollary}

Since $\cZ_h$ is spanned by the $e_\beta$'s, we obtain
a map $\pi_{\la,\mu}:\cZ_h\rightarrow B(H_{\la,\mu})$
by linear extension, and we obtain similarly the following

\begin{corollary}
The map $\pi_{\la,\mu}:\cZ_h\rightarrow B(H_{\la,\mu})$
is a representation of $(\cZ_h,*_v)$.
\end{corollary}

We now tackle the problem of completeness of the system
of matrix units.

\begin{definition}
\label{character}{\rm
For $\la,\mu\in\NXN$ we define the {\sl vertical projector}
$q_{\la,\mu}\in\dta$ by
\be
q_{\la,\mu} =\frac{\sqrt{d_\la d_\mu}}{w^2} \sum_\xi \,
|\Om_\xi^{\la,\mu} \rangle\langle \Om_\xi^{\la,\mu}| \,.
\label{qlmeq}
\ee
}\end{definition}
%
%vertical projector
\thinlines
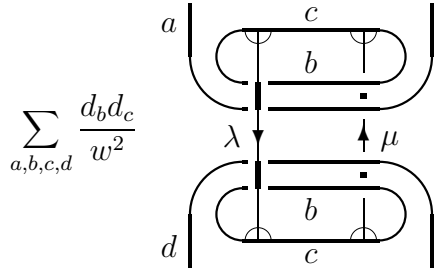
\begin{figure}[htb]
\begin{center}
\unitlength 0.7mm
\begin{picture}(85,50)
\thicklines
\put(15,25){\makebox(0,0){$\displaystyle\sum_{a,b,c,d}
\frac{d_b d_c}{w^2}$}}
\put(37,0){\line(0,1){10}}
\put(83,0){\line(0,1){10}}
\put(47,10){\arc{20}{3.142}{4.712}}
\put(73,10){\arc{20}{4.712}{0}}
\put(47,10){\arc{10}{1.571}{4.712}}
\put(73,10){\arc{10}{4.712}{1.571}}
\put(47,40){\arc{20}{1.571}{3.142}}
\put(73,40){\arc{20}{0}{1.571}}
\put(47,40){\arc{10}{1.571}{4.712}}
\put(73,40){\arc{10}{4.712}{1.571}}
\put(47,5){\line(1,0){26}}
\put(47,15){\line(1,0){1}}
\put(52,15){\line(1,0){21}}
\put(47,20){\line(1,0){1}}
\put(52,20){\line(1,0){21}}
\put(37,50){\line(0,-1){10}}
\put(83,50){\line(0,-1){10}}
\put(47,30){\line(1,0){1}}
\put(52,30){\line(1,0){21}}
\put(47,35){\line(1,0){1}}
\put(52,35){\line(1,0){21}}
\put(47,45){\line(1,0){26}}
\thinlines
\put(50,5){\line(0,1){10}}
\put(50,20){\line(0,1){10}}
\put(70,5){\line(0,1){8}}
\put(70,22){\line(0,1){6}}
\put(50,23){\vector(0,-1){0}}
\put(70,27){\vector(0,1){0}}
\put(50,45){\line(0,-1){10}}
\put(70,45){\line(0,-1){8}}
\put(50,5){\arc{5}{3.142}{0}}
\put(70,5){\arc{5}{3.142}{0}}
\put(50,45){\arc{5}{0}{3.142}}
\put(70,45){\arc{5}{0}{3.142}}
\Thicklines
\put(70,17){\line(0,1){1}}
\put(50,15){\line(0,1){5}}
\put(70,33){\line(0,-1){1}}
\put(50,35){\line(0,-1){5}}
\put(33,4){\makebox(0,0){$d$}}
\put(60,2){\makebox(0,0){$c$}}
\put(60,11){\makebox(0,0){$b$}}
\put(33,46){\makebox(0,0){$a$}}
\put(60,48){\makebox(0,0){$c$}}
\put(60,39){\makebox(0,0){$b$}}
\put(45,25){\makebox(0,0){$\la$}}
\put(75,24){\makebox(0,0){$\mu$}}
\end{picture}
\end{center}
\caption{A vertical projector $q_{\la,\mu}$}
\label{qlamu}
\end{figure}
This is given graphically in Fig.\ \ref{qlamu}. (Clearly, we
can use Lemma \ref{idexp} twice to obtain an equivalent picture
which does not involve pieces of very thick wires
corresponding to $\a_\la^+$ and $\a_\mu^-$.)
We are now ready to prove the main result of this section.

\begin{theorem}
\label{character2}
Under Assumption \ref{set-nondeg},
the vertical projector $q_{\la,\mu}$ is either zero
or a minimal central projection in $(\cZ_h, *_v)$.
We have mutual orthogonality
$q_{\la,\mu} *_v q_{\la',\mu'}=\del \la{\la'}\del \mu{\mu'} q_{\la,\mu}$
and the vertical projectors sum up to the multiplicative
identity of $({\cal Z}_h,*_v)$:
$\sum_{\la,\mu\in\NXN} q_{\la,\mu}=e_0$.
Moreover, $q_{\la,\mu}=0$ whenever $Z_{\la,\mu}=0$
and otherwise the simple summand $q_{\la,\mu} *_v \cZ_h$
is a full $Z_{\la,\mu}\times Z_{\la,\mu}$ matrix algebra,
where $Z_{\la,\mu}$ is the $(\la,\mu)$-entry of the
modular invariant mass matrix of Definition \ref{modular1}.
\end{theorem}

\begin{proof}
It follows from Corollary \ref{matunits} that
$q_{\la,\mu} *_v q_{\la',\mu'}=0$ unless
$\la=\la'$ and $\mu=\mu'$.
We now show that $\sum_{\la,\mu} q_{\la,\mu}=e_0$.
(We denote $e_0\equiv e_{\id_M}$.)
The sum is given graphically by the left-hand side in
Fig.\ \ref{sumqlamu1}.
%
%sum over vertical projectors (1)
\thinlines
\begin{figure}[htb]
\begin{center}
\unitlength 0.6mm
\begin{picture}(205,50)
\thicklines
\put(15,25){\makebox(0,0){$\displaystyle\sum_{a,b,c,d,\la,\mu}
\frac{d_b d_c}{w^2}$}}
\put(37,0){\line(0,1){10}}
\put(83,0){\line(0,1){10}}
\put(47,10){\arc{20}{3.142}{4.712}}
\put(73,10){\arc{20}{4.712}{0}}
\put(47,10){\arc{10}{1.571}{4.712}}
\put(73,10){\arc{10}{4.712}{1.571}}
\put(47,40){\arc{20}{1.571}{3.142}}
\put(73,40){\arc{20}{0}{1.571}}
\put(47,40){\arc{10}{1.571}{4.712}}
\put(73,40){\arc{10}{4.712}{1.571}}
\put(47,5){\line(1,0){26}}
\put(47,15){\line(1,0){1}}
\put(52,15){\line(1,0){21}}
\put(47,20){\line(1,0){1}}
\put(52,20){\line(1,0){21}}
\put(37,50){\line(0,-1){10}}
\put(83,50){\line(0,-1){10}}
\put(47,30){\line(1,0){1}}
\put(52,30){\line(1,0){21}}
\put(47,35){\line(1,0){1}}
\put(52,35){\line(1,0){21}}
\put(47,45){\line(1,0){26}}
\thinlines
\put(50,5){\line(0,1){10}}
\put(50,20){\line(0,1){10}}
\put(70,5){\line(0,1){8}}
\put(70,22){\line(0,1){6}}
\put(50,23){\vector(0,-1){0}}
\put(70,27){\vector(0,1){0}}
\put(50,45){\line(0,-1){10}}
\put(70,45){\line(0,-1){8}}
\put(50,5){\arc{5}{3.142}{0}}
\put(70,5){\arc{5}{3.142}{0}}
\put(50,45){\arc{5}{0}{3.142}}
\put(70,45){\arc{5}{0}{3.142}}
\Thicklines
\put(70,17){\line(0,1){1}}
\put(50,15){\line(0,1){5}}
\put(70,33){\line(0,-1){1}}
\put(50,35){\line(0,-1){5}}
\put(33,4){\makebox(0,0){$d$}}
\put(60,2){\makebox(0,0){$c$}}
\put(60,11){\makebox(0,0){$b$}}
\put(33,46){\makebox(0,0){$a$}}
\put(60,48){\makebox(0,0){$c$}}
\put(60,39){\makebox(0,0){$b$}}
\put(45,25){\makebox(0,0){$\la$}}
\put(75,24){\makebox(0,0){$\mu$}}
\thicklines
\put(115,25){\makebox(0,0){$=\;\displaystyle\sum_{a,b,c,d,\la,\mu,\nu,\rho}
\frac{d_b d_c}{w^2}$}}
\put(152,0){\arc{10}{3.142}{4.712}}
\put(195,0){\arc{10}{4.712}{0}}
%\put(167,10){\line(0,1){30}}
%\put(172,10){\line(0,1){30}}
%\put(162,10){\arc{10}{0}{1.571}}
%\put(162,40){\arc{10}{4.712}{0}}
\put(152,50){\arc{10}{1.571}{3.142}}
\put(195,50){\arc{10}{0}{1.571}}
\put(152,5){\line(1,0){43}}
\put(152,45){\line(1,0){43}}
%\put(177,10){\arc{10}{1.571}{3.142}}
%\put(177,40){\arc{10}{3.142}{4.712}}
%\put(177,5){\line(1,0){18}}
%\put(177,45){\line(1,0){18}}
\thinlines
\put(155,5){\line(0,1){5}}
\put(155,45){\line(0,-1){5}}
\put(169,20){\line(1,0){13}}
\put(169,30){\line(1,0){13}}
\put(155,5){\arc{5}{3.142}{0}}
\put(155,45){\arc{5}{0}{3.142}}
\put(192,45){\arc{5}{0}{3.142}}
\put(192,5){\arc{5}{3.142}{0}}
\put(182,45){\arc{5}{0}{3.142}}
\put(182,5){\arc{5}{3.142}{0}}
\put(167,45){\arc{5}{0}{3.142}}
\put(167,5){\arc{5}{3.142}{0}}
\put(165,10){\arc{20}{3.142}{4.712}}
\put(165,40){\arc{20}{1.571}{3.142}}
\put(182,10){\arc{20}{4.712}{0}}
\put(182,40){\arc{20}{0}{1.571}}
\put(192,5){\line(0,1){5}}
\put(192,45){\line(0,-1){5}}
\put(167,5){\line(0,1){40}}
\put(182,5){\line(0,1){13}}
\put(182,45){\line(0,-1){13}}
\put(182,22){\line(0,1){6}}
\put(182,27){\vector(0,1){0}}
\put(167,23){\vector(0,-1){0}}
\put(176.5,30){\vector(1,0){0}}
\put(172.5,20){\vector(-1,0){0}}
\put(187,24){\makebox(0,0){$\mu$}}
\put(162,25){\makebox(0,0){$\la$}}
\put(153,19){\makebox(0,0){$\rho$}}
\put(153,31){\makebox(0,0){$\nu$}}
\put(144,5){\makebox(0,0){$d$}}
\put(144,45){\makebox(0,0){$a$}}
\put(203,5){\makebox(0,0){$d$}}
\put(203,45){\makebox(0,0){$a$}}
\put(174.5,2){\makebox(0,0){$c$}}
\put(174.5,48){\makebox(0,0){$c$}}
\put(163,9){\makebox(0,0){$b$}}
\put(163,41){\makebox(0,0){$b$}}
\put(186,9){\makebox(0,0){$b$}}
\put(186,41){\makebox(0,0){$b$}}
\end{picture}
\end{center}
\caption{The sum $\sum_{\la,\mu} q_{\la,\mu}$}
\label{sumqlamu1}
\end{figure}
A twofold application of Lemma \ref{idexp} yields the right-hand side
in Fig.\ \ref{sumqlamu1}. Applying Lemma \ref{idexp} twice again,
we obtain the left-hand side of Fig.\ \ref{sumqlamu2}.
%
%sum over vertical projectors (2)
\thinlines
\begin{figure}[htb]
\begin{center}
\unitlength 0.6mm
\begin{picture}(240,60)
\thicklines
\put(24,30){\makebox(0,0){$\displaystyle
\sum_{a,b,c,d\atop\mu,\nu,\rho,\tau,i,j}
\frac{d_c\sqrt{d_\nu d_\rho}}{w^2\sqrt{d_a d_d}}$}}
\put(57,0){\line(0,1){5}}
\put(110,0){\line(0,1){5}}
\put(57,60){\line(0,-1){5}}
\put(110,60){\line(0,-1){5}}
\put(62,5){\arc{10}{3.142}{4.712}}
\put(105,5){\arc{10}{4.712}{0}}
\put(77,15){\line(0,1){30}}
\put(82,15){\line(0,1){30}}
\put(72,15){\arc{10}{0}{1.571}}
\put(72,45){\arc{10}{4.712}{0}}
\put(62,55){\arc{10}{1.571}{3.142}}
\put(105,55){\arc{10}{0}{1.571}}
\put(62,10){\line(1,0){10}}
\put(62,50){\line(1,0){10}}
\put(87,15){\arc{10}{1.571}{3.142}}
\put(87,45){\arc{10}{3.142}{4.712}}
\put(87,10){\line(1,0){18}}
\put(87,50){\line(1,0){18}}
\thinlines
\put(65,10){\line(0,1){40}}
\put(92,50){\arc{5}{0}{3.142}}
\put(92,10){\arc{5}{3.142}{0}}
\put(65,30){\arc{5}{4.712}{1.571}}
\put(102,30){\arc{5}{1.571}{4.712}}
\put(102,10){\line(0,1){40}}
\put(92,10){\line(0,1){18}}
\put(92,50){\line(0,-1){18}}
\put(65,30){\line(1,0){10}}
\put(82,30){\line(1,0){20}}
\put(65,38){\vector(0,-1){0}}
\put(65,18){\vector(0,-1){0}}
\put(102,42){\vector(0,1){0}}
\put(102,22){\vector(0,1){0}}
\put(89,30){\vector(1,0){0}}
\put(92,22){\vector(0,1){0}}
\Thicklines
\put(79,30){\line(1,0){1}}
\put(97,35){\makebox(0,0){$\mu$}}
\put(60,20){\makebox(0,0){$\rho$}}
\put(60,40){\makebox(0,0){$\nu$}}
\put(107,20){\makebox(0,0){$\rho$}}
\put(107,40){\makebox(0,0){$\nu$}}
\put(71,26){\makebox(0,0){$\tau$}}
\put(52,5){\makebox(0,0){$d$}}
\put(52,55){\makebox(0,0){$a$}}
\put(115,5){\makebox(0,0){$d$}}
\put(115,55){\makebox(0,0){$a$}}
\put(84,6){\makebox(0,0){$c$}}
\put(76,53){\makebox(0,0){$b$}}
\put(97,14){\makebox(0,0){$b$}}
\put(97,46){\makebox(0,0){$b$}}
\put(66.5,56.5){\makebox(0,0){$(t_{a\co b}^{\nu;i})^*$}}
\put(102,56.5){\makebox(0,0){$t_{a\co b}^{\nu;i}$}}
\put(65,4){\makebox(0,0){$t_{d\co b}^{\rho;j}$}}
\put(101,4){\makebox(0,0){$(t_{d\co b}^{\rho;j})^*$}}
\thicklines
\put(144,30){\makebox(0,0){$=
\displaystyle\sum_{a,b,c,d,\atop\mu,\nu,\rho,\tau,i,j}
\frac{d_c\sqrt{d_\nu d_\rho}}{w^2\sqrt{d_a d_d}}$}}
\put(177,0){\line(0,1){5}}
\put(230,0){\line(0,1){5}}
\put(177,60){\line(0,-1){5}}
\put(230,60){\line(0,-1){5}}
\put(182,5){\arc{10}{3.142}{4.712}}
\put(225,5){\arc{10}{4.712}{0}}
\put(197,15){\line(0,1){30}}
\put(202,15){\line(0,1){30}}
\put(192,15){\arc{10}{0}{1.571}}
\put(192,45){\arc{10}{4.712}{0}}
\put(182,55){\arc{10}{1.571}{3.142}}
\put(225,55){\arc{10}{0}{1.571}}
\put(182,10){\line(1,0){10}}
\put(182,50){\line(1,0){10}}
\put(207,15){\arc{10}{1.571}{3.142}}
\put(207,45){\arc{10}{3.142}{4.712}}
\put(207,10){\line(1,0){18}}
\put(207,50){\line(1,0){18}}
\thinlines
\put(185,10){\line(0,1){40}}
\put(210,10){\arc{5}{3.142}{0}}
\put(217,10){\arc{5}{3.142}{0}}
\put(185,30){\arc{5}{4.712}{1.571}}
\put(222,30){\arc{5}{1.571}{4.712}}
\put(213.5,38){\arc{7}{3.142}{0}}
\put(222,10){\line(0,1){40}}
\put(217,10){\line(0,1){18}}
\put(217,32){\line(0,1){6}}
\put(210,10){\line(0,1){28}}
\put(185,30){\line(1,0){10}}
\put(202,30){\line(1,0){6}}
\put(212,30){\line(1,0){10}}
\put(185,38){\vector(0,-1){0}}
\put(185,18){\vector(0,-1){0}}
\put(222,42){\vector(0,1){0}}
\put(222,22){\vector(0,1){0}}
\put(194,30){\vector(1,0){0}}
\put(210,18){\vector(0,-1){0}}
\Thicklines
\put(199,30){\line(1,0){1}}
\put(207,40){\makebox(0,0){$\mu$}}
\put(180,20){\makebox(0,0){$\rho$}}
\put(180,40){\makebox(0,0){$\nu$}}
\put(227,20){\makebox(0,0){$\rho$}}
\put(227,40){\makebox(0,0){$\nu$}}
\put(191,26){\makebox(0,0){$\tau$}}
\put(173,5){\makebox(0,0){$d$}}
\put(173,55){\makebox(0,0){$a$}}
\put(235,5){\makebox(0,0){$d$}}
\put(235,55){\makebox(0,0){$a$}}
\put(213.5,14){\makebox(0,0){$c$}}
\put(193,46){\makebox(0,0){$b$}}
\put(220,15){\makebox(0,0){$b$}}
\put(205,55){\makebox(0,0){$b$}}
\put(186.5,56.5){\makebox(0,0){$(t_{a\co b}^{\nu;i})^*$}}
\put(222,56.5){\makebox(0,0){$t_{a\co b}^{\nu;i}$}}
\put(185,4){\makebox(0,0){$t_{d\co b}^{\rho;j}$}}
\put(221,4){\makebox(0,0){$(t_{d\co b}^{\rho;j})^*$}}
\end{picture}
\end{center}
\caption{The sum $\sum_{\la,\mu} q_{\la,\mu}$}
\label{sumqlamu2}
\end{figure}
We can now slide the upper trivalent vertex of the wire $\mu$
around to obtain the right-hand side of Fig.\ \ref{sumqlamu2}.
Next we can use the trick of Fig.\ \ref{trick} to turn around
the small arcs at the trivalent vertices of the wire $\mu$,
yielding a factor $d_\mu/d_c$. This gives the right-hand
side of Fig.\ \ref{sumqlamu2}. Since we have a summation over
$c$, we can again use Lemma \ref{idexp}, and this gives us the
left-hand side of Fig.\ \ref{sumqlamu3}.
%
%sum over vertical projectors (3)
\thinlines
\begin{figure}[htb]
\begin{center}
\unitlength 0.6mm
\begin{picture}(240,60)
\thicklines
\put(24,30){\makebox(0,0){$\displaystyle
\sum_{a,b,d\atop\mu,\nu,\rho,\tau,i,j}
\frac{d_\mu\sqrt{d_\nu d_\rho}}{w^2\sqrt{d_a d_d}}$}}
\put(57,0){\line(0,1){5}}
\put(110,0){\line(0,1){5}}
\put(57,60){\line(0,-1){5}}
\put(110,60){\line(0,-1){5}}
\put(62,5){\arc{10}{3.142}{4.712}}
\put(105,5){\arc{10}{4.712}{0}}
\put(77,15){\line(0,1){30}}
\put(82,15){\line(0,1){30}}
\put(72,15){\arc{10}{0}{1.571}}
\put(72,45){\arc{10}{4.712}{0}}
\put(62,55){\arc{10}{1.571}{3.142}}
\put(105,55){\arc{10}{0}{1.571}}
\put(62,10){\line(1,0){10}}
\put(62,50){\line(1,0){10}}
\put(87,15){\arc{10}{1.571}{3.142}}
\put(87,45){\arc{10}{3.142}{4.712}}
\put(87,10){\line(1,0){18}}
\put(87,50){\line(1,0){18}}
\thinlines
\put(65,10){\line(0,1){40}}
\put(65,30){\arc{5}{4.712}{1.571}}
\put(102,30){\arc{5}{1.571}{4.712}}
\put(93.5,38){\arc{7}{3.142}{0}}
\put(93.5,22){\arc{7}{0}{3.142}}
\put(102,10){\line(0,1){40}}
\put(97,22){\line(0,1){6}}
\put(97,32){\line(0,1){6}}
\put(90,22){\line(0,1){16}}
\put(65,30){\line(1,0){10}}
\put(82,30){\line(1,0){6}}
\put(92,30){\line(1,0){10}}
\put(65,38){\vector(0,-1){0}}
\put(65,18){\vector(0,-1){0}}
\put(102,42){\vector(0,1){0}}
\put(102,22){\vector(0,1){0}}
\put(74,30){\vector(1,0){0}}
\put(90,32){\vector(0,-1){0}}
\Thicklines
\put(79,30){\line(1,0){1}}
\put(87,20){\makebox(0,0){$\mu$}}
\put(60,20){\makebox(0,0){$\rho$}}
\put(60,40){\makebox(0,0){$\nu$}}
\put(107,20){\makebox(0,0){$\rho$}}
\put(107,40){\makebox(0,0){$\nu$}}
\put(71,26){\makebox(0,0){$\tau$}}
\put(52,5){\makebox(0,0){$d$}}
\put(52,55){\makebox(0,0){$a$}}
\put(115,5){\makebox(0,0){$d$}}
\put(115,55){\makebox(0,0){$a$}}
\put(84,6){\makebox(0,0){$b$}}
\put(73,46){\makebox(0,0){$b$}}
\put(66.5,56.5){\makebox(0,0){$(t_{a\co b}^{\nu;i})^*$}}
\put(102,56.5){\makebox(0,0){$t_{a\co b}^{\nu;i}$}}
\put(65,4){\makebox(0,0){$t_{d\co b}^{\rho;j}$}}
\put(101,4){\makebox(0,0){$(t_{d\co b}^{\rho;j})^*$}}
\thicklines
\put(144,30){\makebox(0,0){$=
\displaystyle\sum_{a,b,\nu,i,j}
\frac 1{w d_a}$}}
\put(177,0){\line(0,1){5}}
\put(230,0){\line(0,1){5}}
\put(177,60){\line(0,-1){5}}
\put(230,60){\line(0,-1){5}}
\put(182,5){\arc{10}{3.142}{4.712}}
\put(225,5){\arc{10}{4.712}{0}}
\put(199,15){\line(0,1){30}}
\put(208,15){\line(0,1){30}}
\put(194,15){\arc{10}{0}{1.571}}
\put(194,45){\arc{10}{4.712}{0}}
\put(182,55){\arc{10}{1.571}{3.142}}
\put(225,55){\arc{10}{0}{1.571}}
\put(182,10){\line(1,0){12}}
\put(182,50){\line(1,0){12}}
\put(213,15){\arc{10}{1.571}{3.142}}
\put(213,45){\arc{10}{3.142}{4.712}}
\put(213,10){\line(1,0){12}}
\put(213,50){\line(1,0){12}}
\thinlines
\put(188,10){\line(0,1){40}}
\put(219,10){\line(0,1){40}}
\put(188,28){\vector(0,-1){0}}
\put(219,32){\vector(0,1){0}}
\put(183,30){\makebox(0,0){$\nu$}}
\put(224,30){\makebox(0,0){$\nu$}}
\put(173,5){\makebox(0,0){$a$}}
\put(173,55){\makebox(0,0){$a$}}
\put(235,5){\makebox(0,0){$a$}}
\put(235,55){\makebox(0,0){$a$}}
\put(195,30){\makebox(0,0){$b$}}
\put(212,30){\makebox(0,0){$b$}}
\put(188.5,56.5){\makebox(0,0){$(t_{a\co b}^{\nu;i})^*$}}
\put(220,56.5){\makebox(0,0){$t_{a\co b}^{\nu;i}$}}
\put(187,4){\makebox(0,0){$t_{a\co b}^{\nu;j}$}}
\put(219,4){\makebox(0,0){$(t_{a\co b}^{\nu;j})^*$}}
\end{picture}
\end{center}
\caption{The sum $\sum_{\la,\mu} q_{\la,\mu}$}
\label{sumqlamu3}
\end{figure}
As we have a prefactor $d_\mu$, the summation over
$\mu$ provides a killing ring, and only $\tau=\id_N$
survives it: We obtain a factor $w\del \tau 0$.
Now our picture starts to collapse. The factor
$\del \tau 0$ yields, with the normalization
convention as in Fig.\ \ref{e-beta}, a factor
$d_\nu^{-1} \del \nu\rho$. Since our picture
is now disconnected into two parts which
represent intertwiners in $\Hom(a,d)$, they are
scalars and we obtain a factor $\del ad$. This gives
us the right-hand side of Fig.\ \ref{sumqlamu3}. Therefore
we are now left with a sum over scalars times two straight
vertical wires labelled by $a$, representing a scalar intertwiner
in $\Hom(a\co a,a\co a)$. The scalar value of each
connected part of the picture is
$\del ij \sqrt{d_\nu d_b / d_a}$, therefore we can
compute the prefactor as
\[ \frac 1{w d_a} \sum_{b,\nu}
\sum_{i,j=1}^{N_{a\co b}^\nu} \left( \sqrt{\frac{d_\nu d_b}{d_a}}
\del ij \right)^2 = \frac 1{w d_a^2} \sum_{b,\nu}
d_b N_{a,\co b}^\nu d_\nu = \frac 1{w d_a} \sum_b d_b^2
= \frac 1{d_a} \,. \]
Thus we are left with a sum over two vertical straight wires
with label $a$ and prefactor $d_a^{-1}$. This is $e_0$.

Next, we can expand each vector
$\Om_\xi^{\la,\mu}\in H_{\la,\mu}$,
in an orthonormal basis as
\[ \Om_\xi^{\la,\mu} = \sum_{i=1}^{{\rm{dim}} H_{\la,\mu}}
\langle E_i^{\la,\mu},\Om_\xi^{\la,\mu}\rangle E_i^{\la,\mu} \,.\]
Inserting this in \erf{qlmeq} yields
\[ q_{\la,\mu} = \frac{\sqrt{d_\la d_\mu}}{w^2}
\sum_{i,j}^{{\rm{dim}} H_{\la,\mu}} \sum_\xi \,
\langle E_i^{\la,\mu},\Om_\xi^{\la,\mu} \rangle
\langle \Om_\xi^{\la,\mu}, E_j^{\la,\mu} \rangle
\; |E_i^{\la,\mu} \rangle\langle E_j^{\la,\mu}| \,. \]
Now using $\sum_{\la,\mu} q_{\la,\mu}=e_0$
and Corollary \ref{matunits} we compute
\[ \bearll
\del ij \; |E_i^{\la,\mu} \rangle\langle E_j^{\la,\mu}|
&= \sum_{\la',\mu'} |E_i^{\la,\mu} \rangle\langle E_i^{\la,\mu}|
*_v q_{\la',\mu'} *_v |E_j^{\la,\mu} \rangle\langle E_j^{\la,\mu}|
\\[.5em]
&= \displaystyle\frac{\sqrt{d_\la d_\mu}}{w^2} \sum_\xi \,
\langle E_i^{\la,\mu},\Om_\xi^{\la,\mu} \rangle
\langle \Om_\xi^{\la,\mu}, E_j^{\la,\mu} \rangle
\; |E_i^{\la,\mu} \rangle\langle E_j^{\la,\mu}| \,,
\eear \]
hence
\[ q_{\la,\mu} = \sum_{i=1}^{{\rm{dim}} H_{\la,\mu}}
 |E_i^{\la,\mu} \rangle\langle E_i^{\la,\mu}| \,. \]
Thus $q_{\la,\mu}$ is a projection and we also have
$e_0=\sum_{\la,\mu}\sum_{i=1}^{{\rm{dim}} H_{\la,\mu}}
|E_i^{\la,\mu} \rangle\langle E_i^{\la,\mu}|$.
Hence for any $\beta\in\MXM$ we find
\[ e_\beta = e_0 *_v e_\beta *_v e_0 =
\sum_{\la,\mu}\sum_{i,j=1}^{{\rm{dim}} H_{\la,\mu}}
\langle E_i^{\la,\mu}, \pi_{\la,\mu} (e_\beta) E_j^{\la,\mu}\rangle
\; | E_i^{\la,\mu}\rangle\langle\ E_j^{\la,\mu}| \]
by Corollary \ref{decompieb}. Thus each $e_\beta$ can
be expanded in our matrix units, and since $\cZ_h$ is
spanned by the $e_\beta$'s we conclude that
$\{ | E_i^{\la,\mu}\rangle\langle\ E_j^{\la,\mu}| \}_{\la,\mu,i,j}$
is a complete system of matrix units. It follows that
the non-zero vertical projectors are minimal central projections
in $(\cZ_h,*_v)$, and that the simple summand
$q_{\la,\mu} *_v \cZ_h$ is a full
$\dim H_{\la,\mu} \times \dim H_{\la,\mu}$ matrix algebra.
It remains to show $\dim H_{\la,\mu}=Z_{\la,\mu}$.
The dimension of $H_{\la,\mu}$ can be counted as
\[ \dim H_{\la,\mu} = \sum_{i=1}^{{\rm{dim}} H_{\la,\mu}}
\langle E_i^{\la,\mu}, E_i^{\la,\mu} \rangle =
\sum_{i=1}^{{\rm{dim}} H_{\la,\mu}} \frac 1{d_\la d_\mu}\,
\tau_v(| E_i^{\la,\mu}\rangle\langle E_i^{\la,\mu}|) =
\frac 1{d_\la d_\mu}\, \tau_v (q_{\la,\mu}) \,. \]
Now $d_\la^{-1} d_\mu^{-1} \tau_v (q_{\la,\mu})$
is given graphically in Fig.\ \ref{tauqlm}.
%
% trace q_{la,\mu}
\thinlines
\begin{figure}[htb]
\begin{center}
\unitlength 0.6mm
\begin{picture}(108,50)
\thicklines
\put(18,25){\makebox(0,0){$\displaystyle\sum_{a,b,c}\,
\frac{d_a d_b d_c}{w^2 d_\la d_\mu}$}}
\put(52,5){\line(0,1){5}}
\put(98,5){\line(0,1){5}}
\put(62,10){\arc{20}{3.142}{4.712}}
\put(88,10){\arc{20}{4.712}{0}}
\put(62,10){\arc{10}{1.571}{4.712}}
\put(88,10){\arc{10}{4.712}{1.571}}
\put(62,40){\arc{20}{1.571}{3.142}}
\put(88,40){\arc{20}{0}{1.571}}
\put(62,40){\arc{10}{1.571}{4.712}}
\put(88,40){\arc{10}{4.712}{1.571}}
\put(47,5){\arc{10}{0}{3.142}}
\put(47,45){\arc{10}{3.142}{0}}
\put(103,5){\arc{10}{0}{3.142}}
\put(103,45){\arc{10}{3.142}{0}}
\put(42,5){\line(0,1){40}}
\put(108,5){\line(0,1){40}}
\put(62,5){\line(1,0){26}}
\put(62,15){\line(1,0){1}}
\put(67,15){\line(1,0){21}}
\put(62,20){\line(1,0){1}}
\put(67,20){\line(1,0){21}}
\put(52,45){\line(0,-1){5}}
\put(98,45){\line(0,-1){5}}
\put(62,30){\line(1,0){1}}
\put(67,30){\line(1,0){21}}
\put(62,35){\line(1,0){1}}
\put(67,35){\line(1,0){21}}
\put(62,45){\line(1,0){26}}
\thinlines
\put(65,5){\line(0,1){10}}
\put(65,20){\line(0,1){10}}
\put(85,5){\line(0,1){8}}
\put(85,22){\line(0,1){6}}
\put(65,23){\vector(0,-1){0}}
\put(85,27){\vector(0,1){0}}
\put(65,45){\line(0,-1){10}}
\put(85,45){\line(0,-1){8}}
\put(65,5){\arc{5}{3.142}{0}}
\put(85,5){\arc{5}{3.142}{0}}
\put(65,45){\arc{5}{0}{3.142}}
\put(85,45){\arc{5}{0}{3.142}}
\Thicklines
\put(85,17){\line(0,1){1}}
\put(65,15){\line(0,1){5}}
\put(85,33){\line(0,-1){1}}
\put(65,35){\line(0,-1){5}}
\put(38,10){\makebox(0,0){$a$}}
\put(75,2){\makebox(0,0){$c$}}
\put(75,11){\makebox(0,0){$b$}}
\put(75,48){\makebox(0,0){$c$}}
\put(75,39){\makebox(0,0){$b$}}
\put(60,25){\makebox(0,0){$\la$}}
\put(90,24){\makebox(0,0){$\mu$}}
\end{picture}
\end{center}
\caption{The number $d_\la^{-1} d_\mu^{-1} \tau_v(q_{\la,\mu})$}
\label{tauqlm}
\end{figure}
By the IBFE's we can pull out the circle with label $a$ which
gives us another factor $d_a$. We can therefore proceed with
the summation over $a$, and this yields a factor $w$, the
global index, and then we are left exactly with the picture
in Fig.\ \ref{Zgraph}.
\end{proof}

Note that we learn from the proof that putting
$\Tr_v(z)=\sum_{\la,\mu}d_\la^{-1}d_\mu^{-1}\tau_v(q_{\la,\mu}*_v z)$
for $z\in\cZ_h$ gives a matrix trace $\Tr_v$ on $(\cZ_h,*_v)$ which
sends the minimal projections to one. Next we have learnt that
for all $\la,\mu$ with $Z_{\la,\mu}\neq 0$, the
$\pi_{\la,\mu}$'s are the irreducible representations of
$(\cZ_h,*_v)$ and hence the $\pi_{\la,\mu}\circ\Phi$'s are the
irreducible representations of the $M$-$M$ fusion rule algebra.

\begin{corollary}
\label{commutativity}
Under Assumption \ref{set-nondeg},
the $M$-$M$ fusion rule algebra is commutative
if and only if $Z_{\la,\mu}\in\{0,1\}$
for all $\la,\mu\in\NXN$.
\end{corollary}

\begin{corollary}
\label{MM-number}
Under Assumption \ref{set-nondeg}, the total number
of morphisms in $\MXM$ is equal to
$\tr(Z \,{}^{{\rm{t}}}\!Z)=\sum_{\la,\mu\in\NXN} Z^2_{\la,\mu}$.
\end{corollary}

\subsection{The left action on $M$-$N$ sectors}
\label{sec-actMN}

The decomposition of $(\cZ_h,*_v)$ into simple matrix algebras
is equivalent to the irreducible
decomposition of the ``regular representation'' 
(up to multiplicities given as the dimensions) of the
$M$-$M$ fusion rule algebra, i.e.\ the representation
obtained by its action on itself as a vector space.
There is another representation of the $M$-$M$ fusion rule
algebra, namely the one obtained by its (left) action on
the $M$-$N$ sectors. This is what we study in the
following.

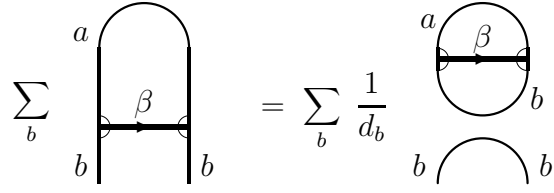
\begin{figure}[htb]
\begin{center}
\unitlength 0.6mm
\begin{picture}(125,40)
\thicklines
\put(5,16){\makebox(0,0){$\displaystyle\sum_b$}}
\put(20,0){\line(0,1){30}}
\put(40,0){\line(0,1){30}}
\put(30,30){\arc{20}{3.142}{0}}
\Thicklines
\put(20,12.5){\line(1,0){20}}
\put(32,12.5){\vector(1,0){0}}
\thinlines
\put(20,12.5){\arc{5}{4.712}{1.571}}
\put(40,12.5){\arc{5}{1.571}{4.712}}
\put(16,33){\makebox(0,0){$a$}}
\put(16,4){\makebox(0,0){$b$}}
\put(44,4){\makebox(0,0){$b$}}
\put(30,17.5){\makebox(0,0){$\beta$}}
\thicklines
\put(70,16){\makebox(0,0){$=\;\displaystyle\sum_b \; \frac 1{d_b}$}}
\put(105,0){\arc{20}{3.142}{0}}
\put(105,30){\arc{20}{3.142}{0}}
\put(105,25){\arc{20}{0}{3.142}}
\put(95,25){\line(0,1){5}}
\put(115,25){\line(0,1){5}}
\Thicklines
\put(95,27.5){\line(1,0){20}}
\put(107,27.5){\vector(1,0){0}}
\thinlines
\put(95,27.5){\arc{5}{4.712}{1.571}}
\put(115,27.5){\arc{5}{1.571}{4.712}}
\put(93,36){\makebox(0,0){$a$}}
\put(91,4){\makebox(0,0){$b$}}
\put(119,4){\makebox(0,0){$b$}}
\put(105,32){\makebox(0,0){$\beta$}}
\put(117,19){\makebox(0,0){$b$}}
\end{picture}
\end{center}
\caption{The element $\varrho(e_\beta) v_{\co a} \in K$}
\label{picrhodef}
\end{figure}

We define the vector space $K$ by
$K=\bigoplus_{a\in\NXM}\Hom(\id_N,a\co a)$. Note that each
block consists just of scalar multiples of the
isometries ${\co r}_a$ but we need the explicit form of $K$.
We define basis vectors $v_{\co a}\in K$ corresponding to
$d_a^{-1/2}{\co r}_a$ in each block $\Hom(\id_N,a\co a)$.
We can display each $v_{\co a}$ graphically by a thick wire
``cap'' with label $a\in\NXM$ together with a prefactor
$1/d_a$. We furnish $K$ with a Hilbert space structure
by putting $\langle v_{\co a},v_{\co b} \rangle=\del ab$.
For each $a\in\NXM$ we define a vector
$\varrho(e_\beta)v_{\co a}$ by putting
\be
\varrho(e_\beta)v_{\co a} = d_\beta \sum_b
N_{\beta,\co a}^{\co b} \, v_{\co b} \,.
\label{rhodef}
\ee
We can display the right-hand side graphically as in
Fig.\ \ref{picrhodef}.
%
% picture for rho
\thinlines
The left- and right-hand side in Fig.\ \ref{picrhodef} are the
same because both sides are scalar multiples of the isometry
${\co r}_a$ in each block $\Hom(\id_N,a\co a)$. The map
$\varrho(e_\beta):v_{\co a}\mapsto\varrho(e_\beta) v_{\co a}$
clearly defines a linear operator on $K$ for each $\beta\in\MXM$,
and we can extend the map $e_\beta\mapsto\varrho(e_\beta)$
linearly to $\cZ_h$. Graphically, this action of $\cZ_h$ is
quite similar to the vertical product. (Note that there also
appears a factor $d_a$ cancelling the $d_a^{-1}$ in the definition
of $v_{\co a}$ when gluing the picture for $v_{\co a}$ on top of
that for $e_\beta$.)

We observe that the map $\varrho:e_\beta\mapsto\varrho(e_\beta)$
extends linearly to a representation of $(\cZ_h,*_v)$ as we can
compute for $\beta,\beta'\in\MXM$ as follows:
\[ \bearll
\varrho(e_\beta) (\varrho(e_{\beta'}) v_{\co a})
&= \varrho(e_\beta) \left(d_\beta \sum_b
N_{\beta',\co a}^{\co b} \, v_{\co b} \right)
= d_\beta d_{\beta'} \sum_{b,c} N_{\beta,\co b}^{\co c}
N_{\beta',\co a}^{\co b} v_{\co c} \\[.4em]
&= d_\beta d_{\beta'} \sum_{\beta'',c} N_{\beta,\beta'}^{\beta''}
N_{\beta'',\co a}^{\co c} v_{\co c}
= d_\beta d_{\beta'} \sum_{\beta'',c} d_{\beta''}^{-1}
N_{\beta,\beta'}^{\beta''} \varrho(e_{\beta''}) v_{\co a} \\[.4em]
&= \varrho (e_\beta *_v e_{\beta'}) v_{\co a} \,,
\eear \]
where we used associativity of the sector product in the
third equality.
Consequently, $\varrho(q_{\la,\mu})$ is a projection onto
a subspace, and $\varrho|_{\varrho(q_{\la,\mu})K}$ is a
subrepresentation.

\begin{lemma}
\label{Kdecom}
We have $K=\bigoplus_{\la\in\NXN}K_\la$,
where $K_\la=\varrho(q_{\la,\la})K$.
\end{lemma}

\begin{proof}
The vector $\varrho(q_{\la,\mu})v_{\co a}\in K$
is given graphically by the left-hand side of
Fig.\ \ref{rhoqlmv}.
%
%the vector \varrho qlamu (v_a)
\thinlines
\begin{figure}[htb]
\begin{center}
\unitlength 0.7mm
\begin{picture}(200,70)
\thicklines
\put(15,30){\makebox(0,0){$\displaystyle\sum_{b,c,d}
\frac{d_b d_c}{w^2}$}}
\put(37,0){\line(0,1){15}}
\put(83,0){\line(0,1){15}}
\put(47,15){\arc{20}{3.142}{4.712}}
\put(73,15){\arc{20}{4.712}{0}}
\put(47,15){\arc{10}{1.571}{4.712}}
\put(73,15){\arc{10}{4.712}{1.571}}
\put(47,55){\arc{20}{1.571}{3.142}}
\put(73,55){\arc{20}{0}{1.571}}
\put(47,55){\arc{10}{1.571}{4.712}}
\put(73,55){\arc{10}{4.712}{1.571}}
\put(47,10){\line(1,0){26}}
\put(47,20){\line(1,0){1}}
\put(52,20){\line(1,0){21}}
\put(47,25){\line(1,0){1}}
\put(52,25){\line(1,0){21}}
\put(37,65){\line(0,-1){10}}
\put(83,65){\line(0,-1){10}}
\put(47,45){\line(1,0){1}}
\put(52,45){\line(1,0){21}}
\put(47,50){\line(1,0){1}}
\put(52,50){\line(1,0){21}}
\put(47,60){\line(1,0){26}}
\put(42,65){\arc{10}{3.142}{4.712}}
\put(78,65){\arc{10}{4.712}{0}}
\put(42,70){\line(1,0){36}}
\thinlines
\put(50,10){\line(0,1){10}}
\put(50,25){\line(0,1){20}}
\put(70,10){\line(0,1){8}}
\put(70,27){\line(0,1){16}}
\put(50,33){\vector(0,-1){0}}
\put(70,37){\vector(0,1){0}}
\put(50,60){\line(0,-1){10}}
\put(70,60){\line(0,-1){8}}
\put(50,10){\arc{5}{3.142}{0}}
\put(70,10){\arc{5}{3.142}{0}}
\put(50,60){\arc{5}{0}{3.142}}
\put(70,60){\arc{5}{0}{3.142}}
\Thicklines
\put(70,22){\line(0,1){1}}
\put(50,20){\line(0,1){5}}
\put(70,48){\line(0,-1){1}}
\put(50,50){\line(0,-1){5}}
\put(33,9){\makebox(0,0){$d$}}
\put(60,7){\makebox(0,0){$c$}}
\put(60,16){\makebox(0,0){$b$}}
\put(33,61){\makebox(0,0){$a$}}
\put(60,63){\makebox(0,0){$c$}}
\put(60,54){\makebox(0,0){$b$}}
\put(45,35){\makebox(0,0){$\la$}}
\put(75,34){\makebox(0,0){$\mu$}}
\put(115,30){\makebox(0,0){$=\;\displaystyle\sum_{b,c,d,i,j}
\;\frac{\del \la\mu}{w^2}$}}
\thicklines
\put(137,0){\line(0,1){15}}
\put(183,0){\line(0,1){15}}
\put(147,15){\arc{20}{3.142}{4.712}}
\put(173,15){\arc{20}{4.712}{0}}
\put(147,15){\arc{10}{1.571}{4.712}}
\put(173,15){\arc{10}{4.712}{1.571}}
\put(147,55){\arc{20}{1.571}{3.142}}
\put(173,55){\arc{20}{0}{1.571}}
\put(147,55){\arc{10}{1.571}{4.712}}
\put(173,55){\arc{10}{4.712}{1.571}}
\put(147,10){\line(1,0){26}}
\put(147,20){\line(1,0){1}}
\put(152,20){\line(1,0){21}}
\put(147,25){\line(1,0){1}}
\put(152,25){\line(1,0){21}}
\put(137,65){\line(0,-1){10}}
\put(183,65){\line(0,-1){10}}
\put(147,45){\line(1,0){1}}
\put(152,45){\line(1,0){21}}
\put(147,50){\line(1,0){1}}
\put(152,50){\line(1,0){21}}
\put(147,60){\line(1,0){26}}
\put(142,65){\arc{10}{3.142}{4.712}}
\put(178,65){\arc{10}{4.712}{0}}
\put(142,70){\line(1,0){36}}
\thinlines
\put(150,10){\line(0,1){10}}
\put(170,10){\line(0,1){8}}
\put(150,25){\line(0,1){2}}
\put(150,45){\line(0,-1){2}}
\put(155,32){\line(1,0){10}}
\put(155,38){\line(1,0){10}}
\put(158,32){\vector(-1,0){0}}
\put(162,38){\vector(1,0){0}}
\put(155,27){\arc{10}{3.142}{4.712}}
\put(165,27){\arc{10}{4.712}{0}}
\put(155,43){\arc{10}{1.571}{3.142}}
\put(165,43){\arc{10}{0}{1.571}}
\put(150,60){\line(0,-1){10}}
\put(170,60){\line(0,-1){8}}
\Thicklines
\put(170,22){\line(0,1){1}}
\put(150,20){\line(0,1){5}}
\put(170,48){\line(0,-1){1}}
\put(150,50){\line(0,-1){5}}
\put(133,9){\makebox(0,0){$d$}}
\put(160,7){\makebox(0,0){$c$}}
\put(160,16){\makebox(0,0){$b$}}
\put(133,61){\makebox(0,0){$a$}}
\put(160,63){\makebox(0,0){$c$}}
\put(160,54){\makebox(0,0){$b$}}
\put(145,31){\makebox(0,0){$\la$}}
\put(175,38){\makebox(0,0){$\la$}}
\put(150,3.5){\makebox(0,0){$t_{b,\co c}^{\la;i}$}}
\put(170,4){\makebox(0,0){$t_{c,\co b}^{\co\la;i}$}}
\put(150,64.5){\makebox(0,0){$(t_{b,\co c}^{\la;i})^*$}}
\put(171,64.5){\makebox(0,0){$(t_{c,\co b}^{\co\la;i})^*$}}
\end{picture}
\end{center}
\caption{The vector $\varrho(q_{\la,\mu})v_{\co a}\in K$}
\label{rhoqlmv}
\end{figure}
Now note that the upper part of the diagram represents
an intertwiner in  $\Hom(\id_N,\la\co\mu)$. Therefore
it vanishes unless $\la=\mu$ and then it must be a scalar
multiple of ${\co r}_\la$. Hence we can insert a term
${\co r}_\la {\co r}_\la^*$ which corresponds
graphically to the disconnection of the wires as on the
right-hand side in Fig.\ \ref{rhoqlmv} and multiplication
by $d_\la^{-1}$. Then the factor $d_b d_c/d_\la$ disappears
because of the normalization convention for trivalent
vertices with small arcs, and we are left
exactly with the right-hand side of Fig.\ \ref{rhoqlmv}.
It follows in particular that $\varrho(q_{\la,\mu})K=0$
unless $\la=\mu$. The claim follows now since the vertical
projectors sum up to $e_0$ and $\varrho(e_0)$ is the
identity on $K$.
\end{proof}

We are now ready to prove the following

\begin{theorem}
\label{varrhodecom}
The representation $\varrho$ of $(\cZ_h,*_v)$
on $K$ obtained by \erf{rhodef}
is unitarily equivalent to the direct sum over the
irreducible representations $\pi_{\la,\la}$:
\be
\varrho \simeq \bigoplus_{\la\in\NXN} \pi_{\la,\la} \,.
\ee
Consequently, the representation $\varrho\circ\Phi$ of the
$M$-$M$ fusion rule algebra which is obtained by the action
on the $M$-$N$ sectors arising from $\MXN$ decomposes
into irreducibles as
$\varrho\circ\Phi\simeq\bigoplus_\la \pi_{\la,\la}\circ\Phi$.
\end{theorem}

\begin{proof}
For $b,c\in\NXM$ and isometries
$t\in\Hom(\la,b\co c)$ and $s\in\Hom(\co\la,c\co b)$ we define
a vector $k_{b,c,t,s}^\la\in K$ by the diagram in
Fig.\ \ref{kbcts}.
%
%The vector k bcts\la
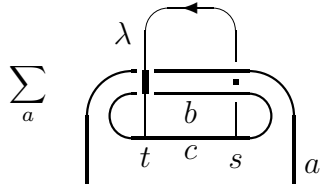
\begin{figure}[htb]
\begin{center}
\unitlength 0.6mm
\begin{picture}(69,40)
\thicklines
\put(17,0){\line(0,1){15}}
\put(63,0){\line(0,1){15}}
\put(27,15){\arc{20}{3.142}{4.712}}
\put(53,15){\arc{20}{4.712}{0}}
\put(27,15){\arc{10}{1.571}{4.712}}
\put(53,15){\arc{10}{4.712}{1.571}}
\put(27,10){\line(1,0){26}}
\put(27,20){\line(1,0){1}}
\put(32,20){\line(1,0){21}}
\put(27,25){\line(1,0){1}}
\put(32,25){\line(1,0){21}}
\thinlines
\put(30,10){\line(0,1){10}}
\put(30,25){\line(0,1){9}}
\put(50,10){\line(0,1){8}}
\put(50,27){\line(0,1){7}}
\put(35,39){\line(1,0){10}}
\put(38,39){\vector(-1,0){0}}
\put(35,34){\arc{10}{3.142}{4.712}}
\put(45,34){\arc{10}{4.712}{0}}
\Thicklines
\put(50,22){\line(0,1){1}}
\put(30,20){\line(0,1){5}}
\put(67,4){\makebox(0,0){$a$}}
\put(40,7){\makebox(0,0){$c$}}
\put(40,16){\makebox(0,0){$b$}}
\put(30,5.8){\makebox(0,0){$t$}}
\put(50,5){\makebox(0,0){$s$}}
\put(25,33.5){\makebox(0,0){$\la$}}
\put(4,20){\makebox(0,0){$\displaystyle\sum_a$}}
\end{picture}
\end{center}
\caption{The vector $k_{b,c,t,s}^\la\in K$}
\label{kbcts}
\end{figure}
Using again intertwiner bases, we also put
$k_\xi^\la=k_{b,c,t_{b,\co c}^{\la;i},t_{c,\co b}^{\co\la;j}}$
with some multi-index $\xi=(b,c,i,j)$. It follows from
the right-hand side in Fig.\ \ref{rhoqlmv} that
$K_\la\subset {\rm span} \{k_\xi^\la\,|\,\xi=(b,c,i,j)\}$.
Conversely, we obtain by Lemma \ref{keyrel} that
$\varrho(q_{\mu,\mu}) k_\xi^\la=0$ unless $\la=\mu$,
hence $K_\la = {\rm span} \{k_\xi^\la\,|\,\xi=(b,c,i,j)\}$.
With $\la=\mu$, closing the wires on the bottom
and on the top on both sides of Fig.\ \ref{keyrelpict} yields
\[ \langle k_\xi^\la,k_{\xi'}^\la \rangle = d_\la
\langle \Om_\xi^{\la,\la},\Om_{\xi'}^{\la,\la} \rangle \,.\]
Hence linear extension of
$\Om_\xi^{\la,\la}\mapsto d_\la^{-1/2} k_\xi^\la$
defines a unitary operator
$U_\la:H_{\la,\la}\rightarrow K_\la$.
Note that $U$ means multiplication by ${\co r}_\la$ from
the right in each block $\Hom(\la\co\la,a\co a)$
and this corresponds graphically to closing the
open ends of the wires $\la$ in Fig.\ \ref{Ombcts}
and multiplying by $d_\la^{-1/2}$. Therefore we find
\[ U \left[\pi_{\la,\la}(e_\beta) \Om_\xi^{\la,\la} \right] =
d_\la^{-1/2} \varrho_\la(e_\beta) k_\xi^\la =
\varrho_\la(e_\beta) U \left[ \Om_\xi^{\la,\la} \right] \,,\]
where $\varrho_\la=\varrho|_{K_\la}$. Thus
$\varrho_\la\simeq\pi_{\la,\la}$.
\end{proof}

Since the dimension of $K$ is the cardinality of $\NXM$ we
immediately obtain the following

\begin{corollary}
\label{NM-number}
Under Assumption \ref{set-nondeg}, the total number of morphisms in
$\NXM$ (or, equivalently, in $\MXN$)
is equal to $\tr(Z)=\sum_{\la\in\NXN} Z_{\la,\la}$.
\end{corollary}

\section{Conclusions and Outlook}
\label{sec-concl}

We have analyzed braided type III subfactors
and shown that in the non-degenerate case the
system of $M$-$M$ system is entirely generated by
$\a$-induction, including in particular the
subsectors of Longo's canonical endomorphism $\can$.
We established that in that case the essential
structural information about the $M$-$M$ fusion rules
is encoded in the modular invariant mass matrix $Z$.
Our setting applies in particular to $\SUn$ loop group
subfactors $\pi^0(\LISUn)''\subset\pi^0(\LIG)''$
of conformal inclusions $\SUn_k\subset G_1$ and
$\pi_0(\LISUn)''\subset\pi_0(\LISUn)''\rtimes_\sigma\bbZ_m$
which were analyzed by $\a$-induction in \cite{BE2,BE3}.
Here $\pi^0$ denotes the level $1$ vacuum representation
of the loop group $\LG$, $\pi_0$ the level $k$ representation
of $\LSUn$, $I\subset S^1$ is an interval,
and $\sigma$ is a ``simple current''.
The braiding here arises from the localized transportable
endomorphisms of the net of local algebras
$A(I)=\pi_0(\LISUn)''$. Since it follows from
Wassermann's work \cite{W2} that these endomorphisms
obey the $\SUn_k$ fusion rules and from the conformal
spin-statistics theorem \cite{GL} that the statistics
phases are given by $\om_\la=\E^{2\pi\I h_\la}$ with
$h_\la$ denoting the $\SUn_k$ conformal dimensions,
it follows that the S- and T-matrices from the braiding coincide
with the well-known S- and T-matrices which transform the
conformal characters. Therefore Theorem \ref{generating}
shows in particular that Condition 4 in Proposition 5.1
in \cite{BE3} holds in the setting of conformal inclusions,
and in turn it proves Conjecture 7.1 in \cite{BE3}.
It also follows that in the setting of Proposition 5.1 in
\cite{BE3}, the sum of $e_\beta$ for ``marked vertices''
$[\beta]$ (the $M$-$M$ sectors arising from the positive
energy representations of the ambient theory) correspond to
the projections appearing in the decomposition
of $\sum_{\la,\mu} p_\la^+ *_h p_\mu^-$, the
``ambichiral projector'' in Ocneanu's language.
Similarly, the results of this paper also prove
Conjecture 7.2 in \cite{BE3}.
Theorem \ref{generating} shows in particular that there are
{\sl no} counter-examples for conformal inclusions where
the $M$-$M$ sectors arising from the conformal inclusion
subfactor are not generated by the mixed $\alpha$-induction
(cf.\ \cite{X2}). Xu made some computation in \cite{X1}
(see also \cite{BE2}) to find an example with non-commutative
fusion rules of ($M$-$M$) sectors generated by the image of
only one ``positive'' induction for subfactors arising from
conformal inclusions.  By Corollary \ref{commutativity},
it is at least very easy to find examples of a
non-commutative entire $M$-$M$ fusion rule algebra.
The D$_4$ case mentioned in
\cite[Subsection 6.1]{BE3} is one such example. In fact,
the whole D$_{2n}$ series arising from simple current extension
of $SU(2)_{4n-4}$ also give examples of non-commutative $M$-$M$
fusion rule algebras. Such non-commutativity for
D$_{{\rm{even}}}$ has been also pointed out in the
setting of \cite{O7} (though not in the context of
conformal inclusions or simple current extensions).

We will present the details and more analysis about $\SUn_k$
loop group subfactors, including the treatment of all $\SUz$
modular invariants, in a forthcoming publication \cite{BEK2}.
Our treatment can now also
incorporate the type II invariants which were not
considered in \cite{BE2,BE3}, because we dropped the
chiral locality condition which automatically forces
the mass matrix $Z$ to be type I, i.e.\ block-diagonal.

Let us remark that we could also have defined $Z_{\la,\mu}$
with exchanged $\pm$-signs in Def.\ \ref{modular1},
and this would correspond to replacing $Z$ by the
transposed mass matrix ${}^{{\rm{t}}}\! Z$.
It is not hard to see that all our calculations
go through with ${}^{{\rm{t}}}\! Z$ as well.
That means $\a$-induction for a (non-degenerately) braided
subfactor determines actually two modular invariant mass
matrices $Z$ and ${}^{{\rm{t}}}\! Z$, and it is not
clear to us at present whether they can in fact be
different in our general setting. (We have
$Z={}^{{\rm{t}}}\! Z$ for all $\SUz$ and $\SUd$
modular invariants).

A notion of subequivalent paragroups was introduced in \cite{K4}.
Since $\NXN$ and $\MXM$ are equivalent 
systems of endomorphisms by definition, $\a$-induction
produces an example of a subequivalent paragroup.  That is,
for $\la\in\NXN$, the subfactors $\a^\pm_\la(M)\subset M$
are subequivalent to $\la(N)\subset N$.  Various examples in
\cite{K4} arise from this construction.
Indeed, the most fundamental example in \cite{K4} comes from
the Goodman-de la Harpe-Jones subfactor \cite[Section 4.5]{GHJ}
with index $3+\sqrt3$.  In our current setting, this example
comes from the conformal inclusion $SU(2)_{10}\subset SO(5)_1$
and shows that the two paragroups with principal graph E$_6$
are subequivalent to the paragroup with principal graph A$_{11}$.

As a corollary of a rigidity theorem presented by Ocneanu in
Madras in January 1997, there are only finitely many
paragroups with global index below a given upper bound.
This implies that for a given paragroup we have only
finitely many subequivalent paragroups
since their global indices are less than or equal to the
global index of the given paragroup. In the context of
modular invariants, a simple argument of Gannon \cite{G1}
shows $\sum_{\la,\mu} Z_{\la,\mu} \le 1/S_{0,0}^2$, which in
turn implies that there are only finitely many modular invariant
mass matrices $Z$ for a given unitary representation of
$\SLZ$, where the S-matrix satisfies the standard relations
$S_{0,\la}\ge S_{0,0}>0$.
As for a non-degenerately braided system of morphisms
this bound coincides with the global index, $w=1/S_{0,0}^2$,
and in view of the relations between modular invariants and
subfactors elaborated in this paper, it is natural to expect
that these two finiteness arguments are not completely
unrelated. We consider a good understanding of the connections
between these two arguments to be highly desirable.

Let us finally remark that in a recent paper of
Rehren \cite{R4} the embedding of left and right chiral
observables in a $2D$ conformal field theory are studied.
Such embeddings give rise to subfactors and in turn to
coupling matrices which are invariant mass matrices if the
Fourier transform matrix of the chiral fusion rules is modular.
As these subfactors are quite different from ours which
appear in a framework considering chiral observables only,
the relation between the two approaches also calls for
a coherent understanding.

%%%%%%%%%%%%%%%%%%%%%%%%%%%%%%%%%%%%%%%%%%%%%%%%%%%%%%%%%%%%%%%%%%%%%%%%%

\vspace{0.5cm}
\begin{footnotesize}
\noindent{\it Acknowledgement.}
Part of this work was done during visits of the third
author to the University of Wales Swansea and the
University of Wales Cardiff, a visit of the
second author to the University of Tokyo, visits of
all the three to Universit\`a di Roma ``Tor Vergata''
and visits of the first two authors to the
Australian National University, Canberra.
We thank R.\ Longo, L.\ Zsido, J.\ E.\ Roberts,
D.\ W.\ Robinson and these institutions for
their hospitality. We would like to thank
S.\ Goto for showing us a preliminary manuscript
of \cite{O7}, M.\ Izumi for explaining \cite{I3},
T.\ Kohno, H.\ Murakami, and T.\ Ohtsuki for helpful
explanations on topological invariants,  and
J.\ E.\ Roberts for his comments.
Y.K. thanks A.\ Ocneanu for various conversations
on \cite{O7} at the Fields Institute in 1995.
We acknowledge the financial support of
the Australian National University,
CNR (Italy), EPSRC (U.K.),
the EU TMR Network in Non-Commutative Geometry,
Grant-in-Aid for Scientific Research, Ministry of
Education (Japan), the Kanagawa Academy of Science and Technology
Research Grants,  the Universit\`a di Roma ``Tor Vergata'',
and the University of Wales.
\end{footnotesize}

\begin{footnotesize}

\end{footnotesize}
\end{document}